\documentclass[11pt]{report}
\usepackage{amsmath}
\usepackage{amssymb}
\usepackage{amsthm}
\usepackage{latexsym}
\usepackage{graphicx}
\usepackage{hyperref}
\usepackage{enumerate}

\newtheorem{thm}{Theorem}[section]
\newtheorem{cor}[thm]{Corollary}
\newtheorem{lem}[thm]{Lemma}
\newtheorem{prop}[thm]{Proposition}

\newtheorem{conj}[thm]{Conjecture}
\newtheorem{cond}[thm]{Condition}
\newtheorem{rem}[thm]{Remark}

\providecommand{\norm}[1]{\left\| #1 \right\|}
\newcommand{\enuma}[1]{\begin{enumerate}[\textup{(}a\textup{)}] {#1} \end{enumerate}}

\newcommand{\mh}{\mathbb}
\newcommand{\mr}{\mathrm}
\newcommand{\mc}{\mathcal}
\newcommand{\mf}{\mathfrak}
\newcommand{\ds}{\displaystyle}

\newcommand{\ts}{\textstyle}

\newcommand{\isom}{\xrightarrow{\;\sim\;}}

\newcommand{\N}{\mathbb N}
\newcommand{\Z}{\mathbb Z}
\newcommand{\Q}{\mathbb Q}
\newcommand{\R}{\mathbb R}
\newcommand{\C}{\mathbb C}
\newcommand{\ep}{\epsilon}

\newcommand{\af}{\mr{aff}}

\newcommand{\pch}{periodic cyclic homology}

\newcommand{\inp}[2]{\langle #1 \,,\, #2 \rangle}
\newcommand{\prefix}[3]{\vphantom{#3}#1#2#3}

\newcommand{\Spr}{\zeta^*}
\newcommand{\WG}{W'_{0}}
\newcommand{\Sc}{\zeta}
\newcommand{\se}{\tilde{\sigma}_\epsilon}
\newcommand{\vv}[2]{\left( \begin{smallmatrix} #1 \\ #2 \end{smallmatrix} \right)}

\begin{document}

\begin{center}
\textbf{\LARGE On the classification of irreducible representations
of affine Hecke algebras\\ with unequal parameters}\\[1cm]
\Large Maarten Solleveld\\[3mm]
\normalsize Mathematisches Institut,
Georg-August-Universit\"at G\"ottingen\\
Bunsenstra\ss e 3-5, 37073 G\"ottingen, Germany\\
email: Maarten.Solleveld@mathematik.uni-goettingen.de \\
March 2011\\[2cm]
\end{center}

\begin{minipage}{13cm}
\textbf{Abstract.}
Let $\mc R$ be a root datum with affine Weyl group $W$, and let $\mc H = \mc H (\mc R,q)$
be an affine Hecke algebra with positive, possibly unequal, parameters~$q$. Then $\mc H$
is a deformation of the group algebra $\C [W]$, so it is natural to compare the representation
theory of $\mc H$ and of $W$.

We define a map from irreducible $\mc H$-representations to $W$-representations and
we show that, when extended to the Grothendieck groups of finite dimensional representations,
this map becomes an isomorphism, modulo torsion. The map can be adjusted to a (nonnatural)
continuous bijection from the dual space of $\mc H$ to that of $W$. We use this to prove 
the affine Hecke algebra version of a conjecture of Aubert, Baum and Plymen, which predicts a 
strong and explicit geometric similarity between the dual spaces of $\mc H$ and $W$.

An important role is played by the Schwartz completion $\mc S = \mc S (\mc R,q)$ of $\mc H$, an
algebra whose representations are precisely the tempered $\mc H$-representations. We construct 
isomorphisms $\Sc_\ep : \mc S (\mc R,q^\ep) \to \mc S (\mc R,q) \; (\ep >0)$ and injection 
$\Sc_0 : \mc S (W) = \mc S (\mc R,q^0) \to \mc S (\mc R,q)$, depending continuously on $\ep$. 

Although $\Sc_0$ is not surjective, it behaves like an algebra isomorphism 
in many ways. Not only does $\Sc_0$ extend to a bijection on Grothendieck groups of finite 
dimensional representations, it also induces isomorphisms on topological $K$-theory and on \pch \
(the first two modulo torsion). This proves a conjecture of Higson and Plymen, which says that 
the $K$-theory of the $C^*$-completion of an affine Hecke algebra $\mc H (\mc R,q)$ does not 
depend on the parameter(s) $q$.\\[2mm]

\end{minipage}

\tableofcontents

\chapter*{Introduction}
\addcontentsline{toc}{chapter}{Introduction}

\newtheorem{thmintro}{Theorem}
\newtheorem{corintro}[thmintro]{Corollary}

Let $\mc R = (X,R_0,Y,R_0^\vee, F_0)$ be a based root datum with finite Weyl group $W_0$ 
and (extended) affine Weyl group $W = W_0 \ltimes X$. For every parameter function 
$q : W \to \C^\times$ there is an affine Hecke algebra $\mc H = \mc H (\mc R,q)$. 

The most important and most studied case is when $q$ takes the same value on all simple
(affine) reflections in $W$. If this value is a prime power, then $\mc H$ is isomorphic to
the convolution algebra of Iwahori-biinvariant functions on a reductive $p$-adic group with root
datum $\mc R^\vee = (Y,R_0^\vee,X,R_0,F_0^\vee)$ \cite{IwMa}. Suppose that the complex
Lie group $G_\C$ with root datum $\mc R$ has simply connected derived group. Then 
the generic Hecke algebra obtained by replacing $q$ with a formal variable $\mathbf q$ is 
known \cite{KaLu1,ChGi} to be isomorphic to the $G_\C \times \C^\times$-equivariant 
$K$-theory of the Steinberg variety of $G_\C$. Via this geometric 
interpretation Kazhdan and Lusztig \cite{KaLu} classified and constructed all irreducible 
representations of $\mc H (\mc R,q)$ (when $q \in \C^\times$ is not a root of unity), 
in accordance with the local Langlands program.

On the other hand, the definition of affine Hecke algebras with generators and relations allows
one to choose the values of $q$ on non-conjugate simple reflections independently. Although
this might appear to be an innocent generalization, much less is known about affine Hecke
algebras with unequal parameters. The reason is that Lusztig's constructions in equivariant
$K$-theory \cite{KaLu1} allow only one deformation parameter. Kato \cite{Kat2}
invented a more complicated variety, called an exotic nilpotent cone, which plays a similar role
for the three parameter affine Hecke algebra of type $C_n^{(1)}$. From this one can extract
a classification of the tempered dual, for an arbitrary parameter function $q$ \cite{CiKa}.

Like equal parameter affine Hecke algebras, those with unequal parameters also arise as
intertwining algebras in the smooth representation theory of reductive $p$-adic groups.
One can encounter them if one looks at non-supercuspidal Bernstein components (in the smooth 
dual) \cite{Lus-Uni,Mor}. Even for split groups unequal parameters occur, albeit not in the
principal series \cite{Roc1}.
It is expected that every Bernstein component of a $p$-adic group can be described 
with an affine Hecke algebra or a slight variation on it. 
In \cite{BaMo1,BaMo2,BaCi} it is shown, in increasing generality, that under certain conditions 
the equivalence between the module category of an affine Hecke algebra and a Bernstein block 
(in the category of smooth modules) respects unitarity. Thus affine Hecke algebras are an 
important tool for the classification of both the smooth dual and the unitary smooth dual of a 
reductive $p$-adic group.

The degenerate versions of affine Hecke algebras are usually called graded Hecke algebras. 
Their role in the representation theory of reductive $p$-adic groups \cite{Lus-Gr,Lus-Uni,BaMo2,Ciu}, 
is related to affine Hecke algebras in the way that Lie algebras stand to Lie groups. They have a 
very simple presentation, which makes is possible to let them act on many function spaces. 
Therefore one encounters graded Hecke algebras (with possibly unequal parameters) also 
independently from affine Hecke algebras, for instance in certain systems of differential equations 
\cite{HeOp} and in the study of the unitary dual of a real reductive group \cite{CiTr}.

In view of the above connections, it is of considerable interest to classify the dual of an affine or
graded Hecke algebra with unequal parameters. Since $\mc H (\mc R,q)$ is a deformation of
the group algebra $\C [W]$, it is natural to expect strong similarities between the duals 
Irr$(\mc H (\mc R,q))$ and Irr$(W)$. Indeed, for equal parameters the 
Deligne--Langlands--Kazhdan--Lusztig parametrization provides a bijection between these duals
\cite{KaLu,Lus-Rep}.
For unequal parameters we approach the issue via harmonic analysis on affine Hecke algebras,
which forces us to consider only parameter functions $q$ with values in $\R_{>0}$. We will assign
to every irreducible $\mc H$-representation $\pi$ in a natural way a representation $\Spr (\pi)$
of the extended affine Weyl group $W$. Although this construction does not always preserve
irreducibility, it has a lot of nice properties, the most important of which is:

\begin{thmintro} \label{thm:0.1} 
\textup{(see Theorem \ref{thm:2.7})} \\
The collection of representations $\{ \Spr (\pi) : \pi \in \mr{Irr}(\mc H (\mc R,q)) \}$ forms a
$\Q$-basis of the representation ring of $W$.
\end{thmintro}

Since the Springer correspondence for finite Weyl groups realizes Irr$(W_0)$, via 
Kazhdan--Lusztig theory, as a specific subset of Irr$(\mc H (\mc R,q))$, Theorem \ref{thm:0.1} 
can be regarded as an affine Springer correspondence without irreducibility. By picking suitable
irreducible subquotients of the $\Spr (\pi)$ one can refine Theorem \ref{thm:0.1} to a bijection 
Irr$(\mc H (\mc R,q)) \to \text{Irr}(W)$. This is related to a conjecture of Aubert, Baum and 
Plymen \cite{ABP1,ABP2} (ABP-conjecture for short) which we sketch here. 

Recall that $W = W_0 \ltimes X$ and let $T$ be the complex torus $\mr{Hom}_\Z (X,\C^\times)$. 
Clifford theory says that the irreducible $W$-representations with an $X$-weight $t \in T$ are
in natural bijection with the irreducible representations of the isotropy group $W_{0,t}$. 
Write $T^w = \{ t \in T : w(t) = t \}$. The extended quotient of $T$ by $W_0$ is defined as
\[
\widetilde{T} / W_0 = \bigsqcup\nolimits_{w \in W_0} \{w\} \times T^w / W_0 ,
\]
with respect to the $W_0$-action $w \cdot (w',t) = (w w' w^{-1},w t)$. We endow $\widetilde{T} / W_0$
with the topology coming from the discrete topology on $W_0$ and the Zariski topology on $T$.
It is a model for Irr$(W)$ with the Jacobson topology, in the sense that there exist continuous 
bijections $\widetilde{T} / W_0 \to \text{Irr}(W)$ which respect the projections to $T / W_0$.

The Bernstein presentation (included as Theorem \ref{thm:1.1}) shows that $\C [X] \cong \mc O (T)$
is naturally to isomorphic to a commutative subalgebra $\mc A \subset \mc H (\mc R,q)$, and that
the center of $\mc H (\mc R,q)$ is $\mc A^{W_0} \cong \mc O (T/W_0)$. Hence we have a natural
map Irr$(\mc H (\mc R,q)) \to T / W_0$, sending a representation to its central character. It is 
continuous with respect to the Jacobson topology on Irr$(\mc H (\mc R,q))$. A simplified
version of the ABP-conjecture for affine Hecke algebras reads:

\begin{thmintro} \label{thm:0.2}
\textup{(see Theorem \ref{thm:5.9})} \\
Let $\mc H (\mc R,q)$ be an affine Hecke algebra with positive, possibly unequal, parameters.
Let $\mc Q (\mc R)$ be the variety of parameter functions $W \to \C^\times$, endowed with
the Zariski topology. There exist a
continuous bijection $\mu : \widetilde{T} / W_0 \to \mr{Irr}(\mc H (\mc R,q))$ and a map
$h : \widetilde{T} / W_0 \times \mc Q (\mc R) \to T$ such that:
\begin{itemize}
 \item $h$ is locally constant in the first argument;
 \item for fixed $c \in  \widetilde{T} / W_0, \; h(c,v)$ is a monomial in the variables $v(s)^{\pm 1}$,
where $s$ runs through all simple affine reflections;
 \item the central character of $\mu (W_0 (w,t))$ is $W_0 \, h(W_0 (w,t), q^{1/2}) t$.
\end{itemize}
\end{thmintro}

The author hopes that Theorem \ref{thm:0.2} will be useful in the local Langlands program. 
This could be the case for a Bernstein component $\mf s$ of a reductive $p$-adic group $G$
for which $\mr{Mod}_{\mf s} (G)$ is equivalent to Mod$(\mc H (\mc R,q))$ (see Section 
\ref{sec:padic} for the notations). Recall that a Langlands parameter for $G$ is a certain kind 
of group homomorphism $\phi : \mc W_{\mathbb F} \ltimes \C \to \prefix{^L}{}{G}$, where 
$\prefix{^L}{}{G}$ is the Langlands dual group and $\mc W_{\mathbb F}$ is the Weil group
of the local field $\mathbb F$ over which $G$ is a variety. Let us try to extract a Langlands 
parameter from $W_0 (w,t) \in \widetilde{T} / W_0$. The image of the distinguished Frobenius
element of $\mc W_{\mathbb F}$ describes the central character, so modulo conjugacy it
should be $h(W_0 (w,t), q^{1/2}) t \in T \subset \prefix{^L}{}{G}$. 

Problematic is that the map $\mu$ from Theorem \ref{thm:0.2} is not canonical, its construction 
involves some arbitrary choices, which could lead to a different $w \in W_0$. Yet it is precisely this 
element $w$ that should determine the unipotent conjugacy class that contains the image of 
$\C \setminus \{0\}$ under $\phi$. Recently Lusztig \cite{Lus-conj} defined a map from conjugacy
classes in a Weyl group to unipotent classes in the associated complex Lie group. Whether 
this map yields a suitable unipotent class for our hypothetical $\phi$ remains to be seen, for this
it is probably necessary to find a more canonical construction of $\mu$. The rest of such a 
Langlands parameter $\phi$ is completely beyond affine Hecke algebras, it will have to depend
on number theoretic properties of the Bernstein component $\mf s$. 

Now we describe the most relevant results needed for Theorems \ref{thm:0.1} and \ref{thm:0.2}.
A large step towards the determination of Irr$(\mc H (\mc R,q))$ is the Langlands classification
(see Theorem \ref{thm:2.4}), which reduces the problem to the tempered duals of parabolic 
subalgebras of $\mc H (\mc R,q)$. This is of course a well-known result, for graded Hecke algebras
due to Evens \cite{Eve}, but the proof for affine Hecke algebras has not been published before.

In line with Harish-Chandra's results for reductive groups, every tempered irreducible 
$\mc H (\mc R,q)$-representation appears as a direct summand of a representation that is induced from a 
discrete series representation of a parabolic subalgebra, a point of view advocated by Opdam \cite{Opd-Sp}. 
We note that in this setting tempered and discrete series representations can be defined very easily via the 
$\mc A$-weights of a representation. For affine Hecke algebras with irreducible root data, Opdam
and the author classified the discrete series in \cite{OpSo2}, but we emphasize that we do not use
that classification in the present paper. 

To every set of simple roots $P \subset F_0$ we associate a root subsystem $R_P \subset R_0$ and a 
parabolic subalgebra $\mc H_P \subset \mc H$, which generated by $\mc H (W (R_P),q)$ and by a
commutative subalgebra $\mc A_P \cong \C [X / X \cap \Q P]$. Induction from $\mc H_P$
allows an induction parameter in $T^P$, the subtorus of $T$ orthogonal to $P^\vee$. So we consider
induction data $\xi = (P,\delta,t)$ where $P \subset F_0, t \in T^P$ and $\delta$ is a discrete series 
representation of $\mc H_P$. Such a triple gives rise to a parabolically induced representation
\[
\pi (\xi) = \pi (P,\delta ,t) = \mr{Ind}_{\mc H^P}^{\mc H} (\delta \circ \phi_t) .
\]
Here $\mc H^P \subset \mc H$ is more or less a central extension of $\mc H_P$ by $\mc O (T^P)$
and $\phi_t : \mc H^P \to \mc H_P$ is a twisted projection. The representation $\pi (\xi)$ is tempered if
and only if $t \in T^P_{un}$, the unitary part of  $T^P$. For the classification of the dual it 
remains to decompose all tempered parabolically induced representations and to determine when they
have constituents in common.

These phenomena are governed by intertwining operators between such representations
\cite[Section 4]{Opd-Sp}. It is already nontrivial to show that these operators are well-defined
on all $\pi (\xi)$ with $t \in T^P_{un}$, and it is even more difficult to see that they span
$\mr{Hom}_{\mc H}(\xi,\xi')$. For reductive groups this is known as Harish-Chandra's completeness
theorem, and for affine Hecke algebras it is a deep result due to Delorme and Opdam \cite{DeOp1}.
The intertwining operators can be collected in a groupoid $\mc G$, which acts on the space $\Xi$
of induction data $(P,\delta,t)$. With these tools one can obtain a partial classification of the
dual of $\mc H (\mc R,q)$:

\begin{thmintro} \label{thm:0.3}
\textup{(see Theorem \ref{thm:3.10})} \\
There exists a natural map $\mr{Irr}(\mc H (\mc R ,q)) \to \Xi / \mc G ,\; 
\rho \mapsto \mc G \xi^+ (\rho)$, such that:
\begin{itemize}
\item the map is surjective and has finite fibers;
\item $\rho$ is a subquotient of $\pi (\xi^+ (\rho))$;
\item $\rho$ does not occur in $\pi (\xi)$ when $\xi$ is larger than $\xi^+ (\rho)$,
where the size of $\xi = (P,\delta,t)$ is the norm of the absolute value of the central 
character of $\delta$.
\end{itemize}
\end{thmintro}

We note that this part of the paper is rather similar to \cite{SolGHA} for graded Hecke algebras. 
The results here are somewhat stronger, and most of them can not be derived quickly from their
counterparts in \cite{SolGHA}.

Yet all this does not suffice for Theorem \ref{thm:0.1}, because we do not have much control over
the number of irreducible constituents of parabolically induced representations. Ultimately the proof 
of Theorem \ref{thm:0.1} is reduced to irreducible tempered representations with central
character in $\mr{Hom}_\Z (X ,\R_{>0})$. The author dealt with this case in \cite{SolHomGHA},
via the \pch \ of graded Hecke algebras.

What we discussed so far corresponds more or less to chapters 1--3 of the article. From Theorem 
\ref{thm:0.1} to Theorem \ref{thm:0.2} is not a long journey, but we put Theorem \ref{thm:0.2}
at the end of the paper because we prove it together with other parts of the ABP-conjecture.

Chapters 4 and 5 are of a more analytic nature. The main object of study is the Schwartz 
algebra $\mc S (\mc R,q)$ of $\mc H (\mc R,q)$ \cite{Opd-Sp,DeOp1}, the analogue of the 
Harish-Chandra--Schwartz algebra of a reductive $p$-adic group. By construction a 
$\mc H (\mc R,q)$-representation extends continuously to $\mc S (\mc R,q)$ if and only if 
it is tempered.
The Schwartz algebra is the ideal tool for the harmonic analysis of affine Hecke algebras,
among others because it admits a very nice Plancherel theorem (due to Delorme and Opdam,
see \cite{DeOp1} or Theorem \ref{thm:3.8}), because the discrete series of $\mc H (\mc R,q)$
is really discrete in the dual of $\mc S (\mc R,q)$, and because the inclusion $\mc H (\mc R,q)
\to \mc S (\mc R,q)$ preserves Ext-groups of tempered representations \cite{OpSo1}.

If we vary the parameter function $q$, we obtain families of algebras $\mc H (\mc R,q)$ and
$\mc S (\mc R,q)$. It is natural to try to connect the representation theory of $\mc H (\mc R,q)$
with that of $\mc H (\mc R,q')$ for $q'$ close to $q$. For general parameter deformations this
is too difficult at present, but we can achieve it for deformations of the form $q \mapsto q^\ep$
with $\ep \in \R$. On central characters of representations, this "scaling" of $q$ fits well 
with the map 
\[
\sigma_\ep : T \to T ,\; t \mapsto t \, |t|^{\ep -1} .
\]
Notice that $\sigma_\ep (t) = t$ for $t \in T_{un} = \mr{Hom}_\Z (X,S^1)$. 
Let $\{ N_w : w \in W \}$ be the standard basis of $\mc H (\mc R,q)$ and
let $\mr{Mod}_{f,W_0 t} (\mc H (\mc R,q))$ be the category of finite dimensional 
$\mc H (\mc R,q)$-representations with central character $W_0 t \in T / W_0$. 

\begin{thmintro}\label{thm:0.4}
\textup{(see Corollary \ref{cor:4.4})} \\
There exists a family of additive functors 
\[
\tilde \sigma_{\ep,t} : \mr{Mod}_{f,W_0 t} (\mc H (\mc R,q)) \to 
\mr{Mod}_{f,W_0 \sigma_\ep (t)} (\mc H (\mc R,q^\ep)) \qquad \ep \in [-1,1] ,
\]
such that:
\begin{itemize}
\item for all $(\pi,V) \in \mr{Mod}_{f,W_0 t} (\mc H (\mc R,q))$ and all $w \in W$, the map
$\ep \mapsto \tilde \sigma_{\ep,t} (\pi) (N_w) \in \mr{End}_\C (V)$ is analytic;
\item for $\ep \neq 0 ,\; \tilde \sigma_{\ep,t}$ is an equivalence of categories;
\item $\tilde \sigma_{\ep,t}$ preserves unitarity.
\end{itemize}
\end{thmintro}

For $\ep > 0$ this can already be found in \cite{Opd-Sp}, the most remarkable part is
precisely that it extends continuously to $\ep = 0$, that is, to the algebra
$\mc H (\mc R,q^0) = \C [W]$. In general the functors $\tilde \sigma_{\ep,t}$ cannot
be constructed if we work only inside the algebras $\mc H (\mc R,q^\ep)$, they are obtained 
via localizations of these algebras at certain sets of central characters. We can do much 
better if we replace the algebras by their Schwartz completions:

\begin{thmintro}\label{thm:0.5}
\textup{(see Theorem \ref{thm:4.8})} \\
For $\ep \in [0,1]$ there exist homomorphisms of Fr\'echet *-algebras $\Sc_\ep :
\mc S (\mc R,q^\ep) \to \mc S (\mc R,q)$, such that:
\begin{itemize}
\item $\Sc_\ep$ is an isomorphism for all $\ep > 0$, and $\Sc_1$ is the identity;
\item for all $w \in W$ the map $\ep \mapsto \Sc_\ep (N_w)$ is piecewise analytic;
\item for every irreducible tempered $\mc H (\mc R,q)$-representation $\pi$ with central
character $W_0 t$, the $\mc S (\mc R,q^\ep)$-representations $\tilde \sigma_{\ep,t}(N_w)$
and $\pi \circ \Sc_\ep$ are equivalent.
\end{itemize}
\end{thmintro}

There are some similarities with the role played by Lusztig's asymptotic Hecke algebra 
\cite{Lus-C2} in \cite{KaLu}. In both settings an algebra is constructed, which contains
$\mc H (\mc R,q)$ for a family of parameter functions $q$. The asymptotic Hecke algebra
is of finite type over $\mc O (T / W_0)$, so it is only a little larger than $\mc H (\mc R,q)$.
So far it has only been constructed for equal parameter functions $q$, but Lusztig \cite{Lus-Une}
conjectures that it also exists for unequal parameter functions. On the hand, the algebra
$\mc S (\mc R,q)$ is of finite type over $C^\infty (T_{un})$, so it is much larger than
$\mc H (\mc R,q)$. Although $\Sc_\ep$ is an isomorphism for $\ep \in (0,1]$, the algebras
$\mc H (\mc R,q^\ep)$ are embedded in $\mc S (\mc R,q)$ in a nontrivial way, in most cases
$\Sc_\ep (\mc H (\mc R,q^\ep))$ is not contained in $\mc H (\mc R,q)$.

Of particular interest is the homomorphism
\begin{equation}\label{eq:0.Sc0}
\Sc_0 : \mc S (W) = \mc S (\mc R,q^0) \to \mc S (\mc R,q) .
\end{equation}
It cannot be an isomorphism, but it is injective and for all irreducible tempered 
$\mc H (\mc R,q)$-representations $\pi$ with central character $W_0 t$ we have 
$\pi \circ \Sc_0 \cong \tilde \sigma_{0,t} (\pi) \cong \Spr (\pi)$. 
Together with Theorem \ref{thm:0.1} this results in:

\begin{corintro}\label{cor:0.6}
\textup{(see Corollary \ref{cor:Sprsigma0})} \\
The functor $\mr{Mod} (\mc S (\mc R,q)) \to \mr{Mod} (\mc S (W)) : \pi \mapsto \pi \circ \Sc_0$
induces an isomorphism between the Grothendieck groups of finite dimensional representations,
tensored with~$\Q$.
\end{corintro}

So Theorem \ref{thm:0.1} does not stand alone, but forms the end of a continuous family of
representations (of a family of algebras). Actually the author first discovered the algebra
homomorphism $\Sc_0$ and only later realized that the corresponding map on representations
can also be obtained in another way, thus gaining in naturality.

Apart from representation theory, the aformentioned results have some interesting consequences
in the noncommutative geometry of affine Hecke algebras. Let $C^* (\mc R,q)$ be the 
$C^*$-completion of $\mc H (\mc R,q)$. It contains $\mc S (\mc R,q)$ and $\Sc_\ep$ extends to
a $C^*$-algebra homomorphism $\Sc_\ep : C^* (\mc R,q^\ep) \to C^* (\mc R,q)$, for which
Theorem \ref{thm:0.5} remains valid. It follows quickly from this and Corollary \ref{cor:0.6}
that $\Sc_0$ induces an isomorphism on topological $K$-theory, see Theorem \ref{thm:5.3}.
More precisely,
\begin{equation}\label{eq:0.KQ}
K_* (\Sc_0) \otimes \mr{id}_\Q : K_* (C^* (W) \rtimes \Gamma) \otimes_\Z \Q \to
K_* (C^* (\mc R,q) \rtimes \Gamma) \otimes_\Z \Q
\end{equation}
is an isomorphism, while for equal parameters the argument also goes through without
$\otimes_\Z \Q$. This solves a conjecture that was posed first by Higson and Plymen \cite{Ply1,BCH}.

Furthermore $C^* (\mc R,q)$ and $\mc S (\mc R,q)$ have the same topological $K$-theory, and
via the Chern character the complexification of the latter is isomorphic to the \pch \ of
$\mc S (\mc R,q)$. As already proved in \cite{SolPadic}, $\mc H (\mc R,q)$ and $\mc S (\mc R,q)$
have the same \pch , so we obtain a commutative diagram 
\[
\begin{array}{*{7}{c}}
\!\! HP_* (\C [W]) & \!\!\! \to \!\!\! & HP_* (\mc S (W))  & \!\!\! \leftarrow \!\!\! 
& K_* (\mc S (W)) \otimes_\Z \C & \!\!\! \to \!\!\! & K_* (C^* (W)) \otimes_\Z \C \\
\downarrow & & \downarrow \scriptstyle{HP_* (\Sc_0)} & & \downarrow  \scriptstyle{K_* (\Sc_0)} & & 
\downarrow  \scriptstyle{K_* (\Sc_0)}  \\
\!\! HP_* (\mc H (\mc R,q))  & \!\!\! \to \!\!\! & HP_* (\mc S (\mc R,q)) & 
\!\!\! \leftarrow \!\!\! &  K_* (\mc S (\mc R,q)) \otimes_\Z \C & \!\!\! \to \!\!\! & 
K_* (C^* (\mc R ,q)) \otimes_\Z \C ,
\end{array}
\]
where all the arrows are natural isomorphisms (see Corollary \ref{cor:5.6}).
Notice that the Schwartz algebra $\mc S (\mc R,q)$ forms a bridge between the purely algebraic 
$\mc H (\mc R,q)$ and the much more analytic $C^* (\mc R,q)$.

For the sake of clarity, the introduction is written in less generality than the actual paper.
Most notably, we can always extend our affine Hecke algebras by a group $\Gamma$ of automorphisms
of the Dynkin diagram of $\mc R$. On the one hand this generality is forced upon us, in particular
by Lusztig's first reduction theorem (see Theorem \ref{thm:2.1}), which necessarily involves
diagram automorphisms. On the other hand, one advantage of having $\mc H (\mc R,q) \rtimes \Gamma$
instead of just $\mc H (\mc R,q)$ is that our proof of the Aubert--Baum--Plymen conjecture
applies to clearly more Bernstein components of reductive $p$-adic groups. 

For most of the results of this paper, the extension from $\mc H (\mc R,q)$ to $\mc H (\mc R,q)
\rtimes \Gamma$ is easy, mainly a matter of some extra notation. An exception is the Langlands
classification, which hitherto was only known for commutative groups of diagram automorphisms
\cite{BaJa1}. In our generalization (see Corollary \ref{cor:2.8}) we add a new ingredient to the
Langlands data, and we show how to save the uniqueness part.

A substantial part of this article is based on the author's PhD-thesis \cite{SolThesis}, which
was written under the supervision of Opdam. We refrain from indicating all the things that
stem from \cite{SolThesis}, among others because some of proofs in \cite{SolThesis} were not
worked out with the accuracy needed for research papers. Moreover, in the years after writing
this thesis many additional insights were obtained, so that in the end actually no part of
\cite{SolThesis} reached this article unscathed. The technical Chapter 4 comes closest.
It should also be mentioned that the conjecture \eqref{eq:0.KQ} formed a central part of the
author's PhD-research. At that time it was still too difficult for the author, mainly because
Theorem \ref{thm:0.1} was not available yet. \\[2mm]

\textbf{Acknowledgements.}
The author learned a lot about affine Hecke algebras from Eric Opdam, first as a PhD-student
a later as co-author. Without that support, this article would not have been possible.
The author also thanks Roger Plymen for providing background information about several conjectures,
and Anne-Marie Aubert for many detailed comments, in particular concerning Section \ref{sec:padic}.
Finally, two anonymous referees read the paper very carefully and made many useful suggestions
for improvements.

\chapter{Preliminaries}

This chapter serves mainly to introduce some definitions and notations that we will use later on.
The results that we recall can be found in several other sources, like \cite{Lus-Gr,Ree1,Opd-Sp}.
By default, our affine Hecke algebras are endowed with unequal parameters and may be extended
with a group of automorphism of the underlying root datum.

In the section dedicated to $p$-adic groups we recall what is known about the (conjectural)
relation between Bernstein components and affine Hecke algebras. On one hand this motivates the
generality that we work in, on the other hand we will use it in Section \ref{sec:ABP} to translate
a conjecture of Aubert, Baum and Plymen to the setting of affine Hecke algebras.

\section{Root systems}

Let $\mf a$ be a finite dimensional real vector space and let $\mf a^*$ be its dual. 
Let $Y \subset \mf a$ be a lattice and $X = \mr{Hom}_\Z (Y,\Z) \subset \mf a^*$ the dual lattice. Let
\[
\mc R = (X, R_0, Y ,R_0^\vee ,F_0) .
\]
be a based root datum. Thus $R_0$ is a reduced root system in $X ,\, R^\vee_0 \subset Y$ 
is the dual root system, $F_0$ is a basis of $R_0$ and the set of positive roots is denoted $R_0^+$.
Furthermore we are given a bijection $R_0 \to R_0^\vee ,\: \alpha \mapsto \alpha^\vee$ such 
that $\inp{\alpha}{\alpha^\vee} = 2$ and such that the corresponding reflections
$s_\alpha : X \to X$ (resp. $s^\vee_\alpha : Y \to Y$) stabilize $R_0$ (resp. $R_0^\vee$).
We do not assume that $R_0$ spans $\mf a^*$. 

The reflections $s_\alpha$ generate the Weyl group $W_0 = W (R_0)$ of $R_0$, 
and $S_0 := \{ s_\alpha : \alpha \in F_0 \}$ is the collection of simple reflections. 
We have the affine Weyl group and its extended version:
\[
\begin{array}{rrrrr}
W^\af & = & W^\af (\mc R) & = & \mh Z R_0 \rtimes W_0 , \\
W & = & W (\mc R) & = & X \rtimes W_0 . 
\end{array}
\]
Both can be considered as groups of affine transformations of $\mf a^*$. 
We denote the translation corresponding to $x \in X$ by $t_x$.
As is well known, $W^\af$ is a Coxeter group, and the basis of $R_0$ gives rise to a set $S^\af$ 
of simple (affine) reflections. More explicitly, let $F_M^\vee$ be the set of maximal elements of
$R_0^\vee$, with respect to the dominance ordering coming from $F_0$. Then
\[
S^\af = S_0 \cup \{ t_\alpha s_\alpha : \alpha \in F_M \} .
\]
We write
\begin{align*}
&X^+ := \{ x \in X : \inp{x}{\alpha^\vee} \geq 0
\; \forall \alpha \in F_0 \} , \\
&X^- := \{ x \in X : \inp{x}{\alpha^\vee} \leq 0
\; \forall \alpha \in F_0 \} = -X^+ .
\end{align*}
It is easily seen that the center of $W$ is the lattice
\[
Z(W) = X^+ \cap X^- .
\]
We say that $\mc R$ is semisimple if $Z (W) = 0$ or equivalently if $R_0$ spans $\mf a^*$.
Thus a root datum is semisimple if and only if the corresponding reductive algebraic group is so.

The length function $\ell$ of the Coxeter system $(W^\af ,S^\af )$ extends naturally to 
$W$, such that \cite[(1.3)]{Opd-Tr}
\begin{equation}\label{eq:ellwtx}
\ell (w t_x) = \ell (w) + \sum\nolimits_{\alpha \in R_0^+} \inp{x}{\alpha^\vee} 
\qquad w \in W_0 , x \in X^+ .
\end{equation}
The elements of length zero form a subgroup $\Omega \subset W$, and 
$W = W^\af \rtimes \Omega$.
With $\mc R$ we also associate some other root systems.
There is the non-reduced root system
\[
R_{nr} := R_0 \cup \{ 2 \alpha : \alpha^\vee \in 2 Y \} .
\]
Obviously we put $(2 \alpha )^\vee = \alpha^\vee / 2$. Let $R_1$
be the reduced root system of long roots in $R_{nr}$:
\[
R_1 := \{ \alpha \in R_{nr} : \alpha^\vee \not\in 2 Y \} .
\]
We denote the collection of positive roots in $R_0$ by $R_0^+$, and similarly
for other root systems.

\section{Affine Hecke algebras}
\label{sec:defAHA}

There are three equivalent ways to introduce a complex parameter function for $\mc R$.
\begin{enumerate}
\item[(1)] A map $q : S^\af \to \mh C^\times$ such that $q(s) = q(s')$ if $s$ and $s'$ 
are conjugate in $W$.
\item[(2)] A function $q : W \to \C^\times$ such that
\begin{equation}\label{eq:2.18}
\begin{array}{lll@{\quad}l}
q (\omega ) & = & 1 & \text{if } \ell (\omega ) = 0 , \\
q (w v) & = & q (w) q(v) & \text{if } w,v \in W \quad
\text{and} \quad \ell (wv) = \ell (w) + \ell (v) .
\end{array}
\end{equation}
\item[(3)] A $W_0$-invariant map $q : R_{nr}^\vee \to \C^\times, \alpha^\vee \mapsto q_{\alpha^\vee}$.
\end{enumerate}
One goes from (2) to (1) by restriction, while the relation between (2) and (3) is given by
\begin{equation}\label{eq:parameterEquivalence}
\begin{array}{lll}
q_{\alpha^\vee} = q(s_\alpha) = q (t_\alpha s_\alpha) & \text{if} & \alpha \in R_0 \cap R_1, \\
q_{\alpha^\vee} = q(t_\alpha s_\alpha) & \text{if} & \alpha \in R_0 \setminus R_1, \\
q_{\alpha^\vee / 2} = q(s_\alpha) q(t_\alpha s_\alpha)^{-1} & \text{if} & \alpha \in R_0 \setminus R_1. 
\end{array}
\end{equation}
We speak of equal parameters if $q(s) = q(s') \; \forall s,s' \in S^\af$ and of positive parameters if 
$q(s) \in \R_{>0} \; \forall s \in S^\af$. 

We fix a square root $q^{1/2} : S^\af \to \mh C^\times$.
The affine Hecke algebra $\mc H = \mc H (\mc R ,q)$ is the unique associative
complex algebra with basis $\{ N_w : w \in W \}$ and multiplication rules
\begin{equation}\label{eq:multrules}
\begin{array}{lll}
N_w \, N_v = N_{w v} & \mr{if} & \ell (w v) = \ell (w) + \ell (v) \,, \\
\big( N_s - q(s)^{1/2} \big) \big( N_s + q(s)^{-1/2} \big) = 0 & \mr{if} & s \in S^\af .
\end{array}
\end{equation}
In the literature one also finds this algebra defined in terms of the
elements $q(s)^{1/2} N_s$, in which case the multiplication can be described without
square roots. This explains why $q^{1/2}$ does not appear in the notation $\mc H (\mc R ,q)$.

Notice that $N_w \mapsto N_{w^{-1}}$ extends to a $\C$-linear anti-automorphism of $\mc H$,
so $\mc H$ is isomorphic to its opposite algebra.
The span of the $N_w$ with $w \in W_0$ is a finite dimensional Iwahori--Hecke algebra,
which we denote by $\mc H (W_0,q)$.

Now we describe the Bernstein presentation of $\mc H$. For $x \in X^+$ we put
$\theta_x := N_{t_x}$. The corresponding semigroup morphism 
$X^+ \to \mc H (\mc R ,q)^\times$ extends to a group homomorphism
\[
X \to \mc H (\mc R ,q)^\times : x \mapsto \theta_x .
\]

\begin{thm}\label{thm:1.1}
\textup{(Bernstein presentation)}
\enuma{
\item The sets $\{ N_w \theta_x : w \in W_0 , x \in X \}$ and
$\{ \theta_x N_w : w \in W_0 , x \in X \}$ are bases of $\mc H$.
\item The subalgebra $\mc A := \mr{span} \{ \theta_x : x \in X \}$
is isomorphic to $\mh C [X]$.
\item The center of $Z(\mc H (\mc R ,q))$ of $\mc H (\mc R ,q)$ is 
$\mc A^{W_0}$, where we define the action of $W_0$ on $\mc A$ by 
$w (\theta_x ) = \theta_{wx}$.
\item For $f \in \mc A$ and $\alpha \in F_0 \cap R_1$
\[
f N_{s_\alpha} - N_{s_\alpha} s_\alpha (f) = 
\big( q (s_\alpha )^{1/2} - q(s_\alpha )^{-1/2} \big) (f - s_\alpha (f)) 
(\theta_0 - \theta_{-\alpha} )^{-1} ,
\]
while for $\alpha \in F_0 \setminus R_1$: 
\[
f N_{s_\alpha} - N_{s_\alpha} s_\alpha (f) = 
\big( q (s_\alpha )^{1/2} - q(s_\alpha )^{-1/2} + ( q_{\alpha^\vee}^{1/2} - 
q_{\alpha^\vee}^{-1/2}) \theta_{-\alpha} \big) {\ds \frac{f - s_\alpha (f)}
{\theta_0 - \theta_{-2\alpha } }} .
\]
}
\end{thm}
\emph{Proof.}
These results are due to Bernstein, see \cite[\S 3]{Lus-Gr}. $\qquad \Box$
\\[3mm]

The following lemma was claimed in the proof of \cite[Lemma 3.1]{Opd-Tr}.
\begin{lem}\label{lem:1.2}
For $x \in X^+$
\begin{equation}\label{eq:spanthetax}
\mr{span} \{ N_u \theta_x N_v : u,v \in W_0 \} = \mr{span} \{ N_w : w \in W_0 t_x W_0 \} .
\end{equation}
Let $W_x$ be the stabilizer of $x$ in $W_0$ and let $W^x$ be a set of representatives 
for $W_0 / W_x$. Then the elements $N_u \theta_x N_v$ with $u \in W^x$ and $v \in W_0$ 
form a basis of \eqref{eq:spanthetax}. 
\end{lem}
\emph{Proof.}
By \eqref{eq:ellwtx} $\ell (u t_x) = \ell (u) + \ell (t_x)$, so $N_u \theta_x = N_{u t_x}$. Recall 
the Bruhat ordering on a Coxeter group, for example from \cite[Sections 5.9 and 5.10]{Hum}.
With induction to $\ell (v)$ it follows from the multiplication rules \eqref{eq:multrules} that
\[
N_w N_v - N_{wv} \in \text{span} \{ N_{w \tilde{v}} : \tilde{v} < v \text{ in the Bruhat ordering} \} .
\]
Hence the sets $\{ N_u \theta_x N_v : v \in W_0 \}$ and $\{ N_w : w \in u t_x W_0 \}$ have the same
span. They have the same cardinality and by definition the latter set is linearly independent, so
the former is linearly independent as well. Clearly $W_0 t_x W_0 = \sqcup_{u \in W^x} u t_x W_0$, so
\begin{multline*}
\text{span} \{ N_w : w \in W_0 t_x W_0 \} = 
\bigoplus\nolimits_{u \in W^x} \text{span} \{ N_w : w \in u t_x W_0 \} \\
= \bigoplus\nolimits_{u \in W^x} \text{span} \{ N_u \theta_x N_v : v \in W_0 \} =
 \text{span} \{ N_u \theta_x N_v : u \in W^x , v \in W_0 \} .
\end{multline*}
The number of generators on the second line side equals the dimension of the first line, 
so they form a basis. $\qquad \Box$
\\[2mm]

Let $T$ be the complex algebraic torus 
\[
T = \mr{Hom}_{\mh Z} (X, \mh C^\times ) \cong Y \otimes_\Z \C^\times ,
\]
so $\mc A \cong \mc O (T)$ and $Z (\mc H ) = \mc A^{W_0} \cong \mc O (T / W_0 )$. From Theorem 
\ref{thm:1.1} we see that $\mc H$ is of finite rank over its center. Let $\mf t = \mr{Lie}(T)$ 
and $\mf t^*$ be the complexifications of $\mf a$ and $\mf a^*$. The direct sum
$\mf t = \mf a \oplus i \mf a$ corresponds to the polar decomposition 
\[
T = T_{rs} \times T_{un} = \mr{Hom}_\Z (X, \R_{>0}) \times \mr{Hom}_\Z (X, S^1) ,
\]
where $T_{rs}$ is the real split (or positive) part of $T$ and $T_{un}$ the unitary part. 
The exponential map $\exp : \mf t \to T$ is bijective on the real parts, and we denote its inverse 
by $\log : T_{rs} \to \mf a$.

An important role in the harmonic analysis of $\mc H (\mc R ,q)$ is played by the Macdonald 
$c$-functions $c_\alpha \in \C (T)$ (cf. \cite[3.8]{Lus-Gr} and \cite[Section 1.7]{Opd-Tr}), defined~as
\begin{equation}\label{eq:calpha}
c_\alpha = \frac{\theta_\alpha + q(s_\alpha )^{-1/2} q_{\alpha^\vee}^{1/2}}{\theta_\alpha + 1} \,
\frac{\theta_\alpha - q(s_\alpha )^{-1/2} q_{\alpha^\vee}^{-1/2}}{\theta_\alpha - 1} .
\end{equation}
Notice that $c_\alpha = 1$ if and only if $q(s_\alpha) = q (t_\alpha s_\alpha) = 1$. 
For $\alpha \in R_0 \cap R_1$ we have $q(s_\alpha) = q_{\alpha^\vee}$, so \eqref{eq:calpha} 
simplifies to $c_\alpha = (\theta_\alpha - q(s_\alpha )^{-1}) (\theta_\alpha - 1)^{-1}$.
With these $c$-functions we can rephrase Theorem \ref{thm:1.1}.d as
\[
f N_{s_\alpha} - N_{s_\alpha} s_\alpha (f) = 
q(s_\alpha)^{-1/2} \big( f - s_\alpha (f) \big)  (q(s_\alpha) c_\alpha - 1) . 
\]
An automorphism of the Dynkin diagram of the based root system $(R_0,F_0 )$ is a 
bijection $\gamma : F_0 \to F_0$ such that
\begin{equation}
\inp{\gamma (\alpha )}{\gamma (\beta )^\vee} = \inp{\alpha}{\beta^\vee} 
\qquad \forall \alpha ,\beta \in F_0 \,.
\end{equation}
Such a $\gamma$ naturally induces automorphisms of $R_0, R_0^\vee, W_0$ and 
$W^\af$. It is easy to classify all diagram automorphisms of $(R_0 ,F_0)$: they permute the
irreducible components of $R_0$ of a given type, and the diagram automorphisms of a
connected Dynkin diagram can be seen immediately.

We will assume that the action of $\gamma$ on $W^{\af}$ has been extended in
some way to $W$. For example, this is the case if $\gamma$ belongs to the Weyl 
group of some larger root system contained in $X$. We regard two diagram automorphisms 
as the same if and only if their actions on $W$ are equal.

Let $\Gamma$ be a finite group of diagram automorphisms of $(R_0,F_0)$ and assume that
$q_{\alpha^\vee} = q_{\gamma (\alpha^\vee)}$ for all $\gamma \in \Gamma, \alpha \in R_{nr}$.
Then $\Gamma$ acts on $\mc H$ by algebra automorphisms $\psi_\gamma$ that satisfy
\begin{equation}\label{eq:1.2}
\begin{array}{lll@{\quad}l}
\psi_\gamma (N_w) & = & N_{\gamma (w)} & w \in W , \\
\psi_\gamma (\theta_x) & = & \theta_{\gamma (x)} & x \in X .
\end{array}
\end{equation}
Hence one can form the crossed product algebra $\Gamma \ltimes \mc H = \mc H \rtimes \Gamma$, 
whose natural basis is indexed by the group $(X \rtimes W_0) \rtimes \Gamma = 
X \rtimes (W_0 \rtimes \Gamma)$. It follows easily from \eqref{eq:1.2} and Theorem \ref{thm:1.1}.c 
that $Z(\mc H \rtimes \Gamma) = {\mc A}^{W_0 \rtimes \Gamma}$. We say that the central 
character of an (irreducible) $\mc H \rtimes \Gamma$-representation is positive if it lies in 
$T^{rs} / (W_0 \rtimes \Gamma)$.

We always assume that we have an $\Gamma \ltimes W_0$-invariant inner product on $\mf a$.
The length function of $W^\af$ also extends to $X \rtimes W_0 \rtimes \Gamma$, and the 
subgroup of elements of length zero becomes
\[
\{ w \in W \rtimes \Gamma : \ell (w) = 0 \} = \Gamma \ltimes \{ w \in W : \ell (w) = 0 \}
= \Gamma \ltimes \Omega .
\]
More generally one can consider a finite group $\Gamma'$ that acts on $\mc R$ by diagram
automorphisms. Then the center of $\mc H \rtimes \Gamma'$ can be larger than 
$\mc A^{W_0 \rtimes \Gamma'}$, but apart from that the structure is the same.

Another variation arises when $\Gamma \to \mr{Aut}(\mc H)$ is not a group homomorphism,
but a homomorphism twisted by a 2-cocycle $\kappa : \Gamma \times \Gamma \to \C^\times$.
Instead of $\mc H \rtimes \Gamma$ one can construct the algebra $\mc H \otimes \C [\Gamma,\kappa]$,
whose multiplication is defined by
\[
\begin{array}{lll}
N_\gamma N_{\gamma'} & = & \kappa (\gamma,\gamma') N_{\gamma \gamma'} , \\
N_\gamma h N_\gamma^{-1} & = & \gamma (h),
\end{array}
\qquad \gamma, \gamma' \in \Gamma, h \in \mc H .
\]
By \cite[Section 7]{Mor} such algebras can appear in relevant examples, although it is no explicit 
nontrivial are known. Let $\Gamma^*$ be the Schur multiplier of $\Gamma$, also known as 
representation group \cite{CuRe1}. 
It is a central extension of $\Gamma$ that classifies projective 
$\Gamma$-representations, and its group algebra $\C [\Gamma^*]$ contains $\C [\Gamma,\kappa]$
as a direct summand. The algebra $\mc H \rtimes \Gamma^*$ is well-defined and contains
$\mc H \otimes \C [\Gamma,\kappa]$ as a direct summand. Thus we can reduce the study of
affine Hecke algebras with twisted group actions the case of honest group actions.

\section{Graded Hecke algebras}

Graded Hecke algebras are also known as degenerate (affine) Hecke algebras.
They were introduced by Lusztig in \cite{Lus-Gr}. We call 
\begin{equation}
\tilde{\mc R} = (\mf a^* ,R_0, \mf a, R_0^\vee, F_0 )
\end{equation}
a degenerate root datum. We pick complex numbers $k_\alpha$ for $\alpha \in F_0$,
such that $k_\alpha = k_\beta$ if $\alpha$ and $\beta$ are in the same $W_0$-orbit.
The graded Hecke algebra associated to these data is the complex vector space
\[
\mh H = \mh H (\tilde{\mc R},k) = S( \mf t^*) \otimes \C [W_0] ,
\]
with multiplication defined by the following rules:
\begin{itemize}
\item $\mh C[W_0]$ and $S (\mf t^* )$ are canonically embedded as subalgebras;
\item for $x \in \mf t^*$ and $s_\alpha \in S$ we have the cross relation
\begin{equation}\label{eq:1.1}
x \cdot s_\alpha - s_\alpha \cdot s_\alpha (x) = 
k_\alpha \inp{x}{\alpha^\vee} .
\end{equation}
\end{itemize}
Multiplication with any $\ep \in \mh C^\times$ defines a bijection $m_\ep : \mf t^* \to \mf t^*$,
which clearly extends to an algebra automorphism of $S(\mf t^* )$. From the cross relation
\eqref{eq:1.1} we see that it extends even further, to an algebra isomorphism
\begin{equation}\label{eq:1.3}
m_\ep : \mh H (\tilde{\mc R},zk) \to \mh H (\tilde{\mc R}, k)
\end{equation}
which is the identity on $\mh C[W_0]$. 

The $c$-functions are considerably easier than for affine Hecke algebras:
\begin{equation}
\tilde c_\alpha = (\alpha + k_\alpha ) \alpha^{-1} = 1 + k_\alpha \alpha^{-1} \in \C (\mf t) .
\end{equation}
They can be used to rewrite \eqref{eq:1.1} as
\[
x s_\alpha - s_\alpha \, s_\alpha (x) = \big( x - s_\alpha (x) \big) (\tilde c_\alpha - 1) 
\qquad \alpha \in F_0, x \in \mf t^* .
\]
Let $\Gamma$ be a group of diagram automorphisms of $\tilde{\mc R}$ and assume that 
$k_{\gamma (\alpha)} = k_\alpha$ for all $\alpha \in R_0 , \gamma \in \Gamma$.
Then $\Gamma$ acts on $\mh H$ by the algebra automorphisms
\begin{equation}
\begin{split}
& \psi_\gamma : \mh H \to \mh H \,, \\
& \psi_\gamma (x s_\alpha ) = \gamma (x) s_{\gamma (\alpha )} 
  \qquad x \in \mf t^* , \alpha \in \Pi \,.
\end{split}
\end{equation}
By \cite[Proposition 5.1.a]{SolHomGHA} the center of the resulting crossed product algebra is
\begin{equation}\label{eq:1.4}
Z (\mh H \rtimes \Gamma) = S(\mf t^*)^{W_0 \rtimes \Gamma} = 
\mc O (\mf t / (W_0 \rtimes \Gamma)) .
\end{equation}
We say that the central character of an $\mh H \rtimes \Gamma$-representation is real if it lies in 
$\mf a / (W_0 \rtimes \Gamma)$.

\section{Parabolic subalgebras}
\label{sec:parabolic}

For a set of simple roots $P \subset F_0$ we introduce the notations
\begin{equation}
\begin{array}{l@{\qquad}l}
R_P = \mh Q P \cap R_0 & R_P^\vee = \mh Q R_P^\vee \cap R_0^\vee , \\
\mf a_P = \R P^\vee & \mf a^P = (\mf a^*_P )^\perp ,\\
\mf a^*_P = \R P & \mf a^{P*} = (\mf a_P )^\perp  ,\\
\mf t_P = \C P^\vee & \mf t^P = (\mf t^*_P )^\perp ,\\
\mf t^*_P = \C P & \mf t^{P*} = (\mf t_P )^\perp  ,\\
X_P = X \big/ \big( X \cap (P^\vee )^\perp \big) &
X^P = X / (X \cap \mh Q P ) , \\
Y_P = Y \cap \mh Q P^\vee & Y^P = Y \cap P^\perp , \\
T_P = \mr{Hom}_{\mh Z} (X_P, \mh C^\times ) &
T^P = \mr{Hom}_{\mh Z} (X^P, \mh C^\times ) , \\
\mc R_P = ( X_P ,R_P ,Y_P ,R_P^\vee ,P) & \mc R^P = (X,R_P ,Y,R_P^\vee ,P) , \\
\tilde{\mc R}_P = ( \mf a_P^* ,R_P ,\mf a_P ,R_P^\vee ,P) & 
\tilde{\mc R}^P = (\mf a^*,R_P ,\mf a,R_P^\vee ,P) .
\end{array}
\end{equation}
We denote the image of $x \in X$ in $X_P$ by $x_P$. Although $T_{rs} = T_{P,rs} 
\times T_{rs}^P$, the product $T_{un} = T_{P,un} T^P_{un}$ is not direct, 
because the intersection 
\[
K_P := T_{P,un} \cap T_{un}^P = T_P \cap T^P
\]
can have more than one element (but only finitely many).

We define parameter functions $q_P$ and $q^P$ on the root
data $\mc R_P$ and $\mc R^P$, as follows. Restrict $q$ to a function on
$(R_P )_{nr}^\vee$ and use \eqref{eq:parameterEquivalence} to extend it to
$W (\mc R_P )$ and $W (\mc R^P )$. Similarly the restriction of $k$ to $P$ is
a parameter function for the degenerate root data $\tilde{\mc R}_P$ and
$\tilde{\mc R}^P$, and we denote it by $k_P$ or $k^P$. Now we can define
the parabolic subalgebras
\[
\begin{array}{l@{\qquad}l}
\mc H_P = \mc H (\mc R_P ,q_P ) & \mc H^P = \mc H (\mc R^P ,q^P ) , \\
\mh H_P = \mh H (\tilde{\mc R}_P ,k_P ) & \mh H^P = \mh H (\tilde{\mc R}^P ,k^P ) .
\end{array}
\]
We notice that $\mh H^P = S (\mf t^{P*} ) \otimes \mh H_P$, a tensor product of 
algebras. Despite our terminology $\mc H^P$ and $\mc H_P$ are not subalgebras of $\mc H$, 
but they are close. Namely, $\mc H (\mc R^P ,q^P )$ is isomorphic to the subalgebra of 
$\mc H (\mc R ,q)$ generated by $\mc A$ and $\mc H (W (R_P) ,q_P)$. 
We denote the image of $x \in X$ in $X_P$ by $x_P$ and we let $\mc A_P \subset \mc H_P$ 
be the commutative subalgebra spanned by $\{ \theta_{x_P} : x_P \in X_P \}$. 
There is natural surjective quotient map
\begin{equation}\label{eq:quotientP}
\mc H^P \to \mc H_P : \theta_x N_w \mapsto \theta_{x_P} N_w .
\end{equation}
Suppose that $\gamma \in \Gamma \ltimes W_0$ satisfies $\gamma (P) = Q \subseteq F_0$.
Then there are algebra isomorphisms
\begin{equation}\label{eq:psigamma}
\begin{array}{llcl}
\psi_\gamma : \mc H_P \to \mc H_Q , &
\theta_{x_P} N_w & \mapsto & \theta_{\gamma (x_P)} N_{\gamma w \gamma^{-1}} , \\
\psi_\gamma : \mc H^P \to \mc H^Q , &
\theta_x N_w & \mapsto & \theta_{\gamma x} N_{\gamma w \gamma^{-1}} , \\
\psi_\gamma : \mh H_P \to \mh H_Q  , &
f_P w & \mapsto & (f_P \circ \gamma^{-1}) w , \\
\psi_\gamma : \mh H_P \to \mh H_Q  , &
f w & \mapsto & (f \circ \gamma^{-1}) w ,
\end{array}
\end{equation}
where $f_P \in \mc O (\mf t_P)$ and $f \in \mc O (\mf t)$. Sometimes we will abbreviate 
$W \rtimes \Gamma$ to $W'$ and $W_0 \rtimes \Gamma$ to $\WG$. For example the group
\begin{equation}\label{eq:GammaP}
W'_P := \{ \gamma \in \Gamma \ltimes W_0 : \gamma (P) = P \} 
\end{equation}
acts on the algebras $\mc H_P$ and $\mc H^P$. Although $W'_{F_0} = \Gamma$, for 
proper subsets $P \subsetneq F_0$ the group $W'_P$ need not be contained in $\Gamma$.
In other words, in general $W'_P$ strictly contains the group
\[
\Gamma_P := \{ \gamma \in \Gamma : \gamma (P) = P \} = W'_P \cap \Gamma .
\]
To avoid confusion we do not use the notation $W_P$. Instead the parabolic subgroup 
of $W_0$ generated by $\{ s_\alpha : \alpha \in P \}$ will be denoted $W (R_P)$. Suppose that 
$\gamma \in \WG$ stabilizes either the root system $R_P$, the lattice $\Z P$ or the vector space
$\Q P \subset \mf a^*$. Then $\gamma (P)$ is a basis of $R_P$, so
$\gamma (P) = w (P)$ and $w^{-1} \gamma \in W'_P$ for a unique $w \in W(R_P)$. Therefore
\begin{equation}\label{eq:WZP}
W'_{\Z P} := \{ \gamma \in \WG : \gamma (\Z P) = \Z P \} \text{ equals } W(R_P) \rtimes W'_P .
\end{equation}
For all $x \in X$ and $\alpha \in P$ we have
\[
x - s_\alpha (x) = \inp{x}{\alpha^\vee} \alpha \in \Z P,
\]
so $t (s_\alpha (x)) = t(x)$ for all $t \in T^P$. Hence $t (w(x)) = t (x)$ for all $w \in W (R_P)$, 
and we can define an algebra automorphism
\begin{equation}
\phi_t : \mc H^P \to \mc H^P, \quad \phi_t (\theta_x N_w) = t (x) \theta_x N_w  \qquad t \in T^P .
\end{equation}
In particular, for $t \in K_P$ this descends to an algebra automorphism 
\begin{equation}\label{eq:twistKP}
\psi_t : \mc H_P \to \mc H_P , \quad \theta_{x_P} N_w \mapsto t(x_P) \theta_{x_P} N_w  
\qquad t \in K_P .
\end{equation}
We can regard any representation $(\sigma ,V_\sigma)$ of $\mc H (\mc R_P ,q_P )$ as a 
representation of $\mc H (\mc R^P ,q^P)$ via the quotient map \eqref{eq:quotientP}. 
Thus we can construct the $\mc H$-representation
\[
\pi (P,\sigma ,t) := \mr{Ind}_{\mc H (\mc R^P ,q^P )}^{\mc H (\mc R ,q)} (\sigma \circ \phi_t ) .
\]
Representations of this form are said to be parabolically induced.
Similarly, for any $\mh H_P$-representation $(\rho, V_\rho)$ and any $\lambda \in \mf t^P$
there is an $\mh H^P$-representation $(\rho_\lambda, V_\rho \otimes \C_\lambda)$.
The corresponding parabolically induced representation is
\[
\pi (P,\rho,\lambda) := \mr{Ind}_{\mh H^P}^{\mh H} (\rho_\lambda) = 
\mr{Ind}_{\mh H^P}^{\mh H} (V_\rho \otimes \C_\lambda) .
\]

\section{Analytic localization}
\label{sec:localiz}

A common technique in the study of Hecke algebras is localization at one or more characters of the center. 
There are several ways to implement this. Lusztig \cite{Lus-Gr} takes a maximal ideal $I$ of $Z (\mc H)$
and completes $\mc H$ with respect to the powers of this ideal. This has the effect of considering only
those $\mc H$-representations which admit the central character corresponding to $I$.

For reasons that will become clear only in Chapter \ref{chapter:scaling}, we prefer to localize with
analytic functions on subvarieties of $T / \WG$. Let $U \subset T$ be a nonempty $\WG$-invariant subset 
and let $C^{an}(U)$ (respectively $C^{me}(U)$) be the algebra of holomorphic (respectively 
meromorphic) functions on $U$. There is a  natural embedding
\[
Z(\mc H \rtimes \Gamma ) = \mc A^{\WG} \cong \mc O (T)^{\WG} \to C^{an}(U)^{\WG}
\]
and isomorphisms of topological algebras
\[
C^{an}(U)^{\WG} \otimes_{\mc A^{\WG}} \mc A \cong C^{an}(U) ,\quad
C^{me}(U)^{\WG} \otimes_{\mc A^{\WG}} \mc A \cong C^{me}(U) .
\]
Thus we can construct the algebras
\begin{equation}\label{eq:defHan}
\begin{array}{ccccc}
\mc H^{an}(U) \rtimes \Gamma & := & 
C^{an}(U)^{\WG} \otimes_{Z (\mc H \rtimes \Gamma)} \mc H \rtimes \Gamma & \cong & 
C^{an}(U) \otimes_\C \C [\WG] , \\
\mc H^{me}(U) \rtimes \Gamma & := & 
C^{me}(U)^{\WG} \otimes_{Z (\mc H \rtimes \Gamma)} \mc H \rtimes \Gamma
& \cong & C^{me}(U) \otimes_\C \C [\WG].
\end{array}
\end{equation}
The isomorphisms with the right hand side are in the category of topological vector spaces, the
algebra structure on the left hand side is determined by 
\begin{equation}\label{eq:multfh}
(f_1 \otimes h_1)(f_2 \otimes h_2) = f_1 f_2 \otimes h_1 h_2 \qquad
f_i \in C^{me} (U)^{\WG}, h_i \in \mc H \rtimes \Gamma . 
\end{equation}
By \cite[Proposition 4.5]{Opd-Sp} $Z (\mc H^{an}(U) \rtimes \Gamma) \cong C^{an}(U)^{\WG}$,
and similarly with meromorphic functions.

For any $T' \subset T$, let $\mr{Mod}_{f,T'} (\mc H \rtimes \Gamma)$ be the category of finite 
dimensional $\mc H \rtimes \Gamma$-modules all whose $\mc A$-weights lie in $T'$. 
By \cite[Proposition 4.3]{Opd-Sp} $\mr{Mod}_f (\mc H^{an}(U) \rtimes \Gamma)$ is naturally equivalent to 
$\mr{Mod}_{f,U} (\mc H \rtimes \Gamma)$. (On the other hand, $\mc H^{me}(U) \rtimes \Gamma$ does not 
have any nonzero finite dimensional representations over $\C$.) 

Of course graded Hecke algebras can localized in exactly the same way, and the resulting algebras 
have analogous properties. By \eqref{eq:1.4} the center of the 
algebra $\mh H \rtimes \Gamma$ is isomorphic to $\mc O (\mf t / \WG) = \mc O (\mf t)^{\WG}$.  
For nonempty open $\WG$-invariant subsets $V$ of $\mf t$ we get the algebras
\begin{equation}\label{eq:defgHan}
\begin{array}{ccccc}
\mh H^{an}(V) \rtimes \Gamma & := & 
C^{an}(V)^{\WG} \otimes_{Z (\mh H \rtimes \Gamma)} \mh H \rtimes \Gamma & \cong & 
C^{an}(V) \otimes_\C \C [\WG] , \\
\mh H^{me}(V) \rtimes \Gamma & := & 
C^{me}(V)^{\WG} \otimes_{Z (\mh H \rtimes \Gamma)} \mh H \rtimes \Gamma
& \cong & C^{me}(V) \otimes_\C \C [\WG].
\end{array}
\end{equation}
Let $\C (T/\WG)$ be the quotient field of $Z (\mc H \rtimes \Gamma) \cong \mc O (T /\WG)$ and consider
\[
\C (T/\WG) \otimes_{Z (\mc H \rtimes \Gamma)} \mc H \rtimes \Gamma .
\]
As a vector space this is $\C (T) \otimes_{\mc A} \mc H \rtimes \Gamma \cong \C (T) \otimes_\C \C [\WG]$, 
while its multiplication is given by \eqref{eq:multfh}.
Similarly, we let $\C (\mf t /\WG)$ be the quotient field of $\mc O (\mf t /\WG)$ and we construct the algebra
\[
\C (\mf t /\WG) \otimes_{Z (\mh H \rtimes \Gamma)} \mh H \rtimes \Gamma .
\]
Its underlying vector space is 
\[
\C (\mf t) \otimes_{\mc O (\mf t)} \mh H \rtimes \Gamma \cong \C (\mf t) \otimes_\C \C [\WG] ,
\]
and its multiplication is the obvious analogue of \eqref{eq:multfh}. Given a simple root $\alpha \in F_0$ 
we define elements $\imath^0_{s_\alpha} \in \C (T/W_0) \otimes_{Z (\mc H)} \mc H$ and 
$\tilde \imath_{s_\alpha} \in \C (\mf t /W_0) \otimes_{Z (\mh H)} \mh H$ by
\begin{equation}\label{eq:imatho}
\begin{array}{rcl}
q(s_\alpha) c_\alpha (1 + \imath^o_{s_\alpha}) & = & 1 + q(s_\alpha)^{1/2} N_{s_\alpha} , \\
\tilde c_\alpha (1 + \tilde \imath_{s_\alpha}) & = & 1 + s_\alpha . 
\end{array}
\end{equation}

\begin{prop}\label{prop:3.4}
The elements $\imath^0_{s_\alpha}$ and $\tilde \imath_{s_\alpha}$ have the following properties:
\enuma{
\item The map $s_\alpha \mapsto \imath^0_{s_\alpha}$ (respectively 
$s_\alpha \mapsto \tilde \imath_{s_\alpha}$) extends to a group homomorphism from $\WG$ to the 
multiplicative group of $\C (T/\WG) \otimes_{Z (\mc H \rtimes \Gamma)} \mc H \rtimes \Gamma$ (respectively 
$\C (\mf t /\WG) \otimes_{Z (\mh H \rtimes \Gamma)} \mh H \rtimes \Gamma$).
\item For $w \in \WG$ and $f \in \C (T ) \cong \C (T/\WG) \otimes_{\mc O (T/\WG)} \mc A$ 
(respectively $\tilde f \in \C (\mf t)$) we have $\imath^0_w f \imath^0_{w^{-1}} = w (f)$ 
(respectively $\tilde \imath_w \tilde f \tilde \imath_{w^{-1}} = w (\tilde f)$).
\item The maps
\[
\begin{array}{rcl@{\; : \;}ccl}
\C (T ) \rtimes \WG & \to & \C (T/\WG) \otimes_{Z (\mc H \rtimes \Gamma)} \mc H \rtimes \Gamma & 
f w & \mapsto & f \imath^0_w , \\
\C (\mf t) \rtimes \WG & \to & \C (\mf t /\WG) \otimes_{Z (\mh H \rtimes \Gamma )} \mh H \rtimes \Gamma &
\tilde f w & \mapsto & \tilde f \tilde \imath_w 
\end{array}
\]
are algebra isomorphisms.
\item Let $P \subset F_0$ and $\gamma \in \WG$ be such that $\gamma (P) \subset F_0$. 
The automorphisms $\psi_\gamma$ from \eqref{eq:psigamma} satisfy
\[
\begin{array}{lll@{\qquad}l}
\psi_\gamma (h) & = & \imath^0_\gamma h \imath^0_{\gamma^{-1}} & h \in \mc H^P \text{ or } h \in \mc H_P, \\
\psi_\gamma (\tilde h) & = & \tilde{\imath}_\gamma \tilde h \tilde{\imath}_{\gamma^{-1}} & \tilde h \in \mh H^P .
\end{array}
\]
}
\end{prop}
\emph{Proof.}
(a), (b) and (c) with $W_0$ instead of $\WG$ can be found in \cite[Section 5]{Lus-Gr}. Notice that Lusztig calls 
these elements $\tau_w$ and $\tilde \tau_w$, while we follow the notation of \cite[Section 4]{Opd-Sp}.
We extend this to $\WG$ by defining, for $\gamma \in \Gamma$ and $w \in W_0$:
\[
\imath^0_{\gamma w} := \gamma \imath^0_w \quad \text{and} 
\quad \tilde \imath_{\gamma w} = \gamma \tilde \imath_w. 
\]
For (d) see \cite[Section 8]{Lus-Gr} or \cite[Lemma 3.2]{SolGHA}. $\qquad \Box$ 
\\[2mm]

We remark that by construction all the $\imath^0_w$ lie in the subalgebra 
$\C (T / T^{F_0}) \mc H (W_0,q) \rtimes \Gamma$ and that the $\tilde \imath_w$ lie in the subalgebra 
$\C (\mf t / \mf t^{F_0}) \C [\WG]$.
As was noticed in \cite[Theorem 4.6]{Opd-Sp}, Proposition \ref{prop:3.4} can easily be generalized:

\begin{cor}\label{cor:1.3}
Proposition \ref{prop:3.4} remains valid under any of the following replacements:
\begin{itemize}
\item $\C (T)$ by $C^{me}(U)$ or, if all the functions $c_\alpha$ are invertible on $U$, by $C^{an}(U)$;
\item $\C (\mf t)$ by $C^{me}(V)$ or, if all the functions $\tilde c_\alpha$ are invertible on $V$, by $C^{an}(V)$. 
\end{itemize}
\end{cor}
In particular
\begin{equation}\label{eq:isoCme}
C^{me}(U) \rtimes \WG \to C^{me}(U)^{\WG} \otimes_{\mc A^{\WG}} \mc H \rtimes \Gamma : f w \to f \imath_w^o
\end{equation}
is an isomorphism of topological algebras.

\section{The relation with reductive $p$-adic groups}
\label{sec:padic}

Here we discuss how affine Hecke algebras arise in the representation theory of reductive
$p$-adic groups. This section is rather sketchy, it mainly serves to provide some motivation and
to prepare for our treatment of the Aubert--Baum--Plymen conjecture in Section \ref{sec:ABP}.
The main sources for this section are \cite{BeDe,BeRu,Roc2,Hei,IwMa,Mor}.

Let $\mathbb F$ be a nonarchimedean local field with a finite residue field.
Let $\mathbf G$ be a connected reductive algebraic group defined over $\mathbb F$ and let 
$G = \mathbf G (\mathbb F)$ be its group of $\mathbb F$-rational points. One briefly calls $G$
a reductive $p$-adic group, even though the characteristic of $\mathbb F$ is allowed to be positive.

An important theme, especially in relation with the arithmetic Langlands program, is the study
of the category Mod$(G)$ of smooth $G$-representations on complex vector spaces. A first step to
simplify this problem is the Bernstein decomposition, which we recall now.
Let $P$ be a parabolic subgroup of $G$ and $M$ a Levi subgroup of $P$. Although $G$ and $M$ are
unimodular, the modular function $\delta_P$ of $P$ is in general not constant. Let $\sigma$ be an
irreducible supercuspidal representation of $M$. In this situation we call $(M,\sigma)$ a
cuspidal pair, and with parabolic induction we construct the $G$-representation
\[
I_P^G (\sigma) := \mr{Ind}_P^G (\delta_P^{1/2} \otimes \sigma) .
\]
This means that first we inflate $\sigma$ to $P$ and then we apply the normalized smooth 
induction functor. The normalization refers to the twist by $\delta_P^{1/2}$, which is useful to 
preserve unitarity. Let Irr$(G)$ be the collection of irreducible representations in Mod$(G)$,
modulo equivalence. For every $\pi \in \mr{Irr}(G)$ there is a cuspidal pair $(M,\sigma)$,
uniquely determined up to $G$-conjugacy, such that $\pi$ is a subquotient of $I_P^G (\sigma)$.

Let $G^0$ be the normal subgroup of $G$ generated by all compact subgroups. Recall that a 
character $\chi : G \to \C^\times$ is called unramified if its kernel contains $G^0$. The group
$X_{ur}(G)$ of unramified characters forms a complex algebraic torus whose character lattice is 
naturally isomorphic to the lattice $X^* (G)$ of algebraic characters $G \to \mathbb F^\times$.
We say that two cuspidal pairs $(M,\sigma)$ and $(M',\sigma')$ are inertially equivalent if there 
exist $g \in G$ and $\chi' \in X_{ur}(M')$ such that 
\[
M' = g M g^{-1} \quad \text{and} \quad \sigma' \otimes \chi' \cong \sigma^g ,
\]
where $\sigma^g (m') = \sigma (g^{-1} m' g)$.
With an inertial equivalence class $\mf s = [M,\sigma]_G$ one associates a full subcategory
$\mr{Mod}_{\mf s} (G)$ of Mod$(G)$. It objects are by definition those smooth representations $\pi$ with
the property that for every irreducible subquotient $\rho$ of $\pi$ there is a $(M,\sigma) \in \mf s$
such that $\rho$ is a subrepresentation of $I_P^G (\sigma)$. The collection $\mf B (G)$ of
inertial equivalence classes is countably infinite (unless $G = \{1\}$).

The Bernstein decomposition \cite[Proposition 2.10]{BeDe} states that
\begin{equation}\label{eq:Bdecomp}
\mr{Mod}(G) = \prod\nolimits_{\mf s \in \mf B (G)} \mr{Mod}_{\mf s} (G) ,
\end{equation}
a direct product of categories. The subcategories $\mr{Mod}_{\mf s} (G)$ (or rather their subsets of
irreducible representations) are also called the Bernstein components of the smooth dual of $G$.

The Hecke algebra $\mc H (G)$ is the vector space of locally
constant, compactly supported functions on $G$, endowed with the convolution product. Mod$(G)$
is naturally equivalent to the category Mod$(\mc H (G))$ of essential $\mc H (G)$-modules.
(A module $V$ is called essential if $\mc H (G) V = V$, which is not automatic because $\mc H (G)$
does not have a unit.) Corresponding to \eqref{eq:Bdecomp} there is a decomposition 
\[
\mc H (G) = \bigoplus\nolimits_{\mf s \in \mf B (G)} \mc H (G)_{\mf s}
\]
of the Hecke algebra of $G$ into two-sided ideals. In several cases $\mc H (\mc G)_{\mf s}$ is
known to be Morita-equivalent to an affine Hecke algebra. 

In the classical case \cite{IwMa,Bor} $G$ is split and $\mr{Mod}_{\mf s} (G)$ is the category of 
Iwahori-spherical representations. That is, those smooth $G$-representations $V$ that are generated 
by $V^I$, where $I$ is an Iwahori-subgroup of $G$. Then $\mc H (G)_{\mf s}$ is Morita equivalent to 
the algebra $\mc H (G,I)$ of $I$-biinvariant functions in $\mc H (G)$, and $\mc H (G,I)$ is isomorphic 
to an affine Hecke algebra $\mc H (\mc R,q)$. The root datum $\mc R = (X,R_0,Y,R_0^\vee,F_0)$ is dual 
to the root datum of $(\mathbf G, \mathbf T)$, where $\mathbf T (\mathbb F)$ is a split maximal torus 
of $G = \mathbf G (\mathbb F)$. More explicitly
\begin{itemize}
\item $X$ is the cocharacter lattice of $\mathbf T$;
\item $Y$ is the character lattice of $\mathbf T$;
\item $R_0^\vee$ is the root system of $(\mathbf G ,\mathbf T)$;
\item $R_0$ is the set coroots of $(\mathbf G ,\mathbf T)$;
\item $F_0$ and $F_0^\vee$ are determined by $I$;
\item $q$ is the cardinality of the residue field of $\mathbb F$.
\end{itemize}
For a general inertial equivalence class $\mf s = [M,\sigma]_G$ it is expected that $\mc H (G)_{\mf s}$
is Morita equivalent to an affine Hecke algebra or to a closely related kind of algebra. We discuss some 
ingredients of this conjectural relation between the representation theory of reductive $p$-adic groups 
and that of affine Hecke algebras.

Let $\sigma^0$ be an irreducible subrepresentation of $\sigma \big|_{M^0}$ and define
$\Sigma = \mr{ind}_{M^0}^M (\sigma^0)$. According to \cite[Theorem 23]{BeRu} $I_P^G (\Sigma)$
is a finitely generated projective generator of the category $\mr{Mod}_{\mf s}(G)$. By 
\cite[Lemma 22]{BeRu} $\mr{Mod}_{\mf s}(G) = \mr{Mod}(\mc H (G)_{\mf s})$ is naturally equivalent 
to the category of right $\mr{End}_G (I_P^G (\Sigma))$-modules. So if $\mr{End}_G (I_P^G (\Sigma))$
would be isomorphic to its opposite algebra (which is likely), then it is Morita equivalent to
$\mc H (G)_{\mf s}$.

Let us describe the center of $\mr{End}_G (I_P^G (\Sigma))$. The map 
\[ 
X_{ur}(M) \to \mr{Irr}_{[M,\sigma]_M} (M) : \chi \mapsto \chi \otimes \sigma
\]
is surjective and its fibers are cosets of a finite subgroup Stab$(\sigma) \subset X_{ur}(M)$. Let 
\[
M_\sigma := \bigcap\nolimits_{\chi \in \mr{Stab}(\sigma)} \! \ker (\chi) \; \subset \; M .
\]
Roche \cite[Proposition 5.3]{Roc2} showed that $\mr{End}_M (\Sigma)$ is a free 
$\mc O (X_{ur}(M_\sigma))$-module of rank $m^2$, where $m$ is the multiplicity of $\sigma^0$ 
in $\sigma \big|_{M^0}$. Moreover the center of $\mr{End}_M (\Sigma)$ is isomorphic to 
$\mc O (X_{ur}(M_\sigma))$, and $\mr{End}_M (\Sigma)$ embeds in $\mr{End}_G (I_P^G (\Sigma))$
by functoriality. The group 
\[
N_G (M,\sigma) := \{ g \in G : g M g^{-1} = M, \sigma^g \cong \chi \otimes \sigma
\text{ for some } \chi \in X_{ur}(M) \}
\]
acts on the family of representations $\{ I_P^G (\chi \otimes \Sigma) : \chi \in X_{ur}(M) \}$,
and via this on $X_{ur} (M_\sigma)$. The subgroup $M \subset N_G (M,\sigma)$ acts trivially, so we get
an action of the finite group 
\[
W_\sigma := N_G (M,\sigma) / M .
\]
According to \cite[Th\'eor\`eme 2.13]{BeDe}, the center of $\mr{End}_G (I_P^G (\Sigma))$ is
isomorphic to $\mc O (X_{ur}(M_\sigma))^{W_\sigma} = \mc O (X_{ur}(M_\sigma) / W_\sigma)$. 
By \cite[Lemma 7.3]{Roc2} $\mr{End}_G (I_P^G (\Sigma))$ is a free $\mr{End}_M (\Sigma)$-module 
of rank $|W_\sigma|$.

Next we indicate how to associate a root datum to $\mr{End}_G (I_P^G (\Sigma))$.
See \cite[Section 6]{Hei} for more details in the case of classical groups.
Let $A$ be the maximal split torus of $Z(M)$, let $X^* (A)$ be its character lattice and
$X_* (A)$ its cocharacter lattice. There are natural isomorphisms 
\[
X_{ur}(M) \cong X^* (A) \otimes_\Z \C^\times \cong \mr{Hom} (X_* (A), \C^\times ).
\]
In $X^* (A)$ we have the root system $R (G,A)$ and in $X_* (A)$ we have the set $R^\vee (G,A)$ of
coroots of $(G,A)$. The parabolic subgroup $P$ determines positive systems $R(P,A)$ and $R^\vee (P,A)$.
Altogether we constructed a (nonreduced) based root datum
\[
\mc R_M := \big( X_* (A),R^\vee (G,A),X^*(A),R(G,A), R^\vee(P,A) \big) ,
\]
from which one can of course deduce a reduced based root datum. 

Yet $\mc R_M$ is not good enough, it does not take $\sigma$ into account. Put 
\[
\begin{array}{ccccccc}
X_\sigma & := & \mr{Hom} (X_{ur}(M_\sigma), \C^\times) & \cong &
\mr{Hom} (X_{ur}(M_\sigma \cap A), \C^\times) & \subset & X_* (A), \\ 
Y_\sigma & := & \mr{Hom}(X_\sigma ,\Z) & \cong & X^* (M_\sigma \cap A) & \supset & X^* (A).
\end{array}
\]
Assume for simplicity that $\sigma \big|_{M^0}$ is multiplicity free, or equivalently that
$\mr{End}_M (\Sigma) = \mc O (X_{ur}(M_\sigma))$.
Then the above says that $\mr{End}_G (I_P^G (\Sigma))$ is a free module of rank $|W_\sigma|$
over $\C [X_\sigma]$. We want to associate a root system to $W_\sigma$. In general $W_\sigma$ 
does not contain $W(G,A) = N_G (M) / M$, so we have to replace $R(G,A)$ by
\[
R_{\sigma,nr} := \{ \alpha \in R (G,A) : s_\alpha \in W_\sigma \} ,
\]
and $R^\vee (G,A)$ by $R^\vee_{\sigma,nr}$. Let $R_\sigma^\vee$ be the reduced root system of 
$R_{\sigma,nr}^\vee$ and $R_\sigma$ the dual root system, which consists of the non-multipliable
roots in $R_{\sigma,nr}$. Let $F_\sigma$ be the unique basis of $R_\sigma$ contained in
$R (P,A)$. Then $W (R_\sigma)$ is a normal subgroup of $W_\sigma$ and
\[
W_\sigma \cong W (R_\sigma) \rtimes \Gamma_\sigma \quad \text{where} \quad 
\Gamma_\sigma = \{ w \in W_\sigma : w (F_\sigma) = F_\sigma \} .
\]
As $\sigma$ is not explicit, it is difficult to say which diagram automorphism groups $\Gamma_\sigma$
can occur here. A priori there does not seem to be any particular restriction.

All this suggests that, if $\mr{End}_G (I_P^G (\Sigma))$ is isomorphic to some affine Hecke algebra,
then to
\begin{equation}\label{eq:AHAsigma}
\mc H_\sigma \rtimes \Gamma_\sigma = \mc H (\mc R_\sigma, q_\sigma) \rtimes \Gamma_\sigma :=
\mc H (X_\sigma, R_\sigma^\vee, Y_\sigma, R_\sigma, F_\sigma^\vee, q_\sigma ) \rtimes \Gamma_\sigma .
\end{equation}
In fact it also possible that the $\Gamma_\sigma$-action is twisted by a cocycle \cite[Section 7]{Mor},
but we ignore this subtlety here. We note that little would change upon replacing $G$ by a 
disconnected group, that would only lead to a larger group of diagram automorphisms.

We note that this description of $\mr{End}_G (I_P^G (\Sigma))$ is compatible with parabolic induction. 
Every parabolic subgroup $Q \subset G$ containing $P$ gives rise to a subalgebra 
\[
\mr{End}_Q (I_P^Q (\Sigma)) \subset \mr{End}_G (I_P^G (\Sigma)) ,
\]
which via \eqref{eq:AHAsigma} and
\[
R_{\sigma,nr}^Q = R(Q,A) \subset R(P,A) = R_{\sigma,nr}
\]
corresponds to a parabolic subalgebra $\mc H_\sigma^Q \rtimes \Gamma_{\sigma,Q} \subset 
\mc H_\sigma \rtimes \Gamma_\sigma$.

By analogy with the Iwahori case the numbers $q_\sigma (w)$ are 
related to the affine Coxeter complex of $X_* (A) \rtimes W (R^\vee (G,A))$.
After fixing a fundamental chamber $C_0$, every $w \in X_* (A) \rtimes W_\sigma$ determines a chamber 
$w (C_0)$. This affine Coxeter complex can be regarded as a subset of the Bruhat--Tits building of
$G$, so $C_0$ has a stabilizer $K \subset G$. In view of \cite[Section 6]{Mor}, a good candidate for
$q_\sigma (w)$ is $[K w K : K]$. In particular, for a simple reflection $s_\alpha \in W(R_\sigma^\vee)$
this works out to $q_\sigma (s_\alpha) = q^{d_\alpha}$, where $q$ is cardinality of the residue field
of $\mathbb F$ and $d_\alpha$ is the dimension of the $\alpha$-weight space in the $A$-representation
Lie$(G)$. Hence $q_\sigma$ is a positive parameter function and, if $A$ is not a split maximal torus of $G$,
$q_\sigma$ tends to be non-constant on the set of simple reflections. 

As said before, most of the above is conjectural. The problem is that in general it is not known whether
one can construct elements $N_w \; (w \in W_\sigma)$ that satisfy the multiplication rules of an
extended affine Hecke algebra. To that end one has to study the intertwining operator between 
parabolically induced representations very carefully. 

Let us list the cases in which it is proven that 
$\mr{End}_G (I_P^G (\Sigma))$ is isomorphic to an (extended) affine Hecke algebra:
\begin{itemize}
\item $G$ split, $\mf s$ the Iwahori-spherical component \cite{IwMa,Bor};
\item $G = GL_n (\mathbb F)\; \mf s$ arbitrary - 
from the work of Bushnell and Kutzko on types \cite{BuKu1,BuKu2,BuKu3};
\item $G = SL_n (\mathbb F)$, many $\mf s$ \cite{GoRo} (for general $\mf s$ the Hecke algebra is known to
have a closely related shape);
\item $G$ a (special) orthogonal group, a symplectic group or an inner form of $GL_n (\mathbb F)$, 
$\mf s$ arbitrary \cite{Hei};
\item $G = GSp_4 (\mathbb F)$ or $G = U(2,1),\; \mf s$ arbitrary \cite{Moy1,Moy2};
\item $G$ classical, certain $\mf s$ \cite{Kim1,Kim2,Blo};
\item $G$ split (with mild restrictions on the residual characteristic), $\mf s$ in the principal series \cite{Roc1};
\item $G$ arbitrary, $\sigma$ induced from a level 0 cuspidal representation of a
parahoric subgroup of $G$ \cite{Mor,MoPr,Lus-Uni};
\end{itemize}
Of course there is a lot of overlap in this list. For $GL_n, SL_n, GSp_4$ and $U(2,1)$ the above references do 
much more, they classify the smooth dual of $G$. In the level 0 case, Morris \cite{Mor} showed that the 
parameters $q_\alpha$ are the same as those for analogous Hecke algebras of finite Chevalley groups. 
Those parameters were determined explicitly in \cite{Lus-fin}, and often they are not equal on all
simple roots.

Apart from this list, there are many inertial equivalence classes $\mf s$
for which $\mr{End}_G (I_P^G (\Sigma))$ is Morita-equivalent a commutative algebra. This is the case
for supercuspidal $G$-representations $\sigma$ such that $\sigma \big|_{G^0}$ is multiplicity-free,
and more generally it tends to happen when $R_\sigma$ is empty.

\chapter{Classification of irreducible representations}

This chapter leads to the main result of the paper (Theorem \ref{thm:2.7}).
It is an affine analogue of the Springer correspondence. Together with Kazhdan--Lusztig-theory,
the classical Springer correspondence parametrizes the irreducible representations of a finite 
Weyl group with certain representations of an affine Hecke algebra with equal parameters 
\cite{BaMo1}. This correspondence is known to remain valid for affine or graded Hecke 
algebras with certain specific unequal parameters \cite{Ciu}.

We construct a natural map from irreducible $\mc H$-representations to representations
of the extended affine Weyl group $W$. Not all representations in the image are 
irreducible, but the image does form a $\Q$-basis of the representation ring of $W$.

The proof proceeds by reduction to a result that the author previously obtained for
graded Hecke algebras \cite{SolHomGHA}. To carry out this reduction, we need variations
on three well-known results in representation theory of Hecke algebras. The first two
are due Lusztig and allow one to descend from affine Hecke algebras to graded Hecke algebras.
We adjust these results to make them more suitable for affine Hecke algebras with
arbitrary positive parameters. 

Thirdly there is the Langlands classification (Theorem \ref{thm:2.4}), 
which comes from reductive groups and reduces 
the classification of irreducible representations to that of irreducible tempered ones. 
For affine Hecke algebras it did not appear in the literature before, although it was of 
course known to experts. Because we want to include diagram automorphisms in our affine 
Hecke algebras, we need a more refined version of the Langlands classification 
(Corollary \ref{cor:2.8}). It turns out that one has to add an extra ingredient to the 
Langlands parameters, and that the unicity claim has to be changed accordingly.

However, these results do not suffice to complete the proof of Theorem \ref{thm:2.7}
for nontempered representations, that will be done in the next chapter.

\section{Two reduction theorems}
\label{sec:reduction}

The study of irreducible representations of $\mc H \rtimes \Gamma$  is simplified by 
two reduction theorems, which are essentially due to Lusztig \cite{Lus-Gr}.
The first one reduces to the case of modules whose central character is positive on
the lattice $\Z R_1$. The second one relates these to modules of an associated graded
Hecke algebra. 

Given $t \in T$ and $\alpha \in R_0$, \cite[Lemma 3.15]{Lus-Gr} tells us that
\begin{equation}\label{eq:sAlphaFixt}
s_\alpha (t) = t \text{ if and only if } \alpha (t) = 
\left\{ \begin{array}{ll}
1 & \text{if } \alpha^\vee \notin 2 Y \\
\pm 1 & \text{if } \alpha^\vee \in 2 Y .
\end{array} \right.
\end{equation}
We define $R_t := \{ \alpha \in R_0 : s_\alpha (t) = t \}$. The collection of long
roots in $R_{t,nr}$ is $\{ \beta \in R_1 : \beta (t) = 1\}$. Let $F_t$ be 
the unique basis of $R_t$ that is contained in $R_0^+$. Then
\[
W'_{F_t ,t} := \{ w \in W_0 \rtimes \Gamma : w(t) = t , w (F_t) = F_t \}
\]
is a group of automorphisms of the Dynkin diagram of $(R_t ,F_t)$. 
Moreover the isotropy group of $t$ in $W_0 \rtimes \Gamma$ is
\[
W'_t = (W_0 \rtimes \Gamma)_t = W (R_t) \rtimes W'_{F_t ,t} .
\]
We can define a parameter function $q_t$ for the based root datum
\[
\mc R_t := (X,R_t,Y,R_t^\vee,F_t)
\]
via restriction from $R_{nr}^\vee$ to $R_{t,nr}^\vee$.

Since $F_t$ does not have to be a subset of $F_0 \,, \mc R_t$ does not always fit
in the setting of Subsection \ref{sec:parabolic}, but this can be fixed without
many problems. For $u \in T_{un}$ we define 
\[
P(u) := F_0 \cap \Q R_u .
\]
Then $R_{P(u)}$ is a parabolic root subsystem of $R_0$ that contains $R_u$ 
as a subsystem of full rank. Although this definition would also make sense for general
elements of $T$, we use it only for $T_{un}$, to avoid a clash with the notation of
\cite[Section 4.1]{Opd-Sp}. We note that the lattice
\[
\Z P(u) = \Z R_0 \cap \Q R_u 
\]
can be strictly larger than $\Z R_u$. 

To study $\mc H$-representations 
with central character $\WG uc$ we need a well-chosen neighborhood of $uc \in T_{un} T_{rs}$.
\begin{cond}\label{cond:ball}
Let $B \subset \mf t$ be such that
\enuma{
\item $B$ is an open ball centred around $0 \in \mf t$;
\item $\Im (\alpha (b)) < \pi$ for all $\alpha \in R_0, b \in B$;
\item $\exp : B \to \exp (B)$ is a diffeomorphism (this follows from (b) if $R_0$ spans $\mf a^*$);
\item if $c_\alpha (t) \in \{0,\infty\}$ for some $\alpha \in R_0 , t \in uc \exp \overline{B}$, then
$c_\alpha (uc) \in \{0,\infty\}$;
\item if $w \in \WG$ and $w (uc \exp \overline{B}) \cap uc \exp \overline{B} \neq \emptyset$,
then $w(uc) = uc$.
}
\end{cond}
Since $\WG$ acts isometrically on $\mf t$, (a) implies that $B$ is $\WG$-invariant. 
There always exist balls satisfiying these conditions, and if we have one such $B$, 
then $\ep B$ with $\ep \in (0,1]$ also works.

We will phrase our first reduction theorem in such a way that it depends mainly on the unitary part of
the central character, it will decompose a representation in subspaces corresponding to the points of
the orbit $\WG u$. We note that $R_{uc} \subset R_u$ and $W'_{uc} \subset W'_u$. 
Given $B$ satisfying the above conditions, we define
\begin{equation*}
U = \WG uc \exp (B) ,\: U_{P(u)} = W'_{\Z P(u)} uc \exp (B) 
\text{ and } U_u = W'_u uc \exp (B) .
\end{equation*}
We are interested in the algebras $\mc H (\mc R,q)^{an}(U) \rtimes \Gamma ,\: 
\mc H (\mc R^{P(u)},q^{P(u)})^{an} (U_{P(u)}) \rtimes W'_{P(u)}$ and
$\mc H (\mc R_u ,q_u)^{an} (U_u)\rtimes W'_{F_u ,u}$. Their respective centers
\[
C^{an}(U)^{\WG} ,\: C^{an}(U_{P(u)})^{W'_{\Z P(u)}} \text{ and }
C^{an}(U_u)^{W'_u} 
\]
are naturally isomorphic, via the embeddings $U_u \subset U_{P(u)} \subset U$.
For any subset $\varpi \subset \WG uc$ we define $1_\varpi \in C^{an}(U)$ by
\[
1_\varpi (t) = \left\{ \begin{array}{ll}
1 & \text{if } \: t \in \varpi \exp (B) \\
0 & \text{if } \: t \in U \setminus \varpi \exp (B) .
\end{array} \right. 
\]

\begin{thm}\label{thm:2.1}\textup{(First reduction theorem)} 
\enuma{
\item There are natural isomorphisms of $C^{an}(U)^{\WG}$-algebras
\[
\begin{array}{lll}
\mc H (\mc R^{P(u)}, q^{P(u)})^{an}(U_{P(u)}) \rtimes W'_{P(u)} & \cong & 1_{W'_{\Z P(u)} uc}
(\mc H^{an}(U) \rtimes \Gamma) \, 1_{W'_{\Z P(u)} uc} , \\
\mc H (\mc R_u ,q_u)^{an} (U_u) \rtimes W'_{F_u ,u} & \cong & 
1_{W'_u uc} (\mc H^{an}(U) \rtimes \Gamma) \, 1_{W'_u uc} .
\end{array}
\]
\item These can be extended (not naturally) to isomorphisms of $C^{an}(U)^{\WG}$-algebras
\[
\begin{array}{lll}
\mc H^{an}(U) \rtimes \Gamma  & \cong & 
M_{[\WG : W'_{\Z P(u)}]} \big( 1_{W'_{\Z P(u)} uc} (\mc H^{an}(U) \rtimes \Gamma) \, 
1_{W'_{\Z P(u)} uc} \big) , \\
\mc H^{an}(U) \rtimes \Gamma  & \cong & 
M_{[\WG : W'_u]} \big( 1_{W'_u uc} (\mc H^{an}(U) \rtimes \Gamma) \, 1_{W'_u uc} \big) ,
\end{array}
\]
where $M_n (A)$ denotes the algebra of $n \times n$-matrices with coefficients in an algebra~$A$. 
\item The following maps are natural equivalences of categories:
}
\[
\begin{array}{ccccc}
\!\!\! \mr{Mod}_{f,U}(\mc H (\mc R ,q) \! \rtimes \! \Gamma) \!\!\! & \leftrightarrow &
\!\!\! \mr{Mod}_{f,U_{P(u)}} (\mc H^{P(u)} \! \rtimes \! W'_{P(u)}) \!\!\! & \leftrightarrow &
\!\!\! \mr{Mod}_{f,U_u} (\mc H (\mc R_u ,q_u) \! \rtimes \! W'_{F_u ,u}) \\
V & \mapsto & 1_{W'_{\Z P(u)} uc} V & \mapsto & 1_{W'_u uc} V \\
\!\!\! \mr{Ind}_{\mc H (\mc R_u,q_u) \rtimes W'_{F_u,u}}^{\mc H \rtimes \Gamma} (\pi) \!\!\! &
\mathrel{\reflectbox{$\mapsto$}} & \mr{Ind}_{\mc H (\mc R_u,q_u) \rtimes 
W'_{F_u,u}}^{\mc H^{P(u)} \rtimes W'_{P(u)}} (\pi) & \mathrel{\reflectbox{$\mapsto$}} & \pi
\end{array}
\]
\end{thm}
\emph{Proof.}
(a) This is a variation on \cite[Theorem 4.10]{Opd-Sp}, which itself varied on 
\cite[Theorem 8.6]{Lus-Gr}. Compared to Lusztig we replaced $Y \otimes \langle v_0 \rangle \subset T$
by $T_{rs} = Y \otimes \R_{>0}$, we substituted his $R_t$ by a larger root system,
we added the group $\Gamma$ and we localized with analytic functions instead of formal completions
at one central character. 

The first change is allowed because we localize only on $T$, and not simultaneously in the $q$-direction 
like Lusztig. Hence the structure of the subgroup $\langle v_0 \rangle \subset \C^\times$ becomes 
irrelevant to the arguments of \cite[Section 8]{Lus-Gr}. With this modification Lusztig's $R_t$ becomes 
our $R_u$ and his equivalence classes $c \in \mathcal P$ correspond to the orbits $W'_u uc \subset \WG uc$.
The root system $R_{P(u)}$ and the orbits $W'_{\Z P(u)} uc \subset \WG uc$ fall into the same
framework, so the constructions from \cite[Section 8]{Lus-Gr} can still be carried out. Since
$\Z R_u \subset \Z R_{P(u)}$ there are fewer orbits in the second case, so it actually becomes somewhat
easier. 

By \cite[Lemma 8.13.b]{Lus-Gr} Lusztig's version of the isomorphisms (a) sends $\gamma \in \Gamma (c)$ 
to $1_c \imath^0_\gamma 1_c$, for some $c \in \mc P$. Translated to our setting this means that we can 
include the appropriate diagram automorphisms by defining
\begin{equation}\label{eq:firstredgamma}
\begin{array}{cccccc}
W'_{P(u)} & \ni & \gamma & \mapsto & 1_{W'_{\Z P(u)} uc} \, \imath^0_\gamma \, 1_{W'_{\Z P(u)} uc} , \\
W'_{F_u,u} & \ni & \gamma & \mapsto & 1_{W'_u uc} \, \imath^0_\gamma \, 1_{W'_u uc} . 
\end{array}
\end{equation}
Finally, that Lusztig's arguments also apply with analytic localization was already checked by
Opdam \cite[Section 4.1]{Opd-Sp}. \\
(b) Knowing (a), this can proved just as in \cite[4.16]{Lus-Gr}.\\
(c) By \cite[Proposition 4.3]{Opd-Sp} the categories in the statement are of the categories
of finite dimensional modules of the algebras figuring in (a) and (b). Therefore the maps in (c)
are just the standard equivalences between the module categories of $B$ and $M_n (B)$,
translated with (a) and (b). $\qquad \Box$

\begin{rem}\label{rem:firstred}
This reduction theorem more or less forces one to consider diagram
automorphisms: the groups $W'_{\Z P(u)}$ and $W'_{F_u ,u}$ can be nontrivial even if 
$\Gamma = \{ \text{id} \}$. 

The notation with induction functors in part (c) is a little sloppy, since $W'_{F_u,u}$ need
not be contained in $\Gamma$ or in $W'_{\Z P(u)}$. In such cases these induction functors are
defined via part (a).

The first reduction theorem also enables us to make sense of $\mr{Ind}_{\mc H (\mc R^P,q^P) 
\rtimes \Gamma'_P}^{\mc H \rtimes \Gamma}$ for any $P \subset F_0$ and any subgroup 
$\Gamma'_P \subset W'_P$. Namely, first induce from $\mc H (\mc R^P,q^P) \rtimes \Gamma'_P$ to 
$\mc H (\mc R^P,q^P) \rtimes W'_{\Z P}$, then choose $u \in T_{un}$ such that $P(u) = P$
and finally use (c).
\end{rem}

By \eqref{eq:sAlphaFixt} we have $\alpha (u) = 1$ for all $\alpha \in R_1 \cap \Q R_t$, so
$\alpha (t) = \alpha (u) \alpha (c) > 0$ for such roots. By definition $u$ is fixed by 
$W'_u$, so Theorem \ref{thm:2.1} allows us to restrict
our attention to $\mc H \rtimes \Gamma$-modules whose central character is positive
on the sublattice $\Z R_1 \subseteq X$.

Next we want to reduce to graded Hecke algebras.
We define a parameter function $k_u$ for the degenerate root datum 
$\tilde{\mc R}_u = (\mf a , R_u, \mf a^*, R_u^\vee, F_u)$ by
\begin{equation}\label{eq:kualpha}
k_{u,\alpha} = \big( \log q(s_\alpha) + \alpha (u) \log q(t_\alpha s_\alpha) \big) / 2 .
\end{equation}
Recall that $\alpha \in R_u$ implies that $s_\alpha (u) = u$ and $\alpha (u)^2 = 1$.
We will see in \eqref{eq:lim0cep} that for this choice of $k_u$ the function $\tilde c_\alpha$
can be regarded as the first order approximation of $q(s_\alpha) c_\alpha$ in a neighborhood of $q=1$
and $u \in T$. 

Let us pick $u \in T_{un}^{\WG}$, so $\alpha (u) = \pm 1$ for all $\alpha \in R_0$. Then the map
\begin{equation}\label{eq:expMap}
\exp_u : \mf t \to T ,\;  \lambda \mapsto u \exp (\lambda)
\end{equation}
is $\WG$-equivariant. To find a relation between $\mc H (\mc R ,q) \rtimes \Gamma$ and 
$\mh H (\tilde {\mc R}_u ,k_u) \rtimes \Gamma$, we first extend these algebras with analytic localization.
For every open nonempty $\WG$-invariant $V \subset \mf t$ we can define an algebra homomorphism
\begin{align}
\nonumber \Phi_u : \mc H^{me}(\exp_u (V)) \rtimes \Gamma & \to
\mh H (\tilde{\mc R}_u,k_u)^{me}(V) \rtimes \Gamma , \\
\label{eq:Phiu} f \imath^0_w & \mapsto (f \circ \exp_u) \tilde \imath_w .
\end{align}

\begin{thm}\label{thm:2.2}\textup{(Second reduction theorem)} \\
Let $u \in T_{un}^{\WG}$ and let $V$ be as above, such that moreover $\exp_u$ is injective on $V$.
\enuma{
\item The map $\exp_u$ induces an isomorphism $C^{an}(\exp_u (V))^{\WG} \to C^{an}(V)^{\WG}$.
\item Suppose that every $\lambda \in V$ satisfies
\begin{equation}\label{eq:condLambda}
\begin{array}{llr@{\quad}l}
\inp{\alpha}{\lambda}, \inp{\alpha}{\lambda} + k_{u,\alpha} & \notin & 2 \pi i \Z \setminus \{0\} & 
\text{for } \alpha \in R_0 \cap R_1 , \\
\inp{\alpha}{\lambda}, \inp{\alpha}{\lambda} + k_{u,\alpha} & \notin & \pi i \Z \setminus \{0\} & 
\text{for } \alpha \in R_0 \setminus R_1 . 
\end{array}
\end{equation}
Then $\Phi_u$ restricts to an isomorphism of $C^{an}(V)^{\WG}$-algebras
\begin{equation*}
\Phi_u : \mc H^{an} (\exp_u (V)) \rtimes \Gamma \to \mh H (\tilde{\mc R}_u, k_u)^{an}(V) \rtimes \Gamma .
\end{equation*}
}
\end{thm}
\emph{Proof.}
(a) This is clear, it serves mainly to formulate (b).\\
(b) The case $\Gamma = \{ \textup{id} \}$ is essentially \cite[Theorem 9.3]{Lus-Gr}. 
The difference is that our conditions on $\lambda$ replace the conditions \cite[9.1]{Lus-Gr}.
The general case follows easily under the assumption that $\Gamma$ fixes $u. \qquad \Box$
\\[2mm]
Given $\mf t' \subset \mf t$ we denote by $\mr{Mod}_{f,\mf t'}(\mh H (\tilde {\mc R},k) \rtimes \Gamma)$
the category of finite dimensional $\mh H (\tilde {\mc R},k) \rtimes \Gamma$-modules all 
whose $\mc O (\mf t)$-weights lie in $\mf t'$.
\begin{cor}\label{cor:2.3}
Let $q$ be positive and let $u c \in T_{un} T_{rs}$. The following categories are equivalent:
\enuma{
\item $\mr{Mod}_{f,\WG uc} (\mc H (\mc R,q) \rtimes \Gamma)$ and
$\mr{Mod}_{f,(W(R_u) \rtimes W'_{F_u,u}) \log (c)} (\mh H (\tilde{\mc R}_u, k_u) \rtimes W'_{F_u,u})$,
\item $\mr{Mod}_{f,\WG u T_{rs}} (\mc H (\mc R,q) \rtimes \Gamma)$ and
$\mr{Mod}_{f,\mf a} (\mh H (\tilde{\mc R}_u, k_u) \rtimes W'_{F_u,u})$. 
}
These equivalences are compatible with parabolic induction.
\end{cor}
\emph{Proof.} (a) follows from Theorems \ref{thm:2.1}.b and \ref{thm:2.2}.b.
Notice that the conditions \eqref{eq:condLambda} are automatically satisfied because $q$ 
is positive and $\log (c) \in \mf a$, so $k_{u,\alpha} \in \R$ and $\inp{\alpha}{\log (c)} \in \R$. 
If we sum that equivalence over all $\WG c \in T_{rs}/\WG$, we find (b).
By \cite[Theorem 6.2]{BaMo2} or \cite[Proposition 5.3.a]{SolGHA} these equivalences of 
categories are compatible with parabolic induction. $\qquad \Box$ 
\\[2mm]

\section{The Langlands classification}

In this section we discuss Langlands' classification of irreducible representations. Basically
it reduces from general representations to tempered ones, and from there to the discrete
series. Actually Langlands proved this only in the setting of real reductive groups, but it holds 
just as well for $p$-adic reductive groups, affine Hecke algebras and graded Hecke 
algebras. We will only write down the results for affine Hecke algebras, the graded Hecke algebra 
case is completely analogous and can be found in \cite{Eve,KrRa,SolGHA}.

An important tool to study $\mc H$-representations is restriction to the commutative 
subalgebra $\mc A \cong \mc O (T)$. We say that $t \in T$ is a weight of $(\pi,V)$ if there
exists a $v \in V \setminus \{ 0 \}$ such that $\pi (a) v = a(t) v$ for all $a \in \mc A$.
It is easy to describe how the collection of $\mc A$-weights behave under parabolic induction.
Recall that 
\begin{equation}
W^P := \{ w \in W_0 : w(P) \subset R_0^+ \} 
\end{equation}
is the set of minimal length representatives of $W_0 / W (R_P)$.

\begin{lem}\label{lem:2.12}
Let $\Gamma'_P$ be a subgroup of $\Gamma_P$ and let $\sigma$ be a representation of
$\mc H^P \rtimes \Gamma'_P$. The $\mc A$-weights of 
$\mr{Ind}_{\mc H^P \rtimes \Gamma'_P}^{\mc H \rtimes \Gamma} (\sigma)$ are the elements
$\gamma w (t) \in T$, where $\gamma \in \Gamma, w \in W^P$ and $t$ is an $\mc A$-weight of $\sigma$.
\end{lem}
\emph{Proof.}
From \cite[Theorem 6.4]{BaMo1} and the proof of \cite[Proposition 4.20]{Opd-Sp} we see that 
this holds in the case $\Gamma = \Gamma'_P = \{ \mathrm{id} \}$. For the general case we only have 
to observe that the operation $\pi \mapsto \pi \circ \psi_\gamma^{-1}$ on $\mc H$-representations 
has the effect $t \mapsto \gamma (t)$ on all $\mc A$-weights $t. \qquad \Box$
\\[2mm]

Temperedness of a representation is defined via its $\mc A$-weights.
Given $P \subseteq F_0$, we have the following positive cones in $\mf a$ and in $T_{rs}$: 
\begin{equation}
\begin{array}{lll@{\qquad}lll}
\mf a^+ & = & \{ \mu \in \mf a : \inp{\alpha}{\mu} \geq 0 \: 
  \forall \alpha \in F_0 \} , & T^+ & = & \exp (\mf a^+) , \\
\mf a_P^+ & = &  \{ \mu \in \mf a_P : \inp{\alpha}{\mu} \geq 0 \: \forall \alpha \in P \} , & 
  T_P^+ & = & \exp (\mf a_P^+) , \\
\mf a^{P+} & = & \{ \mu \in \mf a^P : \inp{\alpha}{\mu} \geq 0 \: 
  \forall \alpha \in F_0 \setminus P \} , & T^{P+} & = & \exp (\mf a^{P+}) , \\
\mf a^{P++} & = & \{ \mu \in \mf a^P : \inp{\alpha}{\mu} > 0 \; 
  \forall \alpha \in F_0 \setminus P \} , & T^{P++} & = & \exp (\mf a^{P++}) .
\end{array}
\end{equation}
The antidual of $\mf a^{*+} :=  \{ x \in \mf a^* :  \inp{x}{\alpha^\vee} \geq 0 
\: \forall \alpha \in F_0 \}$ is
\begin{equation} 
\mf a^- = \{ \lambda \in \mf a : \inp{x}{\lambda} \leq 0 \: \forall x \in \mf a^{*+} \} = 
\big\{ \sum\nolimits_{\alpha \in F_0} \lambda_\alpha \alpha^\vee : \lambda_\alpha \leq 0 \big\} .
\end{equation}
Similarly we define
\begin{equation}
\mf a_P^- = \big\{ \sum\nolimits_{\alpha \in P} \lambda_\alpha \alpha^\vee  \in \mf a_P : 
\lambda_\alpha \leq 0 \big\} .
\end{equation}
The interior $\mf a^{--}$ of $\mf a^-$ equals
$\big\{ {\ts \sum_{\alpha \in F_0}} \lambda_\alpha \alpha^\vee : \lambda_\alpha < 0 \big\}$
if $F_0$ spans $\mf a^*$, and is empty otherwise. We write $T^- = \exp (\mf a^-)$ and
$T^{--} = \exp (\mf a^{--})$. 
Let $t = |t| \cdot t |t|^{-1} \in T_{rs} \times T_{un}$ be the polar decomposition of $t$.

An $\mc H$-representation is called tempered if $|t| \in T^-$ for all its 
$\mc A$-weights $t$, and anti-tempered if $|t|^{-1} \in T^-$ for all such $t$. 
For infinite dimensional representations this is not entirely satisfactory,
but we postpone a more detailed discussion to Section \ref{sec:Schwartz}. Since all irreducible 
$\mc H$-representations have finite dimension, this vagueness does not cause any problems.
Notice that our definition mimics Harish-Chandra's definition of admissible smooth tempered 
representations of reductive $p$-adic groups \cite[Section III.2]{Wal}. In that setting the crucial
condition says that all exponents of such a representation must lie in certain cone.

More restrictively we say that an irreducible  $\mc H$-representation belongs to the discrete 
series (or simply: is discrete series) if $|t| \in T^{--}$, for all its $\mc A$-weights $t$. 
In particular the discrete series is empty if $F_0$ does not span $\mf a^*$. This is the analogue
of Casselman's criterium for square integrable representations of semisimple $p$-adic groups
\cite[Theorem 4.4.6]{Cas}.

The notions tempered and discrete series apply equally well to $\mc H \rtimes \Gamma$,
since that algebra contains $\mc A$ and the action of $\Gamma$ on $T$ preserves $T^-$.
It follows more or less directly from the definitions that the correspondence of Theorem
\ref{thm:2.2} preserves temperedness and provides a bijection between discrete
series representations with the appropriate central characters, see \cite[(2.11)]{Slo}.

It easy to detect temperedness for 
$\mc H (\mc R ,1) \rtimes \Gamma = \C [W']  = \C [X \rtimes \WG]$. 
\begin{lem}\label{lem:2.9}
A finite dimensional $\C [X \rtimes \WG]$-representation is
tempered if and only if all its $\mc A$-weights lie in $T_{un}$.

This algebra has no discrete series representations, unless $X = 0$.
\end{lem}
\emph{Proof.}
Suppose that $V$ is a representation of this algebra, and that $t \in T$ is an $\mc A$-weight 
with weight space $V_t$. For every $g \in \WG ,\: g V_t = V_{g(t)}$ is the $g(t)$-weight space of $V$, 
which shows that every element of the orbit $\WG t$ is an $\mc A$-weight of $V$. But 
$\WG |t|$ can only be contained in $T^-$ if it equals the single element 
$1 \in T_{rs}$. Hence $V$ can only be tempered if $|t| = 1$ for all its weights, or equivalently 
if all its weights lie in $T_{un}$. By definition the latter condition also suffices for temperedness.

Unless $X = 0$, the condition $|t| =1$ implies $|t| \not\in T^{--}$, so 
$\C [X \rtimes \WG]$ has no discrete series representations.
$\qquad \Box$ \\[2mm]

The Langlands classification looks at parabolic subalgebras of $\mc H$ and irreducible representations
of those that are essentially tempered. We will describe such representations with two data:
a tempered representation and a "complementary" part of the central character. This is justified by 
the following result.

\begin{lem}\label{lem:2.11}
Let $P \subset F_0, t_P \in T_P$ and $t^P \in T^P$. 
\enuma{
\item The map $\sigma \mapsto \sigma \circ \phi_{t^P}$ defines an equivalence between the 
categories of $\mc H_P$-representations with central character $W (R_P) t_P \in T_P / W (R_P)$ 
and of $\mc H^P$-representations with central character $W (R_P) t_P t^P \in T / W (R_P)$.
\item Every irreducible $\mc H^P$-representation is of the form $\sigma \circ \phi_{t^P}$, 
where $\sigma$ is an irreducible $\mc H_P$-representation and $t^P \in T^P$. Both these data 
are unique modulo twists coming from $K_P = T^P \cap T_P$, as in \eqref{eq:twistKP}.
}
\end{lem}
\emph{Proof.}
(a) The kernel of $\phi_{t^P}$ followed by the quotient map $\mc H^P \to \mc H_P$ is generated (as an ideal)
by $\{ \theta_x - t^P (x) : x \in X \cap (P^\vee)^\perp \}$. If $\rho$ is an $\mc H^P$-representation with central 
character $W (R_P) t_P t^P$, then the kernel of $\rho$ clearly contains these generators, so $\rho$ factors
via $\phi_{t^P}$ and this quotient map. \\
(b) Let $\rho$ be an irreducible $\mc H^P$-representation with central character $W (R_P) t \in T /W (R_P)$.
Decompose $t = t_P t^P \in T_P T^P$. Then part (a) yields a unique irreducible $\mc H_P$-representation 
$\sigma$ such that $\rho = \sigma \circ \phi_{t^P}$. The only freedom in this constuction comes from elements 
$u \in K_P$. If we replace $t^P$ by $u t^P$, then part (a) again gives a unique $\sigma'$ with 
$\rho = \sigma' \circ \phi_{u t^P}$, and its follows directly that $\sigma' \circ \psi_u = \sigma. \qquad \Box$
\\[2mm]

A Langlands datum for $\mc H$ is a triple $(P,\sigma,t)$ such that
\begin{itemize}
\item $P \subseteq F_0$ and $\sigma$ is an irreducible tempered $\mc H_P$-representation;
\item $t \in T^P$ and $|t| \in T^{P++}$. 
\end{itemize}
We say that two Langlands data $(P,\sigma,t)$ and $(P',\sigma',t')$ are equivalent if $P = P'$ and 
the $\mc H^P$-representations $\sigma \circ \phi_t$ and $\sigma' \circ \phi_{t'}$ are equivalent. 

\begin{thm} \label{thm:2.4} \textup{(Langlands classification)}
\enuma{
\item For every Langlands datum $(P,\sigma,t)$ the $\mc H$-representation 
$\pi (P,\sigma,t) = \mr{Ind}_{\mc H^P}^{\mc H} (\sigma \circ \phi_t)$
has a unique irreducible quotient $L(P,\sigma,t)$. 
\item For every irreducible $\mc H$-representation $\pi$ there exists a Langlands 
datum $(P,\sigma,t)$, unique up to equivalence, such that $\pi \cong L (P,\sigma,t)$.
}
\end{thm}
\emph{Proof.}
The author learned this result from a preliminary version of \cite{DeOp2}, but unfortunately
Delorme and Opdam did not include it in the final version. 
Yet the proof in the setting of affine Hecke algebras is much easier than for reductive 
groups. It is basically the same as the proof of Evens \cite{Eve} for graded Hecke algebras, see 
also \cite[Section 2.4]{KrRa}. For later use we rephrase some parts of that proof in our notation.\\
(a) The dominance ordering on $\mf a$ is defined by
\begin{equation}\label{eq:domorder}
\lambda \leq \mu \text{ if and only if } \inp{\lambda}{\alpha} \leq \inp{\mu}{\alpha} 
\text{ for all } \alpha \in F_0 .
\end{equation}
For $\alpha \in F_0$ we define $\delta_\alpha \in \mf a_{F_0}$ by
\[
\inp{\beta}{\delta_\alpha} = \left\{ \begin{array}{lll}
1 & \mr{if} & \alpha = \beta \\
0 & \mr{if} & \alpha \neq \beta \in F_0 \,.
\end{array} \right.
\]
According to Langlands \cite[Lemma 4.4]{Lan}, for every $\lambda \in \mf a$ there is a 
unique subset $F (\lambda ) \subset F_0$ such that $\lambda$ can be written as
\begin{equation}\label{eq:LanglandsF}
\lambda = \lambda^{F_0} + \sum_{\alpha \in F_0 \setminus F (\lambda )} c_\alpha \delta_\alpha
+ \sum_{\alpha \in F (\lambda )} d_\alpha \alpha^\vee \quad
\mr{with} \; \lambda^{F_0} \in \mf a^{F_0} , c_\alpha > 0, d_\alpha \leq 0 .
\end{equation}
We put $\lambda_0 = \sum_{\alpha \in F_0 \setminus F (\lambda )} c_\alpha \delta_\alpha \in \mf a^+$. 
According to \cite[(2.13)]{KrRa}
\begin{equation}\label{eq:2.3}
(w \mu)_0 < \mu_0 \text{ for all } \mu \in \mf a_P^- \oplus \mf a^{P++}, w \in W^P \setminus \{1\} .
\end{equation}
By the definition of a Langlands datum $\log |s| \in \mf a_P^- \oplus \mf a^{P++}$ for every 
$\mc A$-weight $s$ of $\sigma \circ \phi_t$. Choose $s$ such that $(\log |s|) _0$ is maximal with 
respect to the dominance order. By Lemma \ref{lem:2.12} and \eqref{eq:2.3} $(\log |s|) _0$ is also 
maximal for $s$ regarded as an $\mc A$-weight of $\pi (P,\sigma,t)$. 

Suppose that $\rho$ is an $\mc H$-submodule of $\pi (P,\sigma,t)$ of which $s$ is an $\mc A$-weight. 
By the maximalty of $s$, $\rho$ must contain the $s$-weight space of $\sigma \circ \phi_t$. 
The irreduciblity of $\sigma$ implies that $\rho$ contains the $\mc H^P$-submodule
$1 \otimes_{\mc H^P} V_\sigma \subset \mr{Ind}_{\mc H^P}^{\mc H}(\sigma \circ \phi_t)$, and therefore 
$\rho = \pi (P,\sigma,t)$. Thus the sum of all proper submodules is again proper, which means that 
$\pi (P,\sigma,t)$ has a unique maximal submodule and a unique irreducible quotient.\\
(b) Let $s$ be an $\mc A$-weight of $\pi$ such that $(\log |s|)_0 \in \mf a$ is maximal and put 
$P = F (\log |s|)$. Let $\rho$ be the $\mc H^P$-subrepresentation $\rho$ of $\pi$ generated by the 
$s$-weight space. Then $\log |s| \in \mf a_P^- \oplus \mf a^{P++}$ and according to \cite[p. 38]{KrRa} 
$(\log |s'|)_0 = (\log |s|)_0$
for all $\mc A$-weights $s'$ of $\rho$. By Lemma \ref{lem:2.11} we can write every irreducible 
$\mc H_P$-subrepresentation of $\rho$ as $\sigma \circ \phi_t$, where $\sigma$ is an 
irreducible $\mc H_P$-representation and $\log |t| = (\log |s|)_0$. 
The $\mc A_P$-weights of $\sigma$ are of the form $s' t^{-1}$ and by construction 
\[
\log |s' t^{-1}| = \log |s´| - (\log |s´|)_0 \in \mf a_P^-,
\]
so $\sigma$ is tempered. The inclusion map $\sigma \circ \phi_t \to \pi$ induces a nonzero 
$\mc H$-homomorphism $\pi (P,\sigma,t) \to \pi$. Since $\pi$ is irreducible, this map is surjective.
Together with part (a) this shows that $\pi$ is the unique quotient of $\pi (P,\sigma,t)$.\\
The proof that $(P,\sigma \circ \phi_t)$ is uniquely determined by $\pi$ is easy, and exactly the same 
as in the graded Hecke algebra setting, see \cite[Theorem 2.1.iii]{Eve} or \cite[Theorem 2.4.b]{KrRa}. 
$\qquad \Box$ \\[2mm]

Theorem \ref{thm:2.4} can be regarded as the analogue of the Langlands classification for
connected reductive $p$-adic groups. For disconnected reductive groups the classification is no longer 
valid as such, it has to be modified. In the case that the component group is abelian, this is worked out 
in \cite{BaJa1}, via reduction to cyclic component groups.

We work with a diagram automorphism group $\Gamma$ which is more general than a component group
and does not have to be abelian. For use in Section \ref{sec:Springer} we have to extend Theorem 
\ref{thm:2.4} to this setting.

There is a natural action of $\Gamma$ on Langlands data, by
\begin{equation}\label{eq:gammaLanglands}
\gamma (P,\sigma,t) = (\gamma (P),\sigma \circ \psi_\gamma^{-1}, \gamma (t)) .
\end{equation}
Every Langlands datum yields a packet of irreducible quotients, and all data in one 
$\Gamma$-orbit lead to the same packet. For $\gamma \in \Gamma_P$ the Langlands
classification for $\mc H^P$ shows that the irreducible $\mc H^P$-representations
$\sigma \circ \psi_\gamma \circ \phi_{\gamma (t)}$ and $\sigma \circ \phi_t$ are
equivalent if and only if $\gamma (P,\sigma,t) = (P,\sigma,t)$.

To get a more precise statement one needs Clifford theory, as for example in \cite{RaRa} 
or \cite[Section 53]{CuRe1}.
Let $\Gamma_{P,\sigma,t}$ be the isotropy group of the Langlands datum $(P,\sigma,t)$. 
In \cite[Appendix A]{SolGHA} a 2-cocycle $\kappa$ of $\Gamma_{P,\sigma,t}$ is constructed, 
giving rise to a twisted group algebra $\C [\Gamma_{P,\sigma,t} ,\kappa]$. We define a 
Langlands datum for $\mc H \rtimes \Gamma$ as a quadruple $(P,\sigma,t,\rho)$, where
\begin{itemize}
\item $(P,\sigma,t)$ is a Langlands datum for $\mc H$;
\item $\rho$ is an irreducible representation of $\C [\Gamma_{P,\sigma,t} ,\kappa]$.
\end{itemize}
The action \eqref{eq:gammaLanglands} extends naturally to Langlands data for 
$\mc H \rtimes \Gamma$, since $\psi_\gamma^{-1}$ induces an isomorphism between the
relevant twisted group algebras.

From such a Langlands datum we can construct the $\mc H_P \rtimes \Gamma_{P,\sigma,t}$-representation
$\sigma \otimes \rho$ and the $\mc H \rtimes \Gamma$-representation 
\begin{equation}\label{eq:piLanglands}
\pi^\Gamma (P,\sigma,t,\rho) := \mr{Ind}_{\mc H^P \rtimes \Gamma_{P,\sigma,t}}^{\mc H \rtimes \Gamma}
\big( (\sigma \circ \phi_t) \otimes \rho \big) = \mr{Ind}_{\mc H^P \rtimes \Gamma_{P,\sigma,t}}^{
\mc H \rtimes \Gamma} \big( (\sigma \otimes \rho) \circ \phi_t \big) .
\end{equation}
If $Q \supset P$, then $(P,\sigma,t,\rho)$ can also be considered as a Langlands datum for
$\mc H^Q \rtimes \Gamma_Q$, and we denote the corresponding $\mc H^Q \rtimes \Gamma_Q$-representation
by $\pi^{Q,\Gamma_Q} (P,\sigma,t,\rho)$. In particular $\pi^{P,\Gamma_P} (P,\sigma,t,\rho)$
is an irreducible $\mc H^P \rtimes \Gamma_P$-representation.

\begin{cor} \label{cor:2.8} \textup{(extended Langlands classification)}
\enuma{
\item The $\mc H \rtimes \Gamma$-representation $\pi^\Gamma (P,\sigma,t,\rho)$ has a unique irreducible
quotient $L^\Gamma(P,\sigma,t,\rho)$.
\item For every irreducible $\mc H \rtimes \Gamma$-representation $\pi$ there exists a Langlands 
datum $(P,\sigma,t,\rho)$, unique modulo the action of $\Gamma$, 
such that $\pi \cong L^\Gamma (P,\sigma,t,\rho)$.
\item $L^\Gamma(P,\sigma,t,\rho)$ and $\pi^\Gamma (P,\sigma,t,\rho)$ are tempered if and only if 
$P = F_0$ and $t \in T_{un}^{F_0}$.
}
\end{cor}
\emph{Proof.}
(a) and (b) By \cite[Theorem A.1]{SolGHA} the $\mc H \rtimes \Gamma$-representation
\begin{equation}\label{eq:Lrho}
\mr{Ind}_{\mc H \rtimes \Gamma_{P,\sigma,t}}^{\mc H \rtimes \Gamma} 
(L(P,\sigma,t) \otimes \rho)
\end{equation}
is irreducible, and every irreducible $\mc H \rtimes \Gamma$-representation is of this form,
for a Langlands datum which is unique modulo $\Gamma$. By construction \eqref{eq:Lrho}
is a quotient of $\pi (P,\sigma,t,\rho)$. It is the unique irreducible quotient by Theorem 
\ref{thm:2.4}.a and because $\rho$ is irreducible. \\
(c) If $P \subsetneq F_0$, then $L(P,\sigma,t,\rho)$ and $\pi^\Gamma (P,\sigma,t,\rho)$ are never 
tempered. Indeed $|t| \not\in T^-$, so $|r t| \not\in T^-$ for any $\mc A_P$-weight $r$ of $\sigma$. 
But the construction of $L(P,\sigma,t)$ in the proof of \ref{thm:2.4}.a is precisely such that the 
$\mc A$-weight $rt$ of $\pi(P,\sigma,t)$ survives to the Langlands quotient. Since the group 
$\Gamma$ preserves $T^-$, its presence does not affect temperedness.

Now assume that $P = F_0$. Since $T^{F^0++} \subset T^{F_0}_{rs}$ and 
$T^- \cap T^{F_0}_{rs} = \{1\}$, this representation can only be tempered if $|t| =1$. In that 
case $\sigma$ and $\pi^\Gamma(P,\sigma,t,\rho)$ have the same absolute values of $\mc A$-weights, 
modulo $\Gamma$. But $\Gamma T^- = T^-$, so the temperedness of $\pi^\Gamma(P,\sigma,t,\rho)$ 
and $L^\Gamma(P,\sigma,t,\rho)$ follows from that of $\sigma. \qquad \Box$
\\[2mm]

For connected reductive $p$-adic groups the Langlands quotient always appears with multiplicity 
one in the standard representation of which it is a quotient. Although not stated explicitly
in most sources, that is already part of the proof, see \cite{Kon} or \cite[Theorem 2.15]{SolPadic}.
This also holds for reductive $p$-adic groups with a cyclic component group \cite{BaJa2}.

Closer examination of the proof of Theorem \ref{thm:2.4} allows us to generalize and improve 
upon this in our setting. Let $W (R_P) r_\sigma \in T_P / W (R_P)$ be the central 
character of $\sigma$. Then $|r_\sigma| \in T_{P,rs} = \exp (\mf a_P)$, so we can define
\begin{equation}\label{eq:ccdelta}
cc_P (\sigma) := W (R_P) \log |r_\sigma| \in \mf a_P / W (R_P) .
\end{equation}
Since the inner product on $\mf a$ is $\WG$-invariant, the number $\norm{cc_P (\sigma)}$ 
is well-defined.

\begin{lem}\label{lem:2.10}
Let $(P,\sigma,t,\rho)$ and $(P,\sigma',t,\rho')$ be Langlands data for $\mc H \rtimes \Gamma$.
\enuma{
\item The functor $\mr{Ind}_{\mc H^P \rtimes \Gamma_P}^{\mc H \rtimes \Gamma}$ induces an isomorphism
\begin{multline*}
\mr{Hom}_{\mc H^P \rtimes \Gamma_P} (\pi^{P,\Gamma_P} 
(P,\sigma,t,\rho), \pi^{P,\Gamma_P} (P,\sigma',t,\rho')) \cong \\
\mr{Hom}_{\mc H \rtimes \Gamma} (\pi^\Gamma (P,\sigma,t,\rho), \pi^\Gamma (P,\sigma',t,\rho')) .
\end{multline*}
These spaces are one-dimensional if $(\sigma,t,\rho)$ and $(\sigma',t,\rho')$ are $\Gamma_P$-conjugate, and
zero otherwise.
\item Suppose that $L^\Gamma (Q,\tau,s,\nu)$ is a constituent of $\pi^\Gamma (P,\sigma,t,\rho)$, but
not $L^\Gamma (P,\sigma,t,\rho)$. Then $P \subset Q$ and $\norm{cc_P (\sigma)} < \norm{cc_Q (\tau)}$.
}
\end{lem}
\emph{Proof.}
(a) We use the notation from \eqref{eq:LanglandsF}.
For any weight $s$ of $\sigma \circ \phi_t$ we have $\log |s t^{-1}| \in \mf a_P^-$ and 
$(\log |s| )_0 = (\log |t|)_{F_0}$, where the subscript $F_0$ refers to the decomposition 
of elements of $\mf t$ with respect to $\mf t = \mf t_{F_0} \oplus \mf t^{F_0}$.
Let $s'$ be a weight of $\sigma' \circ \phi_t$. By \eqref{eq:2.3}
\begin{equation}
(w \log |s'| )_0 < (\log |s'| )_0 = (\log |t|)_{F_0} \qquad
\forall w \in W^P \setminus \{ \mathrm{id} \} ,
\end{equation}
with respect to the dominance order on $\mf a_{F_0}^*$. Since $(\gamma \lambda)_0 = \gamma (\lambda_0)$ 
for all $\lambda \in \mf a$ and $\gamma \in \Gamma$, we get
\begin{equation}\label{eq:2.4}
\norm{(\gamma w \log |s'| )_0} < \norm{(\log |t|)_{F_0}} 
\qquad \forall \gamma \in \Gamma, w \in W^P \setminus \{ \mathrm{id} \} .
\end{equation}
In particular $\gamma w (s')$ with $w \in W^P$ can only equal the weight $s$ of $\sigma \circ \phi_t$ 
if $w = 1$. Let $v_s \in V_{\sigma \otimes \rho}$ be a nonzero weight vector. Since 
$\pi^{P,\Gamma_P} (P,\sigma,t,\rho)$ is an irreducible $\mc H^P \rtimes \Gamma_P$-representation,
$1 \otimes v_s \in (\mc H \rtimes \Gamma) \otimes_{\mc H^P \rtimes \Gamma_{P,\sigma,t}} 
V_{\sigma \otimes \rho}$ is cyclic for $\pi^\Gamma (P,\sigma ,t,\rho )$. Therefore the map
\begin{equation}\label{eq:2.5}
\mr{Hom}_{\mc H \rtimes \Gamma} (\pi^\Gamma (P,\sigma,t,\rho) ,\pi^\Gamma (P,\sigma',t,\rho')) \to 
\pi^\Gamma (P,\sigma',t,\rho') : f \mapsto f (1 \otimes v_s )
\end{equation}
is injective. By \eqref{eq:2.4} the $s$-weight space of $\pi^\Gamma (P,\sigma',t,\rho')$
is contained in $1 \otimes V_{\sigma' \otimes \rho'}$.  
So $f (1 \otimes v_s ) \in 1 \otimes V_{\sigma' \otimes \rho'}$ and multiplying by
$\mc H^P \rtimes \Gamma_P$ yields
\[
f (\C [\Gamma_P] \otimes_{\Gamma_{P,\sigma,t}} V_{\sigma \otimes \rho}) \subset 
\mh C [\Gamma_P] \otimes_{\Gamma_{P,\sigma',t}} V_{\sigma' \otimes \rho'} .
\]
Thus any $f \in \mr{Hom}_{\mc H \rtimes \Gamma} (\pi^\Gamma (P,\sigma,t,\rho) 
,\pi^\Gamma (P,\sigma',t,\rho'))$ lies in 
\begin{equation}\label{eq:2.6}
\mr{Ind}_{\mc H^P \rtimes \Gamma_P}^{\mc H \rtimes \Gamma} 
\Big( \mr{Hom}_{\mc H^P \rtimes \Gamma_P} 
(\pi^{P,\Gamma_P} (P,\sigma,t,\rho), \pi^{P,\Gamma_P} (P,\sigma',t,\rho')) \Big).
\end{equation}
From \eqref{eq:2.5} we see that this induction functor is injective on homomorphisms.
The modules in \eqref{eq:2.6} are irreducible, so the dimension of \eqref{eq:2.6} is zero or one.
By Corollary \ref{cor:2.8}.b it is nonzero if and only if $(\sigma,t,\rho)$ and $(\sigma',t,\rho')$
are $\Gamma_P$-conjugate. \\
(b) The proofs of Theorem \ref{thm:2.4}.a and Corollary \ref{cor:2.8}.a show that 
$L^\Gamma (P,\sigma,t,\rho)$ is the unique irreducible subquotient of $\pi^\Gamma (P,\sigma,t,\rho)$ 
which has an $\mc A$-weight $t_L$ with $(\log |t_L| )_0 = (\log |t|)_{F_0}$. Moreover of all 
$\mc A$-weights $s'$ of proper submodules of $\pi^\Gamma (P,\sigma,t,\rho)$ satisfy 
$(\log |s'| )_0 < (\log |t|)_{F_0}$, with the notation of \eqref{eq:LanglandsF}. In particular, for the
subquotient $L^\Gamma (Q, \tau,s,\nu )$ of $\pi^\Gamma (P,\sigma,t,\rho)$ we find that 
\[
(\log |s| )_{F_0} = (\log |s'| )_0 < (\log |t| )_{F_0} .
\]
Since $\log |s| \in \mf a^{Q++}$ and $\log |t| \in \mf a^{P++}$, this implies $P \subset Q$ and
\begin{equation}\label{eq:2.7}
\norm{(\log |s| )_{F_0}}  < \norm{(\log |t| )_{F_0}} .
\end{equation}
According to Lemma \ref{lem:2.12} all constituents of $\pi^\Gamma (P,\sigma,t,\rho)$ have central 
character $\WG( r_\sigma t ) \in T / \WG$. 
The same goes for $(Q,\tau,s,\nu )$, so $r_\sigma t$ and $r_\tau s$ lie in the same 
$W_0 \rtimes \Gamma$-orbit. Thus also
\[
\WG (\log |r_\sigma t| )_{F_0}  = \WG (\log |r_\tau s| )_{F_0} .
\]
By definition $(\log |t|)_{F_0} \perp \mf t_P$ and $(\log |s|)_{F_0} \perp \mf t_Q$, so
\begin{multline}
\| \Re (cc_P (\sigma )) \|^2 + \| (\log |t|)_{F_0} \|^2 = \| (\log |r_\sigma t| )_{F_0} \|^2 \\ 
= \| (\log |r_\tau s| )_{F_0} \|^2 = \| \Re (cc_Q (\tau )) \|^2 + \| (\log |s|)_{F_0}  \|^2 .
\end{multline}
Finally we use \eqref{eq:2.7}. $\qquad \Box$
\\[1mm]

\section{From $\mc H$-representations to $W$-representations}
\label{sec:Springer}

Given an algebra or group $A$, let Irr$(A)$ be the collection of (equivalence classes of)
irreducible complex $A$-representations. Let $G_\Z (A) = G_\Z (\mr{Mod}_f (A))$ denote the 
Grothendieck group of the category of finite length complex representations of $A$, and write 
$G_{\mathbb F} (A) = G_\Z (A) \otimes_\Z \mathbb F$ for any field $\mathbb F$.

The classical Springer correspondence \cite{Spr} realizes all irreducible representation of a 
finite reflection group $W_0$ in the top cohomology of the associated flag variety. 
Kazhdan--Lusztig theory (see \cite{KaLu,Xi}) allows one to interpret this as a bijection between 
Irr$(W_0)$ and a certain collection of irreducible representations of an affine Hecke algebra with 
equal parameters. As such, the finite Springer correspondence is a specialisation of an affine 
Springer correspondence between Irr$(W (\mc R))$ and Irr$(\mc H (\mc R, q))$, see 
\cite[Section 8]{Lus-Rep}. The proof of Kazhdan and Lusztig requires that $\mc R$ is of simply 
connected type and that $q$ is an equal parameter function whose value is either 1 or not
a root of unity. Reeder \cite[Theorem 3.5.4]{Ree2} showed that the resulting parametrization 
of irreducible $\mc H (\mc R ,q)$-modules remains valid without simple connectedness.
We will prove an analogue of this result for all extended affine Hecke algebras with
unequal (but positive) parameters.

For any $\mh H \rtimes \Gamma$-representation $\pi$, let $\pi \big|_{W_0 \rtimes \Gamma} = 
\pi \big|_{\WG}$ be the restriction of $\pi$ to the subalgebra 
$\C [\WG] = \C [W_0 \rtimes \Gamma] \subset \mh H \rtimes \Gamma$. 
Let $\mr{Irr}_0 (\mh H \rtimes \Gamma)$ be the collection of (equivalence classes of) 
irreducible tempered $\mh H \rtimes \Gamma$-representations with real central character. 
In \cite[Theorem 6.5.c]{SolHomGHA} the author proved that the set 
\begin{equation}\label{eq:Irr0}
\{ \pi \big|_{W_0 \rtimes \Gamma} : \pi \in \mr{Irr}_0 (\mh H \rtimes \Gamma) \}
\end{equation}
is a $\Q$-basis of $G_\Q (W_0 \rtimes \Gamma)$. When $\Gamma$ is trivial, this is a kind of Springer 
correspondence for finite Weyl groups. The only problem is that $\pi \big|_{W_0}$ may be reducible,
but that could be solved by picking a suitable (a priori not canonical) irreducible subrepresentation
of $\pi \big|_{W_0}$.

In fact it is known from \cite[Corollary 3.6]{Ciu} that in many cases the matrix that expresses 
\eqref{eq:Irr0} in terms of irreducible $W_0 \rtimes \Gamma$-representations is unipotent and upper 
triangular (with respect to a suitable ordering). That would provide a natural Springer correspondence
for graded Hecke algebras with arbitrary parameters, but unfortunately that improvement is still open 
in our generality.

\begin{thm}\label{thm:2.7}
There exists a unique system of maps 
\[
\Spr : \mr{Irr}(\mc H \rtimes \Gamma) \to \mr{Mod}_f (X \rtimes \WG) ,
\]
for all extended affine Hecke algebras $\mc H \rtimes \Gamma$ with positive parameters, such that:
\enuma{
\item The image of $\Spr$ is a $\Q$-basis of $G_\Q (X \rtimes \WG)$.
\item $\Spr$ preserves the unitary part of the central character.
\item $\Spr (\pi)$ is tempered if and only if $\pi$ is tempered.
\item Let $u \in T_{un}$, let $\tilde \pi \in \mr{Irr}_0 (\mh H (\tilde{\mc R}_u ,k_u) 
\rtimes W'_{F_u ,u})$ and let $\tilde \pi \circ \Phi_u$ be the $\mc H (\mc R_u ,q_u) 
\rtimes W'_{F_u ,u}$-representation associated to it via Theorem \ref{thm:2.2}.b. Then
\[
\Spr \big( \mr{Ind}_{\mc H (\mc R_u ,q_u) \rtimes W'_{F_u ,u}}^{\mc H \rtimes \Gamma} 
(\tilde \pi \circ \Phi_u) \big) = \mr{Ind}_{X \rtimes W'_u}^{X \rtimes \WG}
\big( \C_u \otimes \tilde \pi \big|_{W'_u} \big) ,
\]
where $\C_u$ denotes the one-dimensional $X$-representation with character $u$.
\item If $(P,\sigma,t,\rho)$ is a Langlands datum for $\mc H \rtimes \Gamma$, then
\[
\Spr (L^\Gamma (P,\sigma,t,\rho)) = \mr{Ind}_{X \rtimes (W (R_P) \rtimes \Gamma_{P,\sigma,t})}^{X 
\rtimes \WG} (\Spr (\sigma \otimes \rho) \circ \phi_t) .
\]
}
\end{thm}
\emph{Proof.}
In view of Corollary \ref{cor:2.3}, properties (b) and (d) determine $\Spr$ uniquely for all
irreducible tempered representations. A glance at Lemma \ref{lem:2.9} shows that $\Spr$ 
preserves temperedness.

Next Corollary \ref{cor:2.8}.b and property (e) determine $\Spr$ for all irreducible 
representations. By Corollary \ref{cor:2.8} every nontempered irreducible 
$\mc H \rtimes \Gamma$-representation $\pi$ is of the form $L(P,\sigma,t,\rho)$ for some Langlands
datum with $t \not\in T_{un}$. By construction all $\mc A$-weights of $\Spr (\sigma \otimes \rho)$ lie in 
$T_{un}$, so by property (e) the absolute values of $\mc A$-weights of $\Spr (\pi)$ lie in 
$\WG |t|$. Together with Lemma \ref{lem:2.9} this shows that $\Spr (\pi)$ is not
tempered. 

Now we have (b)--(e), let us turn to (a). By Corollary \ref{cor:2.3} and the result mentioned in
\eqref{eq:Irr0}, (a) holds if we restrict to tempered representations on both sides. The proof that
this restriction is unnecessary is more difficult, we postpone it to Section \ref{sec:dual}.

We will see in Corollary \ref{cor:Sprsigma0} that on tempered representations $\Spr$ is given by composition 
with an algebra homomorphism between suitable completions of $\C [X \rtimes (W_0 \rtimes \Gamma)]$ 
and $\mc H \rtimes \Gamma$. That is not possible for all irreducible representations, since sometimes
$\Spr$ does not preserve the dimensions of representations. More precisely, $\Spr$ preserves the
dimension of an irreducible tempered representation, since property (d) does so. By contrast,
the right hand side of property (e) has the same dimension as $\pi (P,\sigma,t,\rho)$. Thus $\Spr$ 
preserves the dimension of $L (P,\sigma,t,\rho)$ if this Langlands quotient equals $\pi (P,\sigma,t,\rho)$, 
and increases the dimension otherwise.

Since the principal series representations $M(t) = \mr{Ind}_{\mc A}^{\mc H \rtimes \Gamma} \C_t$
and $\mr{Ind}_{X}^{X \rtimes \WG} \C_t$ are irreducible for all $t$ in a Zariski-open dense subset of 
$T, \Spr$ preserves irreducibility on a large part of Irr$(\mc H \rtimes \Gamma)$. To make a nice
affine Springer correspondence out of $\Spr$, one would have to modify it so that it always
respects irreducibilty. This would boil down to refining \eqref{eq:Irr0} (preferably in a canonical way)
to map that preserves irreducibilty.

For some applications the following variation on $\Spr$ is more convenient:

\begin{cor}\label{cor:2.6}
Theorem \ref{thm:2.7} also holds with condition (e) replaced by
\[
(e^\vee) \qquad \zeta^\vee (\pi^\Gamma (P,\sigma,t,\rho)) = \mr{Ind}_{X \rtimes (W (R_P) \rtimes 
\Gamma_{P,\sigma,t})}^{X \rtimes \WG} (\zeta^\vee (\sigma \otimes \rho) \circ \phi_t) .
\] 
The resulting map 
\[
\zeta^\vee : G_\Z (\mc H \rtimes \Gamma) \to G_\Z (X \rtimes \WG) 
\]
commutes with parabolic induction.
\end{cor}
\emph{Proof.} This follows from Theorem \ref{thm:2.7} and Lemma \ref{lem:2.10}.b.
$\qquad \Box$ \\[1mm]

The disadvantage of $\zeta^\vee$ compared to $\Spr$ is that it sends some irreducible 
$\mc H \rtimes \Gamma$-representations to virtual $W \rtimes \Gamma$-representations. 
Suppose for example that a principal series representation $M(t)$ with $t \in T^{++}$ has
only one quotient $\pi$ and only one subrepresentation $\delta$, which is tempered. 
(This occurs already for $\mc R$ of type $A_1^{(1)}$.) Then 
\[
\zeta^\vee (\pi) = \zeta^\vee (M(t)) - \zeta^\vee (\delta) = 
\mr{Ind}_X^{X \rtimes \WG} \C_t - \zeta^\vee (\delta) \; \in G_\Z (X \rtimes \WG) ,
\]
and the right hand side is no ordinary representation because $\zeta^\vee (\delta)$ is 
tempered and $\mr{Ind}_X^{X \rtimes \WG} \C_t$ has no $X$-weights in $T_{un}$.

Of course there also exist versions of Theorem \ref{thm:2.7} and Corollary \ref{cor:2.6} for graded 
Hecke algebras. They can easily be deduced from the above using Theorem \ref{thm:2.2}.b.

\chapter{Parabolically induced representations}

Parabolic induction is a standard tool to create a large supply of interesting representations of
reductive groups, Hecke algebras and related objects.
In line with Harish-Chandra's philosophy of the cusp form, every irreducible tempered representation
of an affine Hecke algebra can be obtained via unitary induction of a discrete series
representation of a parabolic subalgebra. With the Langlands classification we can also reach
irreducible representations that are not tempered. 

Hence we consider induction data $\xi = (P,\delta,t)$, where $\delta$ is a discrete series representation
of $\mc H_P$ and $t \in T^P$ is an induction parameter. With this we associate a representation 
$\pi (\xi) = \mr{Ind}_{\mc H^P}^{\mc H} (\delta \circ \phi_t)$. Among these are the principal 
series representations, which already exhaust the dual space of $\mc H$. But that is not very 
satisfactory, since a principal series representation can have many irreducible subquotients, 
and it is not so easy to determine them, see \cite{Ree1}.

Instead we are mostly interested in induction data $\xi$ for which $|t|$ is positive 
(in an appropriate sense) and in irreducible quotients of $\pi (\xi)$, because the 
Langlands classification applies to these. In Theorem \ref{thm:3.10} we construct, for every
irreducible $\mc H$-representation $\rho$, an essentially unique induction datum $\xi^+ (\rho)$,
such that $\rho$ is a quotient of $\pi (\xi^+ (\rho ))$. However, in general $\pi (\xi^+ (\rho))$
has more than one irreducible quotient.

Another important theme in this chapter are intertwining operators between induced representations 
of the form $\pi (\xi)$. Their definition and most important properties stem from the work of
Opdam and Delorme \cite{Opd-Sp,DeOp1}. Like in the setting of reductive groups, it is already 
nontrivial to show that normalized intertwining operators are regular on unitary induced 
representations. Under favorable circumstances such intertwining operators span $\mr{Hom}_{\mc H}
(\pi (\xi), \pi (\xi'))$. This was already known \cite{DeOp1} for unitary induction data $\xi,\xi'$,
in which case $\pi (\xi)$ and $\pi (\xi')$ are tempered representations. We generalize this to pairs
of positive induction data (Theorem \ref{thm:3.9}).
Crucial in all these considerations is the Schwartz algebra $\mc S$ of $\mc H$, the analogue of
the Harish-Chandra--Schwartz algebra of a reductive $p$-adic group.

For the geometry of the dual space of $\mc H$ it is important to understand the number $n(\xi)$
of irreducible $\mc H$-representations $\rho$ with $\xi^+ (\rho)$ equivalent to $\xi$. This is
governed by a groupoid $\mc G$ that keeps track of all intertwining operators. Indeed, if
$t \mapsto \xi_t$ is a continuous path of induction data such that all $\xi_t$ have the same
isotropy group in $\mc G$, then $n(\xi_t)$ is constant along this path (Proposition \ref{prop:3.11}).

From this we deduce that the dual of $\mc H$ is a kind of complexification of the tempered dual
of $\mc H$. As a topological space, the tempered dual is built from certain algebraic 
subvarieties of compact tori, each with a multiplicity. In this picture the dual of $\mc H$ 
is built from the corresponding complex subvarieties of complex algebraic tori, with the same
multiplicities. This geometric description is used to finish the proof of Theorem \ref{thm:2.7}.

\section{Unitary representations and intertwining operators} 

Like for Lie groups, the classification of the unitary dual of an affine Hecke algebra appears to be 
considerably more difficult than the classification of the full dual or of the tempered dual. 
This is an open problem that we will not discuss in this paper, cf. \cite{BaMo2,BaCi}. 
Nevertheless we will use unitarity arguments, mainly to show that certain representations are 
completely reducible. The algebra $\mc H \rtimes \Gamma$ is endowed with a sesquilinear 
involution * and a trace $\tau$, defined by
\begin{align}
\nonumber & (z N_w \gamma)^* = \bar{z} \gamma^{-1} N_{w^{-1}} 
\qquad z \in \C, w \in W, \gamma \in \Gamma , \\
\label{eq:*tau} & \tau (z N_w \gamma) = \left\{ \begin{array}{ll}
z & \text{if } \gamma = w = e, \\
0 & \text{otherwise} .
\end{array}
\right.
\end{align}
Since $q$ is real-valued, this * is anti-multiplicative and $\tau$ is positive. These give rise to an 
Hermitian inner product on $\mc H \rtimes \Gamma$:
\begin{equation}\label{eq:inptau}
\inp{h}{h'}_\tau = \tau(h^* h') \qquad h,h' \in \mc H \rtimes \Gamma.
\end{equation}
A short calculation using the multiplication rules \ref{eq:multrules} shows that the basis 
$\{ N_w \gamma : w \in W ,\gamma \in \Gamma \}$ of $\mc H \rtimes \Gamma$ 
is orthonormal for this inner product.

We note that $\Gamma$ acts on $\mc H$ by *-automorphisms, and that $\mc H, \mc H (W_0,q)$ 
and $\C [\Gamma]$ are *-subalgebras of $\mc H \rtimes \Gamma$. In general $\mc A$ is not 
a *-subalgebra of $\mc H$. For $x \in X$ \cite[Proposition 1.12]{Opd-Sp} tells us that 
\begin{equation}\label{eq:thetax*}
\theta_x^* = N_{w_0} \theta_{-w_0 (x)} N_{w_0}^{-1} ,
\end{equation}
where $w_0$ is the longest element of the Coxeter group $W_0$.

Let $\Gamma'_P$ be a subgroup of $\Gamma_P$ and let $\tau$ be a 
$\mc H^P \rtimes \Gamma'_P$-representation on an inner product space $V_\tau$. By default we 
will endow the vector space
\[
\mc H \rtimes \Gamma \otimes_{\mc H^P \rtimes \Gamma'_P} V_\tau 
\cong \C [\Gamma W^P] \otimes_{\C [\Gamma'_P]} V_\tau 
\]
with the inner product
\begin{equation}\label{eq:inpind}
\inp{h \otimes v}{h' \otimes v'} = \tau (h^* h') \inp{v}{v'} 
\qquad h,h' \in \C [\Gamma W^P] , v,v' \in V_\tau .
\end{equation}
Recall that a representation $\pi$ of $\mc H \rtimes \Gamma$ on a Hilbert space is unitary if 
$\pi (h^*)$ is the adjoint operator $\pi (h)^*$  of $\pi (h)$, for all $h \in \mc H \rtimes \Gamma$. 
In particular, such representations are completely reducible. 

\begin{lem}\label{lem:3.1}
Let $\Gamma'_P$ be a subgroup of $\Gamma_P$, let $\sigma$ be a finite dimensional 
$\mc H_P \rtimes \Gamma'_P$-representation and let $t \in T^P$.
\enuma{
\item If $\sigma$ is unitary and $t \in T^P_{un}$, then 
$\mr{Ind}_{\mc H^P \rtimes \Gamma'_P}^{\mc H \rtimes \Gamma} (\sigma \circ \phi_t)$ is unitary 
with respect to the inner product \eqref{eq:inpind}.
\item $\mr{Ind}_{\mc H^P \rtimes \Gamma'_P}^{\mc H \rtimes \Gamma} (\sigma \circ \phi_t)$ 
is (anti-)tempered if and only if $\sigma$ is (anti-)tempered and $t \in T^P_{un}$.
}
\end{lem}
\emph{Proof.} 
Since $\Gamma$ acts by *-algebra automorphisms and $\Gamma \cdot T^- = T^-$, 
it does not disturb the properties unitarity and temperedness. Hence it suffices to prove 
the lemma in the case $\Gamma = \Gamma'_P = $\{id\}. Then (a) and the "if"-part of (b) 
are \cite[Propositions 4.19 and 4.20]{Opd-Sp}.

For the "only if"-part of (b), suppose that $t \in T^P \setminus T^P_u$. Since $X \cap (P^\vee)^\perp$ is of 
finite index in $X^P = X / (X \cap \Q P)$, there exists $x \in X \cap (P^\vee)^\perp$ with $|t(x)| \neq 1$.
Possibly replacing $x$ by $-x$, we may assume that $|t(x)| > 1$. But
$\delta(\phi_t (\theta_x)) (v) = x (t) v$ for all $v \in V_\delta$ and
$x \in Z( X \rtimes W_P)$, so the $\mc H^P$-representation $\delta \circ \phi_t$ is not tempered. 
Hence its induction to $\mc H$ cannot be tempered.
Similarly, if $\sigma$ is not tempered, then the restriction of $\sigma \circ \phi_t$ to
$\mc H (X \cap \Q P, R_P, Y / Y \cap P^\perp, R_P^\vee, P,q_P)$ is not tempered. 

The same proof works in the anti-tempered case, we only have to replace $|t(x)| > 1$ by 
$|t(x)| < 1. \qquad \Box$
\\[2mm]
\emph{Remark.} It is possible that $\mr{Ind}_{\mc H^P \rtimes \Gamma'_P}^{\mc H \rtimes \Gamma} 
(\sigma \circ \phi_t)$ is unitary with respect to some inner product other than \eqref{eq:inpind},
if the conditions of part (a) are not met.\\[2mm]

We intend to partition Irr$(\mc H \rtimes \Gamma)$ into finite packets, each of which is obtained
by inducing a discrete series representation of a parabolic subalgebra of $\mc H$.
Thus our induction data are triples $(P,\delta,t)$, where
\begin{itemize}
\item $P \subset F_0$;
\item $(\delta ,V_\delta)$ is a discrete series representation of $\mc H_P$;
\item $t \in T^P$.
\end{itemize} 
Let $\Xi$ be the space of such induction data, where we regard $\delta$ only modulo
equivalence of $\mc H_P$-representations. We say that $\xi = (P,\delta,t)$ is unitary
if $t \in T^P_{un}$, and we denote the space of unitary induction data by $\Xi_{un}$. Similarly
we say that $\xi$ is positive if $|t| \in T^{P+}$, which we write as $\xi \in \Xi^+$. 
Notice that, in contrast to Langlands data, we do not require $|t|$ to be strictly positive.
We have three collections of induction data:
\begin{equation}\label{eq:inductionData}
\Xi_{un} \subseteq \Xi^+ \subseteq \Xi .
\end{equation}
By default we endow these spaces with the topology for which $P$ and $\delta$ are
discrete variables and $T^P$ carries its natural analytic topology.
We will realize every irreducible $\mc H \rtimes \Gamma$-representation as a quotient of a
well-chosen induced representation 
\[
\pi^\Gamma (\xi) := \mr{Ind}_{\mc H}^{\mc H \rtimes \Gamma} \pi (P,\delta,t) = 
\mr{Ind}_{\mc H^P}^{\mc H \rtimes \Gamma} (\delta \circ \phi_t ) .
\] 
We note that for all $s \in T^{W_0 \rtimes \Gamma}$:
\begin{equation}
\pi^\Gamma (P,\delta,ts) = \pi^\Gamma (P,\delta,t) \circ \phi_s .
\end{equation}
As vector space underlying $\pi^\Gamma (\xi)$ we will always take $\C [\Gamma W^P] \otimes
V_\delta$. This space does not depend on $t$, which will allow us to speak of maps that
are continuous, smooth, polynomial or even rational in the parameter $t \in T^P$.

The discrete series representations of affine Hecke algebras with irreducible root data were
classified in \cite{OpSo2}. Here we recall only how their central characters can be determined,
which is related to the singularities of the elements $\imath^0_s$. Consider the $W_0 \rtimes 
\Gamma$-invariant rational function
\[
\eta = \prod\nolimits_{\alpha \in R_0} c_\alpha^{-1} \in \C (T),
\]
where $c_\alpha$ is as in \eqref{eq:calpha}. Notice that $\eta$ depends on the parameter function
$q$, or more precisely on $q^{1/2}$. A coset $L$ of a subtorus of $T$ is said to be 
residual if the pole order of $\eta$ along $L$ equals $\dim_\C (T) - \dim_C (L)$, see \cite{Opd4}.
A residual coset of dimension 0 is also called a residual point. Such points can exist only if
$\mc R$ is semisimple, otherwise all residual cosets have dimension at least rank $Z(W) > 0$.

According to \cite[Lemma 3.31]{Opd-Sp} the collection of central characters of discrete series
representations of $\mc H (\mc R,q)$ is exactly the set of $W_0$-orbits of residual points
for $(\mc R ,q)$. Moreover, if $\delta$ is a discrete series representation of $\mc H_P$ 
with central character $W (R_P) r$, then $r T^P$ is a residual coset for $(\mc R ,q)$
\cite[Proposition 7.4]{Opd-Sp}. Up to multiplication by an element of $W_0$, every residual coset
is of this form \cite[Proposition 7.3.v]{Opd-Sp}.

The map that assigns to $\xi \in \Xi$ the central character of $\pi (\xi) \in \mr{Mod}_f 
(\mc H (\mc R,q))$ is an algebraic morphism $\Xi \to T / W_0$. The above implies that the image of 
\begin{equation}\label{eq:rescosd}
\{ (P,\delta,t) \in \Xi : |F_0 \setminus P| = d \}
\end{equation}
is the union of the $d$-dimensional residual cosets, modulo $W_0$.

Let $\delta_\emptyset$ be the unique onedimensional representation of $\mc H_\emptyset = \C$ 
and consider 
\[
M(t) := \pi^\Gamma (\emptyset, \delta_\emptyset,t) =
\mr{Ind}_{\mc A}^{\mc H \rtimes \Gamma} (\delta_\emptyset \circ \phi_t) \cong
\mr{Ind}_{\mc O (T)}^{\mc H \rtimes \Gamma} \C_t , \qquad t \in T.
\]
The family consisting of these representations is called the principal series of
$\mc H \rtimes \Gamma$ and is of considerable interest. For example, by Frobenius reciprocity
every irreducible $\mc H \rtimes \Gamma$-representation is a quotient of some principal series
representation.

\begin{lem}\label{lem:3.13}
Suppose that $h \in \mc H \rtimes \Gamma$ and that $M(t,h) = 0$
for all $t$ in some Zariski-dense subset of $T$. Then $h = 0$. 
\end{lem}
\emph{Proof.}
Since $M(t,h) \in \mr{End}_\C (\C [\Gamma W_0] )$ depends algebraically on $t$, it is zero for all 
$t \in T$. Write $h = \sum_{\gamma w \in \Gamma \rtimes W_0} a_{\gamma w} N_{\gamma w}$ with 
$a_{\gamma w} \in \mc A$ and suppose that $h \neq 0$. Then we can
find $w' \in W_0$ that $a_{\gamma w'} \neq 0$ for some $\gamma \in \Gamma$, and such that $\ell (w')$
is maximal for this property. From Theorem \ref{thm:1.1}.d we see that 
\[
M(t,h) (N_e) = \sum_{\gamma w \in \Gamma \rtimes W_0} b_{\gamma w} N_{\gamma w}
\]
for some $b_{\gamma w}$ with $b_{\gamma w'} = a_{\gamma w'} \neq 0$. Therefore 
$M(t,h)$ is not identically zero. This contradiction shows that
the assumption $h \neq 0$ is untenable. $\qquad \Box$
\\[2mm]

By \cite[Corollary 2.23]{Opd-Sp} discrete series representations are unitary. 
(Although Opdam only worked in the setting $\Gamma =$ \{id\}, his proof also applies with 
general $\Gamma$.) From this and Lemma \ref{lem:3.1} we observe:

\begin{cor}\label{cor:3.2}
Let $\xi = (P,\delta,t) \in \Xi$. If $t \in T^P_{un}$, then $\pi^\Gamma (\xi)$ is unitary and tempered.
If $t \in T^P \setminus T^P_{un}$, then $\pi^\Gamma (\xi)$ is not tempered.
\end{cor}

For any subset $Q \subset F_0$, let $\Xi^Q, \pi^Q, \ldots$ denote the things $\Xi ,\pi, \ldots$, but for the
algebra $\mc H^Q$ instead of $\mc H$. For $\xi = (P,\delta,t) \in \Xi$ we define 
\begin{equation}
P(\xi) := \{ \alpha \in R_0 : |\alpha (t)| = 1 \} .
\end{equation}

\begin{prop}\label{prop:3.3}
Let $\xi = (P,\delta,t) \in \Xi^+$.
\enuma{
\item The $\mc H^{P(\xi)} \rtimes \Gamma_{P(\xi)}$-representation $\pi^{P(\xi) ,\Gamma_{P(\xi)}}(\xi)$
is completely reducible.
\item Every irreducible summand of $\pi^{P(\xi) ,\Gamma_{P(\xi)}}(\xi)$ is of the form
$\pi^{P(\xi),\Gamma_{P(\xi)}} (P(\xi), \sigma, t^{P(\xi)}, \rho)$, where $(P(\xi), \sigma, t^{P(\xi)}, \rho)$ 
is a Langlands datum for $\mc H \rtimes \Gamma$ and $t^{P(\xi)} t^{-1} \in T_{P(\xi)}$.
\item The irreducible quotients of $\pi^\Gamma (\xi)$ are the representations 
$L^\Gamma (P(\xi), \sigma, t^{P(\xi)}, \rho)$, with $(P(\xi), \sigma, t^{P(\xi)}, \rho)$ coming from (b).
\item Every irreducible $\mc H \rtimes \Gamma$-representation is of the form described in (c).
\item The functor $\mr{Ind}_{\mc H^{P(\xi)} \rtimes \Gamma_{P(\xi)}}^{\mc H \rtimes \Gamma}$
induces an isomorphism
\[
\mr{End}_{\mc H^{P(\xi)} \rtimes \Gamma_{P(\xi)}} \big( \pi^{P(\xi) ,\Gamma_{P(\xi)}}(\xi) \big) 
\cong \mr{End}_{\mc H \rtimes \Gamma} (\pi (\xi)) .
\]
}
\end{prop}
\emph{Remarks.} Part (a) holds for any $\xi \in \Xi$. 
In (b) $t^{P(\xi)}$ is uniquely determined modulo $K_{P(\xi)}$. \\
\emph{Proof.}
(a) By construction there exists $t^{P(\xi)} \in T^{P(\xi)}$ such that 
\begin{equation}\label{eq:ttPxi}
t (t^{P(\xi)})^{-1} \in T_{P(\xi),un} . 
\end{equation}
Then $\pi^{P(\xi)} (P,\delta,t) \circ \phi_{t^{P(\xi)}}^{-1} = 
\pi^{P(\xi)} \big( P,\delta,t (t^{P(\xi)})^{-1} \big)$ 
is unitary by Corollary \ref{cor:3.2}. In particular it is completely reducible, which implies that 
$\pi^{P(\xi)} \big( P,\delta,t (t^{P(\xi)})^{-1}\big)$ is also completely reducible. 
By \cite[Theorem A.1.c]{SolGHA} 
\begin{equation}\label{eq:piPxi}
\pi^{P(\xi) ,\Gamma_{P(\xi)}}(\xi) = \pi^{P(\xi),\Gamma_{P(\xi)}} (P,\delta,t) = 
\mr{Ind}_{\mc H^{P(\xi)}}^{\mc H^{P(\xi)} \rtimes \Gamma_{P(\xi)}} \pi^{P(\xi)} (P,\delta,t)
\end{equation}
remains completely irreducible.\\
(b) By Corollary \ref{cor:3.2} $\pi^{P(\xi)} \big( P,\delta,t (t^{P(\xi)})^{-1}\big)$ is tempered and 
unitary, so by Lemma \ref{lem:2.10} all its irreducible summands are of the form 
\[
\pi^{P(\xi)} (P(\xi),\sigma,k) = L^{P(\xi)} (P(\xi),\sigma,k) \text{, where } k \in T^{P(\xi)}_{un} . 
\]
Moreover $\pi^{P(\xi)} \big( P,\delta,t (t^{P(\xi)})^{-1}\big) \big|_{\C [X \cap (P(\xi)^\vee)^\perp]}$
consists only of copies of the trivial $X \cap (P(\xi)^\vee)^\perp$-representation, so
$k \in K_{P(\xi)} = T^{P(\xi)}_{un} \cap T_{P(\xi),un}$. Together with \eqref{eq:ttPxi} this
implies $k t^{P(\xi)} t^{-1} \in T_{P(\xi),un}$. Hence every irreducible summand of 
\eqref{eq:piPxi} is an irreducible summand of some 
\[
\mr{Ind}_{\mc H^{P(\xi)}}^{\mc H^{P(\xi)} \rtimes \Gamma_{P(\xi)}} 
\pi^{P(\xi)} (P(\xi),\sigma,k t^{P(\xi)}) . 
\]
By Clifford theory (see the proof of Corollary \ref{cor:2.8}) these are of the required form  
$\pi^{P(\xi),\Gamma_{P(\xi)}} (P(\xi),\sigma,k t^{P(\xi)},\rho)$.\\
(c) Follows immediately from (b) and Corollary \ref{cor:2.8}.\\
(d) By Corollary \ref{cor:2.8}.b it suffices to show that every $\mc H^P \rtimes \Gamma_P$-representation
of the form $(\sigma \circ \phi_t) \otimes \rho$ is a direct summand of some 
$\pi^{P,\Gamma_{P,\sigma,t}} (\xi^+)$. Without loss of generality we may assume that $P = F_0$
and that $\Gamma_{P,\sigma,t} = \Gamma$. The $\mc H$-representation $\sigma \circ \phi_{t |t|^{-1}}$ 
is irreducible and tempered, so by \cite[Theorem 3.22]{DeOp1} it is a direct summand of $\pi (\xi')$
for some $\xi' = (P',\rho',t') \in \Xi_{un}$. Then $(P',\rho',t' |t|) \in \Xi^+$ and $\sigma \circ \phi_t$ 
is a direct summand of $\pi (P',\rho',t' |t|)$. By Clifford theory \cite[Theorem A.1.b]{SolGHA} the 
$\mc H \rtimes \Gamma$-representation  $(\sigma \circ \phi_t) \otimes \rho$ is a direct summand 
of $\pi^\Gamma (P',\delta',t' |t|)$. \\
(e) Follows from (a), (b) and Lemma \ref{lem:2.10}.b. $\qquad \Box$
\\[3mm]

The parabolically induced representations $\pi^\Gamma (\xi)$ are by no means all disjoint.
The relations among them are described by certain intertwining operators, whose construction
we recall from \cite{Opd-Tr,Opd-Sp}.

Suppose that $P,Q \subset F_0 , u \in K_P, g \in \Gamma \ltimes W_0$ and $g (P) = Q$.
Let $\delta$ and $\sigma$ be discrete series representations of respectively $\mc H_P$ and
$\mc H_Q$, such that $\sigma$ is equivalent with $\delta \circ \psi_u^{-1} \circ \psi_g^{-1}$.
Choose a unitary map $I_\delta^{g u} : V_\delta \to V_\sigma$ such that
\begin{equation}\label{eq:Idelta}
I_\delta^{g u} (\delta (h) v) = \sigma (\psi_g \circ \psi_u (h)) (I_\delta^{g u} (v))
\qquad \forall v \in V_\delta, h \in \mc H_P .
\end{equation}
Notice that any two choices of $I_\delta^{g u}$ differ only by a complex number of norm 1.
In particular $I_\delta^{g u}$ is a scalar if $\sigma = \delta$. 

We obtain a bijection
\begin{align}
& \nonumber I_{g u} : (\C (T/W_0) \otimes_{Z (\mc H)} \mc H) \rtimes \Gamma \otimes_{\mc H^P} 
V_\delta \to (\C (T/W_0) \otimes_{Z (\mc H)} \mc H) \rtimes \Gamma \otimes_{\mc H^Q} V_\sigma ,\\
& \label{eq:defintop} I_{g u} (h \otimes v) = h \imath^o_{g^{-1}} \otimes I_\delta^{g u} (v) .
\end{align}

\begin{thm}\label{thm:3.5}
\enuma{
\item The map $I_{g u}$ defines an intertwining operator
\[
\pi^\Gamma (g u,P,\delta,t) : \pi^\Gamma (P,\delta,t) \to \pi^\Gamma (Q,\sigma,g (ut)) .
\]
As a map $\C [\Gamma W^P] \otimes_\C V_\delta \to \C [\Gamma W^Q] \otimes_\C V_\sigma$
it is rational in $t \in T^P$ and constant on $T^{F_0}$-cosets.
\item This map is regular and invertible on an open neighborhood of $T^P_{un}$ in $T^P$
(with respect to the analytic topology).
\item $\pi^\Gamma (g u,P,\delta,t)$ is unitary if $t \in T^P_{un}$.
}
\end{thm}
\emph{Remark.} Due to the freedom in the choice of \eqref{eq:Idelta}, for composable 
$g_1, g_2 \in \mc G$ the product $\pi (g_1, g_2 \xi) \pi (g_2,\xi)$ need not be equal to
$\pi (g_1 g_2,\xi)$. The difference is a locally constant function whose absolute value is 1 everywhere.

\emph{Proof.}
If $g = \gamma w \in \Gamma \ltimes W_0$, then $I_{g u} = I_\gamma \circ I_{w u}$,
modulo this locally constant function. It follows directly from the definitions that
the theorem holds for $I_\gamma$, so the difficult part is $I_{w u}$, which is dealt with
in \cite[Theorem 4.33 and Corollary 4.34]{Opd-Sp}. $\qquad \Box$
\\[2mm]

The intertwining operators for reflections acting on the unitary principal series can be 
made reasonably explicit:

\begin{lem}\label{lem:3.15}
Suppose that $\beta \in R_0$ and $t \in T_{un}$. Then $\pi^\Gamma (s_\beta, \emptyset, 
\delta_\emptyset, t)$ is a scalar operator if and only if $c_\beta^{-1} (t) = 0$.
\end{lem}
\emph{Proof.}
Suppose that $\alpha \in F_0, t \in T$ and $c_\alpha^{-1}(t) = 0$. Then \eqref{eq:imatho} implies
that $1 + \imath^0_{s_\alpha} (t) = 0$, regarded as an element of $\mc H (W_0,q)$. 
Hence $\pi^\Gamma (s_\alpha, \emptyset, \delta_\emptyset, t)$ is a scalar operator.
Conversely, if $c_\alpha^{-1}(t) \neq 0$, then \eqref{eq:imatho} shows that 
$1 + \imath^0_{s_\alpha} (t)$ is not scalar, because the action of $1 + q(s_\alpha)^{1/2} N_{s_\alpha}$ 
on $\mc H (W_0,q)$ has two different eigenvalues.

With Theorem \ref{thm:3.5} we can see that this is not specific for simple reflections. 
Find $w \in W_0$ such that $w (\beta) = \alpha$ is a simple root. Then $s_\beta = w^{-1} s_\alpha w$, 
so up to a nonzero scalar
\[
\pi^\Gamma (s_\beta ,\emptyset, \delta_\emptyset, t) = 
\pi^\Gamma (w^{-1} ,\emptyset, \delta_\emptyset, w t) \,
\pi^\Gamma (s_\alpha ,\emptyset, \delta_\emptyset, w t) \,
\pi^\Gamma (w ,\emptyset, \delta_\emptyset, t) .
\]
Now we notice that  $c_\beta^{-1} (t) = 0$ if and only if $c_\alpha^{-1}(w t) = 0$, and that 
$\pi^\Gamma (w^{-1} ,\emptyset, \delta_\emptyset, w t) = 
\pi^\Gamma (w ,\emptyset, \delta_\emptyset, t)^{-1}$ up to a scalar. $\qquad \Box$
\\[2mm]

Thus it is possible to determine the $\mc H$-endomorphisms for unitary principal series 
representations, at least when the isotropy groups of points $t \in T_{un}$ are generated by 
reflections. The reducibility and intertwining operators for nonunitary principal series are more 
complicated, and have been subjected to ample study \cite{Kat1,Rog,Ree1}. 
For other parabolically induced representations the intertwining operators are less explicit.
They can be understood better with the theory of R-groups \cite{DeOp2}. 

The action of these intertwining operators on the induction data space $\Xi$ is described
most conveniently with a groupoid $\mc G$ that includes all pairs $(g,u)$ as above.
The base space of $\mc G$ is the power set of $F_0$, and for
$P,Q \subseteq F_0$ the collection of arrows from $P$ to $Q$ is
\begin{equation}\label{eq:GPQ}
\mc G_{PQ} = \{ (g,u) : g \in \Gamma \ltimes W_0 , u \in K_P , g (P) = Q \} .
\end{equation}
Whenever it is defined, the multiplication in $\mc G$ is 
\[
(g',u') \cdot (g,u) = (g' g, g^{-1} (u') u) .
\]
Usually we will write elements of $\mc G$ simply as $gu$.
This groupoid acts from the left on $\Xi$ by
\[
(g,u) \cdot (P,\delta,t) := (g (P),\delta \circ \psi_u^{-1} \circ \psi_g^{-1},g (ut)) ,
\]
the action being defined if and only if $g (P) \subset F_0$.
Since $T^+ \supset T^{P+}$ is a fundamental domain for the action of $W_0$ on $T$,
every element of $\Xi$ is $\mc G$-associate to an element of $\Xi^+$.

Although $\pi^\Gamma (gu(P,\delta,t))$ and $\pi^\Gamma (P,\delta,t)$ are not always
isomorphic, the existence of rational intertwining operators has the following consequence:

\begin{lem}\label{lem:3.6}
The $\mc H \rtimes \Gamma$-representations $\pi^\Gamma (gu(P,\delta,t))$ and 
$\pi^\Gamma (P,\delta,t)$ have the same irreducible subquotients, counted with multiplicity.
\end{lem}
\emph{Proof.}
This is not hard, the proof in the graded Hecke algebra setting \cite[Lemma 3.4]{SolGHA} 
also works here. $\qquad \Box$
\\[1mm]

\section{The Schwartz algebra}
\label{sec:Schwartz}

We recall the construction of various topological completions of $\mc H$ \cite{Opd-Sp}:
a Hilbert space, a $C^*$-algebra and a Schwartz algebra. The latter is the most relevant
from the representation theoretic point of view. All tempered representations of $\mc H$
extend to its Schwartz completion, and a close study of this Schwartz algebra reveals 
facts about tempered representations for which no purely algebraic proof is known.

Let $L^2 (\mc R ,q)$ be Hilbert space completion of $\mc H$ with respect to the inner product
\eqref{eq:inptau}. By means of the orthonormal basis $\{ N_w : w \in W \}$ we can 
identify $L^2 (\mc R,q)$ with the Hilbert space $L^2 (W)$ of square integrable functions
$W \to \C$.

For any $h \in \mc H$ the map $\mc H (\mc R,q) \to \mc H (\mc R,q) : h' \mapsto h h'$ extends 
to a bounded linear operator on $L^2 (\mc R,q)$. This realizes $\mc H (\mc R,q)$ as a *-subalgebra
of $B (L^2 (\mc R,q))$. Its closure $C^* (\mc R,q)$ is a separable unital $C^*$-algebra, called
the $C^*$-algebra of $\mc H$.

The Schwartz completion of $\mc H$ will, as a topological vector space, consist of all rapidly 
decaying functions on $W$, with respect to some length function. For this purpose the length
function $\ell (w)$ of the Coxeter system $(W^\af ,S^\af)$ is unsatisfactory, because its natural 
extension to $W$ is zero on $Z(W)$. To overcome this inconvenience, recall that 
\[
X \otimes_\Z \R = \mf a^* = \mf a^*_{F_0} \oplus \mf a^{*F_0} = 
\mf a^*_{F_0} \oplus (Z(W) \otimes_\Z \R) .
\]
Thus we can decompose any $x \in X \subset \mf a^*$ uniquely as $x = x_{F_0} + x^{F_0} \in
\mf a^*_{F_0} \oplus \mf a^{*F_0}$. Now we define
\[
\mc N (w) = \ell (w) + \norm{w(0)^{F_0}} \qquad w \in W .
\]
Since $W^\af \oplus Z (W)$ is of finite index in $W$, the set $\{ w \in W : \mc N (w) = 0\}$
is finite. For $n \in \N$ we define the following norm on $\mc H$:
\[
p_n \big( \sum_{w \in W} h_w N_w \big) = \sup_{w \in W} |h_w| (\mc N (w) + 1)^n .
\]
The completion $\mc S = \mc S (\mc R ,q)$ of $\mc H$ with respect to the family of norms
$\{ p_n : n \in \N \}$ is a nuclear Fr\'echet space. It consists of all (possibly infinite) sums
$h = \sum_{w \in W} h_w N_w$ such that $p_n (h) < \infty$ for all $n \in \N$.

\begin{thm}\label{thm:3.7}
There exist $C_q > 0 ,\, d \in \mh N$ such that $\forall h,h' \in \mc S (\mc R ,q), n \in \mh N$
\begin{align*}
& \norm{h}_{B(L^2 (\mc R,q))} \leq C_q p_d (h) , \\
& p_n (h \cdot h') \leq C_q p_{n+d}(h) p_{n+d}(h') .
\end{align*}
In particular $\mc S (\mc R ,q)$ is a unital locally convex
*-algebra, and it is contained in $C^* (\mc R ,q)$.
\end{thm}
\emph{Proof.}
This was proven first with representation theoretic methods in \cite[Section 6.2]{Opd-Sp}.
Later the author found a purely analytic proof \cite[Theorem A.7]{OpSo2}. $\qquad \Box$
\\[3mm]
It is easily seen that the action of $\Gamma$ on $\mc H$ preserves all the above norms.
Hence the crossed product $\mc S \rtimes \Gamma = \mc S (\mc R,q) \rtimes \Gamma$ 
(respectively $C^* (\mc R,q) \rtimes \Gamma$) is a well-defined Fr\'echet algebra 
(respectively $C^*$-algebra). For $q = 1$ we obtain the algebras
\begin{equation}
\begin{array}{rrrrrrr}
\mc S(\mc R ,1) \rtimes \Gamma & = & \mc S (W') & = & 
\mc S (X) \rtimes \WG  & = & C^\infty (T_{un}) \rtimes \WG , \\
C^* (\mc R ,1) \rtimes \Gamma & = & C^* (W') & = & C^*(X) \rtimes \WG & 
= & C (T_{un}) \rtimes \WG ,
\end{array}
\end{equation}
where $\mc S (X)$ denotes the algebra of rapidly decreasing functions on $X$.

We can use these topological completions to characterize discrete series and 
tempered representations. According to \cite[Lemma 2.22]{Opd-Sp}, an irreducible 
$\mc H \rtimes \Gamma$-representation $\pi$ is discrete series if and only if it is contained 
in the left regular representation of $\mc H \rtimes \Gamma$ on $L^2 (\mc R,q) \otimes 
\C [\Gamma]$, or equivalently if its character $\chi_\pi : \mc H \rtimes \Gamma \to \C$ 
extends to a continuous linear functional on $L^2 (\mc R,q) \otimes \C [\Gamma]$.

By \cite[Lemma 2.20]{Opd-Sp} a finite dimensional $\mc H \rtimes \Gamma$-representation is 
tempered if and only if it extends continuously to an $\mc S \rtimes \Gamma$-representation. 
More generally, suppose that $\pi$ is a representation of $\mc H \rtimes \Gamma$ on a Fr\'echet 
space $V$, possibly of infinite dimension. As in \cite[Proposition A.2]{OpSo1}, we define $\pi$ to 
be tempered if it induces a continuous map $(\mc S \rtimes \Gamma) \times V \to V$.

A crucial role in the harmonic analysis on affine Hecke algebra is played by a particular Fourier 
transform, which is based on the induction data space $\Xi$. Let $\mc V_\Xi^\Gamma$ be the
vector bundle over $\Xi$, whose fiber at $(P,\delta,t) \in \Xi$ is the representation space
$\C [\Gamma \times W^P] \otimes V_\delta$ of $\pi^\Gamma (P,\delta,t)$. Let
$\mr{End} (\mc V_\Xi^\Gamma)$ be the algebra bundle with fibers 
$\mr{End}_\C (\C [\Gamma \times W^P] \otimes V_\delta)$. The inner product 
\eqref{eq:inpind} endows $\mr{End}_\C (\C [\Gamma \times W^P] \otimes V_\delta)$ and
$\mr{End} (\mc V_\Xi^\Gamma)$ with a canonical involution *. Of course these vector bundles
are trivial on every connected component of $\Xi$, but globally not even the dimensions need be
constant. Since $\Xi$ has the structure of a complex algebraic variety, we can construct the
algebra of polynomial sections of $\mr{End} (\mc V_\Xi^\Gamma)$:
\[
\mc O \big( \Xi ; \mr{End} (\mc V_\Xi^\Gamma) \big) := \bigoplus_{P,\delta} \mc O (T^P) \otimes
\mr{End}_\C (\C [\Gamma \times W^P] \otimes V_\delta) .
\]
Given a reasonable subset (preferably a submanifold) $\Xi' \subset \Xi$, we define the algebras
$L^2 \big( \Xi' ; \mr{End} (\mc V_\Xi^\Gamma) \big), C \big( \Xi' ; \mr{End} (\mc V_\Xi^\Gamma) \big)$ 
and $C^\infty \big( \Xi' ; \mr{End} (\mc V_\Xi^\Gamma) \big)$ in similar fashion. Furthermore,
if $\mu$ is a sufficiently nice measure on $\Xi$ and $\Xi'$ is compact, then the following formula 
defines a Hermitian form on $L^2 \big( \Xi' ; \mr{End} (\mc V_\Xi^\Gamma) \big)$:
\begin{equation}\label{eq:inpmu}
\inp{f_1}{f_2}_\mu := \int_{\Xi'} \text{tr} (f_1 (\xi)^* f_2 (\xi)) \, \textup{d} \mu .
\end{equation}
The intertwining operators from Theorem \ref{thm:3.5} give rise to an action of the groupoid
$\mc G$ on the algebra of rational sections of $\mr{End} (\mc V_\Xi^\Gamma)$, by
\begin{equation}\label{eq:actSections}
(g \cdot f) (\xi) = \pi^\Gamma (g,g^{-1} \xi ) f (g^{-1} \xi) \pi^\Gamma (g, g^{-1} \xi )^{-1} ,
\end{equation}
whenever $g^{-1} \xi \in \Xi$ is defined. This formula also defines groupoid actions of 
$\mc G$ on $C \big( \Xi' ; \mr{End} (\mc V_\Xi^\Gamma) \big)$ and on 
$C^\infty \big( \Xi' ; \mr{End} (\mc V_\Xi^\Gamma) \big)$, provided that $\Xi'$ is a
$\mc G$-stable submanifold of $\Xi$ on which all the intertwining operators are regular.
Given a suitable collection $\Sigma$ of sections of $(\Xi', \mr{End} (\mc V_\Xi^\Gamma) )$, 
we write
\[
\Sigma^{\mc G} = \{ f \in \Sigma : (g \cdot f) (\xi) = f(\xi) \text{ for all } g \in \mc G,\
\xi \in \Xi' \text{ such that } g^{-1} \xi \text{ is defined} \} .
\]
The Fourier transform for $\mc H \rtimes \Gamma$ is the algebra homomorphism
\begin{align*}\
& \mc F : \mc H \rtimes \Gamma \to \mc O \big( \Xi ; \mr{End} (\mc V_\Xi^\Gamma) \big) , \\
& \mc F (h) (\xi) = \pi (\xi) (h) .
\end{align*}
The very definition of intertwining operators shows that the image of $\mc F$ is contained
in the algebra $\mc O \big( \Xi ; \mr{End} (\mc V_\Xi^\Gamma) \big)^{\mc G}$. 
The Fourier transform also extends continuously to various topological completions of
$\mc H \rtimes \Gamma$:

\begin{thm}\label{thm:3.8} 
\textup{(Plancherel theorem for affine Hecke algebras)} \\
The Fourier transform induces algebra homomorphisms
\[
\begin{array}{rrr}
\mc H (\mc R,q) \rtimes \Gamma & \to & 
  \mc O \big( \Xi ; \mr{End} (\mc V_\Xi^\Gamma) \big)^{\mc G} , \\
\mc S (\mc R,q) \rtimes \Gamma & \to &
 C^\infty \big( \Xi_{un} ; \mr{End} (\mc V_\Xi^\Gamma) \big)^{\mc G} , \\
C^* (\mc R,q) \rtimes \Gamma & \to &
 C \big( \Xi_{un} ; \mr{End} (\mc V_\Xi^\Gamma) \big)^{\mc G} .
\end{array}
\]
The first one is injective, the second is an isomorphism of Fr\'echet *-algebras and the third
is an isomorphism of $C^*$-algebras.

Furthermore there exists a unique Plancherel measure $\mu_{Pl}$ on $\Xi$ such that
\begin{itemize}
\item the support of $\mu_{Pl}$ is $\Xi_{un}$;
\item $\mu_{Pl}$ is $\mc G$-invariant;
\item the restriction of $\mu_{pl}$ to a component $(P,\delta,T^P)$ is absolutely continuous
with respect to the Haar measure of $T^P_{un}$;
\item the Fourier transform extends to a bijective isometry
\end{itemize}
\[
\big( L^2 (\mc R ,q) \otimes \C [\Gamma] ; \inp{}{}_\tau \big) \; \to \; 
\big( L^2 \big( \Xi_{un} ; \mr{End} (\mc V_\Xi^\Gamma) \big)^{\mc G} ; \inp{}{}_{\mu_{Pl}} \big) .
\]
\end{thm}
\emph{Proof.}
Once again the essential case is $\Gamma = \{ \text{id} \}$, which is a very deep result proven by 
Delorme and Opdam, see Theorem 5.3 and Corollary 5.7 of \cite{DeOp1} and \cite[Theorem 4.43]{Opd-Sp}.

To include $\Gamma$ in the picture we need a result of general nature. Let $A$ be any complex 
$\Gamma$-algebra and endow $A \otimes_\C \mr{End} (\C [\Gamma])$ with the $\Gamma$-action
\begin{equation}\label{eq:Agammaact}
\gamma \cdot (a \otimes f) (v) = \gamma \cdot a \otimes f (v \gamma) \gamma^{-1} 
\qquad a \in A, \gamma \in \Gamma, v \in \C [\Gamma], f \in \mr{End} (\C [\Gamma]) .
\end{equation}
There is a natural isomorphism
\begin{equation}\label{eq:isocrossed}
A \rtimes \Gamma \cong \big( A \otimes_\C \mr{End} (\C [\Gamma]) \big)^\Gamma .
\end{equation}
This is easy to show, but it appears to be one of those folklore results whose origins
are hard to retrace. In any case a proof can be found in \cite[Lemma A.3]{SolThesis}.
For $A = C \big( \Xi_{un} ; \mr{End} (\mc V_\Xi^\Gamma) \big)$ the action \eqref{eq:Agammaact}
corresponds to the action of $\Gamma$ on $\mr{End}(\mc V_\Xi^\Gamma)$ described in
\eqref{eq:actSections}. The greater part of the theorem follows from \eqref{eq:isocrossed} 
and the case $\Gamma = \{ \text{id} \}$. It only remains to see how the inner products
$\inp{}{}_\tau$ and $\inp{}{}_{\mu_{Pl}}$ behave when $\Gamma$ is included. Let us 
distinguish the new inner products with a subscript $\Gamma$. On the Hecke
algebra side it is easy, as the formula \eqref{eq:*tau} does not change, so
\[
\inp{N_\gamma h}{N_{\gamma'} h'}_{\Gamma,\tau} = \left\{ \begin{array}{lll}
\inp{h}{h'}_\tau & \text{if} & \gamma = \gamma', \\
0 & \text{if} & \gamma \neq \gamma' . 
\end{array} \right.
\]
On the spectral side the inclusion of $\Gamma$ means that we replace every $\mc H$-representation
$\pi (\xi)$ by $\mr{Ind}_{\mc H}^{\mc H \rtimes \Gamma} \pi (\xi)$. In such an induced 
representation the elements of $\Gamma$ permute the $\mc H$-subrepresentations 
$\gamma \pi (\xi)$, while $h \in \mc H$ acts by $\pi (\gamma^{-1}(h))$ on $\gamma \pi (\xi)$.
The action of $\Gamma$ on $\mc H$ preserves the trace and the *, so 
\[
\mr{tr} \big( \pi^\Gamma (\xi, N_\gamma h)^* \pi^\Gamma (\xi, N_{\gamma'} h') \big) =
\left\{ \begin{array}{lll}
\mr{tr} \big( \pi (\xi, h)^* \pi (\xi, h') \big) & \text{if} & \gamma = \gamma', \\
0 & \text{if} & \gamma \neq \gamma' . 
\end{array} \right.
\]
In view of \eqref{eq:inpmu}, this means that the $L^2$-extension of $\mc F$ is an isometry with 
respect to the Plancherel measure $\mu_{\Gamma,Pl} = |\Gamma |^{-1} \mu_{Pl}. \qquad \Box$
\\[1mm]

\begin{cor}\label{cor:3.14}
The center of $\mc S (\mc R,q) \rtimes \Gamma$ (respectively $C^* (\mc R,q) \rtimes \Gamma$)
is isomorphic to $C^\infty (\Xi_{un})^{\mc G}$ (respectively $C (\Xi_{un})^{\mc G}$).
\end{cor}
\emph{Proof.} This is the obvious generalization of \cite[Corollary 5.5]{DeOp1} to our setting.
$\qquad \Box$ \\[3mm]
Notice that $Z (\mc S \rtimes \Gamma)$ is larger than the closure of $Z (\mc H \rtimes \Gamma)$ in 
$\mc S \rtimes \Gamma$, for example $Z (\mc S \rtimes \Gamma)$ contains a nontrivial
idempotent for every connected component of $\Xi_{un} / \mc G$. Varying on the notation
$\mr{Mod}_{f,U}(\mc H \rtimes \Gamma)$ we will denote by 
\[
\mr{Mod}_{f,\Sigma} (\mc S \rtimes \Gamma)
\] 
the category of finite dimensional $\mc S \rtimes \Gamma$-modules with 
$Z(\mc S \rtimes \Gamma)$-weights in $\Sigma \subset \Xi_{un} / \mc G$.
\\[2mm]

Let us compare Schwartz algebras of affine Hecke algebras with those for reductive $p$-adic groups.
Suppose that $G$ is reductive $p$-adic group and that $\mc H \rtimes \Gamma$ is Morita 
equivalent to $\mc H (G)_{\mf s}$, in the notation of Section \ref{sec:padic}. The (conjectural) 
isomorphism described in \eqref{eq:AHAsigma} is such that $\mf a^* = X \otimes_\Z \R$
corresponds to $X_* (A) \otimes_\Z \R$. The conditions for temperedness of finite length representations 
of $\mc H \rtimes \Gamma$ and $\mc H (G)_{\mf s}$ are formulated in terms of corresponding 
negative cones in $\mf a^*$ and in $X_* (A) \otimes_\Z \R$. Therefore such a Morita equivalence 
would preserve temperedness of representations. Thus $\mr{Mod}_f (\mc S \rtimes \Gamma)$
would be equivalent to the category of finite length modules
\begin{equation}\label{eq:ModSG}
\mr{Mod}_{f,\mf s} (\mc S (G)) = \mr{Mod}_f (\mc S (G)_{\mf s}) ,
\end{equation}
where $\mc S (G)$ is the Harish-Chandra--Schwartz algebra of $G$ and $\mc S (G)_{\mf s}$ is its
two-sided ideal corresponding to the inertial equivalence class $\mf s \in \mf B (G)$.

Moreover, is $I$ is an Iwahori subgroup of a split group $G$, it is shown in \cite[Proposition 10.2]{DeOp1}
that the isomorphism $\mc H (G,I) \cong \mc H (\mc R,q)$ extends to an isomorphism $\mc S (G,I) 
\cong \mc S (\mc R,q)$. Therefore it is reasonable to expect that more generally $\mc S (G)_{\mf s}$ 
will be Morita equivalent to $\mc S (\mc R,q) \rtimes \Gamma$ in case of an isomorphism 
\eqref{eq:AHAsigma}. Further support of this is provided by Theorem \ref{thm:3.8} in 
comparison with the Plancherel theorem for $\mc S (G)$ \cite{Wal} and for $C_r^* (G)$ \cite{Ply2}. 
These show that $\mc S (\mc R,q) \rtimes \Gamma$ and $\mc S (G)_{\mf s}$, as well as
their respective $C^*$-completions, have a very similar shape, 
which can almost entirely be deduced from their categories of finite length modules.

\section{Parametrization of representations with induction data}

Theorem \ref{thm:3.8} is extremely deep and useful, a large part of what follows depends on it. 
It shows a clear advantage of $\mc S$ over $\mc H$, namely that the Fourier transform
consists of \emph{all} smooth sections. In particular one can use any smooth section,
without knowing its preimage under $\mc F$. By Corollary \ref{cor:3.14} the irreducible tempered 
$\mc H \rtimes \Gamma$-representations are partitioned 
in finite packets parametrized by $\Xi_{un} / \mc G$. Moreover, from Theorem \ref{thm:3.8}
Delorme and Opdam also deduce analogues of Harish-Chandra's Completeness Theorem 
\cite[Corollary 5.4]{DeOp1} and of Langlands' Disjointness Theorem \cite[Corollary 5.8]{DeOp1}.

We will generalize these results to all irreducible representations. For that we do not need all 
induction data from $\Xi$, in view of Lemma \ref{lem:3.6} it suffices to 
consider $\xi \in \Xi^+$. At the same time this restriction to positive induction data enables
us to avoid the singularities of the intertwiners $\pi (g,\xi)$.

\begin{thm}\label{thm:3.9}
Let $\xi = (P,\delta,t), \xi' = (P',\delta�,t') \in \Xi^+$.
\enuma{
\item The $\mc H \rtimes \Gamma$-representations $\pi^\Gamma (\xi)$ and $\pi^\Gamma (\xi')$
have a common irreducible quotient if and only if there exists a $g \in \mc G$ such that
$g \xi = \xi'$.
\item The operators $\{ \pi^\Gamma (g,\xi) : g \in \mc G , g \xi = \xi' \}$ are regular and invertible,
and they span $\mr{Hom}_{\mc H \rtimes \Gamma} (\pi^\Gamma (\xi), \pi^\Gamma (\xi'))$.
}
\end{thm}
\emph{Proof.}
(a) Suppose that there exists a $g \in \mc G$ with $g \xi = \xi'$. Since $\pi^\Gamma (\gamma, \xi')$ 
is invertible for $\gamma \in \Gamma$, we may replace $\xi'$ by $\gamma \xi'$, which allows
to assume without loss of generality that $g = (w,u) \in W_0 \times K_P$.

Recall that $T^+$ is a fundamental domain for the action of $W_0$ on $T_{rs}$. Since $|t|$ and
$|t'|$ are both in $T^+$ we must have $|t| = w(|t|) = |t'|$ and hence $P(\xi) = P(\xi')$. 
Thus $w u (P,\delta,t |t|^{-1}) = (P',\delta',t' |t'|^{-1})$ and by Theorem \ref{thm:3.5}.b the
$\mc H^{P(\xi)}$-representations
\begin{align*}
& \pi^{P(\xi)} (P,\delta,t |t|^{-1}) = \pi^{P(\xi)} (P,\delta,t) \circ \phi_{|t|}^{-1} \text{ and } \\
& \pi^{P(\xi)} (P',\delta',t' |t|^{-1}) = \pi^{P(\xi)} (P',\delta',t') \circ \phi_{|t|}^{-1}
\end{align*}
are isomorphic. Hence $\pi^{P(\xi)} (P,\delta,t) \cong \pi^{P(\xi)} (P',\delta',t')$, which 
implies that $\pi^\Gamma (\xi)$ and $\pi^\Gamma (\xi')$ are isomorphic. In particular
$\pi^\Gamma (\xi)$ and $\pi^\Gamma (\xi')$ have the same irreducible quotients.

Conversely, suppose that $\pi^\Gamma (\xi)$ and $\pi^\Gamma (\xi')$ have a common 
irreducible quotient. Again we may replace $\xi'$ by $\gamma \xi'$ for any $\gamma \in \Gamma$.
In view of this, Proposition \ref{prop:3.3}.c and Corollary \ref{cor:2.8}.b we may assume that 
$P(\xi) = P(\xi')$ and that the $\mc H^{P(\xi)} \rtimes \Gamma_{P(\xi)}$-representations
$\pi^{P(\xi),\Gamma_{P(\xi)}} (\xi)$ and $\pi^{P(\xi),\Gamma_{P(\xi)}} (\xi')$ have a common
irreducible summand $\pi^{P(\xi),\Gamma_{P(\xi)}} (P(\xi),\sigma,t^{P(\xi)},\rho)$.

Pick $s \in T^{P(\xi)}_{un} t^{P(\xi)}$ such that $t^{P(\xi)}$ and 
$t^{P(\xi)} s^{-1}$ have the same isotropy group in $W_0 \rtimes \Gamma$. 
This is possible because $T^g$ is a complex algebraic subtorus of $T$ for every 
$g \in W_0 \rtimes \Gamma$. The $\mc H^{P(\xi)} \rtimes \Gamma_{P(\xi)}$-representations 
$\pi^{P(\xi),\Gamma_{P(\xi)}} (P,\delta,t s^{-1})$ and 
$\pi^{P(\xi),\Gamma_{P(\xi)}} (P',\delta',t' s^{-1})$ are completely reducible by 
Proposition \ref{prop:3.3}.a, and they have the common irreducible summand
\[
\pi^{P(\xi),\Gamma_{P(\xi)}} (P(\xi),\sigma,t^{P(\xi)} s^{-1},\rho) =
\mr{Ind}_{\mc H^{P(\xi)} \rtimes \Gamma_{P(\xi),\sigma,t^{P(\xi)}}}^{\mc H^{P(\xi)} \rtimes \Gamma_{P(\xi)}} 
\big( \sigma \circ \phi_{t^{P(\xi)}} \circ \phi_s^{-1} \otimes \rho \big) .
\]
Moreover, because every irreducible summand is of this form,
\begin{multline}\label{eq:3.1}
\mr{Hom}_{\mc H^{P(\xi)} \rtimes \Gamma_{P(\xi)}} \big( \pi^{P(\xi),\Gamma_{P(\xi)}} (P,\delta,t s^{-1}) ,
\pi^{P(\xi),\Gamma_{P(\xi)}} (P',\delta',t' s^{-1}) \big) \cong \\
\mr{Hom}_{\mc H^{P(\xi)} \rtimes \Gamma_{P(\xi)}} \big( \pi^{P(\xi),\Gamma_{P(\xi)}} (P,\delta,t) ,
\pi^{P(\xi),\Gamma_{P(\xi)}} (P',\delta',t') \big) \neq 0.
\end{multline}
Since $t^{P(\xi)} s^{-1} \in T_{un}$, we have $t s^{-1}, t' s^{-1} \in T_{un}$. So $|t| = |t'|$ and
the representations $\pi^{P(\xi),\Gamma_{P(\xi)}} (P,\delta,t s^{-1})$ and 
$\pi^{P(\xi),\Gamma_{P(\xi)}} (P',\delta',t' s^{-1})$ extend continuously to 
$\mc S (\mc R^{P(\xi)},q^{P(\xi)}) \rtimes \Gamma_{P(\xi)}$. Now Theorem \ref{thm:3.8} for
this algebra shows that the left hand side of \eqref{eq:3.1} is spanned by the intertwiners 
$\pi^{P(\xi),\Gamma_{P(\xi)}} (g,P,\delta, t s^{-1})$ with $g (P,\delta, t s^{-1}) = (P',\delta',t' s^{-1})$. 
Since \eqref{eq:3.1} is nonzero, there exists at least one such $g \in \mc G$. The choice of $s$
guarantees that $g (P,\delta, t) = (P',\delta',t')$ as well.\\
(b) By Theorem \ref{thm:3.5} the $\pi^{P(\xi),\Gamma_{P(\xi)}} (g,P,\delta, t s^{-1})$ are invertible and 
constant on $T^{P(\xi)}$-cosets. Hence the $\pi^{P(\xi),\Gamma_{P(\xi)}} (g,P,\delta, t)$ span the right 
hand side of \eqref{eq:3.1}, and they are invertible. $\qquad \Box$
\\[3mm]

It is interesting to compare Theorem \ref{thm:3.9} with \cite{Ree1}, which describes the 
$\mc H$-endomorphisms of principal series representations $M (t)$. It transpires that the results 
of \cite{Ree1} simplify considerably when $|t|$ is in the positive Weyl chamber: then 
$\mr{End}_{\mc H}(M(t))$ is semisimple and all its irreducible quotients occur with multiplicity one.

Now we can prove the desired partition of Irr$(\mc H \rtimes \Gamma)$ in packets:

\begin{thm}\label{thm:3.10}
Let $\pi$ be an irreducible $\mc H \rtimes \Gamma$-representation. There exists a unique
association class $\mc G (P,\delta,t) \in \Xi / \mc G$ such that the following equivalent properties hold:
\enuma{
\item $\pi$ is isomorphic to an irreducible quotient of $\pi^\Gamma (\xi^+)$, for some
$\xi^+ \in \Xi^+ \cap \mc G (P,\delta,t)$;
\item $\pi$ is a constituent of $\pi^\Gamma (P,\delta,t)$, and $\norm{cc_P (\delta)}$ 
is maximal for this property.
}
\end{thm}
\emph{Proof.}
Proposition \ref{prop:3.3}.d says that there exists $\xi^+ = (P',\delta',t') \in \Xi^+$ satisfying (a), 
and by Theorem \ref{thm:3.9} its $\mc G$-association class is unique.

Let $\xi = (P,\delta,t) \in \Xi$ such that $\pi$ is a constituent of $\pi^\Gamma (\xi)$ and
$\norm{cc_P (\delta)}$ is maximal under this condition. By Lemma \ref{lem:3.6} we may assume that
$\xi \in \Xi^+$. Suppose that $\pi$ is not isomorphic to a quotient of $\pi^\Gamma (\xi)$. In view of
Proposition \ref{prop:3.3} this means that there exist Langlands data $(P(\xi+),\sigma',t'^{P(\xi^+)},\rho')$
and $(P(\xi), \sigma, t^{P(\xi)},\rho)$ such that $\pi \cong L^\Gamma (P(\xi+),\sigma',t'^{P(\xi^+)},\rho')$
is a constituent but not a quotient of $\pi^\Gamma (P(\xi), \sigma, t^{P(\xi)},\rho)$. Now 
Lemma \ref{lem:2.10}.b tells us that
\[
\norm{cc_{P(\xi)}(\sigma)} < \norm{cc_{P(\xi^+)}(\sigma')} .
\]
But $cc_{P(\xi)}(\sigma) = W_{P(\xi)} cc_P (\delta)$ and $cc_{P(\xi^+)}(\sigma') = W_{P(\xi^+)}
cc_{P'}(\delta')$, so 
\[
\norm{cc_P (\delta)} < \norm{cc_{P'} (\delta')} ,
\] 
contradicting the maximality of $\norm{cc_P (\delta)}$. Therefore $\pi$ must be a quotient of 
$\pi^\Gamma (\xi)$.

Thus the association class $\mc G \xi$ satisfies not only (b) but also (a), which at the same time
shows that is unique. In particular conditions (a) and (b) turn out to be equivalent.
$\qquad \Box$ \\[3mm]

All these constructs with induction data have direct analogues in the setting of graded Hecke 
algebras \cite[Sections 6 and 8]{SolGHA}. Concretely, $\tilde \Xi$ is the space of all triples 
$\tilde \xi = (Q,\sigma,\lambda)$, where $Q \subset F_0 \,, \sigma$ is a discrete series 
representation of $\mh H_Q$ and $\lambda \in \mf t^Q$. The subsets of unitary (respectively 
positive) induction data are obtained by imposing the restriction $\lambda \in i \mf a^Q$ 
(respectively $\lambda \in \mf a^{Q+} + i \mf a^Q$). The corresponding induced representation is
\[
\pi^\Gamma (\tilde \xi) = \mr{Ind}_{\mh H}^{\mh H \rtimes \Gamma} \pi (Q,\sigma,\lambda) =
\mr{Ind}_{\mh H_Q}^{\mh H \rtimes \Gamma} (\sigma_\lambda ) .
\]
The groupoid $\tilde{\mc G}$ and its action on $\tilde \Xi$ are defined like $\mc G$,
but without the parts $K_P$. We would like to understand the relation between induction data 
for $\mc H \rtimes \Gamma$ and for $\mh H \rtimes \Gamma$. We consider, for every $u \in T_{un}$, 
the induction data for $(\tilde{\mc R}_u,k_u)$ with $\lambda \in \mf a$. Thus we arrive at the space 
$\widehat \Xi$ of quadruples $\hat \xi = (u, \tilde P, \tilde \delta, \lambda)$ such that:
\begin{itemize}
\item $u \in T_{un}$;
\item $\tilde P \subset F_u$;
\item $\tilde \delta$ is a discrete series representation of 
$\mh H (\tilde{\mc R}_{u,\tilde P},k_{u,\tilde P})$;
\item $\lambda \in \mf a^{\tilde P}$.
\end{itemize}
The $\mc H \rtimes \Gamma$-representation associated to $\hat \xi$ is
\[
\pi^\Gamma (u, \tilde P, \tilde \delta, \lambda) = \mr{Ind}_{\mc H (\mc R_u,q_u)}^{\mc H \rtimes \Gamma}
\pi (\tilde P, \tilde \delta, \lambda) ,
\]
where the $\mh H (\tilde{\mc R}_u,k_u)$-representation $\pi (\tilde P, \tilde \delta, \lambda)$ 
is considered as a representation of $\mc H (\tilde{\mc R}_u,q_u)$, via Theorem \ref{thm:2.2}.

For $g \in \mc G$ the map $\psi_g$ from \eqref{eq:twistKP} and \eqref{eq:psigamma} induces
an algebra isomorphism $\mh H (\tilde{\mc R}_u, k_u) \to \mh H (\tilde{\mc R}_{g(u)}, k_{g(u)})$,
and the stabilizer in $\mc G$ of $u \in T_{un}$ is the groupoid $\tilde {\mc G}_u$ associated to
$(\tilde{\mc R}_u, k_u)$. This leads to an action of $\mc G$ on $\widehat \Xi$.

The collections of $\mc H \rtimes \Gamma$-representations corresponding to $\Xi$ and to $\widehat \Xi$
are almost the same, but not entirely:

\begin{lem}\label{lem:3.5}
There exists a natural finite-to-one surjection 
\[
\Xi / \mc G \to \widehat \Xi / \mc G, \: \mc G \xi \mapsto \mc G \hat{\xi} ,
\]
with the following property. Given $\hat \xi \in \widehat \Xi$ one can find $\xi_i \in \Xi$ 
(not necessarily all different) such that $\bigcup_i \mc G \xi_i$ is the preimage of 
$\mc G \hat \xi$ and $\pi^\Gamma (\hat \xi) = \bigoplus_i \pi^\Gamma (\xi_i)$.
\end{lem}
\emph{Proof.}
Given $\xi = (P,\delta,t) \in \Xi$, let $t = u^P c^P \in T^P_u \times T^P_{rs}$ be the polar 
decomposition of $t$ and let $u_P c_P \in T_{P,u} \times T_{P,rs}$ be an 
$\mc A_P$-weight of $\delta$. Put 
\begin{equation}\label{eq:WP+}
W_{P,u_P}^+ = \{ w \in W (R_P) : w(u_P) = u_P, w (R_{P,u_P}^+) = R_{P,u_P}^+ \} . 
\end{equation}
By Theorem \ref{thm:2.1} there exists a unique discrete series representation 
\[
\delta_1 \text{ of } \mc H (\mc R_{P,u_P} ,q_{P,u_P}) \rtimes W_{P,u_P}^+ 
\text{  such that  } \delta \cong \mr{Ind}_{\mc H (\mc R_{P,u_P} ,q_{P,u_P}) 
\rtimes W_{P,u_P}^+}^{\mc H_P}(\delta_1). 
\]
Then automatically
\[
\delta \circ \phi_t \cong \mr{Ind}_{\mc H (\mc R^P_{u_P} ,q^P_{u_P}) \rtimes 
W_{P,u_P}^+}^{\mc H^P} (\delta_1 \circ \phi_t) .
\]
Let $\delta'$ be an irreducible direct summand of the restriction of $\delta_1$ to 
$\mc H (\mc R_{P,u_P} ,q_{P,u_P})$, such that $u_P c_P$ is a weight of $\delta'$. Then 
$\delta_1 \circ \phi_t$ is a direct summand of
\begin{equation}
\mr{Ind}_{\mc H (\mc R^P_{u_P} ,q^P_{u_P})}^{\mc H (\mc R^P_{u_P} ,q^P_{u_P}) 
\rtimes W_{P,u_P}^+} (\delta' \circ \phi_t) , 
\end{equation}
so $\pi^\Gamma (\xi)$ is a direct summand of 
\begin{equation}\label{eq:pi''xi}
\mr{Ind}_{\mc H (\mc R^P_{u_P} ,q^P_{u_P})}^{\mc H \rtimes \Gamma} (\delta' \circ \phi_t) .
\end{equation}
By Theorem \ref{thm:2.2} $\delta'$ can also be regarded as a discrete series representation 
$\tilde \delta$ of $\mh H (\tilde{\mc R}_{P,u_P} ,k_{P,u_P})$ with central character 
$W(R_{P,u_P}) c_P$. Then $\delta' \circ \phi_t$ corresponds to the representation 
$\tilde \delta_{\log (c^P)}$ of $\mh H (\tilde{\mc R}^P_{u_P} ,k^P_{u_P})$. Let $\tilde P$ be 
the unique basis of $R_{P,u_P}$ contained in $R_0^+$.

All in all $(P,\delta,t)$ gives rise to the induction datum 
$\tilde \xi = (\tilde P,\tilde \delta, \log (c^P))$
for the graded Hecke algebra $\mh H (\tilde{\mc R}_{P,u_P} ,k_{P,u_P})$. Since $R_{P,u_P}$ 
is a parabolic root subsystem of $R_u \,, \tilde \xi$ can also be regarded as an induction datum for
$\mh H (\tilde {\mc R}_u ,k_u)$. Let us check the possible freedom in the above construction.
All $\mc A_P$-weights of $\delta$ are in the same $W_P$-orbit, so another choice of $u_P c_P$ 
differs only by an element of $W_P$. All possible choices of $\delta'$ above are conjugated 
by the action of the group $W_{P,u_P}^+$, and $(W_0 \rtimes \Gamma)u$ is the unitary part of the 
central character of $\pi^\Gamma (\xi)$. Therefore $\xi$ determines the quadruple 
\[
\hat{\xi} := (u,\tilde P, \tilde \delta, \log (c^P)) \in \widehat{\Xi}
\]
uniquely modulo the action of $\mc G$. That yields a map 
$\Xi \to \widehat \Xi / \mc G, \: \xi \to \mc G \hat{\xi}$, and since the actions of 
$\mc G$ are defined in the same way on both sides, this map factors via $\Xi / \mc G$.

By reversing the above steps one can reconstruct the representations $\delta' \circ \phi_t$ and
\eqref{eq:pi''xi} from $\tilde \xi \,, u_P$ and $u^P$. In fact it 
suffices to know $\tilde \xi$ and the product $u = u_P u^P \in T_{un}$. 
Namely, the only additional ambiguity comes from the group $K_P$, but this is inessential since
\[
(\delta' \circ \psi_{k^{-1}}) \circ \phi_{k t} \cong \delta' \circ \phi_t
\text{ for } k \in K_P . 
\]
By construction $\delta_1$ is a direct summand of $\mr{Ind}_{\mc H 
(\mc R_{P,u_P} ,q_{P,u_P})}^{\mc H (\mc R_{P,u_P} ,q_{P,u_P}) \rtimes W_{P,u_P}^+}(\delta')$, 
and the other constituents $\delta_j \; (1 < j \leq n)$ are also discrete series representations. 
Hence \eqref{eq:pi''xi} is a direct sum of finitely many parabolically induced representations 
\[
\pi^\Gamma (P,  \mr{Ind}_{\mc H (\mc R_{P,u_P} ,q_{P,u_P}) \rtimes W_{P,u_P}^+}^{\mc H_P}(\delta_j),t) .
\] 
Now Corollary \ref{cor:2.3} assures that our map $\Xi / \mc G \to \widehat \Xi / \mc G$ is surjective 
and that the preimage of $\mc G (u,\tilde P, \tilde \delta, \log (c^P))$ consists precisely of
the association classes 
\[
\mc G (P, \mr{Ind}_{\mc H (\mc R_{P,u_P} ,q_{P,u_P}) \rtimes 
W_{P,u_P}^+}^{\mc H_P}(\delta_j),t) \quad (1 \leq j \leq n) . \qquad \Box
\]
\emph{Remark.} 
Things simplify considerably if the group $W^+_{P,u_P} = \{ \mr{id} \}$ in \eqref{eq:WP+}, 
then the map $\Xi / \mc G \to \widehat{\Xi} / \mc G$ is bijective on $(P,\delta,T^P) / \mc G$. 
In many cases this group is indeed trivial, but not always. See \cite[Section 8]{OpSo2}, 
where $W^+_{P,u_P}$ is denoted $\Gamma_{s(e)}$.

\section{The geometry of the dual space}
\label{sec:dual}

For any algebra $A$ the set Irr$(A)$ has a natural topology, the Jacobson topology. This is the
noncommutative generalization of the Zariski topology, by definition all its closed sets are
of the form
\[
V(S) := \{ \pi \in \text{Irr}(A) : \pi (s) = 0 \: \forall s \in S \} \qquad S \subset A.
\]
In this section we discuss the topology and the geometry of Irr$(\mc H \rtimes \Gamma)$, and we
compare it with the dual of $\mc S \rtimes \Gamma$. This will be useful for the proof of
Theorem \ref{thm:2.7} and for our discussion of periodic cyclic homology in Section \ref{sec:pch}.

Parabolic induction gives, for every discrete series representation $\delta$ of a parabolic
subalgebra~$\mc H_P$, a family of $\mc H \rtimes \Gamma$-representations 
$\pi^\Gamma (P,\delta,t)$, parametrized by $t \in T^P$. The group
\[
\mc G_{P,\delta} := \{ g \in \mc G : g (P) = P, \delta \circ \psi_g^{-1} \cong \delta \}
\]
acts algebraically on $T^P$, and by Lemma \ref{lem:3.6} points in the same orbit lead to 
representations with the same irreducible subquotients. 

Theorem \ref{thm:3.10} allows us to associate to every $\pi \in \text{Irr}(\mc H \rtimes \Gamma)$ 
an induction datum $\xi^+ (\pi) \in \Xi^+$, unique modulo $\mc G$, such that $\pi$ is a quotient of 
$\pi^\Gamma (\xi^+ (\pi))$. For any subset $U \subset T^P$ we define
\[
\mr{Irr}_{P,\delta,U} (\mc H \rtimes \Gamma ) = \{ \pi \in \text{Irr}(\mc H \rtimes \Gamma) :
\mc G \xi^+ (\pi) \cap (P,\delta,U) \neq \emptyset \} .
\]
For $U = T^P$ or $U = \{t\}$ we abbreviate this to $\mr{Irr}_{P,\delta}(\mc H \rtimes \Gamma)$ or
$\mr{Irr}_{P,\delta,t}(\mc H \rtimes \Gamma)$. 

\begin{prop}\label{prop:3.11}
Let $U$ be a subset of $T^{P+} T^P_{un}$ such that every $g \in \mc G_{P,\delta}$ 
with $g U \cap U \neq \emptyset$ fixes $U$ pointwise. For arbitrary $t \in U$ there are 
canonical bijections
\[
\mr{Irr}_{P,\delta} (\mc H_P \rtimes \Gamma_{P,\delta,t} ) \times U  \rightarrow
\mr{Irr}_{P,\delta,U} (\mc H^P \rtimes \Gamma_{P,\delta,t} ) \rightarrow
\mr{Irr}_{P,\delta,U} (\mc H \rtimes \Gamma ) .
\]
\end{prop}
\emph{Remark.} It is not unreasonable to expect that the Jacobson topology of 
$\mc H \rtimes \Gamma$ induces the Zariski topology on $\mr{Irr}_{P,\delta} (\mc H_P 
\rtimes \Gamma_{P,\delta,t} ) \times U$, where $\mr{Irr}_{P,\delta} (\mc H_P 
\rtimes \Gamma_{P,\delta,t} )$ is regarded as a discrete space. However, while it is easy to
see that all $V(h)$ become Zariski-closed in $\mr{Irr}_{P,\delta} (\mc H_P 
\rtimes \Gamma_{P,\delta,t} ) \times U$, it is not clear that one can obtain all
Zariski-closed subsets in this way. That might require some extra conditions on $U$.

\emph{Proof.}
By assumption every $t \in U$ has the same stabilizer $\mc G_{P,\delta,t} \subset 
\mc G_{P,\delta}$. According to Theorem \ref{thm:3.9} the operators 
$\{\pi^\Gamma (g,P,\delta,t) : g \in \mc G_{P,\delta,t} \}$ span 
$\mr{End}_{\mc H \rtimes \Gamma} (\pi^\Gamma (P,\delta,t))$. By definition all elements
of $\mr{Irr}_{P,\delta,t}(\mc H \rtimes \Gamma)$ occur as a quotient of $\pi^\Gamma (P,\delta,t)$, 
but the latter representation also has other constituents if it is not completely reducible.
We have to avoid that situation if we want to find a direct relation between
$\mr{End}_{\mc H \rtimes \Gamma} (\pi^\Gamma (P,\delta,t))$ and 
$\mr{Irr}_{P,\delta,t}(\mc H \rtimes \Gamma)$.

Since $\pi^\Gamma (P,\delta,t)$ and $\pi^\Gamma (g,P,\delta,t)$ are unitary for $t \in T^P_{un}$,
there exists an open $\mc G_{P,\delta}$-stable tubular neighborhood $T^P_\epsilon$ of $T^P_{un}$ 
in $T^P$, such that $\pi^\Gamma (P,\delta,t)$ is completely reducible and $\pi^\Gamma (g,P,\delta,t)$ 
is regular and invertible, for all $t \in T^P_\epsilon$ and $g \in \mc G_{P,\delta}$. For every $t \in U$
we can find $r \in \R_{>0}$ such that $t |t|^{r-1} \in T^P_\epsilon$. Let $U_\epsilon \subset T^P_\epsilon$ 
be the resulting collection. For every $t \in U_\epsilon$ the algebras
\[
\{ \pi^\Gamma (P,\delta,t) (h) : h \in \mc H \rtimes \Gamma \} 
\text{ and span}\{ \pi^\Gamma (g,P,\delta,t) : g \in \mc G_{P,\delta,t} \}
\]
are each others' commutant in 
\[
\mr{End}_\C (\pi^\Gamma (P,\delta,t)) =
\mr{End}_\C \big( \C [\Gamma \Gamma^P] \otimes_\C V_\delta \big) .
\]
Hence there is a natural bijection between
\begin{itemize}
\item isotypical components of $\pi^\Gamma (P,\delta,t)$ as a $\mc H \rtimes \Gamma$-representation;
\item isotypical components of $\pi^\Gamma (P,\delta,t)$ as a $\mc G_{P,\delta,t}$-representation.
\end{itemize}
The operators $\pi^\Gamma (P,\delta,t)$ are rational in $t \in T^P$, and regular on $T^P_\epsilon$. 
As there are only finitely many inequivalent $\mc G_{P,\delta,t}$-representations
of a fixed finite dimension, the isomorphism class of $\pi^\Gamma (P,\delta,t)$ as a 
$\mc G_{P,\delta,t}$-representation does not depend on $t \in U_\epsilon$. This provides a natural bijection
\begin{equation}\label{eq:bijectionU1}
\mr{Irr}_{P,\delta,U_\epsilon}(\mc H \rtimes \Gamma) \longleftrightarrow 
\mr{Irr}_{P,\delta,t}(\mc H \rtimes \Gamma) \times U_\epsilon  \qquad t \in U_\epsilon .
\end{equation}
The extended Langlands classification (Corollary \ref{cor:2.8}) shows that there is a canonical
bijection $\mr{Irr}_{P,\delta,t |t|^{r-1}}(\mc H \rtimes \Gamma) \leftrightarrow
\mr{Irr}_{P,\delta,t}(\mc H \rtimes \Gamma)$ for every $r \in \R_{>0}$, which allows us to extend
\eqref{eq:bijectionU1} uniquely to
\begin{equation}\label{eq:bijectionU}
\mr{Irr}_{P,\delta,U}(\mc H \rtimes \Gamma) \longleftrightarrow 
\mr{Irr}_{P,\delta,t_0}(\mc H \rtimes \Gamma) \times U \qquad t_0 \in U .
\end{equation}
The above also holds with the algebras $\mc H \rtimes \Gamma_{P,\delta,t}$ or
$\mc H_P \rtimes \Gamma_{P,\delta,t}$ in the role of $\mc H \rtimes \Gamma$. Since the construction
of the intertwiners corresponding to $g \in \mc G_{P,\delta,t}$ is the same in all three cases, 
we obtain natural isomorphisms
\[
\mr{End}_{\mc H_P \rtimes \Gamma_{P,\delta,t}} 
\big( \mr{Ind}_{\mc H_P}^{\mc H_P \rtimes \Gamma_{P,\delta,t}} \delta \big) \cong
\mr{End}_{\mc H^P \rtimes \Gamma_{P,\delta,t}} (\pi^{P,\Gamma_{P,\delta,t}} (P,\delta,t)) \cong
\mr{End}_{\mc H \rtimes \Gamma} (\pi^\Gamma (P,\delta,t)) .
\]
Now the above bijection between isotypical components shows that the maps
\[
\begin{array}{cccccc}
\mr{Irr}_{P,\delta} (\mc H_P \rtimes \Gamma_{P,\delta,t} )& \to &
\mr{Irr}_{P,\delta,t} (\mc H^P \rtimes \Gamma_{P,\delta,t} ) & \to &
\mr{Irr}_{P,\delta,t} (\mc H \rtimes \Gamma ) \\
\rho & \mapsto & \rho \circ \phi_t & \mapsto & 
\mr{Ind}_{\mc H^P \rtimes \Gamma_{P,\delta,t}}^{\mc H \rtimes \Gamma} (\rho \circ \phi_t)
\end{array}
\]
are bijective, and \eqref{eq:bijectionU} allows us to extend this from one $t$ to $U. \qquad \Box$
\\[2mm]

Theorem \ref{thm:3.10} shows that $\mr{Irr}_{P,\delta}(\mc H \rtimes \Gamma)$ and 
$\mr{Irr}_{Q,\sigma}(\mc H \rtimes \Gamma)$ are either equal or disjoint, depending on whether 
or not $(P,\delta)$ and $(Q,\sigma)$ are $\mc G$-associate.
The sets $\mr{Irr}_{P,\delta}(\mc H \rtimes \Gamma)$ are usually not closed in
Irr$(\mc H \rtimes \Gamma)$, because we did not include all constituents of the representations 
$\pi^\Gamma (P,\delta,t)$. We can use their closures to define a stratification of 
Irr$(\mc H \rtimes \Gamma)$, and a corresponding stratification of Irr$(\mc S \rtimes \Gamma)$. 
By Corollary \ref{cor:3.2} we may identify $\mr{Irr}_{P,\delta} (\mc S \rtimes \Gamma)$ with the 
tempered part $\mr{Irr}_{P,\delta,T^P_{un}} (\mc H \rtimes \Gamma) $ of 
$\mr{Irr}_{P,\delta}(\mc H \rtimes \Gamma)$.

Let $\Delta$ be a set of representatives for the action of $\mc G$ on pairs $(P,\delta)$. Then the
cardinality of $\Delta$ equals the number of connected components of $\Xi_{un} / \mc G$ and
by Theorem \ref{thm:3.8}
\begin{equation}\label{eq:FourierDelta}
\mc S \rtimes \Gamma \cong \bigoplus\nolimits_{(P,\delta) \in \Delta} \big( C^\infty (T^P_{un}) \otimes_\C 
\mr{End}(\C [\Gamma W^P] \otimes_\C V_\delta ) \big)^{\mc G_{P,\delta}} .
\end{equation}

\begin{lem}\label{lem:3.12}
There exist filtrations by two-sided ideals
\[
\begin{array}{ccccccc}
\mc H \rtimes \Gamma = F_0 (\mc H \rtimes \Gamma) & \supset & F_1 (\mc H \rtimes \Gamma) &
\supset & \cdots & \supset & F_{|\Delta |} (\mc H \rtimes \Gamma) = 0 ,\\
\mc S \rtimes \Gamma = F_0 (\mc S \rtimes \Gamma) & \supset & F_1 (\mc S \rtimes \Gamma) &
\supset & \cdots & \supset & F_{|\Delta |} (\mc S \rtimes \Gamma) = 0 ,
\end{array}
\]
with $F_i (\mc H \rtimes \Gamma) \subset F_i (\mc S \rtimes \Gamma)$, such that for all $i > 0$:
\enuma{
\item $\mr{Irr} \big( F_{i-1} (\mc S \rtimes \Gamma) / F_i (\mc S \rtimes \Gamma) \big) \cong 
\mr{Irr}_{P_i ,\delta_i}(\mc S \rtimes \Gamma)$,
\item $\mr{Irr} \big( F_{i-1} (\mc H \rtimes \Gamma) / F_i (\mc H \rtimes \Gamma) \big) \cong 
\mr{Irr}_{P_i ,\delta_i}(\mc H \rtimes \Gamma)$.
}
\end{lem}
\emph{Remark.} Analogous filtrations of Hecke algebras of reductive $p$-adic groups 
are described in \cite[Lemma 2.17]{SolPadic}. The proof in our setting is basically the same.

\emph{Proof.}
We number the elements of $\Delta$ such that
\begin{equation}\label{eq:numberij}
\norm{cc_{P_i}(\delta_i)} \geq \norm{cc_{P_j}(\delta_j)} \quad \text{if } j \leq i ,
\end{equation}
and we define
\[
\begin{array}{ccc}
F_i (\mc H \rtimes \Gamma) & = & \{ h \in \mc H \rtimes \Gamma : \pi (h) = 0 \text{ for all }
\pi \in \mr{Irr}_{P_j ,\delta_j}(\mc H \rtimes \Gamma) , j \leq i \} , \\
F_i (\mc S \rtimes \Gamma) & = & \{ h \in \mc S \rtimes \Gamma : \pi (h) = 0 \text{ for all }
\pi \in \mr{Irr}_{P_j ,\delta_j}(\mc S \rtimes \Gamma) , j \leq i \} .
\end{array}
\]
Clearly $F_i (\mc H \rtimes \Gamma) \subset F_i (\mc S \rtimes \Gamma)$ and
\[
F_{i-1}(\mc S \rtimes \Gamma) / F_i (\mc S \rtimes \Gamma) \cong 
\big( C^\infty (T^P_{un}) \otimes_\C \mr{End}(\C [\Gamma \times W^P]
\otimes_\C V_\delta ) \big)^{\mc G_{P_i ,\delta_i}} ,
\]
which establishes (a). For (b), we first show that 
\begin{equation}\label{eq:J-closed}
\bigcup\nolimits_{j \leq i} \mr{Irr}_{P_j ,\delta_j} (\mc H \rtimes \Gamma) 
\end{equation}
is closed in the Jacobson topology of Irr$(\mc H \rtimes \Gamma)$. Its Jacobson-closure
consists of all irreducible subquotients $\pi$ of $\pi^\Gamma (P_j,\delta_j,t)$, for $j \leq i$
and $t \in T^{P_j}$. Suppose that $\pi \notin \mr{Irr}_{P_j ,\delta_j}(\mc H \rtimes \Gamma)$.
By Theorem \ref{thm:3.10} $\pi \in \mr{Irr}_{P,\delta}(\mc H \rtimes \Gamma)$ for some
discrete series representation $\delta$ of $\mc H_P$ with $\norm{cc_P (\delta)} > 
\norm{cc_{P_j}(\delta_j)}$. Select $(P_n ,\delta_n)$ and $g \in \mc G$ such that 
$g (P,\delta) = (P_n, \delta_n)$. Then $\pi \in \mr{Irr}_{P_n ,\delta_n}(\mc H \rtimes \Gamma)$ and
\[
\norm{cc_{P_n}(\delta_n)} =  \norm{cc_P (\delta)} > \norm{cc_{P_j}(\delta_j)},
\]
so $n < j \leq i$ by \eqref{eq:numberij}. Therefore \eqref{eq:J-closed} is indeed Jacobson-closed and
\[
\mr{Irr}(F_i (\mc H \rtimes \Gamma)) \cong 
\bigcup\nolimits_{j < i} \mr{Irr}_{P_j ,\delta_j} (\mc H \rtimes \Gamma) ,
\]
which implies (b). $\qquad \Box$
\\[2mm]

The filtrations from Lemma \ref{lem:3.12} help us to compare the dual of $\mc H \rtimes \Gamma$
with its tempered dual, which can be identified with the dual of $\mc S \rtimes \Gamma$. The space
$\mr{Irr}_{P,\delta}(\mc H \rtimes \Gamma)$ comes with a finite-to-one projection onto 
\begin{equation}\label{eq:projPdelta}
T^{P+}T^P_{un} / \mc G_{P,\delta} \cong T^P / (W (R_P) \rtimes \mc G_{P,\delta}) .
\end{equation}
The subspace $\mr{Irr}_{P,\delta}(\mc S \rtimes \Gamma) \subset \mr{Irr}_{P,\delta}(\mc H \rtimes \Gamma)$ 
is the inverse image of $T^P_{un} / (W (R_P) \rtimes \mc G_{P,\delta})$ under this projection.
By Proposition \ref{prop:3.11} the fiber at $t \in T^P$ essentially depends only on the stabilizer 
$\mc G_{P,\delta,t}$. Since $\mc G_{P,\delta}$ acts algebraically on $T^P$, the collection of points
$t \in T^P$ such that the fiber at $(W_P \rtimes \mc G_{P,\delta})t$ has exactly $m$ points (for some fixed
$m \in \N$) is a complex affine variety, say $T^{P,m}$. As the action of $\mc G_{P,\delta}$ stabilizes
$T^P_{un}$, the variety $T^{P,m}$ is already determined by its intersection with $T^P_{un}$. Hence
one can reconstruct the set $\mr{Irr}_{P,\delta}(\mc H \rtimes \Gamma)$ from 
$\mr{Irr}_{P,\delta}(\mc S \rtimes \Gamma)$. 

To complete the proof of Theorem \ref{thm:2.7}.a we will introduce some terminology.
We know from Lemma \ref{lem:3.1} that representations of the form 
\[
\pi^{P,\Gamma_{P,\delta,u}} (P,\delta,u) = \mr{Ind}_{\mc H^P}^{\mc H^P \rtimes \Gamma_{P,\delta,u}} 
(\delta \circ \phi_u) \quad \text{with } (P,\delta,u) \in \Xi_{un}
\]
are tempered and unitary. Let $\sigma$ be an irreducible summand of such an 
$\mc H^P \rtimes \Gamma_{P,\delta,u}$-representation and let $T^\sigma$
be the connected component of $T^{W(R_P) \rtimes \Gamma_{P,\delta,u}}$ that contains $1 \in T$. 
Notice that $T^\sigma \subset T^P$ because $T_P^{W (R_P)}$ is finite. Following \cite{Opd-ICM} we call 
\begin{equation}\label{eq:SmoothFamily}
\big\{ \mr{Ind}_{\mc H^P \rtimes \Gamma_{P,\delta,u}}^{\mc H \rtimes \Gamma} (\sigma \circ \phi_t) : 
t \in T^\sigma \big\} 
\end{equation}
a smooth $d$-dimensional family of representations, where $d$ is the dimension of the complex algebraic 
variety $T^{W(R_P) \rtimes \Gamma_{P,\delta,u}}$. If we restrict the parameter $t$ to $T_{un}^\sigma =
T^\sigma \cap T_{un}$, then we add the adjective tempered to this description. Since there are only
finitely many pairs $(P,\delta)$ and since two $u$'s with the same $u T^\sigma$ give rise to the same
smooth family, there exist only finitely many smooth families of $\mc H \rtimes \Gamma$-representations.

For $t \in T^\sigma$ such that $\mc G_{P,\delta,ut} = \mc G_{P,\delta,u}$ and all the intertwiners
$\pi^\Gamma (g,P,\delta,ut)$ with $g \in \mc G_{P,\delta,u}$ are invertible, the representation 
$\mr{Ind}_{\mc H^P \rtimes \Gamma_{P,\delta,u}}^{\mc H \rtimes \Gamma} (\sigma \circ \phi_t)$
is irreducible by Proposition \ref{prop:3.11}. Hence the family \eqref{eq:SmoothFamily} consists
of representations which are irreducible on a Zariski-open dense subset of the parameter space $T^\sigma$. 

The proof of Proposition \ref{prop:3.3}.b shows that 
\[
\mr{Ind}_{\mc H^P \rtimes \Gamma_{P,\delta,u}}^{\mc H \rtimes \Gamma}  (\sigma \circ \phi_t)
\]
is always a direct sum of representations of the form $\pi^\Gamma (P',\sigma',t',\rho')$, where 
$(P',\sigma',t',\rho')$ is almost a Langlands datum for $\mc H \rtimes \Gamma$, the only difference being 
that $|t'| \in T^{P'}_{rs}$ need not be positive. Nevertheless we can specify a "Langlands constituent"
of $\pi^\Gamma (P',\sigma',t',\rho')$:

\begin{lem}\label{lem:3.16}
Let $Q \subset F_0, t \in T^Q$, let $\tau$ be an irreducible tempered $\mc H_Q$-representation
and let $\rho$ be an irreducible representation of $\C [\Gamma_{Q,\tau,t}, \kappa]$. The representation
$\pi^\Gamma (Q,\tau,t,\rho)$ defined by \eqref{eq:piLanglands} has a unique "Langlands" constituent
$L \big( \pi^\Gamma (Q,\tau,t,\rho) \big)$ which is minimal in the sense of Lemma \ref{lem:2.10}.b.
Moreover $L \big( \pi^\Gamma (Q,\tau,t,\rho) \big)$ occurs with multiplicity one in $\pi^\Gamma (Q,\tau,t,\rho)$.
\end{lem}
\emph{Proof.}
Choose $g \in \mc G$ such that $g(Q,\tau,t) \in \Xi^+$. By Proposition \ref{prop:3.3}.b $g (Q,\tau,t,\rho)$ 
is a Langlands datum for $\mc H \rtimes \Gamma$ and by Theorem \ref{thm:3.9} it does not matter which
$g \in \mc G$ with this property we take. Like in Lemma \ref{lem:3.6}, 
$\pi^\Gamma (g (Q,\tau,t,\rho))$ has the same trace and the same semisimplication as 
$\pi^\Gamma (Q,\tau,t,\rho)$. Using Corollary \ref{cor:2.8}.a we define 
\begin{equation}\label{eq:LanglandsC}
L \big( \pi^\Gamma (Q,\tau,t,\rho) \big) := L^\Gamma (g (Q,\tau,t,\rho)) .
\end{equation}
Lemma \ref{lem:2.10}.b characterizes this as the unique constituent which is minimal in the appropriate
sense, and it shows that $L \big( \pi^\Gamma (Q,\tau,t,\rho) \big)$ appears only once in 
$\pi^\Gamma (Q,\tau,t,\rho) . \qquad \Box$ \\[2mm]

Let $L \big( \mr{Ind}_{\mc H^P \rtimes \Gamma_{P,\delta,u}}^{\mc H \rtimes \Gamma} (\sigma \circ \phi_t) \big)$ 
be the direct sum of the Langlands constituents of the irreducible summands of 
$\mr{Ind}_{\mc H^P \rtimes \Gamma_{P,\delta,u}}^{\mc H \rtimes \Gamma} (\sigma \circ \phi_t)$. The family 
\begin{equation}\label{eq:LSmoothFamily}
\big\{ L \big( \mr{Ind}_{\mc H^P \rtimes \Gamma_{P,\delta,u}}^{\mc H \rtimes \Gamma} (\sigma \circ \phi_t) \big) : 
t \in T^\sigma \big\} 
\end{equation}
cannot be called smooth, because  the traces of these representations do not depend continuously on $t$.
Let us call it an L-smooth $d$-dimensional family of representations. 
\\[2mm]
\emph{Continuation of the proof of Theorem \ref{thm:2.7}}.\\
We have to show that the $\Q$-linear extension
\begin{equation}\label{eq:SprQ}
\Spr_\Q : G_\Q (\mc H \rtimes \Gamma) \to G_\Q (X \rtimes (W_0 \rtimes \Gamma)) .
\end{equation}
of $\Spr$ is a bijection. By properties (d) and (e) of Theorem \ref{thm:2.7}
\begin{equation}\label{eq:SprFamilies}
\Spr_\Q L \big( \mr{Ind}_{\mc H^P \rtimes \Gamma_{P,\delta,u}}^{\mc H \rtimes \Gamma} 
(\sigma \circ \phi_t) \big) = \mr{Ind}_{X \rtimes (W (R_P) \rtimes \Gamma_{P,\delta,u})}^{X \rtimes 
(W_0 \rtimes \Gamma)} (\Spr_\Q (\sigma) \circ \phi_t) . 
\end{equation}
The right hand side is almost a smooth $d$-dimensional family of $X \rtimes (W_0 \rtimes \Gamma)$-representations.
Not entirely, because $\Spr_\Q (\sigma)$ is in general reducible and because a priori this family could
be only a part of a higher dimensional smooth family. 

Let $G_\Q^d (\mc H \rtimes \Gamma)$ be the $\Q$-submodule of $G_\Q (\mc H \rtimes \Gamma)$ 
spanned by the representations \eqref{eq:LSmoothFamily}, for all L-smooth families of dimension at least 
$d \in \Z_{\geq 0}$. This is a decreasing sequence of $\Q$-submodules of $G_\Q (\mc H \rtimes \Gamma)$,
by convention 
\[
G_\Q^0 (\mc H \rtimes \Gamma) = G_\Q (\mc H \rtimes \Gamma)
\]
and $G_\Q^d (\mc H \rtimes \Gamma) = 0$ when $d > \dim_\C (T) = \text{rank}(X)$.
We define $G_\Q^d (\mc S \rtimes \Gamma)$ analogously, with tempered smooth families of
dimension at least $d$. Now \eqref{eq:SprFamilies} says that 
\begin{equation}
\begin{array}{lll}
\Spr_\Q \big( G_\Q^d (\mc H \rtimes \Gamma) \big) & \subset & G_\Q^d (X \rtimes W_0 \rtimes \Gamma) , \\
\Spr_\Q \big( G_\Q^d (\mc S \rtimes \Gamma) \big) & \subset & G_\Q^d (\mc S (X) \rtimes W_0 \rtimes \Gamma) .
\end{array}
\end{equation}
Let us consider the graded vector spaces associated to these filtrations. 
For tempered representations $\Spr_\Q$ induces $\Q$-linear maps
\begin{equation}\label{eq:SprTempGr}
G_\Q^d (\mc S \rtimes \Gamma) / G_\Q^{d+1} (\mc S \rtimes \Gamma) \to 
G_\Q^d (\mc S (W) \rtimes \Gamma) / G_\Q^{d+1} (\mc S (W) \rtimes \Gamma) .
\end{equation}
We proved in Section \ref{sec:Springer} that
\begin{equation}
\Spr_\Q : G_\Q (\mc S \rtimes \Gamma) \to G_\Q (\mc S (W) \rtimes \Gamma)
\end{equation}
is bijective. Since there are only finitely many smooth families, it is impossible to fill up
a $d$-dimensional tempered smooth family with smooth families of dimension smaller than $d$.
Therefore \eqref{eq:SprTempGr} is bijective for all $d \in \Z_{\geq 0}$. We will show that 
\begin{equation}\label{eq:SprGr}
G_\Q^d (\mc H \rtimes \Gamma) / G_\Q^{d+1} (\mc H \rtimes \Gamma) \to 
G_\Q^d (X \rtimes W_0 \rtimes \Gamma) / G_\Q^{d+1} (X \rtimes W_0 \rtimes \Gamma) 
\end{equation}
is bijective as well. For $d > \dim_\C (T)$ there is nothing to prove, so take $d \leq \dim_\C (T)$. Pick smooth 
$d$-dimensional families $\{ \pi_{i,t} : i \in I, t \in V_i \}$ such that 
\begin{equation}\label{eq:basisSmoothFam}
\big\{ \pi_{i,t} : i \in I , t \in (V_i \cap T_{un}) / \sim \big\}
\end{equation}
is a basis of $G_\Q^d (\mc S \rtimes \Gamma) / G_\Q^{d+1} (\mc S \rtimes \Gamma)$, where $t \sim t'$ 
if and only if $\pi_{i,t} \cong \pi_{i,t'}$. The bijectivity of \eqref{eq:SprTempGr} implies that 
\[
\big\{ \Spr_\Q (\pi_{i,t}) : i \in I , t \in (V_i \cap T_{un}) / \sim \big\}
\]
is a basis of $G_\Q^d (\mc S (W) \rtimes \Gamma) / G_\Q^{d+1} (\mc S (W) \rtimes \Gamma)$. 
The discussion following \eqref{eq:projPdelta} shows that L-smooth $d$-dimensional 
families and tempered smooth $d$-dimensional families are in natural bijection. Hence 
\[
\big\{ L (\pi_{i,t}) : i \in I , t \in V_i / \sim \big\}
\]
is a basis of $G_\Q^d (\mc H \rtimes \Gamma) / G_\Q^{d+1} (\mc H \rtimes \Gamma)$ and
\[
\big\{ \Spr_\Q (\pi_{i,t}) : i \in I , t \in V_i / \sim \big\}
\]
is a basis of $G_\Q^d (X \rtimes W_0 \rtimes \Gamma) / G_\Q^{d+1} (X \rtimes W_0 \rtimes \Gamma)$.
Therefore \eqref{eq:SprGr} is indeed bijective. From this we deduce, with some standard applications 
of the five lemma, that \eqref{eq:SprQ} is a bijection. $\qquad \Box$.

\chapter{Parameter deformations}
\label{chapter:scaling}

Let $\mc H = \mc H (\mc R,q)$ be an affine Hecke algebra associated to an equal parameter 
function $q$. Varying the parameter $q \in \C^\times$ yields a family of algebras, whose 
members are specializations of an affine Hecke algebra with a formal variable $\mathbf q$. 
The Kazhdan--Lusztig-parametrization of Irr$(\mc H (\mc R,q))$ \cite{KaLu,Ree2} 
provides a bijection between the irreducible representations of $\mc H (\mc R,q)$ and 
$\mc H (\mc R,q')$, as long as $q,q' \in \C^\times$ are not roots of unity. Moreover, 
every $\pi_q \in \mr{Irr}(\mc H (\mc R,q))$ is part of a family of representations 
$\{ \pi_q : q \in \C^\times \}$ which depends algebraically on $q$.

It is our conviction that a similar structure underlies the representation theory of
affine Hecke algebras with unequal parameters. However, at present a proof seems to be
out of reach. We have more control when we restrict ourselves to positive parameter functions
and to parameter deformations of the form $q \mapsto q^\ep$ with $\ep \in \R$. We call this
scaling of the parameter function, because it corresponds to multiplying the parameters of
a graded Hecke algebra with $\ep$. Notice that $\mc H (\mc R,q^0) = \C [W]$.

We can relate representations of $\mc H (\mc R,q)$ to $\mc H (\mc R,q^\ep)$-representations
by applying a similar scaling operation on suitable subsets of the space 
$T \cong \mr{Irr}(\mc A)$. We construct a family of functors 
\[
\se : \mr{Mod}_f (\mc H (\mc R,q)) \to \mr{Mod}_f (\mc H(\mc R,q^\ep))
\]
which is an equivalence of categories for $\ep \neq 0$, and which preserves many properties 
of representations (Corollary \ref{cor:4.4}). In particular this provides families of 
representations $\{ \se (\pi) : \ep \in \R \}$ that depend analytically on $\ep$.

The Schwartz algebra $\mc S (\mc R,q)$ behaves even better under scaling of the parameter 
function $q$. As $q$ can be varied in several directions, we have a higher dimensional family
of Fr\'echet algebras $\mc S (\mc R,q)$, which is known to depend continuously on $q$ 
\cite[Appendix A]{OpSo2}. This was exploited for the main results of \cite{OpSo2}, but the 
techniques used there to study general deformations of $q$ are specific to discrete series 
representations. 

To get going at the other series, we only scale $q$.
Via a detailed study of the Fourier transform of $\mc S (\mc R,q)$ 
(see Theorem \ref{thm:3.8}) we construct homomorphisms of Fr\'echet *-algebras
\[
\Sc_\ep : \mc S (\mc R,q^\ep ) \to \mc S (\mc R,q) \qquad \ep \in [0,1],
\]
which depend piecewise analytically on $\ep$ and are isomorphisms for $\ep > 0$ 
(Theorem \ref{thm:4.8}). It is not known whether this is possible with $\mc H (\mc R,q)$ 
instead of $\mc S (\mc R,q)$, when $q$ is not an equal parameter function.

The most remarkable part is that these maps extend continuously
to $\ep = 0$, that is, to a map $\Sc_0 : \mc S (W) \to \mc S (\mc R,q)$. Of course $\Sc_0$
cannot be an isomorphism, but it is injective and and some ways behaves like an isomorphism. 
In fact, we show that for irreducible tempered $\mc H (\mc R,q)$-representations the map
$\Spr$ from Section \ref{sec:Springer} agrees with the functors 
$\tilde \sigma_0$ and $\pi \mapsto \pi \circ \Sc_0$ (Corollary \ref{cor:Sprsigma0}).

\section{Scaling Hecke algebras}
\label{sec:scarep}

As we saw in Section \ref{sec:Springer}, there is a correspondence between tempered representations 
of $\mc H \rtimes \Gamma$ and of $W \rtimes \Gamma$. On central characters this correspondence
has the effect of forgetting the real split part and retaining the unitary part. These elements of
$T / (W_0 \rtimes \Gamma)$ are connected by the path $(W_0 \rtimes \Gamma) c^\ep u$, with
$\ep \in [0,1]$. Opdam \cite[Section 5]{Opd-Sp} was the first to realize that one can interpolate
not only the central characters, but also the representations themselves. In this section we will
recall and generalize the relevant results of Opdam. In contrast to the previous sections we will not
include an extra diagram automorphism group $\Gamma$ in our considerations, as the notation
is already heavy enough. However, it can be checked easily that all the results admit obvious 
generalizations with such a $\Gamma$.

First we discuss the situation for graded Hecke algebras, which is considerably easier.
Let $V \subset \mf t$ be a nonempty open $W_0$-invariant subset. Recall the elements
$\tilde \imath_w \in C^{me}(V)^{W_0} \otimes_{Z(\mh H (\tilde{\mc R},k))} \mh H (\tilde{\mc R},k)$
from Proposition \ref{prop:3.4}. Given any $\ep \in \C$ we define a scaling map
\[
\lambda_\ep : V \to \ep V ,\; v \mapsto \ep v .
\]
For $\ep \neq 0$ it induces an algebra isomorphism 
\begin{align}
\nonumber m_\ep : C^{me}(\ep V)^{W_0} \otimes_{Z(\mh H (\tilde{\mc R},\ep k))} \mh H (\tilde{\mc R},k) &
\to C^{me}(V)^{W_0} \otimes_{Z(\mh H (\tilde{\mc R},k))} \mh H (\tilde{\mc R},k) , \\
f \tilde \imath_{w,\ep} & \mapsto (f \circ \lambda_\ep) \tilde \imath_w . \label{eq:defmep}
\end{align}
Let us calculate the image of a simple reflection $s_\alpha \in S_0$:
\begin{align*}
m_\ep (1 + s_\alpha) & = m_\ep \big( \tilde c_{\alpha,\ep}^{-1} (1 + \tilde \imath_{s_\alpha,\ep}) \big) =
m_\ep \Big( \frac{\alpha}{\ep k_\alpha + \alpha} (1 + \tilde \imath_{s_\alpha,\ep}) \Big) 
= \frac{\ep \alpha}{\ep k_\alpha + \ep \alpha} (1 + \tilde \imath_{s_\alpha}) \\
& = \tilde c_\alpha^{-1} (1 + \tilde \imath_{s_\alpha}) = 1 + s_\alpha .
\end{align*}
That is, $m_\ep (w) = w$ for all $w \in W_0 \subset \mh H (\tilde{\mc R}, \ep k)$, so $m_\ep$ is indeed the 
same as \eqref{eq:1.3}. In particular we see that $m_\ep$ can be defined without using meromorphic
functions. These maps have a limit homomorphism
\begin{equation}
\begin{array}{ccccc}
m_0 : & \mh H (\tilde{\mc R},0) = \mc O (\mf t) \rtimes W_0 &
\to & \mh H (\tilde{\mc R},k) , \\
 & f w & \mapsto & f (0) w ,
\end{array}
\end{equation}
with the property that
\[
\C \to C^{me}(V)^{W_0} \otimes_{Z(\mh H (\tilde{\mc R},k))} \mh H (\tilde{\mc R},k) :
\ep \mapsto m_\ep (f w)
\]
is analytic for all $f \in \mc O (\mf t)$ and $w \in W_0$.

From now we assume that $\ep$ is real.
Let $q^\ep$ be the parameter function $q^\ep (w) = q(w)^\ep$, which is well-defined
because $q(w) \in \R_{>0}$ for all $w \in W$. We obtain a family of algebras $\mc H_\ep = 
\mc H (\mc R, q^\ep)$ with $\mc H (\mc R ,q^0) = \C [W]$. To relate representations of these
algebras we use their analytic localizations, as in Section \ref{sec:localiz}.

Let $c_{\alpha,\ep}$ be the $c$-function with respect to $(\mc R ,q^\ep)$, as in \eqref{eq:calpha}. 
For $\ep \in [-1,1] \setminus \{0\}$ the ball $\ep B \subset \mf t$ 
satisfies the conditions \ref{cond:ball} with respect to $u c^\ep \in T$ and the parameter 
function $q^\ep$ (which enters via the function $c_{\alpha,\ep}$). For $\ep = 0$ the point 
$\ep B = \{0\}$ trivially fulfills all the conditions, except that it is not open. We write
\[
U_\ep = W_0 u c^\ep \exp (\ep B) \qquad \ep \in [-1,1] 
\]
and we abbreviate $U := U_1$. We define a $W_0$-equivariant scaling map
\[
\sigma_\ep : U \to U_\ep \,, \; w(uc \exp (b)) \mapsto w (u c^\ep \exp (\ep b)) . 
\]
As was noted in \cite[Lemma 5.1]{Opd-Sp}, $\sigma_\ep$ is an analytic diffeomorphism for $\ep \neq 0$,
while $\sigma_0$ is a locally constant map with range $U_0 = W_0 u$. Let $\imath^0_{w,\ep}$ be the
element constructed in Proposition \ref{prop:3.4}.  By \eqref{eq:isoCme} $\sigma_\ep$ induces an
algebra isomorphism 
\begin{equation}\label{eq:defrhoep}
\begin{array}{ccll}
\rho_\ep : \mc H^{me}_\ep (U_\ep) & \to & \mc H^{me}(U) & \ep \in [-1,1] \setminus \{0\} , \\
f \imath^0_{w,\ep} & \mapsto & (f \circ \sigma_\ep) \imath^0_w ,
\end{array}
\end{equation}
where $f \in C^{me}(U)$ and $w \in W_0$.
We will show that these maps depend continuously on $\ep$ and have a well-defined limit as $\ep \to 0$.

\begin{lem}\label{lem:4.1}
For $\ep \in [-1,1] \setminus \{0\}$ and $\alpha \in R_0$ write 
$d_{\alpha,\ep} = (c_{\alpha, \ep} \circ \sigma_\ep) c_\alpha^{-1} \in C^{me}(U)$. 
This defines a bounded invertible analytic function on $U \times ([-1,1] \setminus \{0\})$, 
which extends to a function on $\overline U \times [-1,1]$ with the same properties.
\end{lem}
\emph{Proof.}
This extends \cite[Lemma 5.2]{Opd-Sp} to $\ep = 0$. By the definition \eqref{eq:calpha}
\begin{multline*}
d_{\alpha, \ep}(t) \quad = \quad {\ds \frac{1 + \theta_{-\alpha} (t)}{1 + \theta_{-\alpha} (\sigma_\ep (t))} \:
\frac{1 + q(s_\alpha)^{-\ep/2} q(t_\alpha s_\alpha)^{\ep /2} \theta_{-\alpha}
(\sigma_\ep (t))}{1 + q(s_\alpha)^{-1/2} q(t_\alpha s_\alpha)^{1/2} \theta_{-\alpha}(t)} \: \times } \\
{\ds \frac{1 - \theta_{-\alpha}(t)}{1 - \theta_{-\alpha}(\sigma_\ep (t))} \:
\frac{1 - q(s_\alpha)^{-\ep/2} q(t_\alpha s_\alpha)^{-\ep /2}
\theta_{-\alpha}(\sigma_\ep (t))}{1 - q(s_\alpha)^{-1/2} q(t_\alpha s_\alpha)^{-1/2}
 \theta_{-\alpha}(t)} } .
\end{multline*}
Let us abbreviate to $d_{\alpha, \ep} = \frac{f_1 f_2}{g_1 g_2} \times \frac{f_3 f_4}{g_3 g_4}$.
We see that $d_{\alpha ,\ep}(t)$ extends to an invertible analytic function on $\overline U \times [-1,1]$ 
if none of the quotients $f_n / g_n$ has a zero or a pole on this domain. By Condition \ref{cond:ball}.c there 
is a unique $b \in w ( \log (c) + \overline B )$ such that $t = w(u ) \exp (b)$. This defines a coordinate system 
on $w ( u c ) \exp (\overline B)$, and $\sigma_\ep (t) = w(u ) \exp (\ep b)$. By Condition \ref{cond:ball}.d, 
if either $f_n (t) = 0$ or $g_n (t) = 0$ for some $t \in w ( uc \exp (\overline B)) \subset \overline U$, then 
$f_n (w(uc)) = g_n (w(uc)) = 0$. One can easily check that in this situation
\[
{\ds \frac{f_n (t)}{g_n (t)} = \left( \frac{1 - e^{-\alpha (b) 
\ep}}{1 - e^{-\alpha (b)}} \right)^{(-1)^n} } .
\]
Again by Condition \ref{cond:ball}.d the only critical points of this function are those for which 
$\alpha (b) = 0$. For $\ep \neq 0$ both the numerator and the denominator have a zero of order 1 at such 
points, so the singularity is removable. For the case $\ep = 0$ we need to have a closer 
look. In our new coordinate system we can write
\begin{multline*}
c_{\alpha ,\ep}(\sigma_\ep (t)) \quad = \quad {\ds \frac{f_2 (t)}{g_1 (t)} \: \frac{f_4 (t)}{g_3 (t)}} \quad =\\
{\ds \frac{u (w^{-1} \alpha) + \big( q(s_\alpha)^{-1/2} 
q(t_\alpha s_\alpha)^{1/2} e^{-\alpha (b)} \big)^\ep}{u (w^{-1} \alpha) + 
e^{-\alpha (b) \ep}} \: \frac{u (w^{-1} \alpha) - \big( q(s_\alpha)^{-1/2} 
q(t_\alpha s_\alpha)^{-1/2} e^{-\alpha (b)} \big)^\ep}{u (w^{-1} \alpha) - e^{-\alpha (b) \ep}} } .
\end{multline*}
What happens when $\ep$ goes to 0 depends on $u (w^{-1} \alpha)$. The rational function $f_2 (t) / g_1 (t)$ 
has no zeros or poles in a neighborhood of $\ep = 0$ if $u (w^{-1} \alpha) \neq -1$, so then
$\lim_{\ep \to 0} f_2 (t) / g_1 (t) = 1$. For $u (w^{-1} \alpha) = -1$ we have $f_2 (t) = g_1 (t) = 0$, 
so by L'Hospital's rule
\begin{multline*}
\lim_{\ep \to 0} \frac{f_2 (t)}{g_1 (t)} \; = \;
\lim_{\ep \to 0} \frac{\frac{\partial}{\partial \ep} f_2 (t)}{\frac{\partial}{\partial \ep}  g_1 (t)} \; = \; 
\lim_{\ep \to 0} \frac{\log \big( q(s_\alpha)^{-1/2} q(t_\alpha s_\alpha)^{1/2} e^{-\alpha (b)} \big)}{-\alpha (b)}
\: \times \\
\frac{\big( q(s_\alpha)^{-1/2} q(t_\alpha s_\alpha)^{1/2} e^{-\alpha (b)} \big)^\ep}{
e^{-\alpha (b) \ep}} \quad = \quad 
\frac{-\log q(s_\alpha)/2 + \log q(t_\alpha s_\alpha)/2 - \alpha (b)}{- \alpha (b)} .
\end{multline*}
Similarly $f_4 (t) / g_3 (t)$ has no zeros or poles in a neighborhood of $\ep = 0$ if $u (w^{-1} \alpha) \neq 1$,
so $\lim_{\ep \to 0} f_4 (t) / g_3 (t) = 1$ in that case. For $u (w^{-1} \alpha) = 1$ we have
$f_4 (t) = g_3 (t) = 0$, so again by L'Hopital:
\begin{multline*}
\lim_{\ep \to 0} \frac{f_4 (t)}{g_3 (t)} \; = \;
\lim_{\ep \to 0} \frac{\frac{\partial}{\partial \ep} f_4 (t)}{\frac{\partial}{\partial \ep} g_3 (t)} \; = \; 
\lim_{\ep \to 0} \frac{- \log \big( q(s_\alpha)^{-1/2} q(t_\alpha s_\alpha)^{-1/2} e^{-\alpha (b)} \big)}{\alpha (b)}
\: \times \\
\frac{\big( q(s_\alpha)^{-1/2} q(t_\alpha s_\alpha)^{-1/2} e^{-\alpha (b)} \big)^\ep}{
e^{-\alpha (b) \ep}} \quad = \quad 
\frac{\log q(s_\alpha)/2 + \log q(t_\alpha s_\alpha)/2 + \alpha (b)}{\alpha (b)} .
\end{multline*}
Summarizing, we have
\[
\lim_{\ep \to 0} c_{\alpha ,\ep}(\sigma_\ep (t)) = 
\left\{ \begin{array}{c@{\qquad\mr{if}\quad}ccc}
1 & u (w^{-1} \alpha )^2 & \neq & 1 \\
{\ds \frac{\alpha (b) + (\log q(s_\alpha) - \log q(t_\alpha s_\alpha) ) / 2}{\alpha (b)} } 
& u (w^{-1} \alpha ) & = & -1 \\
{\ds \frac{\alpha (b) + ( \log q(s_\alpha) + \log q(t_\alpha s_\alpha) ) / 2}{\alpha (b)} } & 
u (w^{-1} \alpha ) & = & 1 .
\end{array} \right.
\]
For later use we record that with \eqref{eq:kualpha} we can interpret this as
\begin{equation}\label{eq:lim0cep}
\lim_{\ep \to 0} c_{\alpha ,\ep}\big( \sigma_\ep (w(u) \exp (b)) \big) = 
(k_{w(u),\alpha} + \alpha (b)) / \alpha (b) = \tilde c_\alpha (b)  \qquad \alpha \in R_{w(u)} .
\end{equation}
Now we know that at least $d_{\alpha ,0} = \lim_{\ep \to 0} d_{\alpha ,\ep}$ exists as a meromorphic 
function on $\overline U$. For $u (w^{-1} \alpha )^2 \neq 1 ,\, d_{\alpha ,0} = c_\alpha^{-1}$ 
is invertible by Condition \ref{cond:ball}.d. For $u (w^{-1} \alpha) = -1$ we have
\[
{\ds d_{\alpha ,0}(t) = \frac{1 - e^{-\alpha (b)}}{\alpha (b)}
\frac{\alpha (b) + \log \big( q(s_\alpha) / q(t_\alpha s_\alpha) \big) / 2}{1 - 
q(s_\alpha)^{-1/2} q(t_\alpha s_\alpha)^{1/2} e^{-\alpha (b)}} \frac{1 + e^{-\alpha 
(b)}}{1 + q(s_\alpha)^{-1/2} q(t_\alpha s_\alpha)^{-1/2} e^{-\alpha (b)}} } ,
\]
while for $u (w^{-1} \alpha ) = 1$
\[
{\ds d_{\alpha ,0}(t) = \frac{1 - e^{-\alpha (b)}}{\alpha (b)}
\frac{1 + e^{-\alpha (b)}}{1 + q(s_\alpha)^{-1/2} 
q(t_\alpha s_\alpha)^{1/2} e^{-\alpha (b)}} 
\frac{\alpha (b) + \log \big( q(s_\alpha) q(t_\alpha s_\alpha) \big) / 2}{1 - 
q(s_\alpha)^{-1/2} q(t_\alpha s_\alpha)^{-1/2} e^{-\alpha (b)}} } .
\]
These expressions define invertible functions by Condition \ref{cond:ball}.c. We conclude 
that $d_{\alpha ,\ep}(t)$ and $d_{\alpha ,\ep}^{-1}(t)$ are indeed analytic functions on 
$\overline U \times [-1,1]$. Since this domain is compact, they are bounded. $\qquad \Box$
\\[2mm]

Lemma \ref{lem:4.1} enables us to show that the maps $\rho_\ep$ preserve analyticity: 

\begin{prop}\label{prop:4.2}
The maps \eqref{eq:defrhoep} restrict to isomorphisms of topological algebras
\[
\rho_\ep : \mc H_\ep^{an}(U_\ep ) \isom \mc H^{an}(U)
\] 
There is a well-defined limit homomorphism
\[
\rho_0 = \lim_{\ep \to 0} \rho_\ep : \mh C [W] \to \mc H^{an}(U)
\] 
such that for every $w \in W$ the map
\[ 
[-1,1] \to \mc H^{an}(U) : \ep \mapsto \rho_\ep (N_w )
\] 
is analytic.
\end{prop}
\emph{Proof.} The first statement is \cite[Theorem 5.3]{Opd-Sp},
but for the remainder we need to prove this anyway. It is clear
that $\rho_\ep$ restricts to an isomorphism between 
$C^{an}(U_\ep )$ and $C^{an}(U)$. For a simple reflection 
$s_\alpha \in S_0$ corresponding to $\alpha \in F_0$ we have
\begin{equation}\label{eq:5.19}
\begin{array}{lll}
N_{s} + q(s)^{-\ep /2} & = & q(s)^{\ep /2} c_{\alpha ,\ep} (\imath^0_{s,\ep} + 1) \\
\rho_\ep (N_s ) & = & q(s)^{\ep /2} (c_{\alpha ,\ep} \circ 
\sigma_\ep )(\imath^0_s + 1) - q(s)^{-\ep /2} \\
& = & q(s)^{(\ep -1)/2}(c_{\alpha ,\ep} \circ \sigma_\ep ) 
c_\alpha^{-1} \big( N_s + q(s)^{-1/2} \big) - q(s)^{-\ep /2} \\
& = & q(s)^{(\ep -1)/2} d_{\alpha ,\ep} \big( N_s + q(s)^{-1/2}
\big) - q(s)^{-\ep /2} \\
\end{array}
\end{equation}
By Lemma \ref{lem:4.1} such elements are analytic in $\ep \in 
[-1,1]$ and $t \in U$, so in particular they belong to 
$\mc H^{an}(U)$. Moreover, since every $d_{\alpha ,\ep}$ is 
invertible and by \eqref{eq:isoCme}, the set $\{ \rho_\ep (N_w ) : w \in W_0 \}$ 
is a basis for $\mc H^{an}(U)$ as a $C^{an}(U)$-module. Therefore 
$\rho_\ep$ restricts to an isomorphism between the topological algebras 
$\mc H_\ep^{an}(U_\ep)$ and $\mc H^{an}(U)$ for $\ep \neq 0$. 

Given any $w \in W$, there is a unique $x \in X^+$ such that $w \in W_0 x W_0$. 
By Lemma \ref{lem:1.2} there exist unique coefficients $c_{w,u,v}(q) \in \C$ such that 
\[
N_w = \sum_{u \in W^x , v \in W_0} c_{w,u,v}(q) N_u \theta_x N_v . 
\]
From \eqref{eq:multrules} we see that in fact 
$c_{w,u,v}(q) \in \Q \big( \{ q(s)^{1/2} : s \in S^\af \} \big)$, 
so in particular these coefficients depend analytically on $q$. Moreover $\rho_\ep (\theta_x ) = 
\theta_x \circ \sigma_\ep$ depends analytically on $\ep$, as a function on $U$, so 
\[
\rho_\ep (N_w ) = \sum_{u \in W^x , v \in W_0} c_{w,u,v}(q^\ep) 
\rho_\ep (N_u) (\theta_x \circ \rho_\ep) \rho_\ep (N_v)
\]
is analytic in $\ep \in [-1,1]$. Thus $\rho_0$ exists as a linear map. But, being 
a limit of algebra homomorphisms, it must also be multiplicative. $\qquad \Box$ 
\\[2mm]

Suppose that $u \in T_{un}$ is $W_0$-invariant, so that $R_u = R_0$ and
\[
\exp_u : W_0 (\log (c) + B) \to U = u W_0 c \exp (B)
\]
is a $W_0$-equivariant bijection. Then clearly
\begin{equation}\label{eq:epexpu}
\sigma_\ep \circ \exp_u = \exp_u \circ \lambda_\ep : \ep W_0 (\log (c) + B) \to U_\ep . 
\end{equation}
Let $\Phi_u$ be as in \eqref{eq:Phiu}, with $V = W_0 \log (c) + B$. It follows from 
\eqref{eq:epexpu}, \eqref{eq:defmep} and \eqref{eq:defrhoep} that
\[
m_\ep \circ \Phi_{u,\ep} = \Phi_u \circ \rho_\ep \qquad \ep \in [-1,1] \setminus \{0\}
\]
as maps $\mc H^{me}_\ep (U_\ep) \to C^{me}(W_0 \log (c) + B)^{W_0} 
\otimes_{Z (\mh H (\tilde {\mc R}_u ,k_u))} \mh H (\tilde {\mc R}_u ,k_u)$.
The maps $\Phi_{u,\ep}$ can also be defined for $\ep \to 0$, simply by 
\[
\Phi_{u,0} (f N_w) = (f \circ \exp_u) w \in C^{an}(\mf t)^{W_0} 
\otimes_{Z (\mh H (\tilde {\mc R}_u ,k_u))} \mh H (\tilde {\mc R}_u ,k_u),
\]
where $f \in C^{an} (T)$ and $w \in W_0$. By \eqref{eq:5.19} for $\alpha \in F_0$
\begin{align*}
m_\ep \circ \Phi_{u,\ep} (N_{s_\alpha}) & = 
m_\ep \circ \Phi_{u,\ep} \big(q(s_\alpha)^{\ep /2} c_{\alpha ,\ep} (\imath^0_{s,\ep} + 1) - q(s_\alpha)^{\ep/2} \big) \\
& = m_\ep \big( q(s_\alpha)^{\ep /2} (c_{\alpha ,\ep} \circ \exp_u) (\tilde \imath_{s_\alpha,\ep} + 1) - q(s_\alpha)^{\ep/2} \big) \\
& =  q(s_\alpha)^{\ep /2} (c_{\alpha ,\ep} \circ \exp_u \circ \lambda_\ep) (\tilde \imath_{s_\alpha} + 1) - q(s_\alpha)^{\ep/2} \\
& = q(s_\alpha)^{\ep /2} (c_{\alpha ,\ep} \circ \exp_u \circ \lambda_\ep) \tilde c_\alpha^{-1} (s_\alpha + 1) - q(s_\alpha)^{\ep/2} 
\end{align*}
Since $R_u = R_0$, \eqref{eq:lim0cep} tells us that 
$\lim\limits_{\ep \to 0}(c_{\alpha ,\ep} \circ \exp_u \circ \lambda_\ep) \tilde c_\alpha^{-1} = 1$. Hence
\[
\lim_{\ep \to 0} m_\ep \circ \Phi_{u,\ep} (N_{s_\alpha}) = s_\alpha = m_0 \circ \Phi_{u,0} (N_{s_\alpha}) .
\]
On the other hand, it is clear that for $f \in \C [X] \cong \mc O (T)$
\[
\lim_{\ep \to 0} m_\ep \circ \Phi_{u,\ep} (f) = \lim_{\ep \to 0} f \circ \exp_u \circ \lambda_\ep = f(u) = m_0 \circ \Phi_{u,0} (f) .
\]
Since $\rho_0 = \lim_{\ep \to 0} \rho_\ep$ exists, we can conclude that
\begin{equation}\label{eq:m0Phi}
m_0 \circ \Phi_{u,0} = \Phi_u \circ \rho_0 : \C [ X \rtimes W_0] \to \mh H (\tilde {\mc R}_u,k_u) .
\end{equation}

\section{Preserving unitarity}

Proposition \ref{prop:4.2} shows that for $\ep \in [-1,1] \setminus \{0\}$ there is an equivalence between 
the categories $\mr{Mod}_{f,U} (\mc H)$ and $\mr{Mod}_{f,U_\ep} (\mc H_\ep)$. It would be nice if this 
equivalence would preserve unitarity, but that is not automatic. In fact these categories are not
always closed under taking duals of $\mc H$-representations.

From \eqref{eq:thetax*} we see that an $\mc H$-representation with central character $W_0 t$ 
can only be unitary if $\overline{t^{-1}} \in W_0 t$, where $\overline{t^{-1}}(x) = \overline{t (-x)}$ 
for $x \in X$. To define a * on $\mc H^{an}(U)$ we must thus replace $U$ by 
\[
U^{\pm 1} := U \cup \{ t^{-1} : t \in U \} .
\] 
Let $\pm W_0$ be the group $\{ \pm 1 \} \times W_0$, which acts on $T$ by $-w (t) = w(t)^{-1}$.
The above means that we need the following strengthening of Condition \ref{cond:ball}.e:
\begin{itemize}
\item[(e')] As Condition \ref{cond:ball}.e, but with $\pm W_0 \rtimes \Gamma$ instead of $\WG$.
\end{itemize}
Lemma \ref{lem:4.1} and Proposition \ref{prop:4.2} remain valid under this condition, with the same proof. 
Equations \eqref{eq:thetax*} and \eqref{eq:defHan} show that the involution from $\mc H$ extends 
naturally to $\mc H^{me}(U^{\pm 1}) \rtimes \Gamma$ by
\begin{equation}\label{eq:def*me}
(N_\gamma f)^* = N_{w_0} (f \circ -w_0) N_{w_0}^{-1} N_{\gamma^{-1}} \qquad  
\gamma \in \WG , f \in C^{me}(U^{\pm 1}) ,
\end{equation}
where $w_0$ is the longest element of $W_0$. According to \cite[(4.56)]{Opd-Sp}
\begin{equation}\label{eq:imath*}
(\imath^0_w )^* = N_{w_0} \prod_{\alpha \in R_0^+ \cap w' R_0^-} \frac{c_\alpha}{c_{-\alpha}}
\imath^0_{w'} N_{w_0}^{-1} ,
\end{equation}
where $w \in W_0$ and $w' = w_0 w^{-1} w_0$. We extend the map from Proposition \ref{prop:4.2} to
$\rho_\ep : \mc H_\ep^{an}(U_\ep ) \rtimes \Gamma \to \mc H^{an}(U) \rtimes \Gamma$ by
defining $\rho_\ep (N_\gamma) = N_\gamma$ for $\gamma \in \Gamma$. Usually the maps $\rho_\ep$ 
do not preserve the *, but this can be fixed. For $\ep \in [-1,1]$ consider the element
\[
M_\ep = \rho_\ep (N_{w_0 ,\ep}^{-1})^* N_{w_0} 
\prod\nolimits_{\alpha \in R_0^+} d_{\alpha ,\ep} \in \mc H^{an}(U)
\]
We will use $M_\ep$ to extend \cite[Corollary 5.7]{Opd-Sp}. 
However, this result contained a small mistake: 
the element $A_\ep$ in \cite{Opd-Sp} is not entirely correct, we replace it by $M_\ep$.

\begin{thm}\label{thm:4.3}
For all $\ep \in [-1,1]$ the element $M_\ep \in \mc H^{an}(U^{\pm 1})$ is invertible, positive 
and bounded. It has a positive square root $M_\ep^{1/2}$ and the map 
$\ep \mapsto M_\ep^{1/2}$ is analytic. The map
\begin{align*}
\tilde \rho_\ep : \mc H_\ep^{an}(U_\ep ) \rtimes \Gamma \to \mc H^{an}(U^{\pm 1}) \rtimes \Gamma \,,\; 
h \mapsto M_\ep^{1/2} \rho_\ep (h) M_\ep^{-1/2}
\end{align*}
is a homomorphism of topological *-algebras, and an isomorphism
if $\ep \neq 0$. For any $w \in W \rtimes \Gamma$ the map
\[ 
[-1,1] \to \mc H^{an}(U^{\pm 1}) \rtimes \Gamma : \ep \mapsto \tilde \rho_\ep (N_w )
\] 
is analytic.
\end{thm}
\emph{Proof.} 
By Lemma \ref{lem:4.1} and Proposition \ref{prop:4.2} the $M_\ep$ are 
invertible, bounded and analytic in $\ep$. Consider, for 
$\ep \neq 0$, the automorphism $\mu_\ep$ of $\mc H^{me}(U^{\pm 1})$ given by
\[
\mu_\ep (h) = \rho_\ep (\rho_\ep^{-1} (h)^* )^* .
\]
We will discuss its effect on three kinds of elements.
Firstly, for $f \in C^{me}(U^{\pm 1})$ we have, by \eqref{eq:def*me} 
and the $W_0$-equivariance of $\sigma_\ep $:
\begin{equation}
\begin{array}{lll}
\mu_\ep (f) & = & \rho_\ep ( (f \circ \sigma_\ep )^* )^* \\ 
& = & \rho_\ep \big( N_{w_0} (f \circ (-w_0) \circ \sigma_\ep )
N_{w_0 ,\ep}^{-1} \big)^* \\
& = & \rho_\ep (N_{w_0 ,\ep}^{-1})^* \big( f \circ -w_0 \big)^*
\rho_\ep (N_{w_0 ,\ep})^* \\
& = & \rho_\ep (N_{w_0 ,\ep}^{-1})^* N_{w_0} f N_{w_0}^{-1}
\rho_\ep (N_{w_0 ,\ep})^* \\
& = & \rho_\ep (N_{w_0 ,\ep}^{-1})^* N_{w_0} \prod\limits_{\alpha \in 
R_0^+} d_{\alpha ,\ep} f \prod\limits_{\alpha \in R_0^+} 
d_{\alpha ,\ep}^{-1} N_{w_0}^{-1} \rho_\ep (N_{w_0 ,\ep})^* \\
& = & M_\ep f M_\ep^{-1} .
\end{array}
\end{equation}
Secondly, suppose that the simple reflections $s$ and $s' = w_0 s w_0 \in S_0$ 
correspond to $\alpha$ and put $-w_0 \alpha \in F_0$. Using Propostion 
\ref{prop:3.4}.b for $\mc H^{me}_\ep (U^{\pm 1})$ and \eqref{eq:imath*} we find
\begin{equation}
\begin{array}{lll}
M_\ep \imath^0_s M_\ep^{-1} & = & \rho_\ep (N_{w_0 ,\ep}^{-1}
)^* N_{w_0} \prod\limits_{\alpha \in R_0^+} d_{\alpha ,\ep} 
\imath^0_s \prod\limits_{\alpha \in R_0^+} d_{\alpha 
,\ep}^{-1} N_{w_0}^{-1} \rho_\ep (N_{w_0 ,\ep})^* \\
& = & \rho_\ep (N_{w_0 ,\ep}^{-1})^* N_{w_0} \imath^0_s
d_{\alpha, \ep}^{-1} d_{-\alpha ,\ep} N_{w_0}^{-1} 
\rho_\ep (N_{w_0 ,\ep})^* \vspace{2mm}\\
& = & \rho_\ep (N_{w_0 ,\ep}^{-1})^* N_{w_0} \imath^0_s {\ds 
\frac{c_{-\alpha}}{c_\alpha}  N_{w_0}^{-1} N_{w_0} 
\frac{c_{\alpha ,\ep} \circ \sigma_\ep}{c_{-\alpha ,\ep} \circ 
\sigma_\ep} } N_{w_0}^{-1} \rho_\ep (N_{w_0 ,\ep})^* \\
& = & \rho_\ep (N_{w_0 ,\ep}^{-1})^* (\imath^0_{s'})^* 
\left({\ds \frac{c_{\alpha' ,\ep} \circ \sigma_\ep}{c_{-\alpha' 
,\ep} \circ \sigma_\ep} }\right)^* \rho_\ep (N_{w_0 ,\ep})^* \\
& = & \left( \rho_\ep (N_{w_0 ,\ep}) {\ds \frac{c_{\alpha' ,\ep} 
\circ \sigma_\ep}{c_{-\alpha' ,\ep} \circ \sigma_\ep} } 
\imath^0_{s'} \rho_\ep (N_{w_0 ,\ep}^{-1}) \right)^* \\
& = & \rho_\ep \left( N_{w_0} {\ds \frac{c_{\alpha'}}{c_{
-\alpha'}} } \imath^0_{s',\ep} N_{w_0,\ep}^{-1} \right)^* \\
& = & \rho_\ep \left( (\imath^0_{s,\ep})^* \right)^*
\;=\; \mu_\ep (\imath^0_s ) .
\end{array}
\end{equation}
Thirdly, for $\gamma \in \Gamma$ by definition
\[
\mu_\ep (N_\gamma) =  \rho_\ep (\rho_\ep^{-1} (N_\gamma)^* )^* =
\rho_\ep (N_\gamma^{-1})^* = (N_\gamma^{-1})^* = N_\gamma .
\]
Since elements of the above three types generate 
$\mc H^{me}(U^{\pm 1}) \rtimes \Gamma$, we conclude that 
\[
\mu_\ep (h) = M_\ep h M_\ep^{-1} \quad 
\text{for all } h \in \mc H^{me}(U^{\pm 1}) \rtimes \Gamma .
\]
Now we can see that
\[ 
\begin{array}{lll}
\rho_\ep \big( N_{w_0 ,\ep}^{-1} \big)^* & 
= & \rho_\ep \big( (N_{w_0 ,\ep}^*)^{-1} \big)^*
\;=\; \rho_\ep \big( (N_{w_0 ,\ep}^{-1})^* \big)^* \\
& = & \mu_\ep \big( \rho_\ep (N_{w_0 ,\ep}^{-1} ) \big) ,
\vspace{2mm}
\;=\; M_\ep \rho_\ep (N_{w_0 ,\ep}^{-1} ) M_\ep^{-1} \\
N_e & = & M_\ep^{-1} \rho_\ep (N_{w_0 ,\ep}^{-1})^* N_{w_0} 
\prod\limits_{\alpha \in R_0^+} d_{\alpha ,\ep} 
\;=\; \rho_\ep (N_{w_0 ,\ep}^{-1}) M_\ep^{-1} N_{w_0} 
\prod\limits_{\alpha \in R_0^+} d_{\alpha ,\ep} , \\
M_\ep & = & N_{w_0} \prod\limits_{\alpha \in R_0^+} d_{\alpha 
,\ep} \rho_\ep (N_{w_0 ,\ep}^{-1}) \;=\; \big( \rho_\ep 
(N_{w_0 ,\ep}^{-1})^* \big( N_{w_0} \prod\limits_{\alpha 
\in R_0^+} d_{\alpha ,\ep} \big)^* \big)^* \\
& = & \big( \rho_\ep (N_{w_0 ,\ep}^{-1})^*  N_{w_0} 
\prod_{\alpha \in R_1^+} d_{\alpha ,\ep} \big)^* 
\;=\; M_\ep^*
\end{array}
\]
Thus the elements $M_\ep$ are Hermitian $\forall \ep \neq 0$. By
continuity in $\ep ,\; M_0$ is also Hermitian. Moreover they are 
all invertible, and $M_1 = N_e$, so they are in fact strictly 
positive. We already knew that the element $\ep \mapsto M_\ep$ of
\[
C^{an} ([-1,1] ; \mc H^{an}(U^{\pm 1})) \cong C^{an} ([-1,1] \times U^{\pm 1})
\otimes_{\mc A} \mc H
\]
is bounded, so we can construct its square root using holomorphic functional calculus 
in the Fr\'echet algebra $C_b^{an}([-1,1] \times U^{\pm 1}) \otimes_{\mc A} \mc H$, 
where the subscript $b$ denotes bounded functions. This ensures that $\ep \to M_\ep^{1/2}$ 
is still analytic. Finally, for $\ep \neq 0$
\begin{equation}
\begin{array}{lll}
\tilde \rho_\ep (h)^* & = & 
\left( M_\ep^{1/2} \rho_\ep (h) M_\ep^{-1/2} \right)^* \\
& = & M_\ep^{-1/2} \rho_\ep (h)^* M_\ep^{1/2} \\
& = & M_\ep^{-1/2} \mu_\ep (\rho_\ep (h^*)) M_\ep^{1/2} \\
& = & M_\ep^{1/2} \rho_\ep (h^* ) M_\ep^{-1/2} 
\;=\; \tilde \rho_\ep (h^* )
\end{array}
\end{equation}
Again this extends to $\ep = 0$ by continuity. $\qquad \Box$
\\[2mm]
\begin{cor}\label{cor:4.4}
For $uc \in U$ and $\ep \in [-1,1]$ there is a family of additive functors
\begin{align*}
& \tilde \sigma_{\ep,uc} : \mr{Mod}_{f,\WG uc} (\mc H \rtimes \Gamma) \to 
\mr{Mod}_{f,\WG \sigma_\ep (uc)} (\mc H_\ep \rtimes \Gamma) , \\
& (\pi ,V) \mapsto (\pi \circ \tilde \rho_\ep ,V)
\end{align*}
with the following properties:
\enuma{
\item for all $w \in W \rtimes \Gamma$ and $(\pi ,V)$ the map $[-1,1] \to \mr{End}_\C (V) : 
\ep \mapsto \tilde \sigma_{\ep,uc}(\pi) (N_w)$ is analytic;
\item for $\ep \neq 0 \,, \tilde \sigma_{\ep,uc}$ is an equivalence of categories;
\item $\tilde \sigma_{\ep,uc}$ preserves unitarity;
\item for $\ep < 0 \,, \tilde \sigma_{\ep,uc}$ exchanges tempered and anti-tempered modules,
where anti-tempered means that $|s|^{-1} \in T^-$ for all $\mc A$-weights $s \in T$;
\item for $\ep \geq 0 \,, \tilde \sigma_{\ep,uc}$ preserves temperedness;
\item for $\ep > 0 \,, \tilde \sigma_{\ep,uc}$ preserves the discrete series.
}
\end{cor}
\emph{Proof.}
Parts (a), (b) and (c) follow immediately from Theorem \ref{thm:4.3}. Let $(\pi,V)$ be a finite 
dimensional $\mc H^{an} (U^{\pm 1}) \rtimes \Gamma$-representation. Conjugation by 
$M_\ep^{1/2}$ does not change the isomorphism class of $\pi$, so $\tilde \sigma_{\ep,uc} (\pi)$ 
has the same $\mc A_\ep$-weights as $\pi \circ \rho_\ep$, which by construction are $\sigma_\ep$ 
of the $\mc A$-weights of $\pi$. Now parts (d), (e) and (f) are obvious consequences of
$|\sigma_\ep (t)| = |t|^\ep. \qquad \Box$
\\[2mm]

As the notation indicates, $\tilde \sigma_{\ep,uc}$ depends on the previously chosen base point $uc$. 
For one $t \in T$ there can exist several possible base points such that $t \in U$, and these could 
in principle give rise to different functors $\tilde \sigma_{\ep,t}$. This ambiguity disappears if 
we restrict to $t = uc$ in Corollary \ref{cor:4.4}. Then
\[
\se := \bigoplus\nolimits_{t \in T / W_0} \tilde \sigma_{\ep,t} :
\mr{Mod}_f (\mc H \rtimes \Gamma) \to \mr{Mod}_f (\mc H_\ep \rtimes \Gamma)
\]
is an additive functor which also has the properties described in Corollary \ref{cor:4.4}.
The functor $\tilde \sigma_\ep$ was already used in \cite[Theorem 1.7]{OpSo1}.

The image of $\tilde \sigma_0$ is contained in $\mr{Mod}_{T_{un}}(\C [W])$, so this map
is certainly not bijective, not even after passing to the associated Grothendieck groups of modules.
Nevertheless $\tilde \sigma_0$ clearly is related to the map $\Spr$ from Section \ref{sec:Springer}, 
in fact these maps agree on irreducible tempered $\mc H \rtimes \Gamma$-modules:

\begin{lem}\label{lem:4.9}
Suppose that $\ep \in [-1,1] ,\; uc \in T_{un} T_{rs}$ and $t \in T^{P(u)}$.
Let $\Gamma'_u$ be a subgroup of $W'_{F_u,u}$ and $\pi \in \mr{Mod}_{f,(W(R_u) \rtimes \Gamma'_u) uc} 
(\mc H (\tilde{\mc R}_u, q_u) \rtimes \Gamma'_u)$. 
\enuma{
\item The following $\mc H_\ep \rtimes \Gamma$-representations are canonically isomorphic:
\[
\begin{array}{l@{\qquad}l}
\se \big( \mr{Ind}_{\mc H (\mc R_u ,q_u) \rtimes \Gamma'_u}^{\mc H \rtimes \Gamma} 
(\pi \circ \phi_t) \big) , &
\mr{Ind}_{\mc H (\mc R_u ,q_u^\ep) \rtimes \Gamma'_u}^{\mc H_\ep \rtimes \Gamma} 
\big( \se (\pi) \circ \phi_{t \, |t|^{\ep -1}}) \big) , \\
\big( \mr{Ind}_{\mc H (\mc R_u ,q_u) \rtimes \Gamma'_u}^{\mc H \rtimes \Gamma} 
(\pi \circ \phi_t) \big) \circ \rho_\ep , &
\mr{Ind}_{\mc H (\mc R_u ,q_u^\ep) \rtimes \Gamma'_u}^{\mc H_\ep \rtimes \Gamma} 
\big( (\pi \circ \phi_t) \circ \rho_\ep \big) .
\end{array}
\]
\item $\tilde \sigma_0 = \Spr$ on $\mr{Irr}(\mc S \rtimes \Gamma)$.
}
\end{lem}
\emph{Proof.} 
(a) By definition $\se$ is given by composition with $\tilde \rho_\ep$, and the difference 
between $\rho_\ep$ and $\tilde \rho_\ep$ is only an inner automorphism. Hence the map
\[
\se \big( \mr{Ind}_{\mc H (\mc R_u ,q_u) \rtimes \Gamma'_u}^{\mc H \rtimes \Gamma} 
(\pi \circ \phi_t) \big) \to 
\big( \mr{Ind}_{\mc H (\mc R_u ,q_u) \rtimes \Gamma'_u}^{\mc H \rtimes \Gamma} 
(\pi \circ \phi_t) \big) \circ \rho_\ep : 
v \mapsto M_\ep^{-1/2} v
\]
is an invertible intertwiner. 

The two representations in the right hand column are naturally isomorphic if and only 
if the $\mc H (\mc R_u ,q_u^\ep) \rtimes \Gamma'_u$-representations
\begin{equation}\label{eq:pitrhoep}
\se (\pi) \circ \phi_{t \, |t|^{\ep -1}} \quad \text{and} \quad
\pi \circ \phi_t \circ \rho_\ep
\end{equation}
are so. Notice that $\phi_t$ and $\phi_{t |t|^{\ep -1}}$ are well-defined because 
$t \in T^{P(u)} \subseteq T^{W(R_u)}$.
As we just showed, the left hand side of \eqref{eq:pitrhoep} is naturally isomorphic to
$\pi \circ \rho_\ep \circ \phi_{t \, |t|^{\ep -1}}$. Applying $\sigma_\ep$ once with 
center $uct$ and once with center $uc$ results in
\[
\sigma_{\ep ,uct} (uct) = u t |t|^{-1} c^\ep |t|^\ep = \sigma_{\ep ,uc} (uc) t |t|^{\ep - 1} .
\]
This implies that $\rho_{\ep ,uc} \circ \phi_{t |t|^{\ep -1}} = \phi_t \circ \rho_{\ep, uct}$,
which shows that the representations \eqref{eq:pitrhoep} can indeed by identified.

Now we turn to the most difficult case, the two representations in the bottom row.
In view of Theorem \ref{thm:2.1}.b 
\[
\big( \mr{Ind}_{\mc H (\mc R_u ,q_u) \rtimes \Gamma'_u}^{\mc H \rtimes \Gamma} 
(\pi \circ \phi_t) \big) \circ \rho_\ep \quad \text{and} \quad
\mr{Ind}_{\mc H (\mc R_u ,q_u^\ep) \rtimes \Gamma'_u}^{\mc H_\ep \rtimes \Gamma} 
\big( \pi \circ \phi_t \circ \rho_\ep \big)
\]
are isomorphic if only if the 
\[
\mc H (\mc R_u ,q_u^\ep)^{an}(U_{\ep,u}) \rtimes W'_{F_u,u} \text{-representation }
\mr{Ind}_{\mc H (\mc R_u ,q_u^\ep) \rtimes \Gamma'_u}^{\mc H (\mc R_u, q_u^\ep ) \rtimes W'_{F_u,u}}  
\big( \pi \circ \phi_t \circ \rho_\ep \big)
\] 
corresponds to the
\[
1_{W'_u uc} (\mc H (\mc R ,q^\ep)^{an}(U_\ep) \rtimes \Gamma) 1_{W'_u uc} \text{-representation }
1_{W'_u uc} \big( \mr{Ind}_{\mc H (\mc R_u ,q_u) \rtimes \Gamma'_u}^{\mc H \rtimes \Gamma} 
(\pi \circ \phi_t) \big) \circ \rho_\ep
\]
via the isomorphism from Theorem \ref{thm:2.1}.a.
It is clear from the definition \eqref{eq:defrhoep} that 
\[
\mr{Ind}_{\mc H (\mc R_u ,q_u^\ep) \rtimes \Gamma'_u}^{\mc H (\mc R_u, q_u^\ep ) 
\rtimes W'_{F_u,u}} \big( (\pi \circ \phi_t) \circ \rho_\ep \big) \cong 
\big( \mr{Ind}_{\mc H (\mc R_u ,q_u^\ep) \rtimes \Gamma'_u}^{\mc H (\mc R_u, q_u^\ep ) 
\rtimes W'_{F_u,u}} (\pi \circ \phi_t) \big) \circ \rho_\ep ,
\]
so it suffices to show that the following diagram commutes:
\[
\begin{array}{ccc}
\mc H (\mc R_u ,q_u^\ep)^{an}(U_{\ep,u}) \rtimes W'_{F_u,u} & \to &
1_{W'_u uc} (\mc H (\mc R ,q^\ep)^{an}(U_\ep) \rtimes \Gamma) 1_{W'_u uc} \\
\downarrow \rho_\ep & & \downarrow \rho_\ep \\
\mc H (\mc R_u ,q_u)^{an}(U_u) \rtimes W'_{F_u,u} & \to &
1_{W'_u uc} (\mc H (\mc R ,q)^{an}(U) \rtimes \Gamma) 1_{W'_u uc}.
\end{array}
\]
For elements of $\mc H (\mc R_u ,q_u^\ep)^{an}(U_{\ep,u})$ this is easy, since the effect
of the vertical arrows is only extension of functions from $U_{\ep,u}$ (resp. $U_u$) to
$U_\ep$ (resp. $U$) by 0. For elements of $W'_{F_u,u}$ the commutativity follows from 
\eqref{eq:firstredgamma} and \eqref{eq:defrhoep}.\\
(b) By Corollary \ref{cor:2.3}.b every irreducible tempered $\mc H \rtimes \Gamma$-representation
$\pi$ is of the form $\mr{Ind}_{\mc H (\mc R_u ,q_u) \rtimes W'_{F_u ,u}}^{\mc H \rtimes \Gamma} 
(\tilde \pi \circ \Phi_u)$ for some $\tilde \pi \in \mr{Irr}_0 \big( \mc H (\mc R_u ,q_u)^{an}(U_u) 
\rtimes W'_{F_u,u} \big)$. By Theorem \ref{thm:2.7}.d
\[
\Spr (\pi) = \mr{Ind}_{X \rtimes W'_u}^{X \rtimes \WG} \big( \C_u \otimes \tilde \pi \big|_{W'_u} \big) .
\]
Using \eqref{eq:m0Phi} we can rewrite
\[
\C_u \otimes \tilde \pi \big|_{W'_u} = \C_u \otimes (\tilde \pi \circ m_0) \big|_{W'_u} =
\tilde \pi \circ m_0 \circ \Phi_{u,0} = \tilde \pi \circ \Phi_u \circ \rho_0 \cong
\tilde \pi \circ \Phi_u \circ \tilde \rho_0 = \tilde \sigma_0 (\tilde \pi \circ \Phi_u) .
\]
Now we can apply part (a):
\[
\Spr (\pi) \cong 
\mr{Ind}_{X \rtimes W'_u}^{X \rtimes \WG} \big( \tilde \sigma_0 (\tilde \pi \circ \Phi_u ) \big) \cong 
\tilde \sigma_0 \big( \mr{Ind}_{\mc H (\mc R_u ,q_u) \rtimes W'_{F_u ,u}}^{\mc H \rtimes \Gamma} 
(\tilde \pi \circ \Phi_u) \big) = \tilde \sigma_0 (\pi) . \qquad \Box
\]
\vspace{1mm}

\section{Scaling intertwining operators}

We will show that the scaling maps $\se$ give rise to scaled versions of the intertwining operators 
$\pi^\Gamma (g,\xi)$. We will use this to study the behaviour of the components of the Fourier
transform of $\mc S \rtimes \Gamma$ under scaling of $q$.

As we remarked at the start of Section \ref{sec:scarep}, the results of that section can easily
be extended to $\mc H \rtimes \Gamma$, and we will use that generality here. Recall that the 
groupoid $\mc G$ from \eqref{eq:GPQ} includes $\Gamma$ and is defined independently of $q$. 
Let us realize the representation 
\[
\pi_\ep^\Gamma (P,\se (\delta ),t) \text{ on } \C [\Gamma W^P] \otimes V_\delta
\text{ as } \mr{Ind}_{\mc H_\ep^P}^{\mc H_\ep \rtimes \Gamma} (\delta \circ \tilde \rho_\ep \circ 
\phi_{t,\ep} ) .
\]
For all $\ep \in [-1,1]$ we obtain algebra homomorphisms
\begin{equation}\label{eq:Fep}
\begin{split}
& \mc F_\ep : \mc H (\mc R,q^\ep) \rtimes \Gamma \to \bigoplus\nolimits_{P,\delta } \mc O (T^P) 
\otimes \mr{End}_\C (\C [\Gamma W^P] \otimes V_\delta) , \\
& \mc F_\ep (h) (P,\delta,t) = \pi_\ep^\Gamma (P,\se (\delta ),t,h) .
\end{split}
\end{equation}
Rational intertwining operators $\pi_\ep^\Gamma (g,P,\delta',t)$ can be defined as in \eqref{eq:defintop} 
for all $\mc H_\ep$-representations of the form $\pi_\ep^\Gamma (P,\delta',t)$,
where $\delta'$ is irreducible but not necessarily discrete series. In particular, for $\ep \neq 0$
the $(P,\delta)$-component of the image of $\mc F_\ep$ is invariant under an action of the group 
\[
\mc G_{P,\se (\delta)} := \{ g \in \mc G : g(P) = P, \se (\delta) \circ \psi_g^{-1} \cong \se (\delta) \}
\]
via such intertwiners. As in \eqref{eq:actSections}, the action is not on polynomial
but on rational sections.

\begin{prop}\label{prop:4.5}
Let $\ep \in [-1,1] \setminus \{0\}$, let $g \in \mc G$ with $g(P) \subset F_0$ and let $\delta'$ be 
a discrete series representation of $\mc H_{g(P)}$ that is equivalent to $\delta \circ \psi_g^{-1}$.
\enuma{
\item The $\mc H_{g(P),\ep}$-representations $\se (\delta')$ and 
$\se (\delta) \circ \psi_{g,\ep}^{-1}$ are unitarily equivalent.
\item $\mc G_{P,\se (\delta)} = \mc G_{P,\delta}$ and $\mc G_{P,\se (\delta),t} = \mc G_{P,\delta,t}$ 
for all $t \in T^P$.
\item The intertwiners $\pi_\ep^\Gamma (g,P,\se (\delta),t) \in \mr{Hom}_\C \big( \C [\Gamma W^P] 
\otimes V_\delta, \C [\Gamma W^{g(P)}] \otimes V_{\delta'} \big)$ depend rationally on $t \in T^P$
and analytically on $\ep$, whenever they are regular.
\item For $t \in T^P_{un}$ the $\pi_\ep^\Gamma (g,P,\se (\delta),t)$ are regular and unitary, and 
$\pi_0^\Gamma (g,P,\tilde \sigma_0 (\delta),t) := 
\lim_{\ep \to 0} \pi_\ep^\Gamma (g,P,\se (\delta),t)$ exists.
}
\end{prop}
\emph{Proof.} (a) First we show that
\begin{equation}\label{eq:psirhoep}
\psi_g \circ \rho_\ep = \rho_\ep \circ \psi_{g,\ep} .
\end{equation}
Write $g = \gamma^{-1} u$ with $u \in K_P$ and $\gamma \in W_0 \rtimes \Gamma$. 
The automorphism $\psi_u$ from \eqref{eq:twistKP} reflects the translation of $T^P$ by 
$u \in T^P_{un}$: that changes $U$ to $u U$, but apart from that it commutes with $\sigma_\ep$.
Hence $\psi_u \circ \rho_\ep = \rho_\ep \circ \psi_{u,\ep}$.

The isomorphism $\psi_\gamma : \mc H^P \to \mc H^{\gamma (P)}$ from \eqref{eq:psigamma} is
more difficult to deal with, because it acts nontrivially on the $N_w$ with $w \in W_P$. 
However, by Propostion \ref{prop:3.4}.d this is the restriction to $\mc H^P$ of the automorphism
\[
\psi_\gamma : h \mapsto \imath^0_\gamma h \imath^0_{\gamma^{-1}}
\text{ of the algebra } \big( \C (T/W_0) \otimes_{Z (\mc H)} \mc H \big) \rtimes \Gamma .
\]
Similarly $\psi_{\gamma,\ep} (h) = \imath^0_{\gamma,\ep} h \imath^0_{\gamma^{-1},\ep}$.
From these formula it is clear that $\psi_\gamma \circ \rho_\ep = \rho_\ep \circ \psi_{\gamma,\ep}$
on $\mc H^{me}_\ep (U_\ep)$. This establishes \eqref{eq:psirhoep}.

Let $I^g_\delta : V_\delta \to V_{\delta'}$ be as in \eqref{eq:Idelta}. We claim that 
\begin{equation}\label{eq:Ideltaep}
v \mapsto I_\delta^g \big( \delta (\psi_g^{-1} (M_\ep^{1/2}) M_\ep^{-1/2}) v \big)
\end{equation}
is an intertwiner between $\se (\delta) \circ \psi_{g,\ep}^{-1}$ and $\se (\delta')$. Indeed,
for $v \in V_\delta$ and $h \in \mc H^{an}_\ep (U_\ep)$:
\begin{align*}
& I_\delta^g \big( \delta (\psi_g^{-1} (M_\ep^{1/2}) M_\ep^{-1/2}) 
\se (\delta) \circ \psi_{g,\ep}(h) v \big) = \\
& I_\delta^g \big( \delta (\psi_g^{-1} (M_\ep^{1/2} M_\ep^{-1/2}) \delta ( M_\ep^{1/2}
\rho_\ep (\psi_{g,\ep}^{-1} h) M_\ep^{-1/2}) v \big) = \\
& I_\delta^g \big( \delta \circ \psi_g^{-1} (M_\ep^{1/2} \rho_\ep (h) M_\ep^{-1/2}) 
\delta (\psi_g^{-1} (M_\ep^{1/2}) M_\ep^{-1/2}) v \big) = \\
& \delta' (M_\ep^{1/2} \rho_\ep (h) M_\ep^{-1/2}) 
I_\delta^g \big( \delta (\psi_g^{-1} (M_\ep^{1/2}) M_\ep^{-1/2}) v \big) = \\
& \se (\delta') (h) I_\delta^g \big( \delta (\psi_g^{-1} (M_\ep^{1/2}) M_\ep^{-1/2}) v \big) .
\end{align*}
Obviously \eqref{eq:Ideltaep} is invertible, so it is an equivalence between the irreducible
representations $\se (\delta) \circ \psi_{g,\ep}^{-1}$ and $\se (\delta')$. Since both are 
unitary, there exists a unitary intertwiner between, which by the irreducibility must be a
scalar multiple of \eqref{eq:Ideltaep}. We define $I^g_{\delta,\ep}$ as the unique positive 
multiple of \eqref{eq:Ideltaep} that is unitary. \\
(b) By part (a)
\[
g (P,\se (\delta),t) = (g(P),\se (\delta) \circ \psi_{g,\ep}^{-1},g(t)) \cong 
(g(P), \se (\delta'),g(t)) \cong (g(P) ,\se (\delta \circ \psi_g^{-1}),g(t)) ,
\]
so the stabilizer of $(P,\se (\delta),t)$ does not depend on $\ep \in [-1,1]$.\\
(c) By Theorem \ref{thm:4.3} $I^g_{\delta,\ep}$ depends analytically on $\ep \in [-1,1]$.
By definition $\imath^0_\gamma$ is rational in $t \in T$ and analytic in $\ep$, away from
the poles. By definition \eqref{eq:defintop}
\[
\pi_\ep^\Gamma (g,P,\se (\delta),t) (h \otimes_{\mc H^P_\ep} v) = 
h \imath^0_\gamma \otimes_{\mc H^{g(P)}_\ep} I^g_{\delta,\ep} (v) ,
\]
so $\pi_\ep^\Gamma (g,P,\se (\delta),t)$ has the required properties.\\
(d) All possible singularities of the intertwining operators come from poles and zero
of the $c$-functions from \eqref{eq:calpha}. By Theorem 
\ref{thm:3.5}.c $\pi_\ep^\Gamma (g,P,\se (\delta),t)$ is unitary for
all $t \in T^P_{un}$ and $\ep \in (0,1]$. In particular all the singularities on this
domain are removable. On the other hand, the explicit formula for $c_\alpha$ shows that
the singularities for $q^\ep \; (\ep \neq 0)$ and $t \in T_{un}$ are the same as those for 
$q$ and $t \in T_{un}$. Therefore $\pi_\ep^\Gamma (g,P,\se (\delta),t)$ is also regular for 
$\ep \in [-1,0)$ and $t \in T^P_{un}$. 

Since these linear maps are analytic in $\ep \in [-1,1] \setminus \{0\}$ and unitary for $\ep > 0$, 
they are also unitary for $\ep < 0$. Hence all the matrix
coefficients of $\pi_\ep^\Gamma (g,P,\se (\delta),t)$ are uniformly bounded on
$T^P_{un} \times [-1,1] \setminus \{0\}$, which implies that every possible singularity
at $\ep = 0$ is removable. In particular $\lim_{\ep \to 0} \pi_\ep^\Gamma (g,P,\se (\delta),t)$
exists. $\qquad \Box$
\\[2mm]
We fix a discrete series representation $\delta$ of $\mc H_P$ and we abbreviate
\[
A_{P,\delta} := C^\infty (T^P_{un}) \otimes_\C (\C [\Gamma W^P] \otimes V_\delta) .
\]
Proposition \ref{prop:4.5} says among others that for $g \in \mc G_{P,\delta}$
\[
[-1,1] \to A_{P,\delta}^\times : \ep \mapsto \pi_\ep^\Gamma (g,P,\se (\delta ), \cdot)
\]
is an analytic map. The group $\mc G_{P,\se (\delta)} = \mc G_{P,\delta}$ acts on $A_{P,\delta}$ by
\begin{equation}\label{eq:Gdeltaact}
(g \cdot_\ep f)(t) = \pi_\ep^\Gamma (g,P,\se (\delta ),g^{-1} t) f (g^{-1} t) 
\pi_\ep^\Gamma (g,P,\se (\delta ), g^{-1}t)^{-1} .
\end{equation}
By construction the $\delta$-component of the image of $\mc F_\ep$ consists of 
$\mc G_{P,\se (\delta)}$-invariant sections for $\ep \neq 0$, and by Proposition \ref{prop:4.5} 
this also goes for $\ep = 0$. We intend to show that the algebras
$A_{P,\delta}^{\mc G_{P,\se (\delta)}}$ for $\ep \in [-1,1]$ are all isomorphic. 
(Although $\mc G_{P,\se (\delta)} = \mc G_{P,\delta}$ we prefer the longer notation here, 
because it indicates which action on $A_{P,\delta}$ we consider.) 
We must be careful when taking invariants, because
\begin{equation}\label{eq:projrepG}
\mc G_{P,\delta} \to A_{P,\delta}^\times : g \mapsto \pi_\ep^\Gamma (g,P,\se (\delta ), \cdot)
\end{equation}
is not necessarily a group homomorphism. However, the lack of multiplicativity is small,
it is only caused by the freedom in the choice of a scalar in \eqref{eq:Idelta}. In other
words, \eqref{eq:projrepG} defines a projective representation of $\mc G_{P,\delta}$ on
$A_{P,\delta}$. Recall \cite[Section 53]{CuRe1} that the Schur multiplier $\mc G^*_{P,\delta}$ is a
finite central extension of $\mc G_{P,\delta}$, with the property that every projective
representation of $\mc G_{P,\delta}$ lifts to a unique linear representation of $\mc G^*_{P,\delta}$.
This means that for every lift $g^* \in \mc G^*_{P,\delta}$ of $g \in \mc G_{P,\delta}$ there is a
unique scalar multiple $\pi_\ep^\Gamma (g^*,P,\se (\delta ),\cdot)$ of 
$\pi_\ep^\Gamma (g,P,\se (\delta ),\cdot)$ such that 
\[
\mc G^*_{P,\delta} \to A_{P,\delta}^\times : g^* \mapsto \pi_\ep^\Gamma (g^*,P,\se (\delta ), \cdot)
\]
becomes multiplicative. Since $\pi_\ep^\Gamma (g,P,\se (\delta ), \cdot)$ is unitary, so
is $\pi_\ep^\Gamma (g^*,P,\se (\delta ), \cdot)$. Notice that $\mc G_{P,\delta}$ and $\mc G^*_{P,\delta}$
fix the same elements of $A_{P,\delta}$, because the action \eqref{eq:Gdeltaact} is defined
via conjugation with $\pi_\ep^\Gamma (g,P,\se (\delta ), \cdot)$. According to 
\cite[Section 53]{CuRe1} the way lift \eqref{eq:projrepG} from $\mc G_{P,\delta}$ to $\mc G^*_{P,\delta}$
is completely determined by the cohomology class of the 2-cocycle
\begin{equation}\label{eq:2cocycle}
\mc G_{P,\delta} \times \mc G_{P,\delta} \to \C^\times : (g_1,g_2) \mapsto
I_{\delta,\ep}^{g_1} I_{\delta,\ep}^{g_2} I^{g_2^{-1} g_1^{-1}}_{\delta,\ep} .
\end{equation}
This cocycle depends analytically on $\ep$ and $\mc G_{P,\delta}$ is a finite group, so the class of 
\eqref{eq:2cocycle} in $H^2 (\mc G_{P,\delta},\C)$ does not depend on $\ep$. (In most cases this cohomology 
class is trivial, but examples are known in which it is nontrivial, see \cite[Section 6.2]{DeOp2}.) 
Hence the ratio between $\pi_\ep^\Gamma (g^*,P,\se (\delta ),\cdot)$ and 
$\pi_\ep^\Gamma (g,P,\se (\delta ),\cdot)$ also depends analytically on $\ep$.

For $g^* \in \mc G^*_{P,\delta}$ we define $\lambda_{g^*} : T^P \to T^P$ by $\lambda_{g^*}(t) = g(t)$. 
In the remainder of this section we will work mostly with $\mc G^*_{P,\delta}$, and to simplify the notation
we will denote its typical elements by $g$ instead of $g^*$.
For $g \in \mc G^*_{P,\delta}$ and $t \in T^P_{un}$ we write
\[
u_{g,\ep}(t) := \pi_\ep^{\Gamma} (g,\lambda_g^{-1} (t)) ,
\]
so that the multiplicativity translates into
\begin{equation}\label{eq:ugepmult}
u_{g g',\ep} = u_{g,\ep} (u_{g',\ep} \circ \lambda_g^{-1}) .
\end{equation}
From the above we know that $u_{g,\ep} \in A_{P,\delta}$ is unitary and analytic in $\ep$.
These elements can be used to identify $A_{P,\delta}^{\mc G_{P,\se (\delta)}}$ with a corner in a larger
algebra. Consider the crossed product $A_{P,\delta} \rtimes_\lambda \mc G^*_{P,\delta}$, where the action of
$\mc G^*_{P,\delta}$ on $A_{P,\delta}$ comes only from the action on $C(T^P_{un})$ induced by the 
$\lambda_g$. In particular this action is independent of $\ep$.
On $A_{P,\delta} \rtimes_\lambda \mc G^*_{P,\delta}$ we can define actions of $\mc G^*_{P,\delta}$ by
\begin{equation}
g \cdot_\ep a = u_{g,\ep} g a g^{-1} u_{g,\ep}^{-1} .
\end{equation}
For $a \in A_{P,\delta}$ this recovers the action \eqref{eq:Gdeltaact}. An advantage of introducing the
Schur multiplier is that, by \eqref{eq:ugepmult}, $g \mapsto u_{g,\ep} g$ is a homomorphism from
$\mc G^*_{P,\delta}$ to the unitary group of $A_{P,\delta} \rtimes_\lambda \mc G^*_{P,\delta}$. Hence
\begin{equation}
[-1,1] \to A_{P,\delta} \rtimes_\lambda \mc G^*_{P,\delta} : \ep \mapsto p_{\delta,\ep} :=
|\mc G^*_{P,\delta} |^{-1} \sum\nolimits_{g \in \mc G^*_{P,\delta}} u_{g,\ep} g
\end{equation}
is a family of projections, depending analytically on $\ep$. This was first observed in \cite{Ros} and
worked out in \cite[Lemma A.2]{SolThesis},
\begin{equation}
A_{P,\delta}^{\mc G_{P,\se (\delta)}} \to p_{\delta,\ep} (A_{P,\delta} \rtimes_\lambda \mc G^*_{P,\delta}) 
p_{\delta,\ep} : a \mapsto p_{\delta,\ep} a p_{\delta,\ep}
\end{equation}
is an isomorphism of Fr\'echet *-algebras. Its inverse is the map
\[
\sum\nolimits_{g \in \mc G^*_{P,\delta}} a_g g \mapsto |\mc G^*_{P,\delta}| a_e .
\]
\begin{lem}\label{lem:4.6}
The Fr\'echet *-algebras $A_{P,\delta}^{\mc G_{P,\se (\delta)}} \cong p_{\delta,\ep} 
(A_{P,\delta} \rtimes_\lambda \mc G^*_{P,\delta}) p_{\delta,\ep}$ are all isomorphic, by 
$C^{\infty} (T^P_{un})^{\mc G_{P,\delta}}$-linear isomorphisms that are piecewise analytic in $\ep$.
\end{lem}
\emph{Proof.}
According to \cite[Lemma 1.15]{Phi2} the projections $p_\delta (u_\ep )$ 
are all conjugate, by elements depending continuously on $\ep$. That already proves that all 
the Fr\'echet algebras $A_{P,\delta}^{\mc G_{P,\se (\delta)}}$ are isomorphic. Since $C^{\infty} 
(T^P_{un})^{\mc G_{P,\delta}}$ is the center of both $A_{P,\delta}^{\mc G_{P,\se (\delta)}}$ and
$p_{\delta,\ep} (A_{P,\delta} \rtimes_\lambda \mc G^*_{P,\delta}) p_{\delta,\ep}$, these
isomorphisms are $C^{\infty} (T^P_{un})^{\mc G_{P,\delta}}$-linear.

To show that the isomorphisms can be made piecewise analytic we construct the conjugating elements 
explicitly, using the recipe of \cite[Proposition 4.32]{Bla}.
For $\ep ,\ep' \in [-1,1]$ consider
\[
z(\delta ,\ep ,\ep' ) := (2 p_{\delta,\ep'} ) - 1)(2 p_{\delta,\ep} ) - 1) + 1 
\in A_{P,\delta} \rtimes_\lambda \mc G^*_{P,\delta} .
\]
Clearly this is analytic in $\ep$ and $\ep'$ and
\[
p_{\delta,\ep'} z(\delta ,\ep ,\ep' ) = 2 p_{\delta,\ep'} 
p_{\delta,\ep} = z(\delta ,\ep ,\ep' ) p_{\delta,\ep} . 
\]
The enveloping $C^*$-algebra of $A_{P,\delta} \rtimes_\lambda \mc G^*_{P,\delta}$ is
\[
C_\delta := \mr{End}_\C \big(\C [\Gamma W^P ] \otimes V_\delta \big) \otimes C (T^P_{un}) 
\rtimes_\lambda \mc G^*_{P,\delta} .
\]
Let $\norm{\cdot}$ be its norm and suppose that $\norm{p_\delta (u_\ep ) - p_\delta (u_\ep' )} < 1$. Then
\[
\begin{array}{lll}
\norm{z(\delta ,\ep ,\ep') - 2} & = & \norm{4 p_{\delta,\ep'} p_{\delta,\ep} - 2 p_{\delta,\ep} - 
2 p_{\delta,\ep'} } \\ 
& = & \norm{-2 ( p_{\delta,\ep} - p_{\delta,\ep'} )^2} \\
& \leq & 2 \norm{p_{\delta,\ep} - p_{\delta,\ep'}}^2 \; < \; 2
\end{array}
\]
so $z(\delta ,\ep ,\ep')$ is invertible in $C_\delta$. But 
$A_{P,\delta} \rtimes_\lambda \mc G^*_{P,\delta}$ is closed under the holomorphic functional 
calculus of $C_\delta$, so $z(\delta ,\ep ,\ep')$ is also invertible in this Fr\'echet algebra. 
Moreover, because the Fr\'echet 
topology on $A_{P,\delta} \rtimes_\lambda \mc G^*_{P,\delta}$ is finer than 
the induced topology from $C_\delta$, there exists an open 
interval $I_\ep$ containing $\ep$ such that $z(\delta ,\ep ,\ep')$ is invertible for all $\ep' \in I_\ep$. 
For such $\ep,\ep'$ we construct the unitary element
\[
u(\delta ,\ep ,\ep') := z(\delta ,\ep ,\ep') |z(\delta ,\ep ,\ep')|^{-1} 
\in A_{P,\delta} \rtimes_\lambda \mc G^*_{P,\delta} .
\]
By construction the map
\begin{equation}\label{eq:conjugateu}
p_{\delta,\ep} (A_{P,\delta} \rtimes_\lambda \mc G^*_{P,\delta} ) p_{\delta,\ep} \to 
p_{\delta,\ep'} (A_{P,\delta} \rtimes_\lambda \mc G^*_{P,\delta} ) p_{\delta,\ep'} :
x \mapsto u(\delta ,\ep ,\ep') x u(\delta ,\ep ,\ep')^{-1}
\end{equation}
is an isomorphism of Fr\'echet *-algebras. Notice that \eqref{eq:conjugateu} also
defines an isomorphism between the respective $C^*$-completions.

Let us pick a finite cover $\{I_{\ep_i}\}_{i=1}^m$ of 
$[-1,1]$. Then for every $\ep , \ep' \in [-1,1]$ an isomorphism
between $p_{\delta,\ep} (A_{P,\delta} \rtimes_\lambda \mc G^*_{P,\delta} ) p_{\delta,\ep}$ and
$p_{\delta,\ep'} (A_{P,\delta} \rtimes_\lambda \mc G^*_{P,\delta} ) p_{\delta,\ep'}$ can be obtained by 
composing at most $m$ isomorphisms of the form \eqref{eq:conjugateu}. $\qquad \Box$
\\[1mm]

\section{Scaling Schwartz algebras}

It follows from Corollary \ref{cor:4.4} that, for $\ep \in (0,1]$, the functor $\se$
is an equivalence between the categories of finite dimensional tempered modules of
$\mc H$ and of $\mc H_\ep = \mc H (\mc R ,q^\ep)$. We will combine this with the explicit
description of the Fourier transform of $\mc S (\mc R,q)$ from Theorem \ref{thm:3.8} to
construct "scaling isomorphisms" 
\[
\Sc_\ep : \mc S (\mc R ,q^\ep) \rtimes \Gamma \to \mc S (\mc R,q) \rtimes \Gamma .
\]
These algebra homomorphisms depend continuously on $\ep$ and they turn out to have a
well-defined limit 
\[
\Sc_0 : \mc S (W) \rtimes \Gamma = \mc S (\mc R,q^0) \rtimes \Gamma 
\to \mc S (\mc R,q) \rtimes \Gamma .
\]
Although $\Sc_0$ is no longer surjective, it provides the connection between the representation
theories of $\mc S (\mc R,q)$ and $\mc S (W)$ that we already discussed in Section
\ref{sec:Springer}. 
Recall from \eqref{eq:FourierDelta} that $\Delta$ is a set of representatives for $\mc G$-association
classes of discrete series representations of parabolic subalgebras of $\mc H$, and that
\begin{equation}\label{eq:Fourier}
\mc F : \mc S (\mc R,q) \rtimes \Gamma \to \bigoplus\nolimits_{(P,\delta ) \in \Delta} \big( C^\infty (T^P_{un}) 
\otimes \mr{End}_\C (\C [\Gamma W^P] \otimes V_\delta) \big)^{\mc G_{P,\delta}}
\end{equation}
is an isomorphism of Fr\'echet *-algebras. Together with \eqref{eq:Fep} and Lemma \ref{lem:4.6}, this  
implies the existence of a continuous family of algebra homomorphisms, with some nice properties:

\begin{prop}\label{prop:4.7}
There exists a collection of injective *-homomorphisms 
\[
\Sc_\ep : \mc H (\mc R ,q^\ep) \rtimes \Gamma \to \mc S (\mc R ,q) \rtimes \Gamma
\qquad \ep \in [-1,1] ,
\]
such that:
\enuma{
\item For $\ep < 0$ (respectively $\ep > 0$) the map $\pi \mapsto \pi \circ \Sc_\ep$ is an equivalence 
between $\mr{Mod}_f (\mc S (\mc R,q) \rtimes \Gamma)$ and the category of finite dimensional 
anti-tempered (respectively tempered) $\mc H (\mc R ,q^\ep) \rtimes \Gamma$-modules.
\item $\Sc_1 : \mc H \rtimes \Gamma \to \mc S \rtimes \Gamma$ is the canonical embedding.
\item $\Sc_\ep (N_w) = N_w$ for all $w \in Z(W)$.
\item For all $w \in W$ the map 
\[
[-1,1] \to \mc S (\mc R ,q) \rtimes \Gamma : \ep \mapsto \Sc_\ep (N_w) 
\]
is piecewise analytic, and in particular analytic at 0.
\item For all $\pi \in \mr{Mod}_f (\mc S (\mc R,q) \rtimes \Gamma)$ the representations 
$\pi \circ \Sc_\ep$ and $\se (\pi)$ are equivalent. In particular $\pi^\Gamma (P,\delta,t) 
\circ \Sc_\ep \cong \pi_\ep^\Gamma (P,\se (\delta ),t)$ for all $(P,\delta,t) \in \Xi_{un}$.
}
\end{prop}
\emph{Proof.}
Let $\Sc_{P,\delta,\ep} : A_{P,\delta}^{\mc G_{P,\se (\delta)}} \to A_{P,\delta}^{\mc G_{P,\delta}}$ 
be the isomorphism from Lemma \ref{lem:4.6}. We already observed that $\delta$-component of the 
image of $\mc F_\ep$ is invariant under $\mc G_{P,\se (\delta)}$, so we can define
\[
\Sc_\ep := \mc F^{-1} \circ \big( \bigoplus\nolimits_{(P,\delta) \in \Delta} 
\Sc_{P,\delta,\ep} \big) \circ \mc F_\ep .
\]
Now (b) holds by construction and (c) follows from the 
$C^{\infty} (T^P_{un})^{\mc G_{P,\delta}}$-linearity
in Lemma \ref{lem:4.6}. For (d) we use Theorem \ref{thm:4.3} and Lemma \ref{lem:4.6}. 
From the last lines of the proof of Lemma \ref{lem:4.6} we see how we can arrange that $\Sc_\ep$ 
is analytic at 0: it suffices to take $\ep_i = 0$ and to use the elements $u(\delta,0,\ep')$
for $\ep'$ in a neighborhood of $\ep = 0$.

Any finite dimensional module decomposes canonically as a direct sum of parts corresponding
to different central characters. Hence in (e) it suffices to consider $(\pi,V) \in 
\mr{Mod}_{f,\mc G (P,\delta,t)} (\mc S (\mc R,q) \rtimes \Gamma)$.That is,
$\pi : \mc S (\mc R ,q) \rtimes \Gamma \to \mr{End}_\C (V)$ factors via
\[
\mr{ev}_{P,\delta,t} : \mc S (\mc R,q) \rtimes \Gamma \to \mr{End}_\C (\C [\Gamma W^P] \otimes V_\delta),
\; h \mapsto \pi^\Gamma (P,\delta,t,h) .
\]
As $\Sc_{\ep,P,\delta}$ is $C^{\infty} (T^P_{un})^{\mc G_{P,\delta}}$-linear, $\pi \circ \Sc_\ep :
\mc S (\mc R ,q^\ep) \rtimes \Gamma \to \mr{End}_\C (V)$ factors via 
\[
\mr{ev}_{P,\se (\delta ),t} : \mc S (\mc R,q^\ep) \rtimes \Gamma \to 
\mr{End}_\C (\C [\Gamma W^P] \otimes V_\delta) .
\]
Moreover, by Lemma \ref{lem:4.6} the finite dimensional $C^*$-algebras $\mr{ev}_{P,\delta,t} 
(\mc S (\mc R,q) \rtimes \Gamma)$ and $\mr{ev}_{P,\se (\delta),t} (\mc S (\mc R,q^\ep) \rtimes \Gamma)$
are isomorphic, by isomorphisms depending continuously on $\ep \in [-1,1]$. We have two families
of representations on $V ,\; \pi \circ \Sc_\ep$ and $\se (\pi)$, which agree at $\ep = 1$ and
all whose matrix coefficients are continuous in $\ep$. Since a finite dimensional semisimple
algebra has only finitely many equivalence classes of representations on $V$, such equivalence
classes are rigid under continuous deformations. Therefore $\pi \circ \Sc_\ep \cong \se (\pi)$
for all $\ep \in [-1,1]$. Now Lemma \ref{lem:4.9} (or a simpler version of the above argument) 
shows that $\pi^\Gamma (P,\delta,t) \circ \Sc_\ep \cong \pi_\ep^\Gamma (P,\se (\delta ),t)$,
concluding the proof of (e).

Property (a) is a consequence of (e), Corollary \ref{cor:4.4} and Lemma \ref{lem:3.1}.b.

As concerns the injectivity of $\Sc_\ep$, suppose that $h \in \ker (\Sc_\ep)$. Then $M(t,h) = 0$ 
for all unitary principal series representations $M(t)$. Since $T_{un}$ is Zariski-dense in $T$, 
Lemma \ref{lem:3.13} for $\mc H_\ep \rtimes \Gamma$ shows that $h = 0. \qquad \Box$
\\[2mm]


In general $\Sc_\ep (\mc H_\ep \rtimes \Gamma)$ is not contained in $\mc H \rtimes \Gamma$, for two reasons:
$\Sc_{\delta,\ep}$ usually does not preserve polynomiality, and not every polynomial section is in the image
of $\mc F$. For $\ep \geq 0$ the $\Sc_\ep$ extend continuously to $\mc S (\mc R,q^\ep) \rtimes \Gamma$:

\begin{thm}\label{thm:4.8}
For $\ep \in [0,1]$ there exist homomorphisms of Fr\'echet *-algebras 
\[
\begin{array}{rrrr}
\Sc_\ep : & \mc S (\mc R ,q^\ep ) \rtimes \Gamma & \to & \mc S (\mc R ,q) \rtimes \Gamma \\
\Sc_\ep : & C^* (\mc R ,q^\ep ) \rtimes \Gamma & \to & C^* (\mc R ,q) \rtimes \Gamma
\end{array}
\]
with the following properties:
\enuma{
\item $\Sc_\ep$ is an isomorphism for $\ep > 0$, and $\Sc_0$ is injective.
\item $\Sc_1$ is the identity.
\item $\Sc_\ep (h) = h$ for all $h \in \mc S (Z(W))$.
\item Let $x \in C^* (W \rtimes \Gamma)$ and let $h = \sum_{w \in W \rtimes \Gamma} h_w N_w$
with $p_n (h) < \infty$ for all $n \in \N$. Then the following maps are continuous:
\[
\begin{array}{llllll}
[0,1] & \to & \mc S (\mc R,q) \rtimes \Gamma , & \ep & \mapsto & \Sc_\ep (h) , \\ 

[0,1] & \to & B (L^2 (\mc R ,q)) , & \ep & \mapsto & \Sc_\ep^{-1} \Sc_0 (x) .
\end{array}
\]
\item  For all $\pi \in \mr{Mod}_f (\mc S (\mc R,q) \rtimes \Gamma)$ the representations 
$\pi \circ \zeta_\ep$ and $\se (\pi)$ are equivalent. 
}
\end{thm}
\emph{Proof.}
For any $(P,\delta ) \in \Delta$ the representation $\tilde 
\sigma_\ep (\delta )$, although not necessarily irreducible if
$\ep = 0$, is certainly completely reducible, being unitary.
Hence by Theorem \ref{thm:3.8} every irreducible constituent 
$\pi_1$ of $\tilde \sigma_\ep (\delta )$ is a direct summand of 
\[
\mr{Ind}_{\mc H_{\ep ,P}^{P_1}}^{\mc H_{\ep ,P}} (\delta_1 \circ \phi_{t_1,\ep} )
\]
for a $P_1 \subset P$, a discrete series representation $\delta_1$ of $\mc H (\mc R_{P_1},q^\ep)$ and a
\[
t_1 \in \mr{Hom}_{\mh Z} \left( (X_P)^{P_1} , S^1 \right) =
\mr{Hom}_{\mh Z} \left( X / (X \cap (P^\vee )^\perp + \mh Q P_1 
), S^1 \right) \subset T_u
\]
Consequently, $\pi_\ep (P, \pi_1 ,t)$ is a direct summand of
\[
\mr{Ind}_{\mc H_\ep^P}^{\mc H_\ep} \Big( 
\mr{Ind}_{\mc H_{\ep ,P}^{P_1}}^{\mc H_{\ep ,P}} (\delta_1 
\circ \phi_{t_1 ,\ep}) 
\circ \phi_{t,\ep} \Big) = \pi_\ep (P_1 ,\delta_1 ,t t_1 )
\]
In particular for $t \in T^P_{un}$ every matrix coefficient of $\pi_\ep^\Gamma (P, \se (\delta ),t)$ 
appears in the Fourier transform of $\mc S_(\mc R , q^\ep) \rtimes \Gamma$, so \eqref{eq:Fep} extends to 
the respective Schwartz and $C^*$-completions, as required. By Corollary \ref{cor:4.4}.b and 
\eqref{eq:Fourier} these maps are isomorphisms if $\ep > 0$. Since \eqref{eq:conjugateu} extends
continuously to the appropriate $C^*$-completions, so does the algebra homomorphism $\Sc_\ep$ from
Proposition \ref{prop:4.7}.

Properties (b), (c) and (e) are direct consequences of the corresponding parts of Proposition 
\ref{prop:4.7}.

To see that $\Sc_0$ remains injective we vary on the proof of Proposition \ref{prop:4.7}. 
By (e) the family of representations 
\[
I_t \circ \Sc_\ep \cong 
\pi_\ep^\Gamma (\emptyset , \se (\delta_\emptyset), t) = \pi_\ep^\Gamma (\emptyset, \delta_\emptyset, t)
\]
with $t \in T_{un}$ becomes precisely the unitary principal series of $W \rtimes \Gamma$ when
$\ep \to 0$. By Lemma \ref{lem:2.9} 
and Frobenius reciprocity every irreducible tempered representation of $\mc H (\mc R ,q^0) \rtimes \Gamma = 
\C [W \rtimes \Gamma]$ is a quotient of a unitary principal series representation. Hence every element of 
$C^* (\mc R ,q^0 ) \rtimes \Gamma = C^* (W \rtimes \Gamma)$ that lies in the kernel of $\Sc_0$ annihilates
all irreducible tempered $W \rtimes \Gamma$-representations, and must be 0.

The assumptions in (d) mean that we can consider $h$ as an element of $\mc S (\mc R,q^\ep) \rtimes \Gamma$
for every $\ep$. Moreover the sum $\sum_{w \in W \rtimes \Gamma} h_w N_w$ converges uniformly
to $h$ in $\mc S (\mc R,q^\ep) \rtimes \Gamma$. For every finite partial sum $h'$ the map 
$\ep \mapsto \phi_\ep (h')$ is continuous by Proposition \ref{prop:4.7}.e, so this also holds for $h$ itself.

For $\ep \in (0,1]$ we consider
\begin{equation}\label{eq:Scep0}
\Sc_\ep^{-1} \Sc_0 (x) - x = \mc F_\ep^{-1}
\Big( \bigoplus_{P,\delta \in \Delta} \Sc_{P,\delta,\ep}^{-1} \Sc_{P,\delta,0} (\mc F_0 (x)) - \mc F_\ep (x) \Big) .
\end{equation}
Since $\Sc_{P,\delta,0}$ is invertible, both $ \Sc_{P,\delta,\ep}^{-1} \Sc_{P,\delta,0} (\mc F_0 (x))$ and 
$\mc F_\ep (x)$ are continuous in $\ep$ and converge to $\mc F_0 (x)$ as $\ep \downarrow 0$. The continuity of 
$\mc F_\ep$ with respect to $\ep$ implies that the $\mc F_\ep^{-1}$  are also continuous with respect to $\ep$, 
so $\ep \mapsto \Sc_\ep^{-1} \Sc_0 (x)$ is continuous.

The expression between the large brackets in \eqref{eq:Scep0} also depends continuously on $\ep$ and 
converges to 0 as $\ep \downarrow 0$. Furthermore
\[
\mc F_\ep^{-1} : \bigoplus\nolimits_{(P,\delta ) \in \Delta} \big( C (T^P_{un}) \otimes 
\mr{End}_\C (\C [\Gamma W^P] \otimes V_\delta) \big)^{\mc G_{P,\se (\delta)}} \to B (L^2 (\mc R ,q))
\]
is a homomorphism of $C^*$-algebras, so it cannot increase the norms. Therefore 
\[
\lim_{\ep \downarrow 0} (\Sc_\ep^{-1} \Sc_0 (x) - x) = 0 \text{ in } B (L^2 (\mc R ,q)). \qquad \Box 
\]

The homomorphisms from Theorem \ref{thm:4.8} are by no means natural, the construction involves 
a lot of arbitrary choices. Nevertheless one can prove \cite[Lemma 5.22]{SolThesis} that it is possible
to interpolate continuously between two such homomorphisms, obtained from different choices. In
other words, the homotopy class of $\Sc_\ep : \mc S (\mc R,q^\ep) \rtimes \Gamma \to
\mc S (\mc R ,q) \rtimes \Gamma$ is canonical.

It is quite remarkable that $\Sc_0$ preserves the trace $\tau$, because the measure space
$(\Xi_{un}, \mu_{Pl})$ differs substantially from $T_{un}$ with the Haar measure (which corresponds
to the algebra $\mc S (X) \rtimes \WG$). For example, the former is disconnected and can have isolated
points, while the latter is a connected manifold. So the preservation of the trace already indicates that 
it is not possible to separate all components of $\Xi_{un} / \mc G$ using only elements from the image
of $\Sc_0$.

\begin{cor}\label{cor:Sprsigma0}
For $\pi \in \mr{Irr}(\mc S (\mc R,q) \rtimes \Gamma)$ we have 
$\Spr (\pi) \cong \tilde \sigma_0 (\pi) \cong \pi \circ \Sc_0$, where $\Spr$ is as in Section \ref{sec:Springer}. 
The map $\Sc_0$ induces a linear bijection
\[
G_\Q (\Sc_0) : G_\Q (\mc S (\mc R,q) \rtimes \Gamma) \to G_\Q (\mc S (W) \rtimes \Gamma) .
\]
\end{cor}
\emph{Proof.}
The first claim follows from Theorem \ref{thm:4.8}.e and Lemma \ref{lem:4.9}.b. The second statement
follows from the first and Theorem \ref{thm:2.7}.a $\qquad \Box$

\chapter{Noncommutative geometry}

Affine Hecke algebras have some clear connections with noncommutative geometry. Already classical is the
isomorphism between an affine Hecke algebra (with one formal parameter $\mathbf q$) and the equivariant
$K$-theory of a Steinberg variety, see \cite{Lus-EqK,KaLu,ChGi}. 
Of a more analytic nature is $K$-theory of the $C^*$-completion $C^* (\mc R,q)$ of $\mc H (\mc R,q)$. 
It is relevant for the representation theory of affine Hecke algebras because topological $K$-theory is built 
from finitely generated projective modules. Since $K$-theory tends to be invariant under small perturbations,
it is expected \cite{Ply1,BCH} that $K_* (C^* (\mc R,q))$ does not depend on $q$. We prove this modulo
torsion (Theorem \ref{thm:5.3}).

For the algebra $\mc H (\mc R,q)$ \pch \ is more suitable than $K$-theory. Although \pch \ is not obviously
related to representation theory, there is a link for certain classes of algebras \cite{SolHomGHA}. 
From \cite{BaNi} it is known that $HP_* (\mc H (\mc R,q)) \cong HP_* (\C [W])$ when $q$ is an equal
parameter function, but the proof is by no means straightforward.

We connect these two theories via the Schwartz completion of $\mc H(\mc R,q)$. For this algebra both
topological $K$-theory and \pch \ are meaningful invariants. Notwithstanding the different nature of
the algebras $\mc H (\mc R,q)$ and $\mc S (\mc R,q)$, they have the same \pch \ (Theorem \ref{thm:5.5}).
We deduce the existence of natural isomorphisms
\[
HP_* (\mc H (\mc R,q)) \cong HP_* (\mc S (\mc R,q)) \cong K_* (\mc S (\mc R,q)) \otimes_\Z \C 
\cong K_* (C^* (\mc R,q)) \otimes_\Z \C .
\]
Moreover the scaling maps from Chapter \ref{chapter:scaling} provide isomorphisms from these
vector spaces to the corresponding invariants of group algebras of $W$ (Corollary \ref{cor:5.6}).
Notice the similarity with the ideas of \cite{BHP,SolPadic}.

Our method of proof actually shows that $\mc S (\mc R,q)$ and $\mc S (W)$ are 
\emph{geometrically equivalent} (Lemma \ref{lem:5.7}), a term coined by Aubert, Baum 
and Plymen \cite{ABP1} to formulate a conjecture for Hecke algebras of $p$-adic groups.
This conjecture (which we call the ABP-conjecture) describes the
structure of Bernstein components in the smooth dual of a reductive $p$-adic group. Translated
to affine Hecke algebras this conjecture says among others that the dual of $\mc H$ can
be parametrized with the extended quotient $\widetilde{T} / W_0$. The topological space
Irr$(\mc H (\mc R,q))$, with its central character map to $T / W_0$, should then be obtained from
$\widetilde{T} / W_0$, with its canonical projection onto $T / W_0$, via translating the 
components of $\widetilde{T} / W_0$ in algebraic dependence on $q$. 
We verify the ABP-conjecture for all affine Hecke algebras with positive parameters, 
possibly extended with a group of diagram automorphisms (Theorem \ref{thm:5.9}). 
Hence the ABP-conjecture holds for all Bernstein components of $p$-adic groups, whose 
Hecke algebras are known to be Morita equivalent to affine Hecke algebras.

In the final section we calculate in detail what happens for root data with $R_0$ of type $B_2 / C_2$.
Interestingly, this also shows that the representation theory of $\mc H (\mc R,q)$ seems to behave
very well under general deformations of the parameter function $q$.

\section{Topological $K$-theory}

By means of the canonical basis $\{ N_w : w \in W \rtimes \Gamma \}$ we can identify the topological
vector spaces underlying the algebras $\mc S (\mc R,q) \rtimes \Gamma$ for all positive parameter functions $q$.
It is clear from the rules \eqref{eq:multrules} that the multiplication in affine Hecke algebras depends 
continuously on $q$, in some sense. To make that precise, one endows finite dimensional subspaces of 
$\mc H (\mc R ,q) \rtimes \Gamma$ with their standard topology, and one defines a topology on the space of 
positive parameter functions by identifying them with tuples of positive real numbers. This can be extended to 
the Schwartz completions: by \cite[Lemma A.8]{OpSo2} the multiplication in $\mc S (\mc R ,q) \rtimes \Gamma$ 
is continuous in $q$, with respect to the Fr\'echet topology on $\mc S (\mc R,q) \rtimes \Gamma$. This opens
the possibility to investigate this field of Fr\'echet algebras with analytic techniques that relate the algebras
for different $q$'s, a strategy that was used at some crucial points in \cite{OpSo2}.

We denote the topological $K$-theory of a Fr\'echet algebra $A$ by $K_* (A) = K_0 (A) \oplus K_1 (A)$.
Since $\mc S (\mc R ,q) \rtimes \Gamma$ is closed under the holomorphic functional calculus of $C^* (\mc R,q) 
\rtimes \Gamma$ \cite[Theorem A.7]{OpSo2}, the density theorem \cite[Th\'eor\`eme A.2.1]{Bos} tells us that
the inclusion $\mc S (\mc R ,q) \rtimes \Gamma \to C^* (\mc R ,q) \rtimes \Gamma$ induces an isomorphism
\begin{equation}\label{eq:Kdensity}
K_* (\mc S (\mc R,q) \rtimes \Gamma) \cong K_* (C^* (\mc R ,q) \rtimes \Gamma) .
\end{equation}
$K$-theory for $C^*$-algebras is homotopy-invariant, so it is natural to expect the following:

\begin{conj} \label{conj:5.K}
For any positive parameter function $q$ the abelian groups\\ $K_* (C^* (W) \rtimes \Gamma)$ and
$K_* (C^* (\mc R ,q) \rtimes \Gamma)$ are naturally isomorphic.
\end{conj}

This conjecture stems from Higson and Plymen (see \cite[6.4]{Ply1} and \cite[6.21]{BCH}), at least 
when $\Gamma = \{\mr{id}\}$ and $q$ is constant on $S^\af$. It is similar to the Connes--Kasparov 
conjecture for Lie groups, see \cite[Sections 4--6]{BCH} for more background. Independently Opdam 
\cite[Section 1.0.1]{Opd-Sp} formulated Conjecture \ref{conj:5.K} for unequal parameters. We will 
discuss its relevance for the representation theory of affine Hecke algebras, and we will prove a 
slightly weaker version, obtained by applying the functor $\otimes_\Z \Q$.

Recall \cite{Phi2} that for any unital Fr\'echet algebra $A ,\; K_0 (A)$ (respectively $K_1 (A)$) is generated
by idempotents (respectively invertible elements) in matrix algebras $M_n (A)$. The $K$-groups are obtained
by taking equivalence classes with respect to the relation generated by stabilization and homotopy equivalence.

\begin{lem}\label{lem:5.1}
Suppose that 
\[
\kappa_\ep : K_* (C^* (W) \rtimes \Gamma) \to K_* (C^* (\mc R,q) \rtimes \Gamma) \qquad \ep \in [0,1]
\]
is a family of group homomorphisms with the following property.\\ 
For every idempotent $e_0 \in M_n (\mc S (\mc R ,q^0) \rtimes \Gamma) = M_n (\mc S (W) \rtimes \Gamma)$ 
(resp. invertible element $x_0 \in M_n (\mc S (\mc R ,q^0) \rtimes \Gamma)$) there exists a 
$C^*$-norm-continuous path $\ep \mapsto e_\ep$ (resp. $\ep \mapsto x_\ep$) in the Fr\'echet space 
underlying $M_n (\mc S (\mc R ,q) \rtimes \Gamma)$, such that $\kappa_\ep [e_0] = [e_\ep]$ 
(resp. $\kappa_\ep [x_0] = [x_\ep]$). \\
Then $\kappa_\ep = K_* (\Sc_\ep^{-1} \Sc_0)$, with  $\Sc_\ep : C^* (\mc R,q^\ep) \rtimes \Gamma \to 
C^* (\mc R,q) \rtimes \Gamma$ as in Theorem \ref{thm:4.8}.
\end{lem}
By "$C^*$-norm-continuous" we mean that the path in $B (L^2 (\mc R,q))$ defined by mapping $\ep$
to the operator of left multiplication by $e_\ep$ (with respect to $q^\ep$), is continuous. It follows from
\cite[Proposition A.5]{OpSo2} that every Fr\'echet-continuous path is $C^*$-norm-continuous.
By Theorem \ref{thm:4.8}.d the maps $K_* (\Sc_\ep^{-1} \Sc_0)$ have the property that the $\kappa_\ep$
are supposed to possess, so at least the statement is meaningful. 

\emph{Proof.} 
By definition $K_* (\Sc_\ep^{-1} \Sc_0) [e_0] = [ \Sc_\ep^{-1} \Sc_0 (e_0) ]$, where we extend the $\Sc_\ep$
to matrix algebras over $C^* (\mc R,q^\ep)$ in the obvious way. The paths $\ep \to e_\ep$ and $\ep \to
\Sc_\ep^{-1} \Sc_0 (e_0)$ are both $C^*$-norm-continuous, so we can find $\ep' > 0$ such that
\[
\norm{\Sc_\ep^{-1} \Sc_0 (e_0) - e_\ep} < \norm{2 e_0 - 1}^{-1} = \norm{2 \Sc_\ep^{-1} \Sc_0 (e_0) - 1}^{-1}
\quad \text{for all } \ep \leq \ep' .
\]
Then by \cite[Proposition 4.3.2]{Bla} $e_\ep$ and $\Sc_\ep^{-1} \Sc_0 (e_0)$ are connected by a path of
idempotents in $M_n (C^* (\mc R,q^\ep) \rtimes \Gamma)$, so 
\[
\kappa_\ep [e_0] = [e_\ep] = [ \Sc_\ep^{-1} \Sc_0 (e_0) ] = K_* (\Sc_\ep^{-1} \Sc_0) [e_0]
\quad \text{for all } \ep \leq \ep' .
\]
For $\ep \geq \ep'$ 
\[
K_* (\Sc_\ep^{-1} \Sc_0) [e_0] = K_* (\Sc_{\ep}^{-1} \Sc_{\ep'}) K_* (\Sc_{\ep'}^{-1} \Sc_0) [e_0] = 
K_* (\Sc_{\ep}^{-1} \Sc_{\ep'}) [e_{\ep'}] .
\]
By parts (a) and (d) of Theorem \ref{thm:4.8}, $K_* (\Sc_{\ep}^{-1} \Sc_{\ep'}) \; (\ep \geq \ep')$ is the only
family of maps $K_0 (C^* (\mc R ,q^{\ep'})) \rtimes \Gamma) \to K_0 (C^* (\mc R,q) \rtimes \Gamma)$
that comes from continuous paths of idempotents. 

Now $K_1$. Choose $\ep' > 0$ such that 
\[
\norm{\Sc_\ep^{-1} \Sc_0 (x_0) x_\ep^{-1} - 1} < 1 \quad \text{for all } \ep \leq \ep' .
\]
Then $\Sc_\ep^{-1} \Sc_0 (x_0) x_\ep^{-1}$ is homotopic to $1$ in $GL_n (C^* ( \mc R,q^\ep) \rtimes \Gamma)$, so
\[
K_1 (\Sc_\ep^{-1} \Sc_0) [x_0] = [ \Sc_\ep^{-1} \Sc_0 (x_0) ] = [ x_\ep ] = \kappa_\ep [x_0]  
\quad \text{for all } \ep \leq \ep' .
\]
The argument for $\ep \geq \ep'$ is just as for $K_0. \qquad \Box$
\\[2mm]

This lemma says that the map 
\[
K_* (\Sc_0 ) : K_* (C^* (W) \rtimes \Gamma) \to K_* (C^* (\mc R,q) \rtimes \Gamma )
\]
is natural: it does not really depend on $\Sc_0$, only the topological properties of idempotents and
invertible elements in matrix algebras over $\mc S (\mc R,q^\ep) \rtimes \Gamma$ with $\ep \in [0,1]$. 
To prove that $K_* (\Sc_0)$ becomes an isomorphism after tensoring with $\Q$, we need some preparations
of a more general nature. 

For topological spaces $Y \subset X$ and a topological algebra $A$ we write
\[
C_0 (X,Y;A) = \{ f : X \to A | f \text{ is continuous and } f \big|_Y = 0 \} .
\]
We omit $Y$ (resp. $A$) from the notation if $Y = \emptyset$ (resp. $A = \C$). A $C(X)$-algebra $B$ is a 
$\C$-algebra endowed with a unital algebra homomorphism from $C(X)$ to the center of the multiplier 
algebra of $B$. A morphism of $C(X)$-algebras is an algebra homomorphism that is also a
$C(X)$-module map.

\begin{lem}\label{lem:5.2}
Let $\Sigma$ be a finite simplicial complex, let $A$ and $B$ be $C(\Sigma)$-Banach-algebras and let
$\phi : A \to B$ a morphism of $C(\Sigma)$-Banach-algebras. Suppose that
\enuma{
\item for every simplex $\sigma$ of $\Sigma$ there are finite
dimensional $\C$-algebras $A_\sigma$ and $B_\sigma$ such that
\[
C_0 (\sigma ,\delta \sigma ) A \cong C_0 (\sigma ,\delta \sigma ;A_\sigma ) \quad \mr{and} 
\quad C_0 (\sigma ,\delta \sigma ) B \cong C_0 (\sigma ,\delta \sigma ; B_\sigma ) ;
\]
\item for every $x \in \sigma \setminus \delta \sigma$ the localized map $\phi (x) : A_\sigma \to B_\sigma$ 
induces an isomorphism on $K$-theory.
}
Then $K_* (\phi ) : K_* (A) \isom K_* (B)$ is an isomorphism.
\end{lem}
\emph{Proof.}
Let $\Sigma^n$ be the $n$-skeleton of $\Sigma$ and consider the ideals
\begin{equation}
I_0 = C (\Sigma ) \supset I_1 = C(\Sigma;\Sigma^0) \supset \cdots \supset 
I_n = C_0 (\Sigma ,\Sigma^{n-1} ) \supset \cdots 
\end{equation}
They give rise to ideals $I_n A$ and $I_n B$. Because $\Sigma$ is finite, all these 
ideals are 0 for large $n$. We can identify
\begin{equation}\label{eq:5.ideals}
\begin{split}
& I_n A / I_{n+1} A \cong C_0 (\Sigma^n ,\Sigma^{n-1}) A \cong \\
& \bigoplus_{\sigma \in \Sigma \,:\, \dim \sigma = n} 
A C_0 (\sigma ,\delta \sigma) := 
\bigoplus_{\sigma \in \Sigma \,:\, \dim \sigma = n} 
C_0 (\sigma ,\delta \sigma ;A_\sigma ) ,
\end{split}
\end{equation}
and similarly for $B$. Because $\phi$ is $C(\Sigma )$-linear, it induces homomorphisms
\[
\phi (\sigma ) : C_0 (\sigma ,\delta \sigma ; A_\sigma ) \to
C_0 (\sigma ,\delta \sigma ; B_\sigma ) .
\]
By the additivity of and the excision property of the $K$-functor (see e.g. \cite[Theorem 9.3.1]{Bla}), 
it suffices to show that every $\phi (\sigma )$ induces an isomorphism on $K$-theory. 
Let $x$ be any interior point of $\sigma$. Because $\sigma \setminus \delta \sigma$ is contractible, 
$\phi_\sigma$ is homotopic to $\mr{id}_{C_0 (\sigma ,\delta \sigma )} \otimes \phi (x_\sigma )$. 
By assumption the latter map induces an isomorphism on $K$-theory. With the homotopy invariance 
of $K$-theory it follows that 
\[
K_* (\phi (\sigma )) = K_* \big( \mr{id}_{C_0 (\sigma ,\delta \sigma )} \otimes \phi (x_\sigma ) \big)
\]
is an isomorphism. $\qquad \Box$ 
\\[2mm]

Obviously this lemma is in no way optimal: one can generalize it to larger classes of algebras and one
can relax the finiteness assumption on $\Sigma$. Because we do not need it, we will not bother
to write down such generalizations. What will need however, is that Lemma \ref{lem:5.2} is also
valid for similar functors, in particular for $A \mapsto K_* (A) \otimes_\Z \C$.

\begin{thm}\label{thm:5.3}
The map 
\[
K_* (\Sc_0) \otimes \mr{id}_\Q : K_* (C^* (W) \rtimes \Gamma) \otimes_\Z \Q \to
K_* (C^* (\mc R,q) \rtimes \Gamma) \otimes_\Z \Q
\]
is a $\Q$-linear bijection.
\end{thm}
\emph{Proof.}
Consider the projection
\[
\text{pr} : \Xi_{un} / \mc G \to T_{un} / \WG , \mc G (P,\delta,t) \mapsto \WG r |r|^{-} t ,
\]
where $W(R_P) r \in T_P / W(R_P)$ is the central character of $\delta$. With this projection
and Theorem \ref{thm:3.8} we make $C^* (\mc R ,q)$ into a $C(T_{un}/\WG)$-algebra. By Theorem
\ref{thm:4.8}.e $\Sc_0 : C^* (W) \rtimes \Gamma \to C^* (\mc R ,q) \rtimes \Gamma$ is a 
morphism of $C(T_{un}/\WG) - C^*$-algebras. Choose a triangulation $\Sigma$ of $T_{un}$, such that:
\begin{itemize}
\item $w (\sigma) \in \Sigma$ for every simplex $\sigma \in \Sigma$ and every $w \in \WG$;
\item $T_{un}^G$ is a subcomplex of $\Sigma$, for every subgroup $G \subset \WG$;
\item the star of any simplex $\sigma$ is $W'_{0,\sigma}$-equivariantly contractible.
\end{itemize}
Then $\Sigma / \WG$ is a triangulation of $T_{un}/\WG$. From Theorem \ref{thm:3.8} and
the proof of \cite[Lemma 7]{SolChern} we see that $A = C^* (W) \rtimes \Gamma$ and 
$B = C^* (\mc R,q) \rtimes \Gamma$ and are of the form required in condition (a) of 
Lemma \ref{lem:5.2}. For any $u \in T_{un}$ we write
\[
\begin{array}{lll}
A_u & := & C^* (W) \rtimes \Gamma / \ker I_u , \\
B_u & := & \bigoplus_{\mc G \xi \in \Xi_{un} / \mc G , \mr{pr}(\xi) = \WG u}
C^* (\mc R,q) \rtimes \Gamma / \ker \pi^\Gamma (\xi) .
\end{array}
\]
Condition (b) of Lemma \ref{lem:5.2} for $K_* (?) \otimes_\Z \Q$ 
means that the map $\Sc_0 (\WG u) : A_u \to B_u$ should induce an isomorphism
\begin{equation}\label{eq:KSc0u}
K_* (\Sc_0 (\WG u)) : K_* (A_u) \otimes_Z \Q \to K_* (B_u) \otimes_\Z \Q .
\end{equation}
As for all finite dimensional semisimple algebras, 
\[
K_* (A_u) = K_0 (A_u) = G_\Z (A) \quad \text{and} \quad K_* (B_u) = K_0 (B_u) = G_\Z (B_u) .
\]
With these identifications $K_0 (\Sc_0 (\WG u))$ sends a projective module
$e M_n (A_u)$ to the projective module $\Sc_0 (\WG u)(e) M_n (B_u)$. The free abelian
groups $G_\Z (A_u)$ and $G_\Z (B_u)$ have natural bases consisting of irreducible modules. 
With respect to these bases the matrix of $K_0 (\Sc_0 (\WG u))$ is the transpose of the
matrix of 
\[
\tilde \sigma_0 : G_\Z (B_u) \to G_\Z (A_u) , \pi \mapsto \pi \circ \Sc_0 (\WG u) .
\]
By Theorem \ref{thm:2.7}.a $\tilde \sigma_0 \otimes \mr{id}_\Q : G_\Q (A_u) \to G_\Q (B_u)$ 
is a bijection, so \eqref{eq:KSc0u} is also a bijection.
Now we can apply Lemma \ref{lem:5.2}, which finishes the proof. $\qquad \Box$
\\[2mm]

So we proved Conjecture \ref{conj:5.K} modulo torsion, which raises the question what 
information is still contained in the torsion part. It is known that $K^* (C^* (\mc R,q) 
\rtimes \Gamma)$ is a finitely generated group. Indeed, by \cite[Theorem 6]{SolChern} 
this is the case for all Fr\'echet algebras if the type described in Theorem \ref{thm:3.8}.
Hence the torsion subgroup of $K^* (C^* (\mc R,q) \rtimes \Gamma)$ is finite. In fact 
the author does not know any examples of nontrivial torsion elements in such 
$K$-groups, but it is conceivable that they exist. It turns out that this is related to the 
multiplicities of $\WG$-representations in Section \ref{sec:Springer}, in particular
\eqref{eq:Irr0}.

\begin{lem}\label{lem:5.4}
The following are equivalent:
\enuma{
\item $K_* (\Sc_0) : K_* (C^* (W) \rtimes \Gamma) \to K_* (C^* (\mc R,q) \rtimes \Gamma)$
is an isomorphism.
\item $K_0 (\Sc_0) : K_0 (C^* (W) \rtimes \Gamma) \to K_0 (C^* (\mc R,q) \rtimes \Gamma)$
is surjective.
\item For every $u \in T_{un}$ the map $\Spr$ induces a surjection from the Grothendieck
group of $\mr{Mod}_{f,\WG u T_{rs}} (\mc S (\mc R,q) \rtimes \Gamma)$ to that of $\mr{Mod}_{f,\WG u}
(W \rtimes \Gamma)$.
\item The map $\Spr$ induces a bijection $\Spr_\Z : G_\Z (\mc S (\mc R,q) \rtimes \Gamma) \to 
G_\Z (\mc S (W) \rtimes \Gamma)$.
}
\end{lem}
\emph{Proof.} \\
(a) $\Rightarrow$ (b) Obvious.\\
(b) $\Rightarrow$ (c) We use the notation from the proof of Theorem \ref{thm:5.3}, in particular
$K_0 (A_u)$ and $K_0 (B_u)$ are the Grothendieck groups referred to in (c).  
We claim that the canonical map
\begin{equation}\label{eq:claimK}
K_0 (C^* (\mc R,q) \rtimes \Gamma) \to K_0 (B_u)
\end{equation}
is surjective. Recall that $K_0 (B_u)$ is built from 
idempotents. Given any idempotent $e_u \in M_n (B_u)$ we want to find an idempotent
$e \in M_n (C^* (\mc R ,q) \rtimes \Gamma)$ that maps to it. By Theorem \ref{thm:3.8} this 
means that on every connected component $(P,\delta,T^P_{un}) / \mc G_{P,\delta}$ of 
$\Xi_{un} / \mc G$ we have to find an idempotent $e_{P,\delta}$ in 
\[
M_n \big( C(T^P_{un}) \otimes \mr{End}_\C (\C [\Gamma W^P] \otimes V_\delta ) \big)^{\mc G_{P,\delta}} ,
\]
which in every point of $\mr{pr}^{-1}(\WG u) \cap (P,\delta,T^P_{un}) / \mc G_{P,\delta}$ takes
the value prescribed by $e_u$. Recall that the groupoid $\mc G$ was built from elements of
$\WG$ and from the groups $K_P = T^P \cap T_P$. The latter elements permute the components of
$\Xi_{un}$ freely, so $\mr{pr}^{-1}(\WG u)$ intersects every component of $\Xi_{un}$ in at most
one $\mc G$-association class. Therefore we can always find such a $e_{P,\delta}$, proving
the claim \eqref{eq:claimK}.

Together with assumption (b) this implies that
\[
K_0 (C^* (W) \rtimes \Gamma) \xrightarrow{K_0 (\Sc_0)} K_0 (C^* (\mc R,q) \rtimes \Gamma)
\to K_0 (B_u)
\]
is surjective. The underlying $C^*$-algebra homomorphism factors via $C^* (W) \rtimes \Gamma
\to A_u$, so 
\begin{equation}\label{eq:K0AB}
K_0 (\Sc_0 (\WG u)) : K_0 (A_u) \to K_0 (B_u) .
\end{equation} 
is also surjective.\\ 
(c) $\Rightarrow$ (d) By Corollary \ref{cor:Sprsigma0} $\Spr_\Z (\pi) = \pi \circ \Sc_0$
for all $\pi \in \mr{Mod}_f (\mc S (\mc R,q) \rtimes \Gamma)$. So in the notation of \eqref{eq:KSc0u} 
$\Spr_\Z$ is the direct sum, over all $\WG u \in T_{un}$, of the maps 
\begin{equation}\label{eq:RBAu}
G_\Z (B_u) \to G_\Z (A_u) : \pi \mapsto \pi \circ \Sc_0 (\WG u) .
\end{equation}
As we noticed in the proof of Theorem \ref{thm:5.3}, the matrix of this map is the transpose of
the matrix of \eqref{eq:K0AB}.
We showed in the aforementioned proof that the latter map becomes an isomorphism
after applying $\otimes_\Z \Q$. As $K_0 (A_u)$ and $K_0 (B_u)$ are free abelian groups,
this implies that $K_0 (\Sc_0 (\WG u))$ is injective. So under assumption (c) \eqref{eq:K0AB} 
is in fact an isomorphism. Hence, with respect to the natural bases it is given by an integral 
matrix with determinant $\pm 1$. Then the same goes for \eqref{eq:RBAu}, so that map is also 
bijective. Therefore $\Spr_\Z$ is bijective. \\
(d) $\Rightarrow$ (a) The above shows that under assumption (d) the maps \eqref{eq:RBAu} and
\eqref{eq:K0AB} are bijections. Since $K_1 (A_u) = K_1 (B_u) = 0$, we may apply Lemma \ref{lem:5.1}. 
$\qquad \Box$ \\[2mm]

By Corollary \ref{cor:2.3}.b and property (d) of Theorem \ref{thm:2.7}, condition (c) of Lemma
\ref{lem:5.4} can be reformulated as follows: for all $u \in T_{un}$ the map $\pi \mapsto 
\pi \big|_{W'_u}$ induces a bijection from the Grothendieck group of finite dimensional 
tempered $\mh H (\tilde{\mc R}_u, k_u) \rtimes W'_{F_u,u}$-modules with real central 
character to $G_\Z (W'_u)$.

According to \cite[Corollary 3.6]{Ciu} this statement is valid for all graded Hecke algebras 
of "geometric type". Hence Conjecture \ref{conj:5.K} holds, including torsion, 
for many important examples of affine Hecke algebras.

In particular, let $I$ be an Iwahori subgroup of a split reductive $p$-adic group $G$ with root datum 
$\mc R$, as in Section \ref{sec:padic}. By \cite{Ply2} the completion 
$C^*_r (G,I)$ of $\mc H (G,I)$ is isomorphic to $C^* (\mc R,q)$, where $q$ is some prime power.
It is interesting to combine Conjecture \ref{conj:5.K} with the Baum--Connes conjecture. Let 
$\beta G$ be the affine Bruhat--Tits building of $G$ and identify $\mf a^*$ with an apartment.
The Baum--Connes conjecture for groups like $G$ and $W$ was proven by V. Lafforgue \cite{Laf},
see also \cite{SolPadic} (For $W$ it can of course be done more elementarily.)
We obtain a diagram
\begin{equation}\label{eq:K*BC}
\begin{array}{ccccccc}
\hspace{-5mm} K_*^{W} (\mf a^*) & \to & K_* (C_r^* (W)) & 
\to & K_* (C^* (\mc R,q)) & \to & K_* (C_r^* (G,I)) \\
& & & & & & \downarrow \; \uparrow \\
& & K_*^G (\beta G) & \to & K_* (C_r^* (G)) & \to & \bigoplus_{\mf s \in \mf B (G)} K_* (C_r^* (G)_{\mf s})
\end{array}
\end{equation}
in which all the horizontal maps are natural isomorphisms, while the vertical maps pick the factor
of $K_* (C_r^* (G))$ corresponding to the Iwahori-spherical component in $\mf B (G)$. For the
group $G = GL_n (\mathbb F)$ this goes back to \cite{Ply1}. Notice that \eqref{eq:K*BC} realizes
$K_*^{W}(\mf a^*)$ as a direct summand of $K_*^G (\beta G)$, which is by no means
obvious in equivariant $K$-homology.

\section{Periodic cyclic homology}
\label{sec:pch}

Periodic cyclic homology is rather similar to topological $K$-theory, but the former functor
is defined on larger classes of algebras. For example one can take the \pch \ of nontopological
algebras like $\mc H (\mc R,q)$, while it is much more difficult to make sense of the topological
$K$-theory of affine Hecke algebras without completing them. By definition the \pch \
of an algebra over a field $\mathbb F$ is an $\mathbb F$-vector space.
Whereas topological $K$-theory for $C^*$-algebras is the generalization of $K$-theory for
topological spaces, \pch \ for noncommutative algebras can be regarded as the analogue of 
De Rham cohomology for manifolds. 

In \cite[Theorem 3.3]{SolHomGHA} the author proved with homological-algebraic techniques
that the \pch \ of an (extended) graded Hecke algebra $\mh H (\tilde{\mc R},k) \rtimes \Gamma$
does not depend on the parameter function $k$. Subsequently he translated this into a 
representation-theoretic statement, which we already used in \eqref{eq:Irr0}: the collection
of irreducible tempered $\mh H (\tilde{\mc R},k) \rtimes \Gamma$-representations with real 
central character forms a basis of $G_\Q (W_0 \rtimes \Gamma)$.

We will devise a reversed chain of arguments for affine Hecke algebras. Via
topological $K$-theory we will use Theorem \ref{thm:2.7} to show that 
$\mc H (\mc R,q) \rtimes \Gamma$ and $\mc S (\mc R,q) \rtimes \Gamma$ have the same \pch , 
and that it does not depend on the (positive) parameter function $q$. The material in this
section can be compared with \cite{BHP,SolPadic}.

Recall that the Chern character is a natural transformation $K_* \to HP_*$, where we write
$HP_* (A) = HP_0 (A) \oplus HP_1 (A)$.
By \eqref{eq:Kdensity} and \cite[Theorem 6]{SolChern} there are natural isomorphisms
\begin{equation}\label{eq:Cherniso}
K_* (C^* (\mc R,q) \rtimes \Gamma) \otimes_\Z \C \leftarrow K_* (\mc S (\mc R,q) \rtimes \Gamma) 
\otimes_\Z \C \to HP_* (\mc S (\mc R,q)) ,
\end{equation}
the first one is induced by the embedding $\mc S(\mc R,q) \rtimes \Gamma \to 
C^* (\mc R,q) \rtimes \Gamma$, the second one by the Chern character. Here and elsewhere
in this paper the \pch \ of topological algebras is always meant with respect to the completed
projective tensor product. (One needs a tensor product to build the differential complex whose
homology is $HP_*$.) By contrast, in the definition of the \pch \ of nontopological algebras we
simply use the algebraic tensor product over~$\C$.

\begin{thm}\label{thm:5.5}
The inclusion $\mc H (\mc R,q) \rtimes \Gamma \to \mc S (\mc R,q) \rtimes \Gamma$ induces an
isomorphism on \pch .
\end{thm}
\emph{Proof.}
In \cite[Theorem 3.3]{SolPadic} the author proved the corresponding result for Hecke algebras
of reductive $p$-adic groups. The proof from \cite{SolPadic} also applies in our setting, the
important representation-theoretic ingredients being Theorem \ref{thm:3.8}, Proposition 
\ref{prop:3.11} and Lemma \ref{lem:3.12}. A sketch of this proof already appeared in \cite{SolPCH}.
$\qquad \Box$ 

\begin{cor}\label{cor:5.6}
There exists a natural commutative diagram
\[
\begin{array}{*{7}{c}}
\!\! HP_* (\C [W] \rtimes \Gamma) & \!\!\! \to \!\!\! & HP_* (\mc S (W) \rtimes \Gamma)  & \!\!\! \leftarrow \!\!\! 
& K_* (\mc S (W) \rtimes \Gamma) & \!\!\! \to \!\!\! & K_* (C^* (W) \rtimes \Gamma) \\
\downarrow & & \downarrow \scriptstyle{HP_* (\Sc_0)} & & \downarrow  \scriptstyle{K_* (\Sc_0)} & & 
\downarrow  \scriptstyle{K_* (\Sc_0)}  \\
\!\! HP_* (\mc H (\mc R,q) \rtimes \Gamma)  & \!\!\! \to \!\!\! & HP_* (\mc S (\mc R,q) \rtimes \Gamma) & 
\!\!\! \leftarrow \!\!\! &  K_* (\mc S (\mc R,q) \rtimes \Gamma) & \!\!\! \to \!\!\! & K_* (C^* (\mc R ,q) \rtimes \Gamma)
\end{array}
\]
After applying $\otimes_\Z \C$ to the $K$-groups, all these maps are isomorphisms.
\end{cor}
\emph{Proof.}
The horizontal maps are induced the inclusion maps 
\[
\mc H (\mc R,q) \rtimes \Gamma \to \mc S (\mc R,q) \rtimes \Gamma \to C^* (\mc R,q) \rtimes \Gamma 
\]
and by the Chern character $K_* \to HP_*$. The vertical maps (expect the leftmost one) are induced by 
the Fr\'echet algebra homomorphisms $\Sc_0$ from Theorem \ref{thm:4.8}. According to \eqref{eq:Cherniso} 
and Theorem \ref{thm:5.5} all the horizontal maps become isomorphisms after tensoring the $K$-groups 
with $\C$. By Lemma \ref{lem:5.1} the maps $K_* (\Sc_0)$ are natural, and by Theorem \ref{thm:5.3} they 
become isomorphisms after applying $\otimes_\Z \C$. The diagram commutes by functoriality, so
$HP_* (\Sc_0)$ is also a natural isomorphism. Finally, we define $HP_* (\C [W] \rtimes \Gamma) \to
HP_* (\mc H (\mc R,q) \rtimes \Gamma)$ as the unique map that makes the entire diagram commute.
$\qquad \Box$

\begin{rem}
Whether the leftmost vertical map comes from a suitable algebra homomorphism $\C [W] \rtimes \Gamma
\to \mc H (\mc R,q) \rtimes \Gamma$ is doubtful, no such homomorphism is known if $q \neq 1$.
\end{rem}
Suppose that $X$ is the weight lattice of $R_0^\vee$, that $q \in \C \setminus \{0\}$ is any complex number 
which is not a root of unity, and that $q (s) = q$ for all $s \in S^\af$. In this setting an isomorphism 
$HP_* (\C [W]) \cong HP_* (\mc H (\mc R,q))$ was already constructed by Baum and Nistor 
\cite[Theorem 11]{BaNi}. Their proof makes essential use of the Kazhdan--Lusztig classification 
\cite[Theorem 7.12]{KaLu} of irreducible $\mc H (\mc R,q)$-representations, and of Lusztig's asymptotic 
Hecke algebra \cite{Lus-C2,Lus-C3}. 

For graded Hecke algebras things are even better than in Corollary \ref{cor:5.6}: 
in \cite[Theorem 3.4]{SolHomGHA} it was proven that not only $HP_* (\mh H (\tilde{\mc R},k) \rtimes 
\Gamma)$, but also the cyclic homology
and the Hochschild homology of $\mh H (\tilde{\mc R},k) \rtimes \Gamma$ are independent of $k$.
Whether or not this can be transferred to $\mc H(\mc R,q)$ is unclear to the author. The point is that the
comparison of $\mh H (\tilde{\mc R},k) \rtimes \Gamma$ with $\mc H (\mc R,q) \rtimes \Gamma$ goes only
via analytic localizations of these algebras. Since the effect of localization on the dual space is
very easy, we can translate the comparison between localized Hecke algebras to a comparison between
their dual spaces. By \cite[Theorem 4.5]{SolHomGHA} the \pch \ of a finite type
algebra essentially depends only on its dual space, so it is not surprising that the parameter independence
of $HP_*$ can be transferred from graded Hecke algebras to affine Hecke algebras,

On the other hand, the Hochschild homology of an algebra changes in a nontrivial way under localization.
Therefore one would in the first instance only find a comparison between the Hochschild homology of
two localized affine Hecke algebras with the same root datum but different parameters $q$. Possibly, 
provided that one would know enough about $HH_* (\mc H (\mc R,q) \rtimes \Gamma)$, one could 
deduce that also this vector space is independent of $q$. We remark that certainly the 
$Z (\mc H (\mc R,q) \rtimes \Gamma)$-module structure of $HH_* (\mc H (\mc R,q) \rtimes \Gamma)$
will depend on $q$, because that is already the case for graded Hecke algebras, see the remark to 
Theorem 3.4 in \cite{SolHomGHA}.

\section{Weakly spectrum preserving morphisms}
\label{sec:weakly}

For the statement and the proof of the Aubert--Baum--Plymen conjecture we need spectrum 
preserving morphisms and relaxed versions of those. These notions were developed in \cite{BaNi,Nis}.
Baum and Nistor work in the category of finite type $\mathbf k$-algebras, where $\mathbf k$ is
the coordinate ring of some complex affine variety. Since we are also interested in certain Fr\'echet
algebras, we work in a larger class of algebras.

We cannot do without some finiteness assumptions, but it suffices to impose them on representations. 
So, throughout this section we assume that for 
all our complex algebras $A$ there exists a $N \in \N$ such that the dimensions of irreducible 
$A$-modules are uniformly bounded by $N$. In particular $\pi \mapsto \ker \pi$ is a bijection from 
Irr$(A)$ to the collection of primitive ideals of $A$. A homomorphism $\phi : A \to B$ between two such
algebras is called spectrum preserving if
\begin{itemize}
\item for every primitive ideal $J \subset B$, there is exactly one primitive ideal $I \subset A$ containing 
$\phi^{-1} (J)$;
\item the map $J \mapsto I$ induces a bijection $\text{Irr}(\phi) : \text{Irr}(B) \to \text{Irr}(A)$. 
\end{itemize}
We can relax these conditions in the following way. Suppose that there exists filtrations 
\begin{equation}\label{eq:filAB}
\begin{array}{ccccccc}
A = A_0 & \supset & A_1 & \supset & \cdots & \supset & A_n = 0, \\
B = B_0 & \supset & B_1 & \supset & \cdots & \supset & B_n = 0
\end{array}
\end{equation}
by two sided ideals, such that $\phi (A_i) \subset B_i$ for all $i$. We call $\phi : A \to B$ weakly spectrum
preserving if all the induced maps $\phi_i : A_i / A_{i+1} \to B_i / B_{i+1}$ are spectrum preserving.
In this case there are bijections 
\begin{align*}
& \sqcup_i \text{Irr}(A_i / A_{i+1}) \to \text{Irr}(A), \\
& \sqcup_i \text{Irr}(B_i / B_{i+1}) \to \text{Irr}(B), \\
& \text{Irr}(\phi) := \sqcup_i \text{Irr}(\phi_i) : \text{Irr}(B) \to \text{Irr}(A) .
\end{align*}
Notice that Irr$(\phi)$ depends not only on $\phi$, but also on the filtrations of $A$ and $B$.

\begin{lem}\label{lem:5.7}
Let $\phi : A \to B$ be a weakly spectrum preserving morphism, and suppose that the 
dimensions of irreducible $B$-modules are uniformly bounded by $N \in \N$.
Then $\mr{Irr}(\phi)^{-1} (V(I)) = V ( \phi (I)^N )$ for every two-sided ideal $I \subset A$.
In particular the bijection $\mr{Irr}(\phi)$ is continuous with respect to the Jacobson topology
(cf. Section \ref{sec:dual}).
\end{lem}
\emph{Proof.}
We proceed with induction to the length $n$ of the filtration. 
For $n=1$ the morphism $\phi$ is spectrum preserving, so the statement reduces to \cite[Lemma 9]{BaNi}. 
For $n >1$ the induction hypothesis applies to the homomorphisms $\phi : A_1 \to B_1$
and $\phi_1 : A / A_1 \to B / B_1$. So for $\pi \in \text{Irr}(B / B_1) \subset \text{Irr}(B)$ we have
\begin{align*}
& \pi \in \text{Irr}(\phi)^{-1} (V(I)) \subset \text{Irr}(B)  \Longleftrightarrow \\
& \pi \in \text{Irr}(\phi_1)^{-1} (V(I + A_1 / A_1)) \subset \text{Irr}(B / B_1)  \Longleftrightarrow \\
& \pi \in V( \phi (I)^N + B_1 / B_1)) \subset \text{Irr}(B / B_1)  \Longleftrightarrow \\
& \pi \in V( \phi (I)^N)) \subset \text{Irr}(B) . 
\end{align*}
A similar argument applies to $\pi \in \text{Irr}(B_1) \subset \text{Irr}(B). \qquad \Box$
\\[2mm]

The automatic continuity of Irr$(\phi)$ enables us to extract a useful map from
the Fr\'echet algebra morphism $\Sc_0$:

\begin{lem}\label{lem:5.8}
The morphism $\Sc_0 : \mc S(W) \rtimes \Gamma \to \mc S (\mc R ,q) \rtimes \Gamma$ 
is weakly spectrum preserving. 
\end{lem}
\emph{Proof.}
We wil make use of the proofs of Lemma \ref{lem:5.2} and Theorem \ref{thm:5.3}.
There we constructed a $\WG$-equivariant triangulation of $T_{un}$, which lead to two-sided ideals
\[
\begin{array}{lll}
I_n = C_0^\infty (\Sigma, \Sigma^{n-1})^{\WG} & \subset & C^\infty (T_{un})^{\WG} , \\
I_n \mc S (\mc R,q) \rtimes \Gamma & \subset & \mc S (\mc R,q) \rtimes \Gamma , \\
I_n \mc S (W) \rtimes \Gamma & \subset & \mc S (W) \rtimes \Gamma .
\end{array}
\]
(Here and below we regard the $n$-skeleton $\Sigma^n$ both as a simplicial complex and \
as a subset of $T_{un}$). It suffices to show that the induced map
\begin{equation} \label{eq:Sc0n}
\Sc_{0,n} : I_n \mc S (W) \rtimes \Gamma / I_{n+1} \mc S (W) \rtimes \Gamma \to
I_n \mc S (\mc R,q) \rtimes \Gamma / I_{n+1} \mc S (\mc R,q) \rtimes \Gamma
\end{equation}
is spectrum preserving, for every $n$. Fortunately the dual spaces of these quotient algebras are
rather simple,  by \eqref{eq:5.ideals} 
\begin{equation}\label{eq:IrrSigma}
\text{Irr} \big( I_n \mc S (\mc R,q) \rtimes \Gamma / I_{n+1} \mc S (\mc R,q) \rtimes \Gamma \big) 
\cong \bigsqcup_{\sigma \in \Sigma / \WG , \dim \sigma = n} (\sigma \setminus \delta \sigma) \times 
\mr{Irr}_{x_\sigma} (\mc S (\mc R,q) \rtimes \Gamma) ,
\end{equation}
where $x_\sigma \in \sigma \setminus \delta \sigma$ and 
$\mr{Irr}_{x_\sigma} (\mc S (\mc R ,q) \rtimes \Gamma)$ denotes the dual space of the algebra
\begin{equation}\label{eq:Sxsigma}
\bigoplus_{\mc G \xi \in \Xi_{un} / \mc G , \mr{pr}(\xi) = \WG x_\sigma}
\mc S (\mc R,q) \rtimes \Gamma / \ker \pi^\Gamma (\xi) .
\end{equation}
By construction $\Sc_{0,n}$ is $C_0^\infty$-linear, so in particular it is linear over
\[
C_0^\infty (\Sigma^n , \Sigma^{n-1}) := I_n / I_{n+1} .
\]
We know from Theorem \ref{thm:2.7}.a and Corollary \ref{cor:Sprsigma0} that $G_\Q (\Sc_{0,n})$ 
is a bijection, so in particular 
\begin{equation}\label{eq:GQSc0}
G_Q (\Sc_{0,n}) : G_\Q ( \mr{Mod}_{x_\sigma} (\mc S (\mc R ,q) \rtimes \Gamma) ) \to
G_\Q ( \mr{Mod}_{x_\sigma} (\mc S (W) \rtimes \Gamma) )
\end{equation}
is a bijection. Any ordering $(\pi_1 ,\pi_2 ,\ldots, \pi_k )$ of $\mr{Irr}_{x_\sigma} (\mc S (\mc R,q) 
\rtimes \Gamma)$ gives rise to a filtration of \eqref{eq:Sxsigma} by ideals
\[
B_i := \bigcap\nolimits_{j=1}^i \ker \pi_j \qquad i=0,1,\ldots,k .
\]
Since we are dealing with two finite dimensional semisimple algebras of the same rank $k$, 
\eqref{eq:GQSc0} can be described completely with a matrix $M \in GL_k (\Z)$. Order
$\mr{Irr}_{x_\sigma} (\mc S (\mc R,q) \rtimes \Gamma)$ and $\mr{Irr}_{x_\sigma} (\mc S (W) 
\rtimes \Gamma)$ such that all the principal minors of $M$ are nonsingular. Then the corresponding
ideals $B_i$ of \eqref{eq:Sxsigma} and $A_i \subset \mc S (W) \rtimes \Gamma / \ker I_{x_\sigma}$ 
are such that $\Sc_{0,n} (x_\sigma)$ induces spectrum
preserving morphisms $A_i / A_{i+1} \to B_i / B_{i+1}$. Hence $\Sc_{0,n}(x_\sigma)$ is weakly
spectrum preserving. 

It follows from this and \eqref{eq:IrrSigma} that for any $n$-dimensional simplex $\sigma \in \Sigma / \WG$
we can construct filtrations by two-sided ideals in
\[
C_0^\infty (\Sigma , \delta \sigma)^{\WG} \mc S (\mc R,q) \rtimes \Gamma / 
C_0^\infty (\Sigma , \sigma)^{\WG} \mc S (\mc R,q) \rtimes \Gamma
\]
and in
\[
C_0^\infty (\Sigma , \delta \sigma)^{\WG} \mc S (W) \rtimes \Gamma / 
C_0^\infty (\Sigma , \sigma)^{\WG} \mc S (W) \rtimes \Gamma ,
\]
with respect to which the map induced by $\Sc_{0,n}$ is weakly spectrum preserving. We do this for all
such simplices $\sigma$, and then \eqref{eq:5.ideals} show that \eqref{eq:Sc0n} is weakly spectrum 
preserving. $\qquad \Box$
\\[2mm]

A related notion that we will use in the next section is geometric equivalence of algebras, as defined
in \cite[Section 4]{ABP1}. The basic idea is to call $A$ and $B$ geometrically equivalent if they are 
Morita-equivalent or if there exists a weakly spectrum preserving morphism $\phi : A \to B$. 
Furthermore two finite type $\mathbf k$-algebras are geometrically equivalent if they only differ by
an algebraic deformation of the $\mathbf k$-module structure. Now one defines geometric 
equivalence to be the equivalence relation (on the category of finite type $\mathbf k$-algebras)
generated by these three elementary moves. 

So whenever two algebras $A$ and $B$ are geometrically equivalent, they are so by various 
sequences of elementary moves. Every such sequence induces a bijection between the dual spaces
of $A$ and $B$, which however need not be continuous, since the map Irr$(\phi)$ from Lemma
\ref{lem:5.7} is usually not a homeomorphism. 
Nevertheless, by \cite[Theorem 8]{BaNi} every weakly spectrum preserving morphism of finite
type algebras $\phi : A \to B$ induces an isomorphism $HP_* (\phi) : HP_* (A) \to HP_* (B)$. The 
other two moves are easily seen to respect \pch , so geometric equivalence implies 
$HP$-equivalence.

\section{The Aubert--Baum--Plymen conjecture}
\label{sec:ABP}

In a series of papers \cite{ABP1,ABP2,ABP3,ABP4} Aubert, Baum and Plymen developed a conjecture 
that describes the structure of Bernstein components in the smooth dual of a reductive $p$-adic
group. We will rephrase this conjecture for affine Hecke algebras, and prove it in that setting.

A central role is played by extended quotients. Let $G$ be a finite group acting continuously
on a topological space $T$. We endow 
\[
\widetilde{T} := \{ (g,t) \in G \times T : g \cdot t = t \}
\]
with the subspace topology from $G \times T$. Then $G$ also acts continuously on $\widetilde T$, by
\[
g \cdot (g',t) = (g g' g^{-1}, g \cdot t) .
\]
The extended quotient of $T$ by $G$ is defined as $\widetilde{T} / G$. It comes with a projection
onto the normal quotient:
\[
\widetilde{T} / G \to T / G : G (g,t) \mapsto G t.
\]
The fiber over $G t \in T /G$ can be identified with the collection $\langle G_t \rangle$ of
conjugacy classes in the isotropy group $G_t$.

The relevance of the extended quotient for representation theory comes from crossed product algebras.
Suppose that $F(T)$ is an algebra of continuous complex valued functions on $T$, which separates
the points of $T$ and is stable under the action of $G$ on $C(T)$. These conditions ensure that
the crossed product $F(T) \rtimes G$ is well-defined. The dual space of this algebra was determined
in classical results that go back to Frobenius and Clifford (see \cite[Section 49]{CuRe1}). Recall that
a $F(T)$-weight of a representation $(\pi,V)$ is an element $t \in T$ such that there exists a
$v \in V \setminus \{0\}$ with $\pi (f) v = f(t) v$ for all $f \in F(T)$. The collection
of irreducible representations with a $F(T)$-weight $t \in T$ is in natural bijection with 
Irr$(G_t)$, via the map
\[
\pi \mapsto \mr{Ind}_{F(T) \rtimes G_t}^{F(T) \rtimes G} \C_t \otimes \pi .
\]
Since $| \text{Irr}(G_t) | = | \langle G_t \rangle |$, there exists a bijection
\begin{equation}\label{eq:T/G}
\widetilde{T} / G \to \text{Irr} (F(T) \rtimes G) 
\end{equation}
which maps $G (g,t)$ to representation with a $F(T)$-weight $t$. With a little more work one can
find a continuous bijection. However, it is not natural and a not a homeomorphism, except
in very simple cases.

We return to an extended affine Hecke algebra $\mc H (\mc R,q) \rtimes \Gamma$. As described in
Section \ref{sec:defAHA}, the parameter function $q$ is completely determined by its values on
the quotient set $R_{nr}^\vee / W_0 \rtimes \Gamma$. Let $\mc Q (\mc R)$ be the complex variety
of all maps $R_{nr}^\vee / W_0 \rtimes \Gamma \to \C^\times$. To every $v \in \mc Q (\mc R)$ we 
associate the parameter function $q_{\alpha^\vee} = v (\alpha^\vee )^2$.

\begin{conj}\label{conj:5.ABP}
\textup{(ABP-conjecture for affine Hecke algebras)} \\
Suppose that the subgroup of $\C^\times$ generated by $\{ q(s)^{\pm 1/2} \}$ contains 
no roots of unity (except 1).
\enuma{
\item The algebras $\C [W] \rtimes \Gamma$ and $\mc H (\mc R ,q) \rtimes \Gamma$ are
geometrically equivalent.
\item There exists a canonical isomorphism $HP_* (\C [W] \rtimes \Gamma) \cong
HP_* (\mc H (\mc R ,q) \rtimes \Gamma)$.
\item There exists a continuous bijection 
$\mu : \widetilde{T} / \WG \to \mr{Irr}(\mc H (\mc R,q) \rtimes \Gamma)$ such that
$\mu ( \widetilde{T_{un}} / \WG ) = \mr{Irr}(\mc S (\mc R,q) \rtimes \Gamma)$.
\item For every connected component $c$ of $\widetilde{T} / \WG$ there exists a smooth morphism 
of algebraic varieties $h_c : \mc Q (\mc R) \to T$ with the following properties.

For all components $c$ we have $h_c (1) = 1$, and
\[
\mr{pr}_{q^{1/2}} (\widetilde{T} / \WG - T / \WG) = 
\{ t \in T : M(t) \text{ is reducible } \} / \WG ,
\]
where $\mr{pr}_v : \widetilde{T} / \WG \to T / \WG$ is defined by
\[
\mr{pr}_v (\WG (w,t)) = \WG h_c (v) t \quad \text{for} \quad v \in \mc Q (\mc R), \WG (w,t) \in c . 
\]
Moreover $\mu$ can be chosen such that the central character of $\mu (\WG (w,t))$ 
is $\WG h_c (q^{1/2}) t$ for $\WG (w,t) \in c$.
}
\end{conj}

The assumption on the roots of unity is probably stronger than necessary, but it is difficult to 
predict which roots of unity really cause problems. In any case, it is known that the above statements 
are false for some specific of roots unity, for example $q = -1$ if $\mc R$ is of type $A_1^{(1)}$.

We will verify the ABP-conjecture in the following cases:

\begin{thm} \label{thm:5.9}
Parts (b), (c) and (d) of Conjecture \ref{conj:5.ABP} hold for every extended affine Hecke
algebra with a positive parameter function $q$. Part (a) holds for the Schwartz completions
of the algebras in question.
\end{thm}

Let us discuss the different parts of the ABP-conjecture. As we mentioned at the end of Section
\ref{sec:weakly}, every explicit geometric equivalence gives rise to an isomorphism on \pch .
However, this isomorphism need not be natural, so (b) does not yet follow from (a). In \cite{ABP1,ABP3} 
we see that Aubert, Baum and Plymen have a geometric equivalence via Lusztig's asymptotic 
Hecke algebra \cite{Lus-C3} in mind. The corresponding isomorphism \cite[Theorem 11]{BaNi}
\[
HP_* (\mc H (\mc R,q)) \cong HP_* (\C [W])
\]
can be regarded as canonical, albeit in a rather weak sense. 

Unfortunately, for unequal parameter functions this asymptotic Hecke algebra exists only as a 
conjecture, see \cite[Chapter 18]{Lus-Une}. Neither does the author
know any other way to construct a geometric equivalence between $\mc H (\mc R,q) \rtimes \Gamma$
and $\C [W] \rtimes \Gamma$, so this part of the conjecture remains open for unequal 
parameter functions. As a substitute we offer Lemma \ref{lem:5.8}, which has approximately the
same strength. It is weaker because it concerns only topologically completed versions of the
algebras, but it is stronger in the sense that the geometric equivalence consists of only
one weakly spectrum preserving morphism.

Part (b) of Conjecture \ref{conj:5.ABP} was already dealt with in Corollary \ref{cor:5.6}.
\\[1mm]
\emph{Proof of Conjecture \ref{conj:5.ABP}.c.}\\
Recall the definition of a (tempered) smooth 
family of $\mc H \rtimes \Gamma$-representations from \eqref{eq:SmoothFamily}. 
By \eqref{eq:basisSmoothFam} there exist tempered smooth families $\{ \pi_{i,t} : t \in V_i \}$ 
which together form a basis of $G_\Q (\mc S (\mc R,q) \rtimes \Gamma)$. The parameter space of 
such a family is of the form $V_i = u_i \exp (\mf a^{g_i})$ for some $u_i \in T_{un}, g_i \in \WG$. 
By \eqref{eq:SprTempGr} the number of families with parameter space of the form $g V_i$ for some 
$g \in \WG$, is precisely the number of components of $\widetilde{T_{un}} / \WG$ whose projection 
onto $T_{un} / \WG$ is $\WG V_i$. Hence we can find a continuous bijection 
\[
\widetilde{T_{un}} / \WG \to \{ \pi_{i,t} : t \in V_i , i \in I \} .
\]
This will be our map $\mu$ whenever the image $\pi_{i,t}$ is irreducible, which is the case on
a Zariski-dense subset of $\widetilde{T_{un}} / \WG$. To extend $\mu$ continuously to all 
nongeneric points, we need to find irreducible subrepresentations $\pi'_{i,t} \subset \pi_{i,t}$,
such that $\{ \pi'_{i,t} : t \in V_i \} = \mr{Irr} (\mc S (\mc R,q) \rtimes \Gamma)$. For 0-dimensional
components all point are generic, so there is nothing to do. If we have already defined $\mu$ on all
components of $\widetilde{T_{un}} / \WG$ of dimension smaller than $d$, and $(i,t)$ corresponds to
a nongeneric point $\WG (w,t)$ in a component of dimension $d$, than we choose for $\mu (\WG (w,t))$
any irreducible subrepresentation of $\pi_{i,t}$ that we did not have yet in the image of the 
previously handled components. This process can be carried out completely 
\eqref{eq:SprTempGr}, and yields a continuous bijection 
\begin{equation}
\mu : \widetilde{T_{un}} / \WG \to \mr{Irr} (\mc S (\mc R,q) \rtimes \Gamma).
\end{equation}
From this and Lemmas \ref{lem:5.7} and \ref{lem:5.8} we obtain continuous bijections
\[
\widetilde{T_{un}} / \WG \to \mr{Irr} (\mc S (\mc R,q) \rtimes \Gamma) \to
\mr{Irr} (\mc S(W) \rtimes \Gamma) .
\]
As explained after \eqref{eq:projPdelta} and \eqref{eq:basisSmoothFam}, 
these can be extended canonically to continuous bijections 
\begin{equation}
\widetilde{T} / \WG \to \mr{Irr} (\mc H (\mc R,q) \rtimes \Gamma) \to \mr{Irr} (W \rtimes \Gamma).
\end{equation}
\emph{Proof of Conjecture \ref{conj:5.ABP}.d.} \\
Suppose that $\WG (w,t_0) \in \widetilde{T_{un}} / \WG$ is such that the corresponding representation
$\pi_{i,t}$ is irreducible. With Theorem \ref{thm:3.10} we can find an induction datum
$\xi^+ (\pi_{i,t}) = (P,\delta,t_1) \in \Xi_{un}$ such that $\pi_{i,t}$ is a subquotient of 
$\pi^\Gamma (P,\delta,t_1)$. Then $\mu (\WG (w,t_0 t_2))$ is a subquotient of $\pi^\Gamma 
(P,\delta,t_1 t_2)$ for all $t_2 \in \exp (\mf t^w)$, so its central character is $\WG r t_1 t_2$, 
where $W (R_P) r \in T_P / W(R_P)$ is the central character of the $\mc H_P$-representation $\delta$. 
According to \cite[Lemma 3.31]{Opd-Sp} $r \in T_P$ is a residual point for $\mc R_P$, which by 
Proposition 2.63 and Theorem 2.58 of \cite{OpSo2} means that the coordinates of $r$ can be expressed 
as an element of $T_{P,un}$ times a monomial in the variables $\{ q (s)^{\pm 1/2} : s \in S_\af \}$.
Hence we can write $|r| = h_c (q^{1/2})$, where $c$ is the component of $\widetilde{T} / \WG$
containing $\WG (w,t)$ and $h_c : \mc Q (\mc R) \to T_P \subset T$ is a smooth algebraic morphism
with $h_c (1) = 1$. Now $\WG h_c (q^{1/2}) t_0 t_2$ is by construction the central character
of $\mu (\WG (w,t_0 t_2))$. We note that the discrete series representation $\delta_\emptyset$ of
$\mc H_\emptyset = \C$ has central character $1 \in T_\emptyset = \{1\}$, so $h_c = 1$ when
$c$ has dimension rank$(X)$.

Let $\mr{pr}_{q^{1/2}}$ be as in part (d) of Conjecture \ref{conj:5.ABP} and temporarily
denote the difference of two sets by $-$. Then $\mr{pr}_{q^{1/2}} (\widetilde{T} / \WG - T / \WG)$ 
is the set of central characters of $\mu (\widetilde{T} / \WG - T / \WG)$. Since $\mu$
parametrizes irreducible representations, and since every $\pi \in \text{Irr}(\mc H (\mc R,q)
\rtimes \Gamma)$ with central character $t$ is a quotient of the principal series
representation $M(t)$, no element of $\mr{pr}_{q^{1/2}} (\widetilde{T} / \WG - T / \WG)$ can be the 
parameter of an irreducible principal series representation. 

Conversely, suppose that $t \in T$ is not in the aforementioned set. In view of 
Lemma \ref{lem:3.6} we may assume that $t \in T^+$. Then there is, up to isomorphism, 
only one $\pi_t \in \text{Irr}(\mc H (\mc R,q) \rtimes \Gamma)$ with central character $\WG t$. 
In particular all constituents of $M(t)$ are isomorphic to $\pi_t$. Restricted to the finite 
dimensional semisimple algebra $\mc H (W_0,q) \rtimes \Gamma$ this means that the
regular representation $M(t)$ of $\mc H (W_0,q) \rtimes \Gamma$ is a direct sum of copies
of $\pi_t \big|_{\mc H (W_0,q) \rtimes \Gamma}$. But 
$\mc H (W_0,q) \rtimes \Gamma$ has irreducible representations that appear only once
in the regular representation, for example the trivial onedimensional representation 
$N_w \gamma \mapsto q(w)^{1/2}$. Thus there can be only one copy of $\pi_t$ involved and
\[
\pi_t \cong \mc H (W_0,q) \rtimes \Gamma \cong M(t)
\]
as representations of $\mc H (W_0,q) \rtimes \Gamma$. Moreover $\pi_t$ is a subquotient of $M(t)$ 
and both have finite dimension, so they are also isomorphic as $\mc H \rtimes \Gamma$-representations. 
Therefore $\mr{pr}_{q^{1/2}} (\widetilde{T} / \WG - T / \WG)$ is precisely the subset of $t \in T$ for which
the principal series representation $M(t)$ is reducible. 
$\qquad \Box$\\[1mm]

By \eqref{eq:rescosd} $\mr{pr}_{q^{1/2}} (\widetilde{T} / \WG - T / \WG)$ contains all residual 
cosets of dimension smaller than $\dim_\C (T)$. However, in general it is larger, because a
unitary principal series representation can be reducible. In fact there is also a more direct
criterium for irreducibility of principal series representations, in terms of the functions
$c_\alpha$ \cite[Theorem 2.2]{Kat1}.

We note that by Proposition \ref{prop:4.2} the same map $\mr{pr}_v$ also makes part (d) 
valid for the scaled parameter functions $q^\ep$ with $\ep \in \R$. However, for other parameter 
functions changes can occur.

Theorem \ref{thm:5.9} also proves the Aubert--Baum--Plymen conjecture for many Bernstein components
of reductive $p$-adic groups:

\begin{cor}\label{cor:5.10}
Let $\mf s$ be a Bernstein component of a reductive $p$-adic group $G$ such that the algebra 
$\mc H (G)_{\mf s}$ is Morita-equivalent to an extended affine Hecke algebra $\mc H (\mc R ,q) \rtimes 
\Gamma$ in the way described in Section \ref{sec:padic}. (In particular this applies to all Bernstein 
components listed at the end of that section.)

The part of the Aubert--Baum--Plymen conjecture \cite{ABP1,ABP2} corresponding to
parts (b), (c) and (d) of Conjecture \ref{conj:5.ABP} holds for $\mf s$. If moreover $\mc S (G)_{\mf s}$
is Morita equivalent to $\mc S (\mc R ,q) \rtimes \Gamma$, then it is geometrically equivalent to 
$\mc S (\mc R ,1) \rtimes \Gamma$.
\end{cor}
\emph{Proof.}
Let us write $\mf s = [M,\sigma]_G$. As discussed in \eqref{eq:AHAsigma}, the assumed Morita 
equivalence between $\mc H (G)_{\mf s}$ and $\mc H (\mc R ,q) \rtimes \Gamma$ comes from an
isomorphism between $(\mc R ,\Gamma)$ and $(\mc R_\sigma ,\Gamma_\sigma)$, the latter being naturally
associated to $(M,\sigma)$. These yield isomorphisms $W_\sigma \to W_0 \rtimes \Gamma$ and 
$X_{ur}(M_\sigma) \to T$, such that the latter is equivariant with respect 
to the former and restricts to bijections between the unitary and the positive subsets of 
$X_{ur}(M_\sigma)$ and $T$, see \eqref{eq:ModSG}. Moreover we obtain an isomorphism between
$\mc Q (\mc R)$ and the variety of parameter functions for $\mc R_\sigma$.

With these correspondences at hand, Theorem \ref{thm:5.9} this proves parts (b), (c)  and (d) of 
Conjecture \ref{conj:5.ABP} for $\mc H (G)_{\mf s}$. As discussed after \eqref{eq:ModSG}, it is not
unlikely that $\mc S (G)_{\mf s}$ is Morita equivalent to $\mc S (\mc R ,q) \rtimes \Gamma$. 
In that case Lemma \ref{lem:5.8} provides the required geometric equivalence. $\qquad \Box$ 
\\[1mm]

\section{Example: type $C_2^{(1)}$}

In the final section we illustrate what the Aubert--Baum--Plymen conjecture looks like for an 
affine Hecke algebra with $R_0$ of type $B_2 / C_2$ and $X$ the root lattice. More general 
results for type $C_n^{(1)}$ affine Hecke algebras can be found in \cite{Kat2,CiKa}.
For other examples we refer to \cite[Chapter 6]{SolThesis}.

Consider the based root datum $\mc R$ with
\begin{align*}
& X = Y = \Z ,\\
& R_0 = \{ x \in X : \norm{x} = 1 \text{ or } \norm{x} = \sqrt 2 \} ,\\
& R_0^\vee = \{ \ y \in Y : \norm{x} = 2 \text{ or } \norm{x} = \sqrt 2 \} ,\\
& F_0 = \{ \alpha_1 = \vv{1}{-1} , \alpha_2 = \vv{0}{1} \}
\end{align*}
Then $\alpha_4 = \vv{1}{1}$ is the longest root and $\alpha_3^\vee = \vv{2}{0}$ is the
longest coroot, so
\[
S_\af = \{ s_{\alpha_1}, s_{\alpha_2}, t_{\alpha_3} s_{\alpha_3} \} .
\]
We write $s_i = s_{\alpha_i}$ for $1 \leq i \leq 4$ and $s_0 = t_{\alpha_3} s_{\alpha_3}$. 
The Weyl group $W_0$ is isomorphic to $D_4$ and consists of the elements
\[
W_0 = \{ e, \rho_{\pi/2}, \rho_\pi, \rho_{-\pi/2} \} \cup \{ s_1,s_2,s_3,s_4 \} ,
\]
where $\rho_\theta$ denotes the rotation with angle $\theta$. 
The affine Weyl group of $\mc R$ is the Coxeter group
\[
W_\af = W = X \rtimes W_0 = 
\langle s_0, s_1, s_2 | s_i^2 = (s_0 s_2)^2 = (s_1 s_2)^4 = (s_0 s_1)^4 = e \rangle .
\]
Furthermore $R_{nr} = R_0 \cup \{ \pm \alpha_2, \pm \alpha_3 \}$ and
$X^+ = \{ \vv{m}{n} \in X : m \geq n \geq 0 \}$.

We note that $\mc R$ is the root datum of the algebraic group $SO_5$, while
$\mc R^\vee$ corresponds to $Sp_4$. Let $\mathbb F$ be a $p$-adic field whose
residue field has $q$ elements, and let $\mf s$ be the Iwahori-spherical component
of $Sp_4 (\mathbb F)$. Then $\mr{Mod}_{\mf s} (Sp_4 (\mathbb F)) \cong 
\mr{Mod}(\mc H (\mc R,q))$ and Kazhdan--Lusztig theory describes the irreducible
representations in this category with data from $Sp_4 (\mathbb F)$.

But there are many more parameter functions for $\mc R$. Since $s_0, s_1$ and $s_2$
are not conjugate in $W$, we can independently choose three parameters
\[
q_0 = q(s_0) = q_{\alpha_2^\vee / 2}, q_1 = q(s_1) = q_{\alpha_1}, q_2 = q(s_2) = q_{\alpha_2^\vee} .
\]
Several combinations of these parameters occur in Hecke algebras associated to
non-split $p$-adic groups, see \cite{Lus-Uni}. The $c$-functions are
\[
c_{\alpha_1} = \frac{\theta_{\alpha_1} - q_1^{-1}}{\theta_{\alpha_1} - 1} \text{ and }
c_{\alpha_2} = \frac{\theta_{\alpha_2} + q_2^{-1/2} q_0^{1/2}}{\theta_{\alpha_2} + 1}
\frac{\theta_{\alpha_2} - q_1^{-1/2} q_0^{-1/2}}{\theta_{\alpha_2} - 1} .
\]
For $q_0 = q_2$ the relations from Theorem \ref{thm:1.1}.d simplify to
\begin{equation}\label{eq:multq0=q2}
f N_{s_i} = N_{s_i} s_i (f) = (q_i^{1/2} - q_i^{-1/2}) (f - s_i (f)) (1 - \theta_{-\alpha_i})^{-1} 
\qquad i = 1,2 .
\end{equation}
In contrast with graded Hecke algebras, $\mc H (\mc R\!=\!\mc R (SO_5), q_1, q_2 = q_0)$ 
is not isomorphic to $\mc H (\mc R^\vee\!=\!\mc R (Sp_4), q_2, q_1)$. The reason is that in 
$\mc H (\mc R^\vee, q_2, q_1)$ the relation
\[
f N_{s_{\vv{0}{2}}} = N_{s_{\vv{0}{2}}} s_{\vv{0}{2}} (f) = 
(q_2^{1/2} - q_2^{-1/2}) (f - s_{\vv{0}{2}} (f)) (1 - \theta_{\vv{0}{-2}})^{-1} \qquad f \in \mc A
\]
holds, which really differs from \eqref{eq:multq0=q2} because the root lattice 
$\Z \vv{-1}{1} + \Z \vv{0}{2}$ does not equal $X$ for the root datum $\mc R^\vee$.

We will work out the tempered dual of $\mc H (\mc R ,q)$ for almost all positive parameter 
functions $q$. To this end we discuss for each parabolic subalgebra $\mc H_P$ with $P \subset 
\{ \alpha_1, \alpha_2 \}$ separately. Its contribution to Irr$(\mc S (\mc R,q))$ will of course depend 
on $q$, and can even be empty in some cases.
\vspace{3mm}

{\Large $\bullet \qquad P = \emptyset$ } 
\begin{align*}
& X_P = \{0\}, X^P = X, Y_P = \{0\}, Y^P = Y, R_P = \emptyset, R_P^\vee = \emptyset , W(R_P) = \{e\} , \\
& T_P = \{1\}, T^P = T, \mc G_P = W_0, \mc H_P = \C, \mc H^P = \mc A \cong \C [X] .
\end{align*}
We must determine the reducibility of the unitary principal series representations 
\[
M(t) = \mr{Ind}_{\mc A}^{\mc H} \C_t = \pi (\emptyset, \delta_\emptyset,t) \qquad t \in T_{un} .
\]
By Theorem \ref{thm:3.9} $\mr{End}_{\mc H}(M(t))$ is spanned by the intertwining operators
$\pi (w,\emptyset,\delta_\emptyset,t)$ with $w \in W_0$ and $w (t) = t$. For a root $\alpha \in R_0$
with $s_\alpha (t) = t$, Lemma \ref{lem:3.15} tells us that $\pi (s_\alpha,\emptyset,
\delta_\emptyset,t)$ is a scalar if and only if $c_\alpha^{-1}(t) = 0$. \\
\begin{minipage}{11cm}
Let us write $t = (t_3,t_2)$ with $t_i = t(\alpha_i)$. A fundamental domain for the action 
of $W_0$ on $T_{un}$ is $\{ t = (e^{i\phi}, e^{i \psi}) : 0 \leq \psi \leq \phi \leq \pi \} $:
\end{minipage}
\raisebox{-15mm}{\includegraphics[width = 3cm]{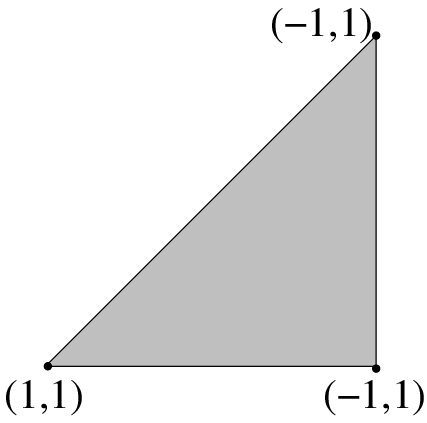}}

The isotropy groups are trivial for all interior points, so $M(t)$ is irreducible for such $t$. 
Below we list the necessary data for all boundary points:
\[
\begin{array}{lllll}
t & W_{0,t} & \text{conditions} & \mr{End}_{\mc H}(M(t)) & \# \text{ irreducibles} \\
\hline
(e^{i \phi},1), \phi \in (0,\pi) & \langle s_2 \rangle & q_0 q_2 \neq 1 & \C & 1 \\
 & & q_0 q_2 = 1 & \C [\langle s_2 \rangle] & 2 \\
(-1, e^{i \psi}), \psi \in (0,\pi) & \langle s_3 \rangle & q_2 \neq q_0 & \C & 1 \\
 & & q_2 = q_0 & \C [\langle s_3 \rangle] & 2 \\
(e^{i \phi},e^{i \phi}), \phi \in (0,\pi) & \langle s_1 \rangle & q_1 \neq 1 & \C & 1 \\
 & & q_1 = 1 & \C [\langle s_1 \rangle] & 2 \\
(-1,1) & \langle s_2,s_3 \rangle & q_0 \neq q_2 , q_0 q_2 \neq 1 & \C & 1 \\
 & & q_0 = q_2 \neq 1 & \C [\langle s_3 \rangle] & 2 \\
 & & q_0 = q_2^{-1} \neq 1 & \C [\langle s_2 \rangle] & 2 \\
 & & q_0 = q_2 = 1 & \C [\langle s_2,s_3 \rangle] & 4 \\
(-1,-1) & W_0 & q_0 \neq q_2, q_1 \neq 1 & \C & 1 \\
 & & q_0 = q_2, q_1 \neq 1 & \C [\langle s_2 \rangle] & 2 \\
 & & q_1 = 1, q_0 \neq q_2 & \C [\langle s_1 \rangle] & 2 \\
 & & q_0 = q_2, q_1 = 1 & \C [W_0] & 5 \\
(1,1) & W_0 & q_0 q_2 \neq 1, q_1 \neq 1 & \C & 1 \\
 & & q_0 = q_2^{-1}, q_1 \neq 1 & \C [\langle s_2 \rangle] & 2 \\
 & & q_1 = 1, q_0 \neq q_2^{-1} & \C [\langle s_3 \rangle] & 2 \\
 & & q_0 = q_2^{-1}, q_1 = 1 & \C [W_0] & 5 \\
\end{array}
\]

{\Large $\bullet \qquad P = \{\alpha_1\}$ } 
\begin{align*}
& X_P = X / \Z \alpha_4 \cong \Z \alpha_1 / 2, X^P = X / \Z \alpha_1, Y_P = \Z \alpha_1^\vee, 
Y^P = \Z \alpha_4^\vee, R_P = \{\pm \alpha_1\}, \\
& R_P^\vee = \{\pm \alpha_1^\vee \}, W(R_P) = \{e, s_1 \}, T^P =\{ t \in T : t_3 = t_2 \}, 
T_P = \{ t \in T : t_2 t_3 = 1 \}, \\
& T^P \cap T_P = \{ (1,1), (-1,-1) \}, \mc G_P = \{e,s_4\} \times T^P \cap T_P , \\
& \mc H_P = \mc H (\mc R_P,q(s_1) = q_1 = q_{\alpha_1^\vee}, \mc H^P = \mc H_P \ltimes \C [X^P] .
\end{align*}
The root datum $\mc R_P$ is of type $C_1^{(1)}$, which means that $R_0^\vee$ is of type 
$C_1 = A_1$ and generates the lattice $Y_P$. For $q_1 = 1$ there are no residual points, for $q_1 \neq 1$
there are two orbits, namely $W (R_P) (q_1^{1/2},q_1^{-1/2})$ and $W (R_P) (-q_1^{1/2}, -q_1^{-1/2})$.
Both orbits carry a unique discrete series representation, which has dimension one. The formulas
for these representations are not difficult, but they depend on whether $q_2 > 1$ or $q_2 < 1$.
So we obtain two families of $\mc H (\mc R,q)$-representations:
\begin{align*}
& \pi \big( \{\alpha_1\}, \delta_1, (t_2,t_2) \big) = 
\mr{Ind}_{\mc H^P}^{\mc H} (\delta_1 \circ \phi_{(t_2,t_2)}) \qquad q_1 \neq 1, t_2 \in S^1 , \\
& \pi \big( \{\alpha_1\}, \delta'_1, (t_2,t_2) \big) = 
\mr{Ind}_{\mc H^P}^{\mc H} (\delta'_1 \circ \phi_{(t_2,t_2)}) \qquad q_1 \neq 1, t_2 \in S^1 .
\end{align*}
The action of $\mc G_P$ on these families is such that $s_4 (t_2,t_2) = (t_2^{-1},t_2^{-1})$, while
$(-1,-1) \in \mc G_P$ simultaneously exchanges $\delta_1$ with $\delta'_1$ and $(t_2,t_2)$ with
$(-t_2,-t_2)$. A fundamental domain for this action is $\{ (\{\alpha_1\},\delta_1,(e^{i \phi},e^{i \phi})) :
\phi \in [0,\pi] \}$. For $\phi \in (0,\pi)$ these points have a trivial stabilizer in $\mc G_P$, so the
corresponding $\mc H$-representations are irreducible. On the other hand, the element $s_4 \in \mc G_P$
fixes the points with $\phi = 0$ or $\phi = \pi$, so the representations $\pi(\{\alpha_1\},\delta_1,(1,1))$ 
and $\pi (\{\alpha_1\},\delta_1,(-1,-1))$ can be reducible. Whether or not this happens depends on more
subtle relations between $q_0,q_1$ and $q_2$.
\vspace{3mm}

{\Large $\bullet \qquad P = \{\alpha_2\}$} 
\begin{align*}
& X_P = X / \Z \alpha_3 \cong \Z \alpha_2, X^P = X / \Z \alpha_2 \cong \Z \alpha_3, 
Y_P = \Z \alpha^\vee_1, Y^P = \Z \alpha_3^\vee / 2, \\
& R_0 = \{ \pm \alpha_2 \}, R_0^\vee = \{\pm \alpha_2^\vee \}, W (R_P) = \{e,s_2\}, \\
& T^P = \{ t \in T : t_2 = 1 \}, T_P = \{ t \in T : t_3 = 1 \}, \mc G_P = \{e,s_3\}, \\
& \mc H_P \cong \mc H ( \mc R_P , q(s_2) = q_2, q_{\alpha_2^\vee} = q_0 ), 
\mc H^P \cong \mc H_P \otimes \C [X^P] .
\end{align*}
The root datum $\mc R_P$ is of type $A_1^{(1)}$, which differs from $C_1^{(1)}$ in the sense that
$X_P$ is the root lattice. There are two orbits of residual points: $W (R_P) (1, q_0^{1/2} q_2^{1/2})$
and $W (R_P) (1, -q_0^{1/2} q_2^{-1/2})$. That is, these points are residual unless they equal
$(1,1)$ or $(1,-1)$. Both orbits admit a unique discrete series representation, of dimension one, which
denote by $\delta_+$ or $\delta_-$. Like for $P = \{\alpha_1\}$, the explicit formulas depend on
which of the $\mc A_P$-characters are in $T_P^{--}$. 
Again we find two families of $\mc H$-representations:
\begin{align*}
& \pi \big( \{\alpha_2\}, \delta_+, (t_3,1) \big) = 
\mr{Ind}_{\mc H^P}^{\mc H} (\delta_+ \circ \phi_{(t_3,1)}) \qquad q_0 q_2 \neq 1, t_3 \in S^1 , \\
& \pi \big( \{\alpha_2\}, \delta_-, (t_3,1) \big) = 
\mr{Ind}_{\mc H^P}^{\mc H} (\delta_- \circ \phi_{(t_3,1)}) \qquad q_0 / q_2 \neq 1, t_3 \in S^1 .
\end{align*}
The group $\mc G_P$ acts on these families by $s_3 (t_3,1) = (t_3^{-1},1)$. A fundamental domain
is, for both families, given by $t \in \{ e^{i\phi} : \phi \in [0,\pi] \}$. For $\phi \in (0,1\pi)$ the 
representations $\pi \big( \{\alpha_2\}, \delta_+, (e^{i \phi},1) \big)$, because the isotropy group
in $\mc G_P$ is trivial. For $\phi \in \{0,\pi \}$ the intertwining operator associated to 
$s_3 \in \mc G_P$ is not necessarily scalar, so we find either one or two irreducible constituents.
Remarkably enough, this depends not only on the parameters $q_0$ and $q_2$ of $\mc H_P$,
but also on $q_1$, as we will see later.
\vspace{3mm}

{\Large $\bullet \qquad P = \{\alpha_2\} = F_0$} \\
Here simply $\mc H^P = \mc H_P = \mc H (\mc R,q)$.
We have to determine all residual points, and how many inequivalent discrete series they carry.
The former can be done by hand, but that is quite elaborate. It is more convenient to use
\cite[Theorem 7.7]{Opd-Sp} and \cite[Proposition 4.2]{HeOp}, which say that for generic $q$
there are 40 residual points. Representatives for the 5 $W_0$-orbits are
\begin{align*}
& r_1 (q) = (q_0^{1/2} q_2^{1/2} q_1^{-1}, q_0^{1/2} q_2^{1/2}) ,\\
& r_2 (q) = (q_0^{-1/2} q_2^{-1/2} q_1^{-1}, q_0^{-1/2} q_2^{-1/2}) ,\\
& r_3 (q) = (-q_0^{1/2} q_2^{-1/2} q_1^{-1}, -q_0^{1/2} q_2^{-1/2}) ,\\
& r_4 (q) = (-q_0^{-1/2} q_2^{1/2} q_1^{-1}, -q_0^{1/2} q_2^{-1/2}) ,\\
& r_5 (q) = (-q_0^{1/2} q_2^{-1/2}, q_0^{-1/2} q_2^{-1/2}) .
\end{align*}
Since every $W_0$-orbit of residual points carries at least one discrete series representation,
\[
\dim_\Q \big( G_\Q^0 (\mc H (\mc R,q)) / G_\Q^1 (\mc H (\mc R,q)) \big) =
\dim_\Q \big( G_\Q^0 (\mc S (\mc R,q)) / G_\Q^1 (\mc S (\mc R,q)) \big) \geq 5.
\]
On the other hand one can easily check, for example with the calculations for $P = \emptyset,
q = 1$, that $\dim_\Q \big( G_\Q^0 (W) / G_\Q^1 (W) \big) = 5$. With \eqref{eq:SprGr}
we deduce that every $W r_i (q)$ carries a unique discrete series representation $\delta (r_i)$.
So far for generic parameter functions.

For nongeneric $q$ the $W r_i (q)$ are still the only possible residual points, but some of them
may cease to be residual for certain $q$. In such cases $r_i (q)$ is absorbed by the tempered
part $r T^P_{un}$ of some onedimensional residual coset $r T^P$. (If $r_i (q)$ is absorbed 
by $T_{un}$, which is the tempered part of the twodimensional residual coset $T$, then
$r_i (q)$ is also absorbed by a onedimensional residual coset.) 
This happens in the following cases:
\[
\begin{array}{lll}
\text{residual point} & q & \text{absorbed by} \\
\hline
r_1 (q) & q_0 q_2 = q_1 & W_0 (q_1^{1/2},q_1^{-1/2}) T^{\alpha_1}_{un} \\
 & q_0 q_2 = q_1^2 & W_0 (1,q_0^{1/2} q_2^{1/2}) T^{\alpha_2}_{un} \\
r_2 (q) & q_0 q_2 = q_1^{-1} & W_0 (q_1^{1/2},q_1^{-1/2}) T^{\alpha_1}_{un} \\
 & q_0 q_2 = q_1^{-2} & W_0 (1,q_0^{1/2} q_2^{1/2}) T^{\alpha_2}_{un} \\
r_3 (q) & q_0 / q_2 = q_1 & W_0 (q_1^{1/2},q_1^{-1/2}) T^{\alpha_1}_{un} \\
 & q_0 / q_2 = q_1^2 & W_0 (1,-q_0^{1/2} q_2^{-1/2}) T^{\alpha_2}_{un} \\
r_4 (q) & q_0 / q_2 = q_1^{-1} & W_0 (q_1^{1/2},q_1^{-1/2}) T^{\alpha_1}_{un} \\
 & q_0 / q_2 = q_1^{-2} & W_0 (1,-q_0^{1/2} q_2^{-1/2}) T^{\alpha_2}_{un} \\
r_5 (q) & q_0 = q_2 & W_0 (1,q_0^{1/2} q_2^{1/2}) T^{\alpha_2}_{un} \\
 & q_0 q_2 = 1 & W_0 (1,-q_0^{1/2} q_2^{-1/2}) T^{\alpha_2}_{un} 
\end{array}
\]
It is also possible that two orbits of residual points confluence, but stay residual.
The deep result \cite[Theorem 3.4]{OpSo2} says that in general situations of this type
the discrete series representations with confluencing central character do not merge
and remain irreducible.

The geometric content of the Aubert--Baum--Plymen conjecture is best illustrated with 
some pictures of the tempered dual of $\mc H (\mc R,q)$, for various $q$. Of course
$T$ has real dimension four, so we cannot draw it. But the unitary principal series can
be parametrized by $T_{un} / W_0$, which is simply a 45-45-90 triangle. 

The other components of Irr$(\mc S (\mc R,q))$ will lie close to $T_{un} / W_0$ if $q$ 
is close to 1, which we will assume in our pictures. We indicate what confluence occurs 
when $q$ is scaled to $1$ by drawing any $\pi \in \text{Irr}(\mc S (\mc R,q))$ close to 
the unitary part of its central character.
To distinguish the three onedimensional components, we denote the series obtained
from inducing $\delta_1 / \delta_+ / \delta_-$ by $L_1 / L_+ / L_-$. 

Finally, we have to represent graphically how many inequivalent irreducible 
representations a given parabolically induced representation $\pi (\xi)$ contains. By
default, $\pi (\xi)$ is itself irreducible. When $\pi (\xi)$ contains two different irreducibles,
we draw the corresponding point fatter. When there are more than two, we write the
number of irreducibles next to it.

\begin{figure}[ht]
\includegraphics[width = 14cm]{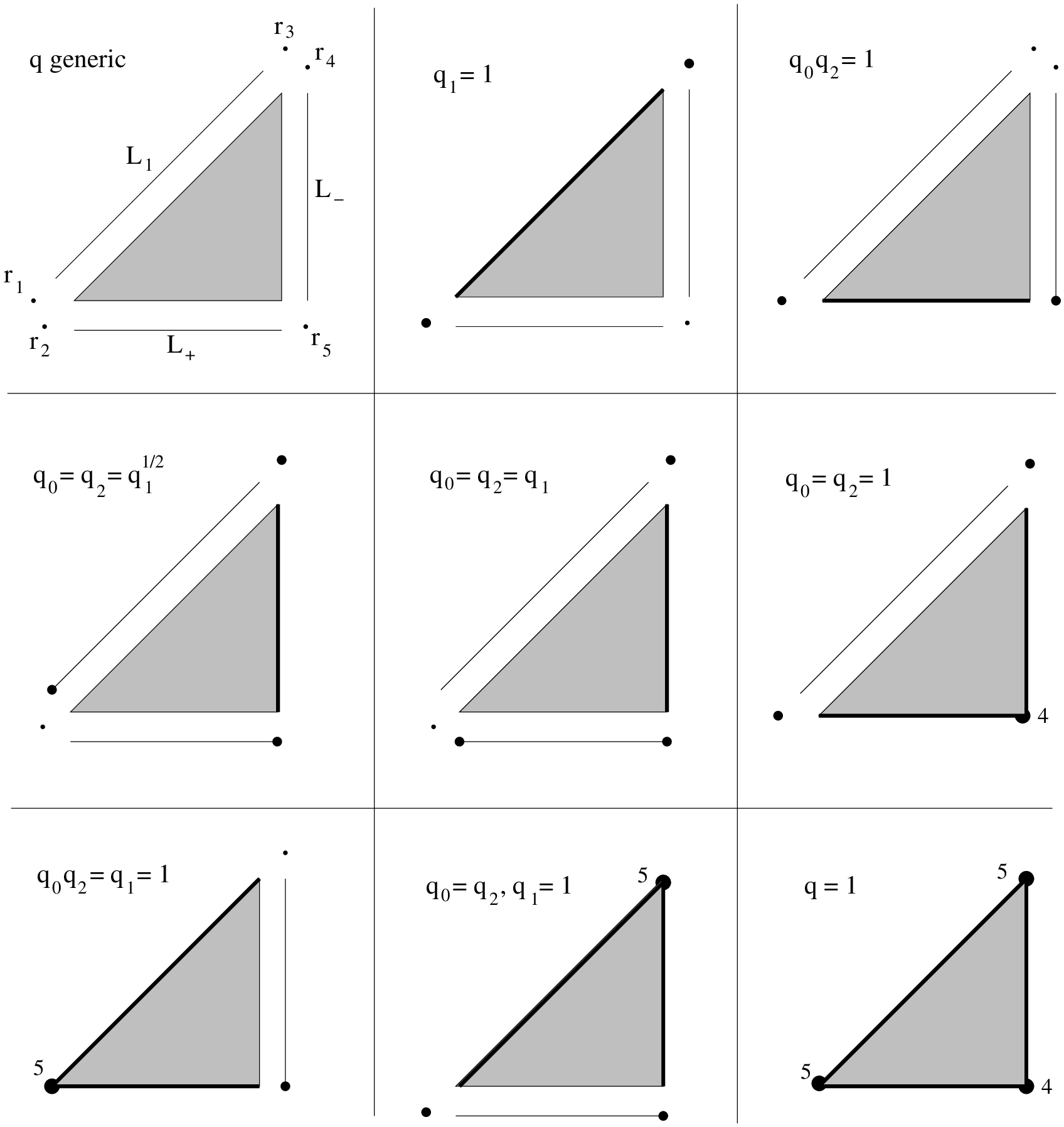}
\caption{the tempered dual for various values of $q$} \label{fig:1}
\end{figure}

Almost everything in Figure \ref{fig:1} can be deduced from the above calculations, the
ABP-conjecture (or rather Theorem \ref{thm:5.9}) and \cite[Theorem 3.4]{OpSo2}.
The only thing that cannot be detected with these methods is what happens at the
confluences $r_1 (q) \to (q_1^{1/2},q_1^{-1/2}) \in L_1$ for $q_0 = q_2 = q_1^{1/2}$
and $r_1 (q) \to (1,q_2) \in L_+$ for $q_0 = q_1 = q_2$. For these $q$ there are four
unitary induction data that give rise to representations with central character in
$T_{rs} / W_0$. So three of them are irreducible, and one contains two different
irreducibles. We can see that the reducibility does not occur in the unitary principal
series or in the discrete series, which leaves the two intermediate series. Fortunately
one can explicitly determine all subrepresentations of $\pi \big( \{\alpha_1\},\delta_1,(1,1) \big)$
(for $q_0 = q_2 = q_1^{1/2}$) and of $\pi \big( \{\alpha_2\},\delta_+,(1,1) \big)$
(for $q_0 = q_1 = q_2$), see \cite[Section 8.1.2]{Slo1}. In fact the graded Hecke
algebras corresponding to these parameter functions are precisely the ones assciated
to $Sp_4$ and to $SO_5$. Hence one may also determine the reducibilty of the 
aforementioned representations via the Deligne--Langlands--Kazhdan--Lusztig parametrization.
Thirdly, it is possible to analyse these parabolically induced representations with R-groups,
as in \cite{DeOp2}.

The comparison of the tempered duals for different parameter functions clearly shows
that this affine Hecke algebra behaves well with respect to general parameter deformations, 
not necessarily of the form $q \mapsto q^\ep$. We see that for small pertubations 
of $q$ it is always possible to find regions $U / W_0 \subset T / W_0$ such that the number
of tempered irreducibles with central character in $U / W_0$ remains stable. 
The $W_0$-types of these representations can change however. It is reasonable to expect
that something similar holds for general affine Hecke algebras, probably that would follow 
from the existence of an appropriate asymptotic Hecke algebra \cite{Lus-Une}.

\end{document}